\documentclass[a4paper,11pt]{report}
\usepackage{amssymb,amsfonts,amsmath}
\usepackage[english,francais]{babel}

\def\d{\delta}

\def\C{\mathbb{C}}

\def\R{\mathbb{R}}
\def\Q{\mathbb{Q}}
\def\Z{\mathbb{Z}}
\def\N{\mathbb{N}}
\def\P{\mathbb{P}}

\def\F{\mathbb{F}}
\def\1{\bold{1}}

\def\a{\alpha}
\def\b{\beta}
\def\e{\varepsilon}
\def\l{\lambda}
\def\f{\varphi}
\def\g{\gamma}
\def\G{\Gamma}
\def\p{\psi}
\def\r{\varrho}
\def\om{\omega}

\def\T{{\mathcal{T}}}

\newtheorem{lem}{Lemme}[chapter]
\newtheorem{pro}[lem]{Proposition}
\newtheorem{defi}[lem]{D\'efinition}
\newtheorem{def/not}[lem]{Definition/Notations}
\newtheorem{Pro/def}[lem]{Proposition/Definition}
\newtheorem{lem/def}[lem]{Lemme/D\'efinition}

\newtheorem{thm}[lem]{Th\'eor\`eme}
\newtheorem{cor}[lem]{Corollaire}
\newtheorem{rqe}[lem]{Remarque}
\newtheorem{rqes}[lem]{Remarques}
\newtheorem{exa}[lem]{Exemple}
\newtheorem{exas}[lem]{Exemples}
\newtheorem{obs}[lem]{Observation}
\newtheorem{ques}[lem]{Question}

\newenvironment{ack}
{\vskip .2cm \noindent {{\bf Remerciements.}}}{\hfill\break}
\newenvironment{preuve}
{\vskip .2cm \noindent {{\it Preuve.}}}{\hfill $\Box$}
\newenvironment{preuve.1.12}
{\vskip .2cm \noindent {{\it Preuve du Lemme 1.12.}}}{\hfill $\Box$}

\newenvironment{preuve.3.1}
{\vskip .2cm \noindent {{\it Preuve du Th\'eor\`eme 3.1.}}}{\hfill $\Box$}

\newenvironment{preuve.3.6}
{\vskip .2cm \noindent {{\it Preuve du Th\'eor\`eme 3.6.}}}{\hfill $\Box$}

\newenvironment{esquisse}
{\vskip .2cm \noindent {{\it Esquisse de preuve.}}}{\hfill $\Box$}

\newenvironment{elements}
{\vskip .2cm \noindent {{\it El\'ements de r\'eponse.}}}{\hfill $\Box$}

\newenvironment{conjG}
{\vskip .2cm \noindent {{\bf Conjecture.}}}{\hfill\break}

\newenvironment{thmA}
{\vskip .1cm \noindent {{\bf Th\'eor\`eme A.}}}{\hfill\break}

\newenvironment{thmD}
{\vskip .1cm \noindent {{\bf Th\'eor\`eme D.}}}{\hfill\break}

\begin{document}

\title{Propri\'et\'es ergodiques des applications rationnelles}
\author{Vincent GUEDJ}
\maketitle

\begin{abstract}

Soit $f:X \rightarrow X$ un endomorphisme rationnel d'une vari\'et\'e projective
complexe compacte. Nous \'etudions la dynamique d'une telle application d'un point de vue
statistique, i.e. nous essayons de d\'ecrire le comportement asymptotique de l'orbite
$O_f(x)=\{f^n(x), \; n \in \N \text{ ou } \Z \}$ d'un point g\'en\'erique.
Nous construisons pour ce faire une mesure de probabilit\'e invariante canonique dont nous \'etudions
les principales propri\'et\'es ergodiques (m\'elange, hyperbolicit\'e, entropie),
et dont nous montrons -- dans certains cas -- qu'elle
refl\`ete une propri\'et\'e d'\'equidistribution des points p\'eriodiques.

\end{abstract}

\chapter*{Pr\'eambule}

Le texte qui suit porte sur l'\'etude statistique de la dynamique des endomorphismes
m\'eromorphes des vari\'et\'es k\"ahl\'eriennes compactes.
Il ne s'agit pas d'un syst\`eme dynamique $(f,X)$ au sens classique du terme:
les endomorphismes $f$ que nous consid\'erons ne sont, en g\'en\'eral, pas bien
d\'efinis sur un sous-ensemble analytique $I_f$ de codimension $\geq 2$ 
-- l'ensemble des points d'ind\'etermination. M\^eme lorsque l'on s'int\'eresse
\`a la dynamique des endomorphismes polynomiaux de $\C^k$, il est important
de consid\'erer la compactification du syst\`eme -- extension m\'eromorphe \`a une
compactification de $\C^k$ -- qui fait, en g\'en\'eral, appara\^itre des points
d'ind\'etermination \`a l'infini.
C'est une source de grande difficult\'e dans l'analyse de la dynamique:
contr\^oler la dynamique pr\`es des points d'ind\'etermination est un enjeu
majeur qui n\'ecessite l'utilisation d'outils relativement sophistiqu\'es 
de G\'eom\'etrie Alg\'ebrique Complexe (probl\`emes de d\'esingularisation
dynamiques) et d'Analyse Complexe (construction, et intersection 
g\'eom\'etrique de courants positifs invariants).

La dynamique des applications de plusieurs variables complexes a connu de 
nombreux d\'eveloppements au cours des quinze derni\`eres ann\'ees, suite aux travaux
fondateurs de E.Bedford, J.E.Fornaess, J.H.Hubbard, M.Lyubich, N.Sibony et J.Smillie.
Ce m\'emoire se propose de faire le point sur une question particuli\`ere, celle
de l'existence d'une mesure d'entropie maximale.

Nous renvoyons le lecteur \`a [CFS], [KH], [W] 
pour une pr\'esentation g\'en\'erale des Syt\`emes
Dynamiques et de la Th\'eorie Ergodiques, 
et \`a [Ci], [De], [GH], [Laz]
pour quelques \'el\'ements d'Analyse et de
G\'eom\'etrie Alg\'ebrique Complexes.
On trouvera dans [BeM], [CG], [Mi 1] une introduction \`a la dynamique holomorphe
en une variable, et dans [F], [M2], [MNTU], [S] une introduction \`a la dynamique holomorphe
en plusieurs variables.

\begin{ack}
C'est un plaisir de remercier tous ceux qui m'ont aid\'e durant l'\'elaboration de 
ce m\'emoire, notamment mes co-auteurs, mes coll\`egues toulousains,
les membres de l'A.C.I. ``Dynamique des applications polynomiales'',
sans oublier ma famille. Une mention sp\'eciale pour mes fr\`eres spirituels,
Charles Favre et Romain Dujardin, ainsi que pour Isabelle et mes parents,
\`a qui je d\'edie ce travail.
\end{ack}

\chapter*{Introduction}

Soit $X$ une vari\'et\'e complexe k\"ahl\'erienne (connexe) compacte et $f:X \rightarrow X$
un endomorphisme m\'eromorphe {\it dominant} (i.e. dont le jacobien ne s'annule pas identiquement).
Nous souhaitons \'etudier la {\it dynamique} de l'application $f$, i.e. d\'ecrire 
statistiquement le comportement
des orbites $O_f(x):=\{f^n(x) \, / \, n \in \N \text{ ou } \Z \}$.
Il s'agit de r\'epondre aussi bien \`a des questions (en apparence) \'el\'ementaires,
\begin{itemize}
\item existe-t-il des points p\'eriodiques ? combien ? de quel type ?
\item comment se r\'epartissent-ils dans l'espace ?
\end{itemize}
qu'\`a des questions plus fines de Th\'eorie Ergodique,
\begin{itemize}
\item quelle est la complexit\'e (entropie topologique) du syt\`eme $(f,X)$ ?
\item existe-t-il une (unique ?) mesure d'entropie maximale ?
\item quelles sont ses propri\'et\'es ergodiques (m\'elange, hyperbolicit\'e,...) ?
\end{itemize}
\vskip.5cm

Lorsque $X$ est une surface de Riemann, la r\'eponse \`a ces questions d\'epend uniquement du degr\'e 
topologique $\deg (f)$ de l'application $f$: lorsque $\deg (f)=1$, la dynamique est pauvre
et se calcule \`a la main.
Lorsque $\deg (f) \geq 2$, on a le r\'esultat suivant 
d\^u \`a H.Brolin [Bro] (cas des polyn\^omes) et M.Lyubich [Ly] (voir \'egalement [FLM]):

\begin{thmA}
{\it 
Si $\dim_{\C} X=1$ et $\deg(f) \geq 2$, alors il existe une mesure de probabilit\'e 
invariante canonique $\mu_f$ qui v\'erifie
\begin{enumerate}
\item $\mu_f$ est l'unique mesure d'entropie maximale
$$
h_{top}(f)=h_{\mu_f}(f)=\log \deg(f)>0 ;
$$
\item $\mu_f$ est m\'elangeante, d'exposant de Lyapunov 
$$
\chi(\mu_f) \geq \frac{1}{2} \log \deg(f)>0;
$$
\item Il y a $(\deg(f))^n+1$ points p\'eriodiques d'ordre $n$. 
Tous --sauf un nombre fini-- sont r\'epulsifs et
s'\'equidistribuent selon la mesure $\mu_f$.
\end{enumerate}
}
\end{thmA}

Le but de ce m\'emoire est d'\'etablir un r\'esultat similaire
au Th\'eor\`eme A lorsque $X$ est de dimension $k \geq 2$: on cherche une
condition num\'erique qui garantisse l'existence d'une mesure canonique
dynamiquement int\'eressante.

A la fin du $20^e$ si\`ecle, E.Bedford, J.E.Fornaess, J.H.Hubbard, M.Lyubich,
N.Sibony et J.Smillie ont construit et \'etudi\'e une telle mesure $\mu_f$
pour deux familles particuli\`eres (mais d\'ecisives) d'endomorphismes:
les applications de H\'enon complexes (i.e. les automorphismes polynomiaux de
$\C^2$ d'entropie positive), et les endomorphismes {\it holomorphes} de
l'espace projectif complexe $\P^k$.
Les dynamiques correspondantes sont tr\`es diff\'erentes: dans le premier cas
les points p\'eriodiques selles s'\'equidistribuent selon $\mu_f$,
alors que dans le deuxi\`eme cas, ce sont les points r\'epulsifs 
qui jouent un r\^ole d\'eterminant (un r\'esultat
plus r\'ecent de J.Y.Briend et J.Duval [BrD 1,2]).

Cette diff\'erence est en partie expliqu\'ee par les valeurs respectives
des degr\'es dynamiques de ces applications.

\vskip.2cm
\noindent {\bf Degr\'es dynamiques.}
Soit $f:\P^k \rightarrow \P^k$ un endomorphisme rationnel
de l'espace projectif complexe $\P^k$.
On d\'efinit, pour $0 \leq j \leq k$, le $j^{\text{i\`eme}}$ degr\'e dynamique de $f$,
$$
\l_j(f):=\liminf_{n \rightarrow +\infty} \left[ \deg f^{-n}L \right]^{1/n},
$$
o\`u $L$ d\'esigne un sous-espace lin\'eaire g\'en\'erique de $\P^k$
de codimension $j$.

Observons que $\l_k(f)$ est le degr\'e topologique de $f$ et que $\l_0(f)=1$.
Lorsque $f$ est un endomorphisme {\it holomorphe}, on v\'erifie
que $\l_j(f)=\l_1(f)^j$, ainsi $\l_k(f)$ est le plus grand degr\'e dynamique
de $f$. A l'inverse, pour une application de H\'enon complexe, on
a $\l_1(f)>\l_2(f)=1$ ($k=2$).

Nous rappelons la d\'efinition 
des degr\'es dynamiques d'un endomorphisme m\'eromorphe
$f:X \rightarrow X$ d'une vari\'et\'e k\"ahl\'erienne compacte quelconque
dans la section 2.1. Ce sont les rapports $\l_j(f)/\l_{j+1}(f)$ qui jouent un
r\^ole crucial dans les ph\'enom\`enes d'\'equidistribution: ceux-ci ont lieu
lorsque $\l_j/\l_{j+1} \neq 1$.
On dira que {\bf f est cohomologiquement hyperbolique} lorsque cette condition
est satisfaite pour tout $0 \leq j \leq k-1$.
Il r\'esulte des propri\'et\'es de concavit\'e de la fonction
$j \mapsto \log \l_j(f)$ (voir Th\'eor\`eme 2.4 de ce texte)
que cette condition est \'equivalente \`a l'existence
d'un entier $l \in [1,k]$ tel que $\l_l(f)$ domine strictement tous les autres degr\'es 
dynamiques.

En dimension $k=1$, $f$ est cohomologiquement hyperbolique si et seulement si
$\l_1(f)>\l_0(f)=1$, i.e. lorsque $\l_1(f)=\deg(f) \geq 2$.
On retrouve la condition du Th\'eor\`eme A. 

En dimension $k=2$, la condition se r\'esume \`a $\l_1(f) \neq \l_2(f)$.
Il y a donc deux cas \`a consid\'erer, selon que
$\l_1(f) <\l_2(f)$ (grand degr\'e topologique, e.g. les endomorphismes holomorphes
de $\P^2$), ou $\l_1(f)>\l_2(f)$ (petit degr\'e topologique, e.g. les applications de
H\'enon).

\vskip.2cm
Les travaux pr\'esent\'es dans ce m\'emoire s'articulent autour de la conjecture
suivante --formul\'ee dans [G 3,7]-- et y apportent des r\'eponses partielles.

\begin{conjG}
{\it 
Supposons que $f$ est {\bf cohomologiquement hyperbolique}, i.e. 
que $\l_l(f)>\max_{j \neq l} \l_j(f)$ pour un entier $l \in [1,k ]$.
Alors il existe une mesure de probabilit\'e invariante canonique $\mu_f$ qui ne
charge pas les hypersurfaces complexes et v\'erifie
\begin{enumerate}
\item $\mu_f$ est l'unique mesure d'entropie maximale
$$
h_{\mu_f}(f)=h_{\text{top}}(f)=\log \l_l(f).
$$

\item $\mu_f$ est m\'elangeante et hyperbolique. Ses exposants de Lyapunov 
v\'erifient
$$
\chi_1 \geq \cdots \geq \chi_l \geq \frac{1}{2} 
\log \left(\l_l(f)/\l_{l-1}(f) \right) >0
$$
et
$$
0>-\frac{1}{2} \log \left(\l_l(f)/\l_{l+1}(f) \right) \geq \chi_{l+1}
\geq \ldots \geq \chi_k.
$$

\item Il y a environ $\l_l(f)^n$ points p\'eriodiques selles de type $(k-l,l)$. 
Ceux-ci s'\'equidistribuent
  selon $\mu_f$.
\end{enumerate}
}
\end{conjG}

\noindent {\bf Les r\'esultats.}
La conjecture est motiv\'ee par les travaux de Bedford-Lyubich-Smillie [BLS 1,2]
qui l'ont \'etablie lorsque $f$ est une application de H\'enon complexe,
ainsi que par ceux de Fornaess-Sibony [FS 2,3,4] et Briend-Duval [BrD 1,2] qui
ont r\'egl\'e le cas des endomorphismes holomorphes non lin\'eaires de $\P^k$.
Nous esquisserons la d\'emonstration de ces r\'esultats ainsi que certaines
de leur g\'en\'eralisations:

\vskip.1cm
\noindent {\bf Th\'eor\`eme B.}
{\it 
La conjecture est vraie lorsque $l=k$, i.e. lorsque le degr\'e topologique 
domine strictement tous les autres degr\'es dynamiques.}
\vskip.1cm

Ce r\'esultat -- qui contient le Th\'eor\`eme A et g\'en\'eralise [BrD 1,2] -- 
est expliqu\'e au Chapitre 3
(dans un souci p\'edagogique, nous d\'etaillons la preuve du Th\'eor\`eme A
au Chapitre 1).

\vskip.1cm
\noindent {\bf Th\'eor\`eme C.}
{\it 
La conjecture est vraie lorsque $f:X \rightarrow X$ est un endomorphisme
birationnel ($\l_2(f)=1$) d'une surface projective ($k=\dim_{\C} X=2$),
moyennant une hypoth\`ese technique.}
\vskip.1cm

Nous pr\'ecisons au Chapitre 4 cette hypoth\`ese
(D\'efinition 4.21), qui porte sur
le contr\^ole quantitatif de la dynamique pr\`es des 
{\it points d'ind\'etermination}
(voir rubrique suivante).
Mentionnons simplement ici qu'elle est v\'erifi\'ee
lorsque $f$ est une application de H\'enon complexe [BLS 1,2],
un automorphisme d'entropie positive [Ca 1],
ou une application birationnelle g\'en\'erique de $\P^2$ [BDi 2], [Du 5].

Notons que tout endomorphisme polynomial cohomologiquement
hyperbolique $f:(z,w) \in \C^2 \mapsto (P(z,w),Q(z,w)) \in \C^2$ avec
$\max (\deg P, \deg Q) \leq 2$ rel\`eve soit du Th\'eor\`eme B,
soit du Th\'eor\`eme C et v\'erifie donc la conjecture [G 3].
De m\^eme, tout endomorphisme {\it holomorphe} cohomologiquement hyperbolique
d'une surface projective v\'erifie la conjecture (Th\'eor\`eme 3.6).

\vskip.2cm
\noindent {\bf Points d'ind\'etermination.}
Les applications $f$ que nous consid\'erons ne sont, en g\'en\'eral, pas partout bien d\'efinies.
Il existe un sous-ensemble analytique $I_f$ de codimension $\geq 2$ constitu\'e de points
en lesquels $f$ n'est pas m\^eme continue. Si $p \in I_f$, son image
$$
f(p):=\bigcap_{\e >0} \overline{ f( B(p,\e) \setminus I_f) }
$$
est un sous-ensemble analytique de dimension positive. On cherche \`a d\'ecrire la dynamique
des points qui se situent hors de l'ensemble d'ind\'etermination it\'er\'e
$$
I_f^{\infty}:= \bigcup_{n \in \Z} f^n(I_f).
$$
En g\'en\'eral cet ensemble peut contenir une hypersurface de $X$ (bien que 
$codim_{\C} I_f \geq 2$). C'est pourquoi il est crucial, dans la conjecture \'enonc\'ee plus haut,
que la mesure $\mu_f$ ne charge pas les hypersurfaces de $X$.
Par ailleurs les degr\'es dynamiques \'etant invariants par un changement birationnel de coordonn\'ees
$\pi$, on veut pouvoir transporter $\mu_f$ par $\pi$: il faut donc que $\mu_f$ ne charge pas
une hypersurface qui pourrait \^etre contract\'ee par $\pi$.

Nous expliquons dans la section 3.3.1 qu'on ne peut pas
se contenter d'\'etudier les endomorphismes
holomorphes: ce sont des objets beaucoup trop rares.
D\`es lors, le contr\^ole de la dynamique pr\`es des points d'ind\'etermi-nation 
constitue une des difficult\'es  majeures de cette \'etude.

\vskip.2cm
\noindent {\bf La strat\'egie.}
L'endomorphisme $f:X \rightarrow X$ induit des actions lin\'eaires $f^*,f_*$ 
par image inverse et directe
sur les espaces vectoriels de cohomologie de deRham $H^q(X,\R)$.
Comme $X$ est k\"ahl\'erienne, on peut tirer profit de la d\'ecomposition de Hodge. L'hypoth\`ese
num\'erique $\l_l(f)>\max_{j \neq l} \l_j(f)$ permet de montrer que les actions ``dominantes''
sont
$$
f^*:H^{l,l}(X,\R) \rightarrow H^{l,l}(X,\R)
\text{ et }
f_*:H^{k-l,k-l}(X,\R) \rightarrow H^{k-l,k-l}(X,\R).
$$
Cela se traduit par exemple
dans la formule des points fixes de Lefschetz (voir paragraphe 2.3.1) qui montre
que le taux asymptotique de croissance des points p\'eriodiques est contr\^ol\'e par $\l_l(f)$, ou
encore dans la majoration de l'entropie topologique due \`a M.Gromov (voir paragraphe 2.2.1)
qui donne ici 
$$
h_{top}(f) \leq \log \l_l(f).
$$
Pour d\'emontrer la conjecture il est naturel de vouloir
utiliser le principe variationnel pour exhiber une mesure d'entropie maximale. La pr\'esence de
points d'ind\'etermination implique malheureusement qu'il n'y a pas toujours \'egalit\'e dans
le principe variationnel (cf [G 7]), du moins lorsque $f$ n'est pas cohomologiquement hyperbolique.
Notre strat\'egie consiste \`a construire une mesure invariante canonique en utilisant l'Analyse
Spectrale des actions lin\'eaires $f^*,f_*$ et des propri\'et\'es fines d'Analyse Complexe 
(notamment la th\'eorie des courants positifs), puis \`a \'etudier les propri\'et\'es ergodiques
de cette mesure en entrem\^elant des outils classiques de Syst\`emes Dynamiques (notamment la th\'eorie
de Pesin) et de G\'eom\'etrie Alg\'ebrique Complexe (notamment le Th\'eor\`eme de Bezout).
La premi\`ere partie repose sur les observations suivantes:

$\bullet$ l'hypoth\`ese $\l:=\l_l(f)>\l_{l-1}(f)$ garantit des ph\'enom\`enes d'\'equidistribu-tion: si
$\om_l$ et $\om_l'$ sont deux formes lisses ferm\'ees cohomologues de bidegr\'e $(l,l)$, alors
$\l^{-n}(f^n)^*(\om_l-\om_l') \rightarrow 0$. On esp\`ere qu'elles s'\'equidistribuent selon un
courant positif ferm\'e invariant $T_l^+$ de bidegr\'e $(l,l)$, i.e.
$$
\frac{1}{\l^n} (f^n)^* \om_l \longrightarrow T_l^+;
$$

$\bullet$l'hypoth\`ese $\l=\l_l(f)>\l_{l+1}(f)$ se traduit par dualit\'e en 
$\l_{k-l}(f_*) >\l_{k-l-1}(f_*)$, donnant de mani\`ere analogue des ph\'enom\`enes 
d'\'equidistribution par image directe des formes lisses ferm\'ees 
$\om_{k-l}$ de bidegr\'e $(k-l,k-l)$.
On esp\`ere de m\^eme que
$$
\frac{1}{\l^n} (f^n)_* \om_{k-l} \longrightarrow T_{k-l}^-,
$$
o\`u $T_{k-l}^-$ d\'esigne un courant positif ferm\'e de bidegr\'e $(k-l,k-l)$;

$\bullet$ on tente alors de d\'efinir la mesure canonique
$$
\mu_f:=T_l^+ \wedge T_{k-l}^-.
$$

La seconde partie consiste \`a donner du sens au slogan suivant: {\it les propri\'et\'es 
g\'eom\'etriques d'extr\'emalit\'e des courants invariants refl\`etent les propri\'et\'es ergodiques
de la mesure $\mu_f$.}

Il est raisonnable de penser que cette strat\'egie va aboutir prochainement
en dimension deux, notamment dans le cas des endomorphismes polynomiaux
(voir Chapitre 4 et [FaJ 3], [DDG]).
La dimension sup\'erieure en est \`a ses premiers balbutiements.
Nous indiquons au Chapitre 5 quelques unes des pistes explor\'ees
jusqu'\`a pr\'esent.

\vskip.2cm
\noindent {\bf Le cas non hyperbolique.}
Nous nous int\'eressons dans ce m\'emoire uniquement aux endomorphismes
m\'eromorphes cohomologiquement hyperboliques. Cela nous permet 
d'\'eviter le cas peu int\'eressant des
produits directs dont un facteur est lin\'eaire. 

On s'attend plus g\'en\'eralement 
\`a ce que les endomorphismes non cohomologiquement hyperboliques
pr\'eservent une fibration: c'est le cas des endomorphismes bim\'eromorphes
des surfaces tels que $\l_1(f)=\l_2(f)=1$, 
comme l'ont montr\'e J.Diller et C.Favre [DF].

Notons que certains endomorphismes m\'eromorphes non cohomologiquement
hyperboliques interviennent dans l'analyse spectrale des op\'erateurs
diff\'erentiels (Laplace, Schr\"odinger)
sur des structues mod\`eles self-similaires, en tant
qu'op\'erateurs de renormalisation (voir [Sa]).

\vskip.2cm
\noindent {\bf Quelles vari\'et\'es ?}
La conjecture sous-tend la question naturelle de savoir quelles sont les vari\'et\'es $X$ qui
admettent des endomorphismes m\'eromorphes cohomologiquement hyperboliques. 
Lorsque $\dim_{\C} X=1$ la r\'eponse 
est bien connue: seules les courbes elliptiques et la sph\`ere de Riemann $\P^1$ admettent des
endomorphismes holomorphes non inversibles. Ce sont les surfaces de Riemann $X$ dont la dimension
de Kodaira $kod(X)$ est n\'egative ou nulle. Nous montrons dans la section 2.4 le

\begin{thmD}
{\it 
Si $\dim_{\C} X=2$ et $f:X \rightarrow X$ est tel que $\l_1(f) \neq \l_2(f)$, alors 
$kod(X) \leq 0$. Plus pr\'ecis\'ement,
\begin{itemize}
\item soit $kod(X)=0$,
\item soit $X$ est rationnelle,
\item soit $X$ est une surface r\'egl\'ee au dessus d'une courbe elliptique; dans ce cas 
 $f$ pr\'eserve la fibration rationnelle et $\l_2(f)>\l_1(f)$.
\end{itemize}
}
\end{thmD}

On peut construire beaucoup d'exemples d'endomorphismes m\'eromorphes
$f:X \rightarrow X$ tels que $\l_1(f)>\l_2(f)$ (resp. $\l_1(f)<\l_2(f)$) lorsque 
$X$ est rationnelle. La situation est beaucoup plus rigide lorsque $kod(X)=0$
(et moins int\'eressante lorsque $X \simeq \P^1 \times E$, $E$ courbe elliptique). 
Le cas le plus important --et le plus d\'elicat \`a \'etudier-- est donc celui d'une surface
$X$ rationnelle: il s'agit d'\'etudier la dynamique d'une application rationnelle dominante
$f:\P^2 \rightarrow \P^2$ \`a conjugaison birationnelle pr\`es.
Il est cependant int\'eressant de consid\'erer \'egalement le cas des surfaces
de dimension de Kodaira nulle. La situation est plus riche qu'en dimension $1$
-- notamment sur les surfaces $K3$ -- et certains endomorphismes ayant un groupe fini de sym\'etries
induisent des endomorphismes sur une surface rationnelle 
(g\'en\'eralisant la construction de S.Latt\`es).

Nous conjecturons, en dimension sup\'erieure, que seules les vari\'et\'es 
de dimension de Kodaira $\leq 0$
admettent des endomorphismes m\'eromorphes cohomologiquement hyperboliques.

\vskip1cm
Pr\'ecisons \`a pr\'esent le contenu du m\'emoire.
Nous esquissons dans la {\it premi\`ere partie} la d\'emonstration du Th\'eor\`eme 
de Brolin-Lyubich. Les \'el\'ements de preuve que nous indiquons
sont en partie originaux et se g\'en\'eralisent en plusieurs variables.

Nous introduisons les degr\'es dynamiques $\l_j(f)$ dans la {\it deuxi\`eme partie}.
Nous \'etablissons quelques unes de leurs propri\'et\'es dans la {\it section 2.1},
leur lien avec l'entropie topologique dans la {\it section 2.2}, puis
avec le nombre de points p\'eriodiques dans la {\it section 2.3}.
Nous donnons dans la {\it section 2.4} plusieurs exemples d'endomorphismes 
cohomologiquement hyperboliques.

Nous \'etudions dans la {\it troisi\`eme partie} le cas des endomorphismes
de grand degr\'e topologique. Nous d\'emontrons la conjecture dans 
les {\it sections 3.1, 3.2} et donnons quelques exemples dans 
la {\it section 3.3}.
Dans la {\it section 3.4} nous \'etudions plus avant
certains invariants num\'eriques (dimension de la mesure d'entropie maximale,
minimalit\'e des exposants de Lyapunov), puis nous nous int\'eressons 
dans la {\it section 3.5} \`a des g\'en\'eralisations de cette situation
dynamique (applications d'allure polynomiale).

Dans la {\it quatri\`eme partie} nous tentons de mettre en place notre strat\'egie
pour d\'emontrer la conjecture lorsque la vari\'et\'e $X$ est de 
dimension deux. Il faut tout d'abord trouver un bon mod\`ele $X$ -- quitte \`a 
effectuer un changement bim\'eromorphe de coordonn\'ees -- sur lequel
l'action lin\'eaire induite en cohomologie est compatible avec la dynamique: nous
abordons ce probl\`eme dans la {\it section 4.1}.
Nous pr\'ecisons l'Analyse Spectrale des op\'erateurs $f^*,f_*$ 
dans la {\it section 4.2}. Nous construisons des courants invariants canoniques
$T^+,T^-$ dans la {\it section 4.3}, et \'etudions leurs propri\'et\'es
g\'eom\'etriques d'extr\'emalit\'e. La mesure canonique $\mu_f=T^+ \wedge T^-$
est construite dans la {\it section 4.4}, puis nous \'etablissons certaines 
de ses propri\'et\'es ergodiques dans la {\it section 4.5}.
Nous illustrons par des exemples ({\it section 4.6}) les difficult\'es rencontr\'ees.

La {\it cinqui\`eme partie} aborde le cas des endomorphismes de petit degr\'e topologique
en dimension sup\'erieure ou \'egale \`a trois.
Nous construisons une mesure canonique d'entropie maximale pour deux familles 
d'exemples, les automorphismes d'entropie positive ({\it section 5.2})
et certains automorphismes polynomiaux de $\C^k$ ({\it section 5.3}).
Plusieurs exemples sont donn\'es dans la {\it section 5.4}.

\tableofcontents

\chapter{La dimension un revisit\'ee}

Nous pr\'esentons dans cette section la d\'emonstration 
du Th\'eor\`eme de Brolin-Lyubich. Nous supposons donc dans toute cette partie
que $X$ est une surface de Riemann compacte munie d'un endomorphisme
$f:X \rightarrow X$ de degr\'e topologique $\l:=\deg(f) \geq 2$. 
Notons que seules les courbes elliptiques $\C/\Lambda$ et la sph\`ere de Riemann $\P^1$
admettent de tels endomorphismes: cela r\'esulte par exemple de la formule
de Hurwitz $f^*K_X+R_f=K_X$ que nous utiliserons \'egalement en dimension sup\'erieure.

 Nous allons donner une construction r\'ecente de la mesure
d'entropie maximale $\mu_f$
inspir\'ee du $dd^c-$lemma de la th\'eorie
de Hodge. Cette technique d\'evelopp\'ee dans [G 1] a depuis \'et\'e appliqu\'ee
avec succ\`es par plusieurs auteurs \`a des probl\`emes de dynamique
\`a plusieurs variables (voir notamment [G 5], [DS 6,8]).
Elle nous permet ici d'\'etablir de fa\c{con} \'el\'ementaire plusieurs propri\'et\'es
ergodiques de $\mu_f$ (m\'elange, d\'ecroissance des corr\'elations,etc).

 Nous donnons \'egalement une d\'emonstration originale --en suivant [G2]--
de l'\'equidistribution des pr\'eimages de points qui utilise des
estim\'ees de volume. Cette approche permet, en plusieurs variables,
d'obtenir des r\'esultats d'\'equidistribution des pr\'eimages de diviseurs.
 Les autres aspects (\'equidistribu-tion des points r\'epulsifs, 
unicit\'e de la mesure d'entropie maximale)
suivent de pr\`es la preuve de M.Lyubich. Ils ne seront qu'esquiss\'es.

\section{La mesure invariante canonique}

H.Brolin a donn\'e dans [Bro] une construction simple de la fonction de Green d'un polyn\^ome.
Trente ans plus tard (soit une dizaine d'ann\'ees apr\`es les travaux de M.Lyubich),
J.H.Hubbard et P.Papadopol [HP 1] ont obtenu une construction \'el\'ementaire de la fonction de Green
d'une fraction rationnelle $f$. Leur construction consiste \`a utiliser la structure complexe homog\`ene
de $\P^1$, en relevant l'application $f$ en un endomorphisme polynomial homog\`ene de $\C^2$.
La technique revient ainsi \`a faire de la dynamique en deux variables, puis \`a utiliser
l'homog\'eneit\'e pour en d\'eduire des informations unidimensionnelles.
Cette construction a le m\'erite de donner des informations sur la r\'egularit\'e du potentiel
de la mesure d'entropie maximale, ainsi que de se g\'en\'eraliser au cas des endomorphismes
de $\P^k$. On peut \'egalement la g\'en\'eraliser \`a des endomorphismes d'autres vari\'et\'es 
homog\`enes (voir [FaG] pour les cas des espaces multiprojectifs). Elle ne s'applique
malheureusement pas au cas de vari\'et\'es non-homog\`enes (qui interviennent
fr\'equemment en dimension sup\'erieure). Elle impose par ailleurs d'inconfortables va-et-vient
entre la dimension 1 et la dimension 2 homog\`ene.
Nous proposons ici une approche plus canonique mise au point dans [G 1].

\subsection{La construction}

Soit $\omega$ une mesure de probabilit\'e lisse sur $X$. Comme $f$ est de degr\'e topologique
$\l$, $\l^{-1}f^*\om$ est encore une mesure de probabilit\'e lisse sur $X$. On peut
donc trouver $\g:X \rightarrow \R$, une fonction lisse, telle que
\begin{equation}
\frac{1}{\l}f^* \om=\om+\Delta \g.
\end{equation}
Nous supposons l'op\'erateur $\Delta$ normalis\'e de sorte 
qu'en coordonn\'ees locales,
$\Delta \log|z|=\e_0$ soit la masse de Dirac en $0$.
Observons que $\g$ est uniquement d\'etermin\'ee \`a une constante additive pr\`es
(principe du maximum). Lorsque $\om$ est la mesure de Haar sur une courbe elliptique $X$,
on obtient en fait $\l^{-1}f^*\om=\om$ et on peut prendre $\g=0$.
Lorsque $\om$ est la forme de Fubini-Study sur $X=\P^1$, on peut prendre
$$
\g(z)=\frac{1}{2 \l}\log [|P(z)|^2+|Q(z)|^2]-\frac{1}{2}\log[1+|z|^2],
$$
o\`u $z$ d\'esigne la coordonn\'ee d'une carte affine $\C$
dans laquelle $f=P/Q$ est donn\'ee par le quotient de deux polyn\^omes
$P,Q$ premiers entre eux, avec $\l=\max(\deg P,\deg Q)$.

Nous consid\'erons \`a pr\'esent les images inverses de l'\'equation (1.1) par les it\'er\'es $f^n$.
Une r\'ecurrence imm\'ediate donne
$$
\frac{1}{\l^n}(f^n)^* \om=\om+ \Delta g_n,
\text{ o\`u } g_n:=\sum_{j=0}^{n-1} \frac{1}{\l^j} \g \circ f^j.
$$
La suite de fonctions $(g_n)$ converge normalement sur $X$ puisque les fonctions
$\g \circ f^j$ sont de m\^eme norme uniforme, \'egale \`a celle de $\g$.
Il s'ensuit que les mesures de probabilit\'e du terme de gauche convergent,
au sens faible des mesures, vers une mesure limite
$$
\mu_f:=\om+\Delta g_f, \text{ o\`u } g_f:=\sum_{j \geq 0} \frac{1}{\l^j} \g \circ f^j.
$$

\subsection{Ensemble de Julia}
Observons que la mesure $\mu_f$ ne d\'epend pas du choix de la forme $\om$: si $\tilde{\om}$ est une autre
mesure lisse de probabilit\'e, alors $\tilde{\om}=\om+\Delta u$ pour une fonction
$u$ lisse donc born\'ee, d'o\`u
$$
\frac{1}{\l^n}(f^n)^* \tilde{\om}=\frac{1}{\l^n}(f^n)^* \om+ 
\Delta \left( \frac{1}{\l^n} u \circ f^n \right) \rightarrow \mu_f.
$$
En particulier, si la suite des it\'er\'es $(f^n)$ est normale dans un ouvert $U$
alors $\mu_f=0$ dans $U$. En effet on peut supposer
que les images $f^{n_i}(U)$, le long d'une sous-suite $n_i \rightarrow +\infty$,
sont toutes incluses dans un petit ouvert $V \subset X$. Quitte \`a restreindre
$U$ et $V$, on peut supposer qu'il existe une mesure de probabilit\'e lisse $\tilde{\om}$
qui est nulle dans $V$. Il s'ensuit que
$$
\mu_f=\lim_{i \rightarrow +\infty} \frac{1}{\l^{n_i}} (f^{n_i})^* \tilde{\om}
=0 \text{ dans } U.
$$
R\'eciproquement supposons que la mesure $\mu_f$ s'annule dans un ouvert $U$.
On obtient, dans cet ouvert,
$$
|(f^n)'|^2 \om=(f^n)^* \om =\l^n \left[ \frac{1}{\l^n} (f^n)^*\om-\mu_f \right]
=\Delta \p_n,
$$
o\`u $\p_n:=\l^n [g_n-g_f]$ est uniform\'ement born\'ee. Soit
$\f \geq 0$ une fonction lisse \`a support compact dans $U$ et qui vaut
1 sur un ouvert $U'$ l\'eg\`erement plus petit. Il vient
$$
\int_{U'} |(f^n)'|^2 \om \leq \int_X \f (f^n)^* \om=\int_X \p_n dd^c \f \leq C_U,
$$
pour une constante $C_U$ ind\'ependante de $n$
qui tend vers $0$ lorsque $\text{diam} \, U \rightarrow 0$. La norme $L^2$ des d\'eriv\'ees
des fonctions $f^n$ est donc uniform\'ement contr\^ol\'ee dans tout ouvert relativement
compact de $U$. Il s'ensuit que la famille $(f^n)_{n \in \N}$ est normale dans $U$.

Rappelons que l'{\it ensemble de Fatou} de $f$ est le plus grand ouvert de $X$ dans lequel
la suite des it\'er\'es est normale. L'{\it ensemble de Julia} est le compl\'ementaire de 
l'ensemble de Fatou. Nous venons donc de montrer la
\begin{pro}
Le support de $\mu_f$ co\"{\i}ncide avec l'ensemble de  Julia de f.
\end{pro}

\subsection{R\'egularit\'e du potentiel}

La r\'egularit\'e du {\it potentiel} $g_f$ de la mesure $\mu_f$ ne d\'epend pas de la forme
$\om$: si $\tilde{\om}=\om+\Delta u$ est une autre mesure de probabilit\'e lisse, il vient
$$
\mu_f=\tilde{\om}+ \Delta \tilde{g_f} \text{ avec } \tilde{g_f}=g_f-u.
$$
Lorsque $f$ est un polyn\^ome de degr\'e $\l \geq 2$ et $\om$ est la forme volume de Fubini-Study
sur $X=\P^1$, le lien entre la fonction d'\'echappement $G_f$ de $f$
(introduite par H.Brolin [Bro]) et le potentiel $g_f$ est donn\'e, pour $z \in \C$, par
$$
G_f(z)=\frac{1}{2}\log[1+|z|^2]+g_f(z).
$$
Nous reviendrons sur le cas des polyn\^omes dans la section 1.5.
N.Sibony a d\'emontr\'e (cf [CG]) que $G_f$ est une fonction H\"old\'erienne
Nous en donnons ci-dessous une preuve \'el\'ementaire qui s'applique \'egalement aux fractions rationnelles
et illustre l'avantage de travailler directement sur $X$.

\begin{pro}
La fonction $g_f$ est H\"old\'erienne d'exposant $\a>0$, pour 
$$
\a< \frac{\log \l}{\chi_{top}(f)}, \text{ o\`u } 
\chi_{top}(f)=\lim_{n \rightarrow +\infty} \frac{1}{n} \log \sup_{x \in X} ||D_xf^n||
$$
d\'esigne l'exposant de Lyapunov topologique de $f$.
\end{pro}

\begin{preuve}
Soit $d$ une distance sur $X$ et $M=\sup_{x \in X} ||D_xf||$. Une r\'ecurrence imm\'ediate
donne pour tous $(x,y) \in X^2$ et $j \in \N$,
$$
d(f^jx,f^jy) \leq M^j d(x,y).
$$
Comme $\g$ est lisse, elle est en particulier H\"old\'erienne d'exposant $\a>0$ pour tout $\a \leq 1$.
Fixons $\a<\log \l/\log M$. Alors
$$
|g_f(x)-g_f(y)| \leq \sum_{j \geq 0} \frac{1}{\l^j} |\g \circ f^j(x)-\g \circ f^j(y)|
\leq C_{\a} d(x,y)^{\a},
$$
o\`u $C_{\a}=\sum_{j \geq 0} M^{\a j}/\l^j <+\infty$. On peut raffiner l'argument pr\'ec\'edent
pour faire intervenir le taux de croissance asymptotique des d\'eriv\'ees en lieu et place de $M$
qui n'est dynamiquement pas significatif.
\end{preuve}
\vskip.2cm

\begin{rqes}
Le caract\`ere H\"old\'erien de $g_f$ a \'et\'e \'etabli pour les endomorphismes
holomorphes de $\P^k$ par J.-Y.Briend et M.Kosek,
puis g\'en\'eralis\'e par de nombreux auteurs (voir [S], [DS 1]). 
La preuve ci-dessus simplifie les approches
pr\'ec\'edentes et permet d'obtenir un contr\^ole du module de continuit\'e de $g_f$ pour
de nombreux endomorphismes {\it m\'eromorphes} (voir [DG]).
Notons que T.C.Dinh et N.Sibony ont \'egalement fait cette observation
dans [DS 6] (Proposition 2.4).

Lorsque $f=f_t$ d\'epend holomorphiquement d'un param\`etre $t$, il suffit de modifier l\'eg\`erement
la preuve ci-dessus pour montrer que la d\'ependance
$(x,t) \mapsto g_{f_t}(x)$ est h\"old\'erienne.
\end{rqes}

Une cons\'equence int\'eressante
est que $\mu_f$ ne charge pas les ensembles de dimension de Hausdorff plus petite que 
$\log \l/\chi_{top}(f)$:

\begin{cor}
Soit $B(p,r)$ un disque centr\'e en $p \in X$ de rayon $r>0$. Alors
$$
\mu_f( B(p,r)) \lesssim r^{\a}.
$$
\end{cor}

\begin{preuve}
Soit $0 \leq \chi \leq 1$ une fonction test telle que $\chi \equiv 1$ au voisinage de
$\overline{B}(p,r)$ et \`a support dans $B(p,2r)$.
Observons que l'on ne change pas la mesure $\mu_f$ si on retranche une constante 
\`a son potentiel $g_f$. On peut donc supposer que $g_f(p)=0$. Il 
r\'esulte alors du th\'eor\`eme de Stokes que,
$$
\mu_f( B(p,r)) \leq \int \chi d\mu_f =\int \chi \om+ \int g_f \Delta \chi
\lesssim r^2+\sup_{B(p,2r)} |g_f|
\lesssim r^{\a},
$$
en observant que l'explosion en $r^{-2}$ des d\'eriv\'ees secondes de $\chi$ est
compens\'ee par l'aire de $B(p,r)$.
\end{preuve}
\vskip.2cm

Un r\'esultat plus fin de A.Manning [Ma] et F.Przytycki [P]
stipule que la
dimension de la mesure $\mu_f$ est \'egal \`a $\log \l/\chi(\mu_f)$, o\`u
$\chi(\mu_f)$ d\'esigne l'exposant de Lyapunov de $\mu_f$ (voir section 1.3).

\section{Equidistribution des pr\'eimages}

La mesure de probabilit\'e limite $\mu_f$ poss\`ede des propri\'et\'es ergodiques remarquables.
Il r\'esulte de sa d\'efinition que $f^* \mu_f=\l \mu_f$. Il s'ensuit que
$\mu_f$ est invariante ($f_* \mu_f=\mu_f$) et d'entropie maximale $=\log \l$, 
car de jacobien constant. Nous reviendrons plus loin sur ces questions d'entropie
(voir section 2.2).
Nous nous int\'eressons ici au comportement de la suite
$\nu_n=\l^{-n}(f^n)^* \nu$, lorsque $\nu$ est une mesure de probabilit\'e quelconque.

\subsection{D\'ecroissance des corr\'elations}

Lorsque $\nu$ a un potentiel continu, le comportement de $\nu_n$ est facile \`a \'etudier.
Nous illustrons cette observation en donnant une preuve \'el\'ementaire de la d\'ecroissance 
exponentielle des corr\'elations.

\begin{thm}
La mesure $\mu_f$ est m\'elangeante. Plus pr\'ecis\'ement il existe $C>0$ telle que pour toutes
fonctions test $\f,\p$, 
$$
\left| \int_X \p \cdot \f \circ f^n d\mu_f-\int_X \f d\mu_f \int_X \p d\mu_f \right|
\leq C ||\p||_{{\mathcal C}^2} ||\f ||_{L^{\infty}} \l^{-n}.
$$
\end{thm}

\begin{preuve}
 Posons $c_{\f}:=\int_X \f d\mu_f$ et $c_{\p}=\int_X \p d\mu_f$. 
Quitte \`a translater et dilater $\f$, on peut supposer $\f \geq 0$ et $c_{\f}=1$.
Comme $\mu_f$ a un potentiel continu, il en va de m\^eme de $\f \mu_f$.
Ces mesures positives ont m\^eme masse,
donc
$$
\f \mu_f=\mu_f+\Delta u_{\f},
$$
o\`u $u_{\f} \in {\mathcal C}^0(X,\R)$ est uniquement d\'etermin\'ee
si on impose $\sup_X u_{\f}=0$.
Il r\'esulte de l'ellipticit\'e du Laplacien que la norme ${\mathcal C}^0$
de $u_{\f}$ est contr\^ol\'ee par celle de $\f$: $||u_{\f}||_{L^{\infty}} \leq C ||\f||_{L^{\infty}}$.
En observant que $\f \circ f^n \mu_f=\l^{-n}(f^n)^*(\f \mu_f)$, il vient donc
\begin{eqnarray*}
\lefteqn{ \! \! \! \! \! \!
\left|\int_X \p \cdot \f \circ f^n d\mu_f-c_{\f}c_{\p}\right|
=\left| \langle \p, \l^{-n} (f^n)^*(\f \mu_f)-c_{\f} \mu_f\rangle  \right|} \\
&&=\left|\langle \p,\l^{-n} (f^n)^* (\f\mu_f-c_{\f} \mu_f) \rangle  \right| 
= \left| \langle \Delta \p, \l^{-n} u_{\f} \circ f^n \rangle  \right| \\
&&\leq \frac{C'}{\l^n}  ||\p||_{{\mathcal C}^2} ||\f ||_{L^{\infty}}.
\end{eqnarray*}
\end{preuve}

\begin{rqe}
La vitesse de m\'elange pour les endomorphismes de la sph\`ere de Riemann a \'et\'e 
estim\'ee dans [FS 4]. 
\end{rqe}

\subsection{Estim\'ees de volume}

Nous nous int\'eressons \`a pr\'esent au cas plus difficile o\`u $\nu=\e_a$ est une masse de Dirac
en un point $x \in X$. Nous utilisons pour ce faire des estim\'ees de volume.
Notre outil principal est l'estimation suivante
qui assure que le volume d'un bor\'elien quelconque ne d\'ecro\^it pas trop vite sous it\'eration.

\begin{thm}
Il existe $C>0$ tel que pour tout $j \in  \N$, $\Omega \subset X$,
$$
\text{Vol}(f^j \Omega) \geq \exp \left( -\frac{C}{\text{Vol}(\Omega)} \l^j \right).
$$
\end{thm}

Les volumes sont calcul\'es ici par rapport \`a une forme volume
$\om$ fix\'ee.

\begin{preuve}
On commence par observer que si $\Omega$ est un bor\'elien, alors
$f_* \1_{\Omega} \leq \l \1_{f(\Omega)}$, donc pour tout $j \in \N$,
$\1_{f^j(\Omega)} \geq \l^{-j} (f^j)_* \1_{\Omega}$.
Il s'ensuit
$$
\text{Vol} (f^j \Omega) \geq \frac{1}{\l^j} \int_{\Omega} (f^j)^* \om
=\frac{1}{\l^j} \int_{\Omega} |(f^j)'|_{\om}^2 \, \om,
$$
o\`u $|f'|_{\om}^2$ d\'esigne le jacobien (r\'eel) de $f$ par rapport
\`a la forme volume $\om$. La concavit\'e du logarithme implique alors
$$
\text{Vol} (f^j \Omega) \geq \frac{\text{Vol}(\Omega)}{\l^j} 
\exp \left[ \frac{2}{\text{Vol} (\Omega)} \int_{\Omega}
\log |(f^j)'|_{\om} \, \om \right].
$$
Pour minorer cette derni\`ere int\'egrale, on observe que $\log |f'|_{\om}$
est une diff\'erence de fonctions quasisousharmoniques, en particulier
$$
\log |f'|_{\om} \geq u, \; \text{ avec } u \leq 0 
\; \text{ et } \; dd^c u \geq -C_1 \om,
$$
pour une constante $C_1$ suffisamment grande. Comme 
$(f^j)'=\Pi_{l=0}^{j-1} f' \circ f^l$, on obtient
$$
\int_{\Omega} \log |(f^j)'|_{\om} \, \om \geq 
\sum_{l=0}^{j-1} \l^l \left[ \int_X \frac{1}{\l^l} u \circ f^l \om \right].
$$
Or la suite de fonctions $(\l^{-l} u \circ f^l)$ est relativement compacte dans 
$L^1(X)$. En effet $u_l:=C_1 g_l+\l^{-l} u \circ f^l$ est une suite uniform\'ement
major\'ee de fonctions $C_1\om$-sousharmoniques telles que
$$
\int_X u_l d\mu_f=\sum_{p=0}^{l-1} \frac{1}{\l^p} \int_X \g d\mu_f
+\frac{1}{\l^l} \int_X u d\mu_f \geq \frac{\l}{\l-1} \int_X \g d\mu_f+\int_X u d\mu_f>-\infty.
$$
On en d\'eduit qu'il existe $C_2>0$ telle que
$$
\text{Vol} (f^j \Omega) \geq \frac{\text{Vol}(\Omega)}{\l^j} 
\exp \left( -\frac{C_2}{\text{Vol}(\Omega)} \l^j \right).
$$
L'estim\'ee annonc\'ee r\'esulte enfin d'in\'egalit\'es \'el\'ementaires.
\end{preuve}
\vskip.2cm

Observons que lorsque $f$ est (localement) conjug\'e \`a $z^{\l}$ pr\`es
de $0$, alors $f^j \sim z^{\l^j}$ et le volume de petits disques $D$ 
centr\'es \`a l'origine de rayon $r>0$, d\'ecro\^it sous it\'eration 
comme $r^{2 \l^j}=\exp( -2[-\log r] \l^j) $. L'estimation du Th\'eor\`eme 1.7 est
donc essentiellement optimale.

\begin{rqe}
Les premi\`eres estim\'ees de volume ont \'et\'e mises
au point dans le cadre des applications birationnelles de $\C^2$
pour construire et caract\'eriser les courants invariants par la dynamique
[BS 1], [FS 1], [Dil 1], [FaG]. Ces applications sont presque inversibles.
J.E.Fornaess et N.Sibony ont \'etabli dans [FS 4] des estim\'ees
de volume pour les endomorphismes holomorphes de $\P^2$, dont les
ramifications sont plus importantes. C.Favre et M.Jonsson ont 
men\'e \`a bien une \'etude syst\'ematique de ces estim\'ees [FJ 1].

L'estimation du Th\'eor\`eme 1.7 est tir\'ee de [G2], [G4], o\`u nous
d\'eveloppons des estim\'ees moins fines, mais
valables dans un contexte plus g\'en\'eral. 
Elles sont suffisantes pour un certain nombre d'applications.
Nous l'illustrons ici en d\'emontrant l'\'equidistribution des
pr\'eimages de points non-exceptionnels. 
\end{rqe}

Soit $\nu$ une mesure de probabilit\'e sur $X$. Nous cherchons \`a savoir sous quelle condition
la suite de mesures de probabilit\'e $\nu_n=\l^{-n} (f^n)^* \nu$ converge vers $\mu_f$.
Soit $u \in L^1(X)$ un potentiel de $\nu$, i.e.
$$
\nu=\om+\Delta u.
$$
La fonction $u$ est une fonction {\it quasisousharmonique}, i.e. localement diff\'erence d'une fonction lisse
et d'une fonction sousharmonique. En particulier $u$ est s.c.s. donc major\'ee sur $X$: quitte \`a
translater $u$, on supposera dans la suite $u \leq 0$.
Il r\'esulte du Th\'eor\`eme de repr\'esentation de Riesz que
$e^{-\a u} \in L^1(X)$, pour $\a>0$ suffisamment petit: il faut choisir 
$\a<2 [\sup_{x \in X} \nu(\{x\}) ]^{-1}$. Notons en particulier que
$e^{-\a u } \in L^1(X)$ pour tout $\a>0$, lorsque $\nu$ {\it n'a pas d'atome}. Observons que
$$
\nu_n:=\frac{1}{\l^n} (f^n)^* \nu=\frac{1}{\l^n} (f^n)^*\om+\Delta u_n, 
\text{ o\`u } u_n:=\frac{1}{\l^n} u \circ f^n.
$$
Montrer que $\nu_n$ converge vers $\mu_f$ revient ainsi \`a montrer que $u_n$ converge vers $0$
dans $L^1(X)$. Il r\'esulte du carat\`ere quasisousharmonique de $u$ que
la suite $(u_n)$ est relativement compacte dans $L^1(X)$
(voir [GZ 1] pour des crit\`eres de compacit\'e des fonctions qpsh). 
Consid\'erons, pour $\e>0$ fix\'e,
$$
\Omega_n^{\e}:=\left\{ x \in X \, / \, u_n<-\e \right\}.
$$
Il suffit de savoir montrer que $\text{Vol}(\Omega_n^{\e}) \rightarrow 0$ lorsque $n \rightarrow +\infty$
pour conclure \`a $u_n \rightarrow 0$. On observe que
$$
f^n(\Omega_n^{\e})=\left\{ x \in X \, \/ \, u(x)<-\e \l^n \right\}.
$$
Or le volume de l'ensemble de gauche ne d\'ecro\^it pas trop vite gr\^ace aux estim\'ees de volume
(Th\'eor\`eme 1.7), et celui du membre de droite
d\'ecro\^it exponentiellement vite: comme $e^{-\a u} \in L^1(X)$, il r\'esulte de l'in\'egalit\'e
de Chebyshev que pour tout $A<2 [\sup_{x \in X} \nu(\{x\})]^{-1}$, il existe $C_A>0$ tel que
$$
\text{Vol}(u<-\e\l^n) \leq C_A \exp (-A\e\l^n).
$$
En mettant bout \`a bout ces in\'egalit\'es, on obtient
$$
\limsup_{n \rightarrow +\infty} \text{Vol}(\Omega_n^{\e}) \leq \frac{1}{A \e}.
$$
Lorsque $\nu$ n'a pas d'atome, on peut prendre $A$ arbitrairement grand et conclure \`a
$\nu_n \rightarrow \mu_f$. En rempla\c{c}ant $u$ par $u_{n_0}$ dans le raisonnement pr\'ec\'edent,
on obtient la version plus pr\'ecise:
$$
 \nu_n \rightarrow \mu_f  \; \text{ ssi } \; 
\sup_{x \in X} \nu_n(\{x\}) \rightarrow 0 .
$$
Il reste donc \`a analyser la partie atomique de la suite $\nu_n$. Elle est contr\^ol\'ee
par le degr\'e topologique local $d(f^n,x)$ de $f^n$ au point $x$ via la relation
$$
\nu_n(\{x\})=\frac{d(f^n,x) \nu(\{f^nx\})}{\l^n}.
$$
Il faut donc contr\^oler le comportement asymptotique de la suite $(d(f^n,x))$
par rapport \`a $\l^n$.
C'est l'objectif du lemme suivant dont la preuve est laiss\'ee au lecteur.

\begin{lem/def}
On d\'efinit l'ensemble des points exceptionnels par
$$
{\mathcal E}_f:=\{ x \in X \, / \, d(f,x)=d(f,f(x))=d(f,f^2x)=\l \}.
$$
Alors $\text{Card}({\mathcal E}_f) \leq 2$. Plus pr\'ecis\'ement,
\begin{itemize}
\item soit ${\mathcal E}_f$ est vide;
\item soit ${\mathcal E}_f$ est constitu\'e d'un point; dans ce cas $X=\P^1$ et $f$ est conjugu\'e
\`a un polyn\^ome par une homographie qui envoie ${\mathcal E}_f$ sur $\{\infty\}$;
\item soit ${\mathcal E}_f$ est consitu\'e de deux points; dans ce cas $X=\P^1$ et $f$ est conjugu\'e
\`a $z^{\l}$ ou $z^{-\l}$ par une homographie qui envoie ${\mathcal E}_f$ sur 
$\{0,\infty\}$.
\end{itemize}
\end{lem/def}

Il s'av\`ere que ${\mathcal E}_f$ est le plus grand ensemble fini totalement invariant,
$f^{-1}({\mathcal E}_f)={\mathcal E}_f$ (voir [BeM], [CG], [Mi 1]
pour ce point de vue plus traditionnel).
Il r\'esulte de ce qui pr\'ec\`ede qu'on a le r\'esultat
d'\'equidistribution suivant:

\begin{thm}
Soit $\nu$ une mesure de probabilit\'e sur $X$. Alors
$$
\frac{1}{\l^n}(f^n)^* \nu \rightarrow \mu_f \; \text{ ssi } \; \nu({\mathcal E}_f)=0.
$$
En particulier $\l^{-n}(f^n)^* \e_a \rightarrow \mu_f$ ssi $a$
n'est pas un point exceptionnel.
\end{thm}

\begin{preuve}
L'ensemble ${\mathcal E}_f$ est totalement invariant: les pr\'eimages d'une mesure 
qui a un atome
dans ${\mathcal E}_f$ ne peuvent donc pas s'\'equidistribuer selon $\mu_f$.

Supposons r\'eciproquement qu'une mesure de probabilit\'e
$\nu$ n'a pas d'atome dans ${\mathcal E}_f$.
Le degr\'e topologique local $d(f^n,x)=\Pi_{j=0}^{n-1} d(f,f^jx)$
de $f^n$ en un point $x$ en lequel $\nu$ a un atome est donc major\'e par $\g^n$,
avec $\g<\l$.
Il s'ensuit que $\sup_{x \in X} \nu_n(\{x\}) \rightarrow 0$,
donc $\nu_n \rightarrow \mu_f$ d'apr\`es la discussion qui pr\'ec\`ede
le Lemme 1.9.
\end{preuve}

\section{Exposant de Lyapunov}

Nous avons construit ci-dessus une mesure $\mu_f$ ergodique dont les potentiels sont continus. 
En particulier $\log^+ ||Df|| \in L^1(\mu_f)$.
Il r\'esulte du Th\'eor\`eme ergodique de Birkhoff que pour $\mu_f$-presque tout $x$,
la suite $n \mapsto$ $n^{-1} \log ||Df^n_x||$ converge vers un nombre 
$\chi(\mu_f)$ ind\'ependant de $x$: c'est l'exposant de Lyapunov de $\mu_f$.
Comme $\mu_f$ est d'entropie $\log \l$, l'in\'egalit\'e de Margulis-Ruelle
assure [R] que
$$
\chi(\mu_f) \geq \frac{1}{2} \log \l >0.
$$
Nous indiquons un peu plus loin une preuve \'el\'ementaire et ``complexe'' de cette in\'egalit\'e due \`a
M.Lyubich. Nous cherchons \`a savoir pour l'instant s'il est possible d'am\'eliorer
cette in\'egalit\'e.

Lorsque $X=\C/\Lambda$ est une courbe elliptique,
$f$ est un rev\^etement \'etale dont on montre
ais\'ement qu'il est induit par une application affine $z \mapsto az+b$
qui pr\'eserve le r\'eseau $\Lambda$.
Il s'ensuit que $f^* dz=a dz$ et
la mesure de Lebesgue $\mu=\frac{i}{2}dz \wedge \overline{dz}$ est de jacobien
constant, $f^* \mu= |a|^2\mu$, avec $|a|^2=\l$: on a donc $\mu=\mu_f$. Observons \'egalement
que $f$ est de diff\'erentielle constante, $D_xf=a Id$, donc
$$
\chi(\mu_f)=\frac{1}{2} \log \l
$$
dans ce cas. Comme $f$ commute avec l'involution
$\sigma:x \mapsto \!-x$, elle induit un endomorphisme de degr\'e $\l$ sur 
le quotient $X/ \!\langle  \! \sigma \! \rangle  \simeq \P^1$. 
Un tel endomorphisme de la sph\`ere de Riemann est un
{\it exemple de Latt\`es}.

Les endomorphismes de la sph\`ere de Riemann $X=\P^1$ sont les fractions rationnelles
$f=P/Q$ o\`u $P,Q$ sont des polyn\^omes sans facteur commun
dont les degr\'es v\'erifient
$\max( \deg P,\deg Q)= \l$. Lorsque $f$ est un exemple de Latt\`es, la mesure $\mu_f$
est absolument continue par rapport \`a la mesure de Lebesgue et on obtient
$\chi(\mu_f)=\frac{1}{2} \log \l$. On ne peut donc pas, en g\'en\'eral, 
am\'eliorer la minoration de l'exposant de Lyapunov de $\mu_f$. 
Il r\'esulte cependant des travaux de F.Ledrappier [L] et A.Zdunik [Z] que cette \'egalit\'e
ne se produit que lorsque $f$ est un exemple de Latt\`es. 
Il est donc naturel d'\'etudier comment $\chi(\mu_f)$ varie lorsque l'on fait bouger
$f$ dans l'espace des param\`etres.

Soit $(f_t)_{t \in M}$ une famille holomorphe d'endomorphismes de $\P^1$ de degr\'e $\l \geq 2$,
o\`u $M$ d\'esigne une vari\'et\'e complexe. Rappelons que la fonction de Green
dynamique de $f_t$ est
$$
g_t(x):=\sum_{j \geq 0} \frac{1}{\l^j} \g_t \circ f_t^j(x),
$$
o\`u $\g_t \in {\mathcal C}^1(X,\R)$ est d\'etermin\'e de fa\c{c}on unique par
$$
\frac{1}{\l}f_t^* \om=\om+\Delta \g_t \; \; \text{ et } \; \int_{\P^1} \g_t \om=0.
$$
Notons que $(x,t) \mapsto g_t(x)$ est une fonction H\"old\'erienne (voir Remarques 1.3).
On note $\mu_t=\om+\Delta g_t$ la mesure d'entropie maximale de $f_t$, et
$$
\chi(t)=\int_{\P^1} \log|f_t'| d\mu_t
$$
l'exposant de Lyapunov de $(f_t,\mu_t)$.

Le point de vue potentialiste permet de donner une preuve \'el\'ementaire du
r\'esultat suivant --qui am\'eliore le Th\'eor\`eme B de R.Man\~e [Mn].

\begin{pro}
L'exposant $\chi(t)$ est une fonction plurisousharmonique et h\"old\'erienne de $t$.
\end{pro}

Notons que la r\'egularit\'e h\"old\'erienne de l'exposant de Lyapunov 
peut \^etre \'egalement vue comme une cons\'equence d'une formule
explicite de G.Bassanelli et F.Berteloot [BaBe],
valable en toute dimension. Le caract\`ere psh est d\'emontr\'e dans un cadre
plus g\'en\'eral par T.C.Dinh et N.Sibony (voir Proposition 3.8.3 dans [DS 1]),
en suivant l'approche de E.Bedford et J.Smillie (voir Th\'eor\`eme 5.5 dans [BS 2]).

\begin{preuve}
On int\`egre par parties pour obtenir
$$
\chi(t)=\int \log |f_t'| \om+\sum_{f_t'(c_t)=0} g_t(c_t),
$$
o\`u les $c_t$ d\'esignent les points critiques de $f_t$.
La premi\`ere somme est une fonction lisse de $t$. Comme $(x,t) \mapsto g_t(x)$
est h\"old\'erienne, et comme la seconde somme est une fonction sym\'etrique des points critiques,
racines de l'\'equation alg\'ebrique $f_t'(x)=0$, on obtient la continuit\'e annonc\'ee.

La plurisousharmonicit\'e provient de ce que
$\chi$ est limite uniforme de la suite de fonctions psh
$$
\chi_n(t)=\int_{\P^1} \log |f'_t| \frac{1}{\l^n} (f_t^n)^* \om
=\frac{1}{\l^n} \int_{\P^1} (f_t^n)_* (\log |f'_t| ) \cdot \om.
$$
La convergence r\'esulte de ce que $\Delta \log |f_t'|$ peut s'exprimer comme
une diff\'erence de deux mesures de Radon positives, ainsi
$$
|\chi(t)-\chi_n(t)|=\left| \int_{\P^1} \Delta \log |f_t'| \cdot (g_t-g_{n,t}) \right|
\leq C ||g_t-g_{n,t}||_{L^{\infty}(\P^1)}.
$$
\end{preuve}
\vskip.2cm

Il est int\'eressant d'\'etudier les propri\'et\'es des courants positifs $(dd^c \chi)^j$,
$1 \leq j \leq \dim_{\C} M$. Le fait que $t \mapsto \chi(t)$ soit plurisousharmonique
et atteigne un minimum aux valeurs de $t$ pour lesquelles $f_t$ est un exemple de
Latt\`es montre que ceux-ci jouent un r\^ole sp\'ecial dans l'espace des param\`etres.
Nous renvoyons le lecteur \`a [DeM], [BaBe], [DuF] pour quelques r\'esultats dans cette direction.

\section{L'approche de Lyubich}

M.Lyubich a le premier donn\'e une construction de la mesure d'entropie maximale pour les
endomorphismes rationnels $f$ de $\P^1$ ([Ly], voir \'egalement [FLM]). Son approche consiste
\`a \'etudier les branches inverses de $f^n$ pour en d\'eduire des propri\'et\'es de
l'op\'erateur de Ruelle
$$
\begin{array}{rcl}
R_f: {\mathcal C}^0(\P^1) & \rightarrow & {\mathcal C}^0(\P^1) \\
\f & \mapsto & \frac{1}{\l}f_* \f.
\end{array}
$$
Etudier le comportement de la suite d'op\'erateurs $\l^{-n}(f^n)^*(\cdot)$ sur les mesures
revient en effet, par dualit\'e, \`a it\'erer l'op\'erateur $R_f$. Bien que la construction
de $\mu_f$ due \`a M.Lyubich soit techniquement plus compliqu\'ee que celle donn\'ee plus haut,
son \'etude des branches inverses lui permet de montrer que
\begin{itemize}
\item $\chi(\mu_f) \geq \frac{1}{2} \log \l>0$;
\item les points r\'epulsifs s'\'equidistribuent selon $\mu_f$;
\item $\mu_f$ est l'unique mesure d'entropie maximale.
\end{itemize}
Nous allons esquisser les preuves de ces r\'esultats qui reposent toutes sur le lemme suivant:

\begin{lem} Soit ${\mathcal C}_f$ l'ensemble critique de $f$.
Soit $\e>0$. Il existe $l=l_{\e} \in \N$ tel que pour tout $n \geq l$ et
pour tout disque $D$ qui \'evite $\cup_{j=1}^l f^j({\mathcal C}_f)$, on peut
construire $(1-\e)\l^n$ branches inverses $f_i^{-n}$ de $f^n$ sur $D$ avec
$$
\text{diam}(f_i^{-n} D) \leq C_{\e} \l^{-n/2}.
$$
\end{lem}

\vskip.1cm
Nous d\'emontrons le lemme un peu plus loin et commen\c{c}ons par en tirer les cons\'equences
annonc\'ees.
\vskip.1cm

\noindent {\bf Exposant de Lyapunov}.
Soit $x \in X \setminus \cup_{j \geq 0} f^j({\mathcal C}_f)$. C'est un point g\'en\'erique
pour $\mu_f$. On peut consid\'erer des branches inverses $f_i^{-n}$ de $f^n$ sur un petit disque
$D(x,\e)$. Comme le diam\`etre de $f_i^{-n} D(x,\e)$ d\'ecro\^it en $\l^{-n/2}$,
il r\'esulte des in\'egalit\'es de Cauchy que
$$
|(f^n)'(x^{-n})|=\left| \frac{1}{(f^{-n})'(x)} \right| \geq \frac{1}{C_{\e}} \l^{n/2},
$$
o\`u $x^{-n}$ d\'esigne l'une des pr\'eimages de $x$. Notons $x^{-1},\ldots,x^{-n}$ une
pr\'ehistoire g\'en\'erique de $x$ choisie parmi les
$(1-\e) \l^n$ pr\'eimages construites dans le Lemme 1.12. Il r\'esulte du Th\'eor\`eme de Birkhoff
que
$$
\chi(\mu_f)=\lim \frac{1}{n} \sum_{i=0}^{n-1} \log |f'(x^{-i})|=
\lim \frac{1}{n} \log |(f^n)'(x^{-n})| \geq \frac{1}{2} \log \l.
$$
\vskip.1cm

\noindent
{\bf Points r\'epulsifs}. 
Rappelons qu'il y a $\l^n+1$ points p\'eriodiques d'ordre $\leq n$, compt\'es avec multiplicit\'e.
Cela r\'esulte de la formule des points fixes de Lefschetz (voir section 2.3).
On peut \'egalement s'en convaincre ici par un calcul \'el\'ementaire: si par exemple
$f=P/Q$ est un endomorphisme de $\P^1$, on peut supposer que l'infini n'est pas un point fixe
--quitte \`a conjuguer par une transformation de Moebius-- pour une carte affine $\C$. Cela
revient \`a dire que le polyn\^ome $Q$ est de degr\'e $\l$. Les points p\'eriodiques
d'ordre $\leq n$ sont alors les solutions dans $\C$ de l'\'equation polynomiale
$f^n(z)=z$ qui est de degr\'e $\l^n+1$.
On montre (voir [CG]) que tous les points p\'eriodiques --sauf un nombre fini
d'entre eux-- sont de type r\'epulsif. Consid\'erons
$$
\nu_n:=\frac{1}{\l^n} \sum_{f^nx=x,\; x \text{ r\'epulsif}} \e_x,
$$
et soit $\nu$ une valeur d'adh\'erence de $\nu$. Alors $\nu$ est une mesure de
probabilit\'e sur $X$ dont nous allons montrer qu'elle co\"{\i}ncide avec $\mu_f$.
Il suffit de montrer que $\mu_f(D) \leq \nu(D)$, pour tout petit disque $D=D(x,r)$.

Soit $D$ un tel disque. Comme $\mu_f$ ne charge pas l'ensemble postcritique
$\cup_{j \geq 0} f^j({\mathcal C}_f)$, on peut supposer que 
$D$ ne rencontre pas $\cup_{j=1}^l f^j({\mathcal C}_f)$ et appliquer le Lemme 1.12.
On note $D_i^{-n}$ les petits disques images de $D$ par les branches inverses de $f^n$.
Soit $D'$ un disque l\'eg\`erement plus petit que $D$.
Comme $\mu_f$ est m\'elangeante, on obtient pour $n$ assez grand,
$$
\mu_f(D)^2 \simeq \mu_f(D) \mu_f(D') \simeq \mu_f(f^{-n}D \cap D') \simeq
\sum_{i=1}^{(1-\e)\l^n} \mu_f(D_i^{-n} \cap D').
$$
Or soit $D_i^{-n} \cap D'=\emptyset$, 
soit $D_i^{-n}$ est relativement compact dans $D$
car $D_i^{-n}$ a un petit diam\`etre. Dans ce dernier cas
$f_i^{-n}$ est une contraction de $D$ dans lui m\^eme et produit donc
un point p\'eriodique r\'epulsif d'ordre $\leq n$. On obtient ainsi
$$
\mu_f(D_i^{-n} \cap D') \leq \mu_f(D_i^{-n}) \cdot \l^n \nu_n(D_i^{-n})
=\mu_f(D) \nu_n(D_i^{-n})
$$
car $\mu_f$ est de jacobien constant et $f^n$ est injective sur $D_i^{-n}$.
Comme les $D_i^{-n}$ que nous consid\'erons sont deux \`a deux disjoints et tous 
inclus dans $D$, il vient
$$
\mu_f(D)^2 \lesssim \mu_f(D) \nu_n(D),
$$
d'o\`u le r\'esultat.
\vskip.1cm

\noindent {\bf Unicit\'e}.
Soit $\nu$ une mesure invariante ergodique d'entropie positive. Observons que $\nu$ ne charge pas
les points, en particulier $\nu({\mathcal E}_f)=0$.
Il r\'esulte du Th\'eor\`eme 1.8 que soit $\nu=\mu_f$, soit $f^* \nu \neq \l \nu$.

Supposons que $\nu$ n'est pas de jacobien constant. On construit alors un ouvert $U$
plein (i.e. $\nu(U)=\text{Vol}(U)=1$),
simplement connexe de $X \setminus f({\mathcal C}_f)$ dont les pr\'eimages 
$U_1,\ldots,U_{\l}=f^{-1}(U)$ v\'erifient
$\nu(U_1)>\cdots>\nu(U_{\l})$.

On peut alors coder la dynamique de $f$ sur son graphe it\'er\'e gr\^ace \`a la
partition $U_1,\ldots, U_{\l}$ et montrer que
$h_{\nu}(f)<\log \l$.
Nous renvoyons le lecteur \`a [Ly], [BrD 2] pour une preuve d\'etaill\'ee de ce dernier point
qui assure que $\mu_f$ est l'unique mesure d'entropie maximale.
\vskip.2cm

\begin{preuve.1.12}
On note $\tau$ le nombre de valeurs critiques de $f$ et $V_l:=\cup_{j=1}^l f^j({\mathcal C_f})$.
Soit $D$ un disque qui \'evite $V_l$. L'entier $l$ sera sp\'ecifi\'e plus loin
en fonction de $\e>0$ qui est fix\'e par l'\'enonc\'e du lemme.

On construit des branches inverses $f_i^{-n}$ de $f^n$ par r\'ecurrence: comme
$D$ \'evite les valeurs critiques $V_l$ de $f^l$, 
$f^l$ est un rev\^etement \'etale de degr\'e $\l^l$ en restriction \`a $D$
et on peut donc construire des branches inverses $f_1^{-1},\ldots,f_{\l^l}^{-l}$ 
de $f^l$ sur $D$ d'images disjointes $D_i^{-l}$.
Si $D_i^{-l}$ ne rencontre pas les valeurs critiques de $f$, on peut fabriquer
$\l$ branches inverses de $f$ sur $D_i^{-l}$ qui d\'efinissent $\l$
branches inverses de $f^{l+1}$ sur $D$. Or il y a au plus $\tau$ disques
$D_i^{-l}$ qui rencontrent $V_1$. Si on jette ces $\tau$ disques (\'eventuels),
on obtient ainsi $\l^{l+1}(1-\tau/\l^{l+1})$ branches inverses de $f^{l+1}$
sur $D$. Une r\'ecurrence montre alors que pour tout $n \geq l+1$, on peut construire
au moins
$$
\l^n \left[1-\frac{\tau}{\l^{l+1}}-\cdots-\frac{\tau}{\l^n}\right] \geq \l^n[1-\e/2]
$$
branches inverses de $f^n$ sur $D$ si l'on fixe $l=l_{\e}$ assez
grand pour que $\tau \l^{-l} /(1-1/\l) <\e/2$.

On veut \`a pr\'esent contr\^oler le diam\`etre des images $D_i^{-n}$. Observons que
$$
\sum_{i=1}^{(1-\e/2)\l^n} \text{Aire}(D_i^{-n}) \leq \text{Aire}(X)=1,
$$
si on calcule l'aire par rapport \`a une forme volume normalis\'ee. Il s'ensuit
qu'au moins $(1-\e)\l^n$ des disques $D_i^{-n}$ ont une aire qui v\'erifie
$$
\text{Aire}(D_i^{-n}) \leq \frac{2}{\e} \l^{-n}.
$$
On en d\'eduit un contr\^ole des diam\`etres: quitte \`a r\'eduire l\'eg\`erement $D$, on peut
supposer qu'on contr\^ole en fait l'aire de disques l\'eg\`erement plus grands
$\tilde{D}_i^{-n}=f_i^{-n}(\tilde{D})$, o\`u $D$ est relativement compact dans $\tilde{D}$.
Le module des anneaux $\tilde{D}_i^{-n} \setminus {D}_i^{-n}$ \'etant constant,
\'egal \`a $mod(\tilde{D} \setminus D)$, on obtient
$$
\text{diam}(D_i^{-n}) \lesssim [\text{Aire}(\tilde{D}_i^{-n})]^{1/2} \lesssim \l^{-n/2}.
$$
Nous renvoyons le lecteur \`a l'Appendice de [BrD 2] pour plus de d\'etails sur cette
estimation aire-diam\`etre.
\end{preuve.1.12}
\vskip.2cm

Notons qu'on peut se passer de cette estimation aire-diam\`etre en dimension 1: il r\'esulte
en effet du Th\'eor\`eme de Koebe que le contr\^ole de l'aire des $D_i^{-n}$ suffit
\`a contr\^oler les d\'eriv\'ees $(f^{-n})'(x)$ et conduit ainsi \`a la
minoration de l'exposant de Lyapunov. Pour le m\'elange, on a juste besoin de savoir que
le diam\`etre des $D_i^{-n}$  converge vers $0$, ce qui se voit facilement
par un argument de familles normales. Nous avons n\'eanmoins indiqu\'e cette estimation
due \`a J.-Y.Briend et J.Duval [BrD 2] (et qui ne figure pas dans [Ly]), car c'est un argument
crucial en dimension sup\'erieure o\`u l'on ne dispose pas du Th\'eor\`eme de Koebe.
Comme on l'aura compris, toutes les preuves propos\'ees dans cette premi\`ere partie
sont model\'ees pour pouvoir s'adapter directement en dimension sup\'erieure.

\section{Le cas des polyn\^omes}

Nous interpr\'etons ici les objets introduits jusqu'\`a pr\'esent lorsque $f$
est un {\it polyn\^ome} de degr\'e $\l \geq 2$. Quitte \`a conjuguer $f$ par
un automorphisme $z \mapsto az+b$ convenablement choisi, on peut toujours supposer
-- et nous le ferons -- que $f$ est {\it centr\'e}, {\it unitaire}, i.e.
$$
f(z)=z^{\l}+a_{\l-2}z^{\l-2}+\cdots+a_0,
$$
o\`u $(a_0,\ldots,a_{\l-2}) \in \C^{\l-1}$. 

Observons que le point $\infty$ \`a l'infini joue un r\^ole sp\'ecial ici: c'est un 
point fixe superattractif ($f'(\infty)=0$) qui est totalement invariant:
$\infty \in {\mathcal E}_f$.
Il s'ensuit que si $z \in \C$, soit l'orbite $(f^n z)_{n \in \N}$ de $z$ est une suite
born\'ee, soit elle converge vers $\infty$. Notons 
${\mathcal B}(\infty)$ le bassin d'attraction (dans $\C$) de l'infini et
$$
K_f:=\C \setminus {\mathcal B}(\infty)=\{ z \in \C \, / \, (f^n z)_{n \in \N}
\text{ est born\'ee } \}.
$$
L'ensemble $K_f$ est appel\'e ensemble de Julia rempli. Cette terminologie
est justifi\'ee par l'observation suivante.

\begin{lem}
L'ensemble $K_f$ est un compact simplement connexe dont le bord
co\"{\i}ncide avec l'ensemble de Julia $J_f$.
\end{lem}

\begin{preuve}
Rappelons que $\mu_f=\om+\Delta g_f$, o\`u $g_f=\sum \l^{-j} \g \circ f^j$
et $\g$ est solution de l'\'equation $\l^{-1} f^*\om=\om+\Delta \g$.
La forme de Fubini-Study $\om$ est donn\'ee par $\Delta \frac{1}{2} \log [1+|z|^2]$
dans $\C$. Posons
$$
G_f(z)=\frac{1}{2} \log [1+|z|^2]+g_f(z).
$$
C'est une fonction sousharmonique et continue dans $\C$ telle que
$\Delta G_f=\mu_f$. Si l'on choisit -- comme nous l'avons implicitement fait 
dans le paragraphe 1.1.1 -- $\g(z)=\frac{1}{2\l} \log [1+|f(z)|^2]-\frac{1}{2} \log [1+|z|^2]$,
on obtient
$$
G_f(z)=\lim_{n \rightarrow +\infty} \frac{1}{2\l^n} \log \left[ 1+|f^n(z)|^2 \right].
$$
En particulier $G_f \geq 0$,
$$
(G_f=0)=K_f \; \text{ et } \; (G_f>0)={\mathcal B}(\infty).
$$
Il s'ensuit que $K_f$ est un compact et $J_f=\text{Supp} \mu_f \subset \partial K_f$
(par la Proposition 1.1, car $G_f \equiv 0$ dans l'int\'erieur de $K_f$ et $G_f=\lim \l^{-n} \log |f^n|$
est harmonique dans ${\mathcal B}(\infty)$).
Or $G_f$ ne peut pas \^etre harmonique en un point $p \in \partial K_f$ par le principe
du minimum, donc $J_f=\partial K_f$.

Enfin il r\'esulte du principe du maximum que ${\mathcal B}(\infty)$ est connexe, donc
$K_f$ est simplement connexe.
\end{preuve}
\vskip.2cm

L'ensemble $K_f$ est donc constitu\'e de $J_f$ et de ses trous (composantes
connexes born\'ees de $\C \setminus J_f$).
La fonction $G_f(z)=g_f(z)+\frac{1}{2}\log[1+|z|^2]$, introduite par H.Brolin [Bro],
poss\`ede de nombreuses propri\'et\'es. Son invariance $G_f \circ f=\l G_f$
refl\`ete l'invariance de la mesure $\mu_f$, $f^* \mu_f=\l \mu_f$, qu'on
appelle parfois mesure de Green en raison de l'interpr\'etation suivante.

\begin{pro}
La fonction $G_f$ est la fonction de Green du compact $K_f$ avec p\^ole \`a l'infini.
Elle admet le d\'eveloppement asymptotique 
$$
G_f(z)=\log|z|+o(1).
$$
En particulier $K_f$ est de capacit\'e logarithmique \'egale \`a $1$.
\end{pro}

La preuve est une cons\'equence imm\'ediate de ce que $G_f \geq 0$ est harmonique dans 
$\C \setminus K_f={\mathcal B}(\infty)$ et nulle sur $K_f$.
Le d\'eveloppement asymptotique provient de ce que nous avons pris $f$ {\it unitaire}.
Il assure \'egalement que $G_f$ co\"{\i}ncide avec le potentiel logarithmique de $\mu_f$,
$$
V_{\mu_f}(z)=\int_{\C} \log|z-w| d\mu_f(w).
$$
En effet $\Delta G_f=\Delta V_{\mu_f}=\mu_f$, donc $G_f-V_{\mu_f}$ est une fonction
harmonique \`a croissance logarithmique dans $\C$: c'est donc une constante
qui vaut z\'ero puisque ces deux fonctions ont le m\^eme
d\'eveloppement asymptotique, $G_f(z)=\log |z|+o(1)=V_{\mu_f}(z)$.

\begin{exas}\text{ }

1) Si $f(z)=z^{\l}$ alors $G_f(z)=\log^+ |z|:=\max(\log|z|,1)$.
L'ensemble $K_f$ est le disque unit\'e ferm\'e et $\mu_f$ est la mesure de Lebesgue
normalis\'ee sur le cercle unit\'e $S^1$.

2) Le polyn\^ome de Tchebychev $g(z)=z^2-2$ est semi-conjugu\'e au polyn\^ome $f(z)=z^2$
par l'application $\Phi(z)=z+1/z$,  via $g \circ \Phi=\Phi \circ f$. On en d\'eduit
que $K_g=J_g$ est l'intervalle $[-2,2]$ et
$$
G_g(z)=\max \left( \log \left| \frac{z+\sqrt{z^2-4}}{2} \right|;
\log \left| \frac{z-\sqrt{z^2-4}}{2} \right| \right).
$$
Notons en particulier que $G_g$ est H\"old\'erienne d'exposant $1/2$
(mais pas plus) et $\chi_{top}(g)=\log 4$ puisque
$g'(2)=4$, $g(2)=2$ et $\sup_{J_g} |g'|=4$. La r\'egularit\'e obtenue \`a 
la Proposition 1.2 est donc optimale.
\end{exas}

L'exposant de Lyapunov $\chi(\mu_f)$ d'un polyn\^ome $f$ par rapport
\`a la mesure d'entropie maximale $\mu_f$ s'exprime simplement en fonction
de la fonction $G_f$, comme l'ont observ\'e A.Manning [Ma]
et F.Przytycki [P].

\begin{pro}
Soit $f$ un polyn\^ome unitaire de degr\'e $\l \geq 2$. Alors
$$
\chi(\mu_f)=\log \l+\sum_{f'(c)=0} G_f(c).
$$
\end{pro}

\begin{preuve}
Il r\'esulte du Th\'eor\`eme ergodique de Birkhoff que
$$
\chi(\mu_f)=\int_{\C} \log|f'| d\mu_f.
$$
Or $f'(z)=\l \Pi_{i=1}^{\l-1} (z-c_i)$, o\`u les $c_i$ sont les points critiques
de $f$, i.e. les points $c \in \C$ tels que $f'(c)=0$. 
Comme $G_f$ co\"{\i}ncide avec le potentiel logarithmique de $\mu_f$
(voir la discussion qui pr\'ec\`ede les Exemples 1.15), il vient
$$
\chi(\mu_f)=\log \l+\sum_{f'(c)=0} \int \log|z-c|d\mu_f(z)=\log \l+\sum_{f'(c)=0} G_f(c).
$$
\end{preuve}

Lorsque l'ensemble de Julia est connexe (i.e. lorsque tous les points critiques sont dans
$K_f$, voir [CG]), on obtient $\chi(\mu_f)=\log \l$.
Il y a donc en g\'en\'eral un d\'ecalage entre l'exposant de Lyapunov mesurable
$\chi(\mu_f)$ et l'exposant de Lyapunov topologique
$\chi_{top}(f)$, qui oscille dans ce cas entre $\log \l$ et $2 \log \l$
(voir [Buf 1]).

\chapter{Invariants num\'eriques}

Dans toute la suite $X$ d\'esigne une vari\'et\'e k\"ahl\'erienne compacte
de dimension $k$ munie d'une forme de K\"ahler $\om$, 
$f:X \rightarrow X$ d\'esigne un endomorphisme m\'eromorphe dominant,
et on note $I_f$ l'ensemble d'ind\'etermination de $f$: c'est un sous-ensemble 
analytique de $X$ de codimension au moins deux, constitu\'e des points en lesquels
$f$ n'est pas continu.

\section{Degr\'es dynamiques}

Etant donn\'ee une forme diff\'erentielle lisse $\theta$ sur $X$, on d\'efinit
son image inverse $f^*\theta$ par $f$ de la fa\c{c}on suivante: soit
$\Gamma_f \subset X \times X$ le graphe de $f$, et soit $\tilde{\Gamma}_f$
une d\'esingularisation de $\Gamma_f$. On a un diagramme commutatif
$$
\begin{array}{ccccc}
\text{ } & \text{ } & \tilde{\G}_f & \text{ } & \text{ } \\
\text{ } & \stackrel{\pi_1}{\swarrow} & \text{ } &
\stackrel{\pi_2}{\searrow} & \text{ } \\
X & \text{ } & \stackrel{f}{\longrightarrow} & \text{ } & X
\end{array}
$$
o\`u $\pi_1,\pi_2$ sont des applications holomorphes. On pose
$f^* \theta:=(\pi_1)_* (\pi_2^* \theta)$, o\`u l'on pousse la forme lisse $\pi_2^* \theta$
par $\pi_1$ au sens des courants. De fa\c{c}on analogue, on
d\'efinit l'image directe de $\theta$ par
$$
f_* \theta:=(\pi_2)_* (\pi_1^* \theta).
$$
On v\'erifie ais\'ement que ces d\'efinitions ne d\'ependent pas du choix de la 
d\'esingularisation de $\G_f$. Observons que les op\'erations
$f^*,f_*$ pr\'eservent les degr\'es,
les bidegr\'es, le caract\`ere r\'eel et les bords; 
elles induisent donc des actions lin\'eaires sur les espaces de cohomologie
$H^{l,l}(X,\R)$.

\begin{defi} Pour $0 \leq l \leq k$,
on note $r_l(f)$ le rayon spectral de l'action lin\'eaire
$f^*:H^{l,l}(X,\R) \rightarrow H^{l,l}(X,\R)$ et on pose
$$
\rho_l(f):=\liminf_{n \rightarrow +\infty} \left[r_l(f^n)\right]^{1/n}.
$$
C'est le rayon spectral dynamique d'ordre $l$ de $f$.
\end{defi}

Notons que $r_l(f^n) \neq r_l(f)^n$ en g\'en\'eral
--sauf bien s\^ur pour les endomorphismes holomorphes.
C'est cette observation qui a pouss\'e S.Friedland \`a introduire les degr\'es 
$\rho_l(f)$ dans [Fr]: ceux-ci v\'erifient
$\rho_l(f^n)=[\rho_l(f)]^n$, pour tout $n \in \N$.

On peut bien s\^ur introduire des quantit\'es analogues --notons
les $r_l(f_*)$, $\rho_l(f_*)$-- correspondant aux actions lin\'eaires
induites par $f_*$.
Il r\'esulte de la d\'efinition que ces actions $f^*,f_*$ sont adjointes l'une de 
l'autre pour la forme d'intersection sur $X$, i.e.
$$
\langle f^*\theta,\eta\rangle =\langle \pi_2^* \theta,\pi_1^* \eta\rangle =\langle \theta,f_* \eta\rangle ,
$$
pour toutes formes lisses ferm\'ees $\theta,\eta$ de degr\'es compl\'ementaires.
On en d\'eduit que $r_l(f_*)=r_{k-l}(f)$ et
$\rho_l(f_*)=\rho_{k-l}(f)$.

Observons que $H^{0,0}(X,\R) \simeq \R$ (fonctions constantes) et
$H^{k,k}(X,\R) \simeq \R$. L'op\'erateur $f^*$ agit sur le premier par multiplication
par la constante $1=r_0(f)=\rho_0(f)$, et sur le deuxi\`eme par multiplication
par le degr\'e topologique (nombre de pr\'eimages par $f$ d'un point g\'en\'erique).
On a donc toujours $r_k(f)=\rho_k(f)$ et
$\rho_k(f)$ est le {\bf degr\'e topologique de } $f$.

\begin{exa}
Lorsque $X=\P^k$, tous les espaces de cohomologie $H^{l,l}(\P^k,\R)$
sont de dimension $1$ et $f^*$ agit par multiplication par la constante
$r_l(f)$ qui admet dans ce cas une interpr\'etation alg\'ebro-g\'eom\'etrique,
$$
r_l(f)=\deg (f^{-1}L), \; L
\text{ un sous-espace lin\'eaire g\'en\'erique de codimension } l,
$$
comme l'ont montr\'e A.Russakovskii et B.Shiffman [RS].
La constante $r_1(f)$ est dans ce cas le degr\'e des polyn\^omes homog\`enes
premiers entre eux 
intervenant dans l'\'ecriture de $f$ en coordonn\'ees homog\`enes.
\end{exa}

Il y a d'autres actions lin\'eaires auxquelles on peut s'int\'eresser, par exemple
l'action lin\'eaire induite par $f^*,f_*$ sur les espaces
$H^{p,q}(X,\R)$, $p \neq q$. C'est l\'egitime si l'on souhaite estimer le nombre de
points p\'eriodiques de $f$ en utilisant la formule de Lefschetz (voir section 2.3.1).
Celle-ci fait intervenir la trace des op\'erateurs $f^*$ sur tous les espaces de deRham
$H^i(X,\R)$. Puisque $X$ est k\"ahl\'erienne, ils se d\'ecomposent
canoniquement (Th\'eor\`eme de Hodge),
$$
H^i(X,\R)=\bigoplus_{p+q=i} H^{p,q}(X,\R).
$$

Lorsque $X$ est {\it projective}, on peut au contraire se restreindre \`a la partie
alg\'ebrique des espaces en ne consid\'erant que le sous-espace 
$H_{alg}^{l,l}(X,\R)$ (stable par $f^*,f_*$)
engendr\'e par les cycles alg\'ebriques,
On notera $r_l^{alg}(f)$ le rayon spectral de la restriction de $f^*$ \`a ce sous-espace
invariant, et $\rho_l^{alg}(f)$ le degr\'e dynamique correspondant.

Il est tr\`es utile d'avoir une description {\it analytique} du comportement asymptotique
des op\'erateurs $(f^n)^*,(f^n)_*$. Nous avons introduit
dans [G 5] les ``degr\'es'' $\d_l(f,\om):=\int_{X \setminus I_f} f^* \om^l \wedge \om^{k-l}$,
et leur version dynamique:

\begin{defi}
Le degr\'e dynamique d'ordre $l$ est 
$$
\l_l(f)=\liminf_{n \rightarrow +\infty} \left[
\int_{X \setminus I_{f^n}} (f^n)^* \om^l \wedge \om^{k-l} \right]^{1/n}.
$$
\end{defi}

Il est clair que $\l_l(f)$ ne d\'epend pas du choix de 
la forme de K\"ahler $\om$.
Le lien entre ces diff\'erentes notions est fourni par le
r\'esultat suivant.

\begin{thm}  
Pour tout $0 \leq l \leq k$, on a $\l_l(f)=\rho_l(f)$. 
En particulier $\rho_l(f)=\rho_l^{alg}(f)$ lorsque $X$ est projective. De plus,

(a) Pour tout $ 1 \leq l \leq k$, on  a
$$
1 \leq \l_{l+1}(f) \l_{l-1}(f) \leq [ \l_l(f) ]^2.
$$ 

(b) Si $g:X' \rightarrow X'$ est un endomorphisme m\'eromorphe dominant
d'une vari\'et\'e k\" ahl\'erienne compacte $X'$ de dimension $k'$, alors
le produit direct $f \times g:X \times X' \rightarrow X \times X'$ v\'erifie
$$
\l_l(f \times g) =\max_{i+j=l} \l_i(f) \l_j(g),
$$
pour tout $l \in [0,k+k']$. 

(c) Il existe $C>0$ telle que pour tout $n \in \N$ et $p,q \in [0,k]$, on a 
$$
\left| Tr\left( (f^n)^*_{|H^{p,q}(X,\R)} \right) \right| \leq C 
\max_{0 \leq j \leq k} \d_j(f^n,\om).
$$
\end{thm}

Ces assertions sont pour la plupart d\'emontr\'ees dans [G 5]. Nous en donnons ici une preuve
l\'eg\`erement simplifi\'ee.

\begin{preuve}
On va montrer qu'il existe $C \geq 1$ telle que pour tout $n \in \N$,
\begin{equation}
\frac{1}{C} ||(f^n)^*|| \leq \d_l(f^n,\om) \leq C ||(f^n)^*||
\end{equation}
pour une norme $||\cdot||$ arbitraire sur $H^{l,l}(X,\R)$. L'\'egalit\'e
$\rho_l(f)=\l_l(f)$ en r\'esulte imm\'ediatement.
Soit $(\a_i)$ une base de $H^{l,l}(X,\R)$. On peut la choisir de sorte que les 
classes $\a_i$ soient repr\'esent\'ees par des $(l,l)$-formes lisses ferm\'ees
{\it positives}, et avec $\a_1=\{\om^l\}$. Fixons $C_1>0$ une constante
telle que la classe $C_1 \{\om^l\}$ domine toutes les classes $\a_i$,
c'est \`a dire que la diff\'erence peut \^etre repr\'esent\'ee par la classe
d'un courant positif. 
 
Soit \`a pr\'esent $\b_j \in H^{k-l,k-l}(X,\R)$ 
les \'el\'ements de la base duale de $(\a_j)$ pour la forme d'intersection
(dualit\'e de Serre). Il vient, pour tout $n \in \N$,
$$
\langle (f^n)^* \a_i,\b_j\rangle =a_{ij}^{(n)},
$$
o\`u les $a_{ij}^{(n)}$ sont les coefficients de la matrice repr\'esentant
$(f^n)^*$ dans la base $(\a_i)$. Soit $C_2>0$ une constante telle que
la classe $C_2\{\om^{k-l}\}$ domine toutes les classes $\pm \b_j$.
Il vient
$$
\langle (f^n)^*\a_i,\pm \b_j\rangle  \leq C_2 \langle (f^n)^* \a_i, \om^{k-l}\rangle  
\leq C_1C_2 \d_l(f^n,\om),
$$
d'o\`u l'in\'egalit\'e $||(f^n)^*|| \leq C \d_l(f^n,\om)$. Pour l'in\'egalit\'e 
r\'eciproque, on observe, puisque $\a_1=\{\om^l\}$, que
$$
\d_l(f^n,\om) =\sum_j a_{1j}^{(n)} \langle \a_j,\om^{k-l}\rangle  \leq 
\sum_j |a_{1j}^{(n)} | \cdot |\langle \a_j,\om^{k-l}\rangle | \leq C ||(f^n)^*||.
$$

Lorsque $X$ est projective, on peut choisir une forme $\om$ dont la classe
est enti\`ere (classe de Hodge). Le raisonnement
pr\'ec\'edent s'applique au sous-espace stable engendr\'e par les cycles alg\'ebriques
et montre que $\l_l(f)=\rho_l^{alg}(f)$, donc 
$\rho_l(f)=\rho_l^{alg}(f)$.

Les in\'egalit\'es (a) sont des cons\'equences des in\'egalit\'es 
mixtes de Teissier-Hovanskii qui assurent
$\d_{l+1}(f,\om) \d_{l-1}(f,\om) \leq [\d_l(f,\om)]^2$ (voir [G 5]).
Ces in\'egalit\'es peuvent s'interpr\'eter comme des propri\'et\'es de
concavit\'e de la fonction $l \mapsto \log \d_l(f,\om)$
(voir [Gr 2]). Comme $\l_k(f) \geq 1$, on en d\'eduit que $\l_j(f) \geq 1$ 
pour tout $j \in [1,k]$.

Pour \'etablir (b), on fixe une forme de K\"ahler $\Omega(x,x')=\om(x)+\om'(x')$
sur le produit $X \times X'$. Un calcul \'el\'ementaire donne
$$
\d_l(h^n,\Omega)=\sum_{i+j=l} \d_i(f^n,\om) \d_j(g^n,\om'),
$$
o\`u $h$ d\'esigne l'endomorphisme produit $(f,g)$. L'\'egalit\'e s'ensuit.

Observons que (c) est une cons\'equence de (2.1) lorsque $p=q$.
Lorsque $p \neq q$, on peut s'y ramener en travaillant dans $X^2$
avec l'endomorphisme produit $h=(f,f)$.
En effet si $\a \in H^{p,q}(X,\R)$ alors $\a(x) \wedge \overline{\a(x')}$
d\'efinit une classe de bidegr\'e $(p+q,p+q)$ dans $X^2$. Estimer la norme
de $(h^n)^*$ sur une telle classe revient \`a estimer le carr\'e de la norme
de $(f^n)^*$ agissant sur $H^{p,q}(X,\R)$. Le calcul fait pour d\'emontrer
(b) donne ainsi
$$
||(f^n)^*_{H^{p,q}(X,\R)}||^2 \leq C_1 \sum_{i+j=p+q} \d_i(f^n,\om) \d_j(f^n,\om)
\leq \left[ C \max_j \d_j(f^n,\om) \right]^2.
$$
\end{preuve}
\vskip.3cm

Notons qu'on peut raffiner l'estimation 2.4.c et montrer
$$
||(f^n)^*_{H^{p,q}(X,\R)}||^2 \leq C \d_p(f^n,\om) \d_q(f^n,\om).
$$
Ce contr\^ole plus pr\'ecis est d\^u \`a  T.C.Dinh (Proposition 5.8 dans [Di 1]).

\vskip.2cm
On note $H^{l,l}_{psef}(X,\R)$ le c\^one des
classes {\it pseudoeffectives} (i.e. repr\'esentables par un courant positif ferm\'e),
et $H^{l,l}_{nef}(X,\R)$ le c\^one des classes {\it nef}, adh\'erence du c\^one
des classes strictement positives (i.e. repr\'esentables par une forme lisse
ferm\'ee qui domine un multiple de $\om^l$): c'est le c\^one dual
de $H^{k-l,k-l}_{psef}(X,\R)$. Ce sont des c\^ones stricts, i.e. ils ne 
contiennent aucune droite r\'eelle.

Lorsque l'un de ces c\^ones est invariant par $f^*$,
il r\'esulte du Th\'eor\`eme de Perron-Frobenius qu'il contient une classe
qui est un vecteur propre associ\'e \`a la plus grande valeur propre
de $f^*:H^{l,l}(X,\R) \rightarrow H^{l,l}(X,\R)$. En particulier
$r_l(f)$ est une valeur propre de $f^*$. Cela a des cons\'equences
particuli\`erement int\'eressantes sur le spectre 
et la dynamique de $f^*$, comme nous le verrons dans la section 4.2.

Le c\^one $H^{1,1}_{psef}(X,\R)$ est toujours pr\'eserv\'e par $f^*$
et $f_*$. Cela r\'esulte de ce qu'on a une bonne d\'efinition de
l'image inverse par $f$
(et de l'image directe) de n'importe quel courant positif ferm\'e de bidegr\'e
$(1,1)$: tout se passe comme si $f$ \'etait holomorphe, car l'ensemble d'ind\'etermination
est de codimension $\geq 2$: c'est donc un ensemble n\'egligeable pour les
courants de bidegr\'e $(1,1)$.
Nous verrons dans la section 4.2 que le c\^one $H^{1,1}_{nef}(X,\R)$ est 
\'egalement pr\'eserv\'e
par $f^*,f_*$ lorsque $X$ est de dimension $2$. Ce n'est pas vrai en dimension
sup\'erieure pour des applications m\'eromorphes, comme le montre l'exemple
suivant qui nous a \'et\'e communiqu\'e par L.Bonavero.

\begin{exa}
Soit $g: \P^3 \rightarrow \P^3$ l'involution donn\'ee en coordonn\'ees homog\`enes
par
$
g[z_0:z_1:z_2:z_3]=[1/z_0:1/z_1:1/z_2:1/z_3]=[z_1z_2z_3:\cdots:z_0z_1z_2].
$
C'est un endomorphisme birationnel (i.e. de degr\'e topologique $r_3(f)=\l_3(f)=1$)
de $\P^3$. 
Son ensemble d'ind\'etermination est
$$
I_g:=\bigcup_{i \neq j} \{ z_i=z_j=0 \}.
$$
Soit $\om$ la forme de Fubini-Study sur $\P^3$, on v\'erifie ais\'ement que 
$$
r_1(f)=\d_1(f,\om)=3 \text{ et } r_1(f^2)=\d_1(f^2)=1,
\text{ donc } \l_1(f)=1.
$$
Comme $f=f^{-1}$, il vient \'egalement $\l_2(f)=\l_1(f^{-1})=1$.

Soit $\pi:X \rightarrow \P^3$ l'\'eclatement de $\P^3$ aux quatre points toriques
$p_0=[1:0:0:0],\, p_1=[0:1:0:0], \, p_2=[0:0:1:0], \, p_3=[0:0:0:1]$.
On note $E_i=\pi^{-1}(p_i)$ les diviseurs exceptionnels. Soit
$f:X \rightarrow X$ l'endomorphisme m\'eromorphe induit par $g$.
On note 
$H_i=(z_i=0)$ l'hyperplan de $\P^3$ passant par les trois points $p_j$, $j \neq i$,
et $H_i'$ sa transform\'ee stricte par $\pi$, 
$$
H_i'=\pi^*H_i-\sum_{j \neq i} E_j.
$$
Comme $g(p_i)=H_i$, on obtient
$$
f^*H_i'=E_i \; \text{ et } \; f^*E_i=H_i'.
$$
Soit $D=\pi^*H_0$. C'est un diviseur nef (il est m\^eme semi-positif), et 
pourtant
$$
f^*D=E_0+H_1'+H_2'+H_3'
$$
n'est pas nef. Soit en effet $L$ la droite passant par les points $p_3,p_4$
et $L'$ sa transform\'ee stricte. On calcule 
$$
f^*D \cdot L'=(3 \pi^*H-3E_0-2\sum_{j \geq 1} E_j) \cdot L'=3-2-2=-1.
$$
Le c\^one $H^{1,1}_{nef}(X,\R)$ n'est donc pas pr\'eserv\'e par $f^*$.
Comme c'est le c\^one dual de $H^{2,2}_{psef}(X,\R)$, cela
montre \'egalement que ce dernier n'est pas pr\'eserv\'e par 
$f_*=(f^{-1})^*$.
\end{exa}

Nous verrons plus loin qu'il est utile d'effectuer des changements 
de coordonn\'ees bim\'eromorphes 
pour analyser la dynamique. Il faut donc s'int\'eresser \`a 
l'invariance des degr\'es dynamiques sous conjugaison bim\'eromorphe.

Observons tout d'abord que l'on peut d\'efinir des degr\'es
$\d_l(f,\om_X,\om_Y)$ lorsque $f:X \rightarrow Y$ est une application
m\'eromorphe entre deux vari\'et\'ess k\"ahl\'eriennes compactes de
m\^eme dimension $k$ par
$$
\d_l(f,\om_X,\om_Y):=\int_{X \setminus I_f} f^* \om_Y^l \wedge \om_X^{k-l}.
$$
Ici $\om_X, \om_Y$ d\'esignent deux formes de K\"ahler sur $X,Y$ respectivement.
On a alors l'estimation suivante due \`a  
T.C.Dinh et N.Sibony [DS 3,4].

\begin{pro}
Soit $X,Y,Z$ trois vari\'et\'es k\"ahl\'eriennes compactes de m\^eme dimension $k$,
munies de formes de K\"ahler $\om_X,\om_Y,\om_Z$. Il existe $C>0$ telle que
pour tout endomorphisme $f:X \rightarrow Y$, $g:Y \rightarrow Z$ et pour tout 
$0 \leq l \leq k$,
$$
\d_l(g \circ f,\om_X,\om_Z) \leq C \d_l(g,\om_Y,\om_Z) \d_l(f,\om_X,\om_Y).
$$
\end{pro}

\begin{esquisse}
Lorsque $l=1$, la d\'emonstration est simplifi\'ee gr\^ace \`a l'observation suivante:
$f^* (g^* \om_Z)$ est un $(1,1)$-courant positif ferm\'e bien d\'efini qui co\"{\i}ncide
avec $(g \circ f)^* \om_Z$ hors d'une hypersurface. Or $(g \circ f)^* \om_Z$ ne charge pas
les hypersurfaces (c'est localement une forme \`a coefficients $L_{loc}^1$), on a 
donc $(g \circ f)^* \om_Z \leq f^* (g^* \om_Z)$. Il s'ensuit que
\begin{eqnarray*}
\d_1(g \circ f,\om_X,\om_Z) & \leq &\int_X f^* (g^* \om_Z) \wedge \om_X^{k-1}
=\int_Y g^* \om_Z \wedge f_*(\om_X^{k-1}) \\
& \leq &  C ||g^*_{H^{1,1}}|| \cdot ||f^*_{H^{1,1}}||
\leq C' \d_1(g,\om_Y,\om_Z) \d_1(f,\om_X,\om_Y),
\end{eqnarray*}
en utilisant les in\'egalit\'es (2.1).

Lorsque $l \geq 2$, la d\'efinition de $f^* (g^* \om_Z^l)$ pose probl\`eme
et l'argument pr\'ec\'edent doit \^etre modifi\'e.
Observons que $\d_k$ est le degr\'e topologique qui se comporte tr\`es bien
par composition. Cela r\`egle le cas de la dimension deux. Lorsque
$X$ est de dimension 3, on peut utiliser la dualit\'e $f^*,f_*$
et travailler en bidegr\'e $(1,1)$ avec $f_*$ (voir paragraphe 1 dans [G 5]).
Pour traiter le cas g\'en\'eral ($2 \leq l \leq k-2$, 
donc $k=\dim_{\C} X \geq 4$), il faut savoir r\'egulariser 
convenablement le courant positif $g^* \om_Z$ comme l'ont observ\'e
A.Russakovskii et B.Shiffman [RS] qui d\'emontrent la Proposition 2.6
lorsque $X=\P^k$: dans ce cas on sait r\'egulariser sans perte
de positivit\'e car $Aut(\P^k)$ agit transitivement sur $\P^k$.
Lorsque la vari\'et\'e $X$ n'est pas homog\`ene, 
T.C.Dinh et N.Sibony ont montr\'e [DS 3,4] qu'on peut r\'egulariser avec
perte uniforme de positivit\'e: cela permet de conclure
(voir e.g. Lemme 4 dans [DS 3]).
\end{esquisse}
\vskip.2cm

Le point important est que la constante $C$ ne d\'epend ni de $f$, ni de $g$.
Si on applique cette estimation \`a $X=Y=Z$, $\om=\om_X=\om_Y=\om_Z$, lorsque
$f$ et $g$ sont les it\'er\'es d'un m\^eme endomorphisme (que nous notons
\`a nouveau $f$), on obtient pour tout $n,m \in \N$,
$$
\d_l(f^{n+m},\om) \leq C \d_l(f^n,\om) \d_l(f^m,\om).
$$
Autrement dit la suite $n \mapsto \d_l(f^n,\om)$ est quasi-sousmultiplicative:
la limite inf\'erieure dans la d\'efinition 2.3 est donc en fait une limite
(et un infimum). Voici une autre application de la Proposition 2.6.

\begin{cor}
Les degr\'es dynamiques sont invariants par conjugaison bi-m\'eromorphe.
\end{cor}

\begin{preuve}
On consid\`ere $\Phi:\tilde{X} \rightarrow X$ une application bim\'eromophe
et $\tilde{f}:\tilde{X} \rightarrow \tilde{X}$ l'endomorphisme induit
par $f, \Phi$ sur $\tilde{X}$. Fixons $\om,\tilde{\om}$ deux formes de 
K\"ahler sur $X,\tilde{X}$. Il r\'esulte d'une double application de la
Proposition 2.6 que
$$
\d_l(\tilde{f},\tilde{\om}) \leq C^2 \d_l(\Phi,\tilde{\om},\om) \d_l(\Phi^{-1},\om,\tilde{\om})
\d_l(f,\om).
$$
On peut remplacer $f,\tilde{f}$ par $f^n,\tilde{f}^n$ dans l'in\'egalit\'e pr\'ec\'edente
sans changer la valeur des constantes, d'o\`u $\l_l(\tilde{f}) \leq \l_l(f)$.
Il suffit enfin d'interchanger les r\^oles de $f,f^{-1}$ (resp. $\Phi,\Phi^{-1}$)
pour conclure.
\end{preuve}

\section{Entropies}

\subsection{Entropie topologique}

Soit $f:X \rightarrow X$ une application que l'on suppose dans un premier temps {\it holomorphe}.
Alors $f$ est en particulier un endomorphisme {\it continu} de $X$.
On peut d\'efinir son {\it entropie topologique} \`a la Bowen 
(voir [Bo] et [Dina]) de la fa\c{c}on suivante,
$$
h_{\text{top}}^{\text{Bow}}(f):=\sup_{\e >0} \limsup_{N \rightarrow +\infty} 
\frac{1}{N} \log \max \left\{ \sharp F
    \,
/ \, F \; \text{ ensemble } (N,\e)\text{-s\'epar\'e} \right\}.
$$
Rappelons qu'un ensemble $F$ est dit $(N,\e)$-s\'epar\'e si 
$d_N(x,y) \geq \e$, pour tout
couple de points distincts (x,y) de $F$, o\`u 
$$
d_N(x,y):=\max_{0 \leq j \leq N-1} d(f^j(x),f^j(y)),
$$
pour une distance $d$ fix\'ee sur $X$. La d\'efinition ne d\'epend du choix de la distance que dans une 
moindre mesure: deux distances \'equivalentes conduisent \`a la m\^eme valeur de l'entropie.
Dans notre contexte, nous choisissons la distance associ\'ee \`a une m\'etrique 
k\"ahl\'erienne sur $X$; la valeur de l'entropie ne d\'epend donc pas 
du choix particulier de cette m\'etrique.

Une approche alternative --inspir\'ee des travaux de M.Gromov [Gr 1]-- est de
consid\'erer 
l'entropie du graphe it\'er\'e de $f$,
$$
\Gamma_f^{\infty}:=\left\{ \hat{x}=(x_n)_{n \in \N} \in X^{\N} \, / \,
x_{n+1}=f(x_n) \text{ pour tout } n \in \N \right\},
$$
sur lequel $f$ agit comme un d\'ecalage unilat\'eral $\hat{f}$. On obtient ainsi
un endomorphisme continu $\hat{f}$ de l'espace topologique 
$\Gamma_f^{\infty}$ qui est compact (pour la topologie produit), et on pose
$$
h_{\text{top}}^{\text{Gr}}(f):=h_{\text{top}}(\hat{f}).
$$
M.Gromov a montr\'e dans [Gr 1] que ces deux notions co\"{\i}ncident et que l'entropie du
graphe it\'er\'e peut s'estimer en cohomologie gr\^ace au taux de croissance
asymptotique du volume des graphes it\'er\'es $\Gamma_f^N$,
$$
\text{lov}(f):=\limsup_{N \rightarrow +\infty} \frac{1}{N} \log \text{Vol}(\Gamma_f^N).
$$

\begin{thm}[Gromov]
$$
h_{\text{top}}(f) \leq \text{lov}(f)=\max_j \log \l_j(f).
$$
\end{thm}

\begin{esquisse}
Rappelons que $\Gamma_{f}^N$, le graphe d'ordre $N$ de $f$, est le sous-ensemble
analytique irr\'eductible de dimension $k$ dans $X^N$, d\'efini par
$$
\Gamma_f^N:=
\left\{ (x_0,\ldots,x_{N-1}) \in X^N \, / \, x_i=f^i(x), \;
0 \leq i \leq N-1 \right\}.
$$
Notons $\pi_i:X^N \rightarrow X$ la projection sur le $i^e$ facteur et
$\om_N:=\sum_{i=0}^{N-1} \pi_i^* \om$ la forme de K\"ahler induite sur $X^N$.
Alors
$$
\text{Vol}(\Gamma_f^N)=\int_{X^N} [\Gamma_f^N] \wedge \om_N^k.
$$
La preuve repose sur les deux observations suivantes:
\begin{enumerate}
\item Si $F$ est un ensemble $(N,\e)$-s\'epar\'e dans $\Omega_f$ pour la
  distance $d$ associ\'ee \`a $\om$, alors 
  $F_N:=\{(x,f(x),\ldots,f^{N-1}(x)) \in X^N \, / \, x \in F \}$ est un
  ensemble $(1,\e)$-s\'epar\'e dans $\Gamma_f^N$ pour la distance $d_N$
  associ\'ee \`a $\om_N$. Les boules $B_{d_N}(y,\e/2)$ sont donc disjointes
  pour $y \in F_N$, ce qui garantit
$$
\sharp F \cdot \min_{y \in F_N}
\int_{B_{d_N}(y,\e/2)} [\Gamma_f^N] \wedge \om_N^k 
\leq \int_{X^N} [\Gamma_f^N] \wedge \om_N^k.
$$
\\

\item Le volume d'une boule centr\'ee en un point $y$ 
  de $\Gamma_f^N$ est minor\'e
  ind\'epen-damment de $N$ et de $y$. 
  Plus pr\'ecis\'ement il existe $C=C(X,\om)>0$ et
  $\e_0>0$ tels que si $0 <\e<\e_0$ et $y \in \Gamma_f^N$, alors
$$
\int_{B_{d_N}(y,\e/2)} [\Gamma_f^N] \wedge \om_N^k \geq C \e^{2k}.
$$
C'est l'observation d\'esormais classique de P.Lelong qui permet de d\'efinir
le nombre de Lelong d'un courant positif ferm\'e de bidimension $(k,k)$
(ici le courant d'int\'egration sur $\Gamma_f^N$). 
\end{enumerate}
Il reste \`a calculer effectivement l'invariant $\text{lov}(f)$. Or
$$
\text{Vol}(\Gamma_f^N)=\sum_{0 \leq i_1,\ldots i_k \leq N-1} 
\int_X (f^{i_1})^* \om \wedge \cdots \wedge (f^{i_k})^* \om.
$$
En particulier $\text{lov}(f) \geq \max_j \log \l_j(f)$ (prendre
$i_1=\cdots =i_j=N-1$ et $i_{j+1}=\cdots =i_k=0$). L'in\'egalit\'e r\'eciproque
s'obtient en observant, lorsque $i_1 \leq i_2 \leq \cdots \leq i_k$, que
$$
\int_X (f^{i_1})^* \om \wedge \cdots \wedge (f^{i_k})^* \om
\leq C_{\e} \l_k(f)^{i_1} [\l_{k-1}(f)+\e]^{i_2-i_1} \cdots 
[\l_1(f)+\e]^{i_k-i_{k-1}},
$$
avec $\e>0$ arbitrairement petit.
Nous renvoyons le lecteur \`a [Gr 1], [Fr] pour plus de d\'etails.

\end{esquisse}
\vskip.2cm

Il r\'esulte par ailleurs des travaux de Y.Yomdin [Y] qu'on a \'egalement
une minoration
$h_{top}(f) \geq \max_j \log \l_j(f)$. Il s'ensuit,
lorsque $f$ est {\bf holomorphe}, qu'on a l'\'egalit\'e 
$$
h_{top}(f)=\max_{1 \leq j \leq k} \log \l_j(f).
$$

Lorsque $f$ est seulement {\it m\'eromorphe}, on peut adopter des d\'efinitions tout \`a fait
similaires --mais il faut \^etre soigneux. Dans la d\'efinition \`a la Bowen,
on consid\`ere uniquement des ensembles $(n,\e)$-s\'epar\'es qui se situent dans le plus grand
ensemble totalement invariant sur lequel $f$ et ses it\'er\'es sont holomorphes,
$$
\Omega_f:= X \setminus \bigcup_{n \in \Z} f^n(I_f).
$$
Dans la d\'efinition via le graphe, on consid\`ere l'adh\'erence dans $X^{\N}$ du graphe 
it\'er\'e holomorphe qui est inclus dans $\Omega_f^{\N}$ (il ne faut pas consid\'erer
le graphe infini fabriqu\'e \`a partir de l'adh\'erence du graphe holomorphe de $f$).

Nous renvoyons le lecteur \`a [G 7] pour des d\'efinitions pr\'ecises,
ainsi que pour la justification des assertions suivantes: lorsque 
$f$ est m\'eromorphe
\begin{itemize}
\item on a encore $h_{\text{top}}^{\text{Bow}}(f)=h_{\text{top}}^{\text{Gr}}(f)
=:h_{\text{top}}(f)$;
\item on a encore $h_{\text{top}}(f) \leq \text{lov}(f)$;
\item on a encore $\text{lov}(f)=\max_{1 \leq j \leq k} \log \l_j(f)$;
\item on peut avoir $h_{\text{top}}(f)<\text{lov}(f)$.
\end{itemize}
Cependant les seuls exemples que nous connaissons pour lesquels
$h_{\text{top}}(f)$ diff\`ere de $\text{lov}(f)$ correspondent \`a des endomorphismes
m\'eromorphes qui ne sont pas cohomologiquement hyperboliques.

Notons que l'\'egalit\'e $\text{lov}(f)=\max_j \log \l_j(f)$ est d\'elicate \`a 
justifier dans le cas m\'eromorphe. Elle est \'etablie dans [G 5] lorsque
$\dim_{\C} X \leq 3$ ou lorsque $X$ est une vari\'et\'e homog\`ene
(par exemple $X=\P^k$), dans [DS 3] en toute dimension lorsque $X$ est projective,
puis dans [DS 4] dans le cas g\'en\'eral: il s'agit, comme dans la preuve
de la Proposition 2.6, de savoir r\'egulariser convenablement les courants positifs ferm\'es.

\subsection{Entropie m\'etrique}

L'image directe d'une masse de Dirac $\e_x$ par un endomorphisme m\'eromorphe 
$f:X \rightarrow X$ est bien d\'efinie lorsque $x \notin I_f$, par
$$
f_* \e_x=\e_{f(x)}.
$$
On peut donc d\'efinir de m\^eme l'image directe $f_*\nu$ de toute mesure de probabilit\'e 
$\nu$ qui ne charge pas l'ensemble d'ind\'etermination $I_f$.
On va s'int\'eresser \`a des mesures invariantes, $f_* \nu=\nu$, et, parfois, \`a des mesures 
invariantes par image inverse. Il faut donc se restreindre et consid\'erer des mesures
de probabilit\'e $\nu$ telles que $\nu(\Omega_f)=1$.
On peut alors d\'efinir l'entropie m\'etrique $h_{\nu}(f)$ de $\nu$ comme on le fait habituellement
(voir e.g. [KH]). Nous montrons dans [G 7] les assertions suivantes:
\begin{itemize}
\item on a toujours $\sup_{\nu} h_{\nu}(f) \leq h_{top}(f)$ (principe variationnel faible);
\item il peut y avoir in\'egalit\'e stricte.
\end{itemize}
Cependant les seuls exemples que nous connaissons pour lesquels il y a in\'egalit\'e stricte
correspondent \`a nouveau
\`a des endomorphismes m\'eromorphes qui ne sont pas cohomologiquement hyperboliques.

Rappelons \`a pr\'esent que M.Misiurewicz et F.Przytycki ont montr\'e [MiP] que si 
$f:M \rightarrow M$ est un endomorphisme de classe ${\mathcal C}^1$ d'une vari\'et\'e lisse
orientable compacte, 
alors son entropie topologique est minor\'ee par le logarithme de son degr\'e topologique.
Ce n'est plus vrai si l'endomorphisme est seulement continu (cf e.g. [KH] p 317).
C'est \'egalement faux en g\'en\'eral pour un endomorphisme m\'eromorphe d'une vari\'et\'e
k\"ahl\'erienne compacte (voir Exemple 2.10). On peut toutefois esp\'erer
obtenir une telle minoration en utilisant la strat\'egie suivante: soit
$\Theta$ une mesure de probabilit\'e lisse sur $X$. Alors la suite
$$
\nu_n:=\frac{1}{n} \sum_{j=0}^{n-1} \frac{1}{\l_k(f)^j} (f^j)^* \Theta
$$
d\'efinit une suite de mesures de probabilit\'e sur $X$. Soit
$\mu=\lim \nu_{n_i}$ une valeur d'adh\'erence de la suite $\nu_n$.
Il se pourrait malheureusement que $\mu$ charge les points d'ind\'etermination $I_f$.
On peut aussi esp\'erer que ce ne sera pas le cas si on fait un choix intelligent
de la mesure $\Theta$ de d\'epart (voir Exemples 2.10 et 2.16 pour des choix heureux et malheureux).
Supposons que $\mu$ ne charge pas $I_f$. Alors la mesure $f_* \mu$ est bien d\'efinie
et v\'erifie
$$
f_* \mu =\lim_{i \rightarrow +\infty} f_* \nu_{n_i}=
\lim_{i \rightarrow +\infty} \nu_{-1+n_i}=\mu,
$$
donc $\mu$ est une mesure invariante.

Supposons \`a pr\'esent que $\mu$ ne charge pas les valeurs critiques $V_f=f({\mathcal C}_f)$
de $f$. C'est automatique sur une vari\'et\'e de dimension de Kodaira positive
(voir section 2.4.2). Cela n'est pas n\'ecessairement vrai lorsque $kod(X)=-\infty$,
mais ce sera toujours le cas lorsque $f$ a un grand degr\'e topologique (voir section 3.1).
On peut alors d\'efinir $f^* \mu:=\sum_{i=1}^{\l_k(f)} (f_i^{-1})_* \mu$, o\`u 
$f_i^{-1}$ d\'esigne les branches inverses (locales) de $f$ qui sont bien d\'efinies
hors de $V_f$. L'op\'erateur $\nu \mapsto f^* \nu$ est continu sur l'espace
des mesures qui ne chargent pas $V_f$. Comme $f^* \nu_n \simeq \l_k(f) \nu_{n+1}$, on
obtient dans ce cas
$$
f^* \mu=\l_k(f) \mu,
$$
autrement dit $\mu$ est de jacobien constant \'egal \`a $\l_k(f)$.
Dans une telle situation on r\'ecup\`ere la minoration de Misiurewicz-Przytycki
gr\^ace \`a l'argument suivant de J.Y.Briend et J.Duval [BrD 2].

\begin{pro}
Soit $\mu$ une mesure de probabilit\'e invariante qui ne charge pas les valeurs critiques de
$f$ et v\'erifie $f^* \mu=\l_k(f) \mu$.
Alors
$$
h_{\text{top}}(f) \geq h_{\mu}(f) \geq \log \l_k(f).
$$
\end{pro}

Remarquons qu'il peut y avoir in\'egalit\'e stricte comme le montre l'exemple
de la mesure de Lebesgue sur un tore complexe qui admet
un endomorphisme de type Anosov.

\begin{preuve}
Nous supposons $\mu$ ergodique pour simplifier --ce sera d'ailleurs dans ce cadre que nous l'utiliserons.
Pour $x \in X, \e>0$ fix\'es, on note
$$
B_n(x,\e):=\{ y \in X \, / \, d_n(x,y)<\e\}
$$
la boule dynamique associ\'ee \`a 
$d_n(x,y)=\max_{0 \leq j \leq n-1} d(f^jx,f^jy)$,
la distance dynamique. Il r\'esulte du Th\'eor\`eme de
Shannon-McMillan (voir [BrKa]) que pour $\mu$ presque tout $x$, on a 
$$
h_{\mu}(f)=\sup_{\e >0} \liminf_{n \rightarrow +\infty} -\frac{1}{n} \log \mu(B_n(x,\e)).
$$
Nous allons minorer $h_{\mu}(f)$ en majorant $\mu(B_n(x,\e))$.
L'id\'ee est la suivante: si $B$ est un Bor\'elien sur lequel $f$ est injective, alors
$\mu(B)=\l_k^{-1} \mu(fB)$ car $\mu$ est de jacobien constant$=\l_k$.
En particulier si $f$ est injective sur $B_n(x,\e)$, alors
$$
\mu(B_n(x,\e))=\frac{1}{\l_k} \mu(f B_n(x,\e)) \leq 
\frac{1}{\l_k} \mu (B_{n-1}(fx,\e))
$$
car $f B_n(x,\e) \subset B_{n-1}(fx,\e)$. Si $f$ est \`a nouveau injective sur 
$B_{n-1}(x,\e)$,  puis sur $B_{n-j}(f^jx,\e)$ pour tout $j \leq n$, on obtient
ainsi
$\mu(B_n(x,\e)) \leq \l_k^{-n} \mu(B(f^nx,\e)) \leq \l_k^{-n}$, ce qui donne
$h_{\mu}(f) \geq \log \l_k$.

Pour rendre rigoureux ce raisonnement simpliste, nous allons montrer que $f$ est
souvent injective sur $B_{n-j}(f^{n-j}(x,\e))$, pour un choix g\'en\'erique de $x$.
Plus pr\'ecis\'ement fixons $0<\d<<1$ et $U$ un petit voisinage de l'ensemble 
$V_f$ des valeurs critiques de $f$, de sorte que $\mu(U)< \d/2$. Posons
$$
X_n(\d):=\left\{ x \in X \, / \, \sharp\{j \in [0,n-1]/ f^jx \in U \} \leq n \d \right\}.
$$
Il r\'esulte du Th\'eor\`eme de Birkhoff que pour $\mu$ presque tout $x \in X$, on a
$$
\lim_{n \rightarrow +\infty} \frac{\sharp\{ j \in [0,n-1]/ f^j(x) \in U \}}{n}
=\mu(U)<\d/2.
$$
Il s'ensuit que $X_n(\d)$ est de mesure positive pour $n$ assez grand. Pour $x \in X_n(\d)$
$\mu$-g\'en\'erique, on applique l'estimation de jacobien constant lorsque
$f^jx \notin U$: il suffit de fixer $\e$ assez petit pour que la boule dynamique
$B_{n-j}(f^jx,\e)$ ne rencontre pas $V_f$ -prendre par exemple $\e<d(V_f,\partial U)$.
Lorsque $f^jx \in U$, on majore brutalement
$$
\mu(B_{n-j}(f^jx,\e)) \leq \mu(B_{n-j-1}(f^{j+1}x,\e))
$$
en utilisant l'invariance $f_*\mu=\mu$. Comme l'orbite de $x$ est souvent 
hors de $U$, on obtient l'estimation
$$
\mu(B_n(x,\e)) \leq \frac{1}{\l_k^{n(1-\d)}}.
$$
Enfin $\d$ \'etant arbitrairement petit, le principe variationnel donne
$$
h_{top}(f) \geq h_{\mu}(f) \geq \log \l_k.
$$
\end{preuve}

\section{Points p\'eriodiques}

\subsection{La formule de Lefschetz}
Nous cherchons ici \`a estimer le nombre de points p\'eriodiques d'ordre $n$ 
de l'endomorphisme $f$.
Quitte \`a changer $f$ en $f^n$, il s'agit d'estimer le nombre de points fixes.
Nous verrons plus loin (2.3.3) qu'il peut exister des courbes de points fixes.
Il s'agit bien s\^ur de points fixes non-hyperboliques, or ce sont 
les points fixes hyperboliques 
que nous souhaitons d\'enombrer. Comme ils sont stables par petite perturbation et comme
l'existence d'une courbe de points fixes est une situation exceptionnelle, 
nous ne consid\'erons que le cas o\`u {\it il n'y a pas de courbe de points fixes.}
Consid\'erons
$$
\Gamma_f:=\overline{ \{ (x,f(x)) \in X^2 \, / \, x \in X \setminus I_f \}}
$$
l'adh\'erence du graphe holomorphe de $f$ dans $X^2$ et 
soit $\Delta:=\{ (x,x) \in X^2 \, / \, x \in X \}$ la diagonale de $X$.
Le produit d'intersection  $\Gamma_f \cdot \Delta$ est bien d\'efini 
au sens des courants positifs --car
il n'y a pas de courbe de points fixes-- et compte, avec multiplicit\'e, le nombre de points fixes
plus le nombre de points d'ind\'etermination. On a donc
$\Gamma_f \cdot \Delta \geq \sharp Fix(f)$.
Par ailleurs ce produit d'intersection se calcule en cohomologie: c'est le contenu de
la formule des points fixes de Lefschetz (voir [GH] p 423) qui assure
$$
\Gamma_f \cdot \Delta=\sum_{0 \leq p,q \leq k} (-1)^{p+q}
\text{Trace} \left( f^*_{|H^{p,q}(X,\R)} \right).
$$

Lorsque $f$ est cohomologiquement hyperbolique, il r\'esulte du Th\'eor\`eme 2.4
qu'on obtient la majoration 
\begin{equation}
\sharp Fix(f^n) \lesssim \d_l(f^n,\om), \text{ pour } n >>1,
\end{equation}
o\`u $\l_l(f)$ d\'esigne le plus grand des degr\'es dynamiques.
On en d\'eduit 
\begin{equation}
\limsup \frac{1}{n}\log \sharp Fix(f^n) \leq \log \l_l(f).
\end{equation}
Notons que l'estimation (2.2) est plus pr\'ecise que (2.3). Lorsque 
$l=k$, on peut remplacer $\d_k$ par $\l_k$ dans (2.2).
Lorsque $l \leq k-1$, on esp\`ere \'egalement pouvoir remplacer $\d_l$
par $\l_l$ dans (2.2). C'est probablement possible si on sait trouver un mod\`ele
$\tilde{X}$ bim\'eromorphe \`a $X$ sur lequel l'action
lin\'eaire induite par $f^*$ sur $H^{l,l}(\tilde{X},\R)$ est
compatible avec la dynamique (voir section 4.2 et l'introduction du chapitre 5;
voir \'egalement [BFJ] pour une r\'eponse en dimension deux).

\subsection{Quels points p\'eriodiques ?}

Lorsque $\dim_{\C}X =1$, $f$ admet une infinit\'e de points p\'eriodiques et tous 
--sauf un nombre fini d'entre eux-- sont r\'epulsifs (voir [CG]).
La situation est plus vari\'ee en dimension sup\'erieure.
Il r\'esulte des travaux de S.Newhouse [N] (voir \'egalement 
[Ga], [Bu]) qu'il peut coexister
une infinit\'e de cycles attractifs et une infinit\'e de points selles
(resp. points r\'epulsifs). Nous allons donner quelques exemples de dimension $2$
qui illustrent d'autres diff\'erences avec la dimension $1$.

\begin{exa}
Soit $f:(z,w) \in \C^2 \mapsto (z^{\l},w+1) \in \C^2$,
avec $\l \geq 2$. Alors $f$ induit un endomorphisme
m\'eromorphe de $\P^2$ qui a un point fixe $q_{\infty}$ et un point d'ind\'etermination $I_f$
\`a l'infini.
On v\'erifie ais\'ement qu'il n'y a pas d'autre point p\'eriodique que $q_{\infty}$ et
que $\l_1(f)=\l_2(f)=\l \geq 2$. On obtient \'egalement (voir [G 7])
$$
0=h_{\text{top}}(f)<\text{lov}(f)=\log \l,
$$
bien que $f$ soit de degr\'e topologique $\l \geq 2$.
En particulier la strat\'egie propos\'ee en 2.2.2 \'echoue dans ce cas: on v\'erifie
en effet ais\'ement que
la suite $\l^{-n}(f^n)^* \Theta$ converge vers le point d'ind\'etermination $I_f$,
quel que soit le choix de la mesure lisse de probabilit\'e $\Theta$.
\end{exa}

En dimension 1, tous les cycles r\'epulsifs sont dans le support de
la mesure d'entropie maximale $\mu_f$ que nous avons construite
dans la section 1.1. Ce n'est plus le cas en dimension
sup\'erieure comme l'ont observ\'e J.H.Hubbard et P.Papadopol [HP 1]:

\begin{exa} 
Soit $f_{\e}:(z,w) \in \C^2 \mapsto (P(w)+\e z^r, Q(z)+R(w)) \in \C^2$,
o\`u $P,Q,R$ sont des polyn\^omes de degr\'e $p,q,r$ avec $pq <r$.
On consid\`ere les extensions m\'eromorphes de ces applications \`a $\P^2$
--encore not\'ees $f_{\e}$.
Lorsque $\e=0$, on obtient un endomorphisme m\'eromorphe de $\P^2$ tel que
$$
\l_1(f_0)=r \text{ et } \l_2(f_0)=pq<\l_1(f_0).
$$
L'application $f_0$ a un unique point d'ind\'etermination $I_0$ \`a 
l'infini. Il est facile d'expliciter des polyn\^omes
$P,Q,R$ tels que $f_0$ ait un point fixe r\'epulsif loin de $I_0$.

Les autres applications ($\e \neq 0$) d\'efinissent des endomorphismes holomorphes
de $\P^2$, pour lesquels on a 
$\l_1(f_{\e})=r$ et $\l_2(f_{\e})=r^2$. On construira un peu plus loin (Chapitre 3)
une mesure d'entropie maximale $\mu_{\e}$ pour les endomorphismes $f_{\e}$,
$\e \neq 0$. Les mesures $\mu_{\e}$
sont localis\'ees pr\`es de $I_0$ et d\'eg\'en\`erent en une masse de Dirac 
au point $I_0$
lorsque $\e \rightarrow 0$. En particulier, les points p\'eriodiques r\'epulsifs de
$f_0$ situ\'es loin de $I_0$ g\'en\`erent des points p\'eriodiques r\'epulsifs
de $f_{\e}$ ($\e$ petit) qui se situent hors du support de $\mu_{\e}$:
il est donc possible que certains points p\'eriodiques r\'epulsifs se situent
loin du support de la mesure d'entropie maximale, bien que la plupart en soient
proches, puisqu'ils
s'\'equidistribuent selon elle (voir chapitre 3.2). Cette observation est due \`a
J.H.Hubbard et P.Papadopol [HP 1] qui consid\`erent des perturbations d'applications
de H\'enon complexes (voir \'egalement [FS 5] pour une autre g\'en\'eralisation).
\end{exa}

\subsection{Courbes de points p\'eriodiques}

Les endomorphismes holomorphes cohomologiquement hyperboliques 
n'ont essentiellement pas de courbes de points p\'eriodiques:

\begin{pro}
Soit $f:X \rightarrow X$ un endomorphisme {\bf holomorphe} d'une surface 
k\"ahl\'erienne compacte tel que $\l_1(f) \neq \l_2(f) \geq 2$.
Alors $f$ est induit par un endomorphisme holomorphe $g:Y \rightarrow Y$ sur un mod\`ele
minimal $Y$ de $X$, tel que $g$ n'admet pas de courbe de points fixes.
\end{pro}

Notons qu'il est facile de produire un endomorphisme holomorphe ``non-minimal''
qui admet une courbe de points fixes: il suffit d'\'eclater en un point
fixe en lequel l'application a une diff\'erentielle \'egale \`a l'identit\'e.

\begin{preuve}
Supposons qu'il existe une courbe irr\'eductible ${\mathcal C}$ de points fixes de $f$.
Alors $f_* {\mathcal C}={\mathcal C}$.
Comme $f_*$ est un op\'erateur inversible sur l'espace $NS(X)$ de N\'eron-Severi
r\'eel (car $f_*f^*D=\l_2(f) D$),
on en d\'eduit que $f^* {\mathcal C}=\l_2(f) {\mathcal C}$.

Soit $\theta \in H^{1,1}_{nef}(X,\R)$ une classe non-nulle telle que $f^* \theta=r_1(f) \theta$.
Il vient
$$
r_1(f) \l_2(f) \langle {\mathcal C},\theta\rangle =\langle f^* {\mathcal C},f^* \theta\rangle =
\l_2(f) \langle {\mathcal C},\theta\rangle .
$$
Si $\langle {\mathcal C},\theta\rangle >0$ alors $r_1(f)=1$, donc $\l_1(f)=1$, donc $\l_2(f)=\l_1(f)=1$
(car $1 \leq \l_2(f) \leq \l_1(f)^2$ par 2.4.a). 
Supposons donc $\langle {\mathcal C},\theta\rangle =0$. Comme $\theta^2 \geq 0$, il r\'esulte du
Th\'eor\`eme de l'indice de Hodge que soit ${\mathcal C}$ est proportionnelle \`a $\theta$
--mais alors $\l_1(f)=\l_2(f)$--, soit ${\mathcal C}^2<0$.

Lorsque ${\mathcal C}^2<0$, ${\mathcal C}$ fait partie du nombre fini de courbes 
d'auto-intersection n\'egative, dont on montrera dans la section 3.3.1 qu'elles sont obtenues par
\'eclatements successifs \`a partir d'un endomorphisme sur un mod\`ele minimal.
On v\'erifie ais\'ement qu'il n'y a pas de courbe de points fixes
lorsque $X$ est minimale et $\l_1(f) \neq \l_2(f)$
(voir Th\'or\`eme 3.6).
\end{preuve}
\vskip.2cm

Ce r\'esultat n'est plus valable si l'on permet des
points d'ind\'etermination:

\begin{exa}
Soit $f:(z,w) \in \C^2 \mapsto (zw^b,w+z^cw^d) \in \C^2$, o\`u $b,c,d \in\N$.
Alors $f$ induit un endomorphisme
m\'eromorphe de $\P^2$ tel que $f(0,w)=(0,w)$. On v\'erifie 
--en travaillant en fait dans $\P^1 \times \P^1$, cf [FaG]-- que
$\l_1(f)$ est le rayon spectral de la matrice 
$\left[ \begin{array}{cc} 1 & b \\ c & d \end{array} \right]$ et que
$$
\l_2(f)=\max(d-bc, [bc-d]+1).
$$
En faisant varier les valeurs des entiers $b,c,d$, on peut obtenir aussi bien
$\l_1(f)>\l_2(f)$ --par exemple pour $b=c=d=1$-- que $\l_1(f)<\l_2(f)$
--par exemple pour $d=0$ et $b=c=2$.
\end{exa}

Cette situation reste n\'eanmoins tr\`es exceptionnelle.
Les endomorphismes birationnels qui admettent une
courbe de points fixes ont \'et\'e partiellement
classifi\'es par D.Jackson [Ja].
Pour des r\'esultats sur la g\'eom\'etrie des courbes invariantes
par des endomorphismes rationnels, nous renvoyons le lecteur \`a
[BDM] et [DJS]. On pourra \'egalement consulter [ABT]
pour une \'etude locale des courbes invariantes.

\section{Quelles vari\'et\'es ?} 

L'objectif de cette section est de fournir des exemples d'endomorphismes 
m\'eromorphes $f:X \rightarrow X$
cohomologiquement hyperboliques, i.e. tels que $\l_l(f)>\max_{j \neq l} \l_j(f)$
pour un indice $l \in [1,k]$.
Nous commen\c{c}ons par justifier pourquoi il faut les chercher sur
les vari\'et\'es dont la dimension de Kodaira vaut $0$ ou $-\infty$.

Soit $s \in H^0(X,K_X)$ une section holomorphe du fibr\'e canonique $K_X$, i.e. une $(k,0)$-forme
holomorphe. Si $f:X \rightarrow X$ est un endomorphisme m\'eromorphe, alors
$f^*s$ est une forme holomorphe sur $X \setminus I_f$ qui s'\'etend holomorphiquement
\`a travers $I_f$, car $codim_{\C} I_f \geq 2$ (ph\'enom\`ene de Bochner).
L'application $f$ induit donc une application {\it lin\'eaire} sur $H^0(X,K_X)$ et,
plus g\'en\'eralement, sur $H^0(X,K_X^m)$, pour tout $m \in \N$.
Soit $s_0,\ldots,s_{N_m}$ une base de $H^0(X,K_X^m)$, alors
$$
\Phi_m:x \in X \mapsto [s_0(x):\cdots:s_{N_m}(x)] \in \P^{N_m}
$$
est une application rationnelle de $X$ sur son image $X_m=\Phi_m(X)$
qui est de dimension $kod(X)$ pour $m>>1$.
On a donc un diagramme commutatif
$$
\begin{array}{ccc}
X & \stackrel{f}{\rightarrow} & X \\
\Phi_m \downarrow & \text{ } & \downarrow \Phi_m \\
X_m & \stackrel{F_m}{\rightarrow } & X_m
\end{array}
$$
o\`u $F_m$ est une application projective lin\'eaire: on dit que $f$ pr\'eserve la
fibration d'Itaka. Notons que $F_m$ est un isomorphisme car $f$ est dominante.

Cela impose des relations entre les degr\'es dynamiques lorsque $kod(X)=\dim_{\C} X_m \geq 1$.
On obtient par exemple $\l_1(f)=\cdots=\l_k(f)=1$ lorsque $X$ est de type g\'en\'eral
($kod(X)=k$): dans ce cas $\Phi_m$ est une application birationnelle
et $\l_j(f)=\l_j(F_m)=1$ (cf corollaire 2.7): c'est une version simple
et faible d'un r\'esultat de S.Kobayashi et T.Ochiai [KO], qui montrent en fait qu'un it\'er\'e
de $f$ est \'egal \`a l'identit\'e.
Nous conjecturons plus g\'en\'eralement que $f$ ne peut pas \^etre
cohomologiquement hyperbolique lorsque $kod(X) \geq 1$. Nous \'etablissons ce fait en dimension deux.

\begin{thm}
Soit $f:X \rightarrow X$ un endomorphisme m\'eromorphe d'une surface k\"ahl\'erienne compacte
tel que $\l_1(f) \neq \l_2(f)$. Alors
\begin{itemize}
\item soit $kod(X)=0$;
\item soit $X$ est rationnelle;
\item soit $X$ est birationnelle \`a $\P^1 \times E$, $\text{genre}(E)=1$, et
$\l_2(f) >\l_1(f)$.
\end{itemize}
\end{thm}

Nous donnons ci-apr\`es des exemples dans chacune de ces cat\'egories.

\begin{preuve}
La d\'emonstration repose sur le Lemme 2.15 ci-dessous.
Lorsque $kod(X) \geq 1$, on applique le lemme 
en utilisant la fibration d'Itaka, avec $\pi=\Phi_m$, $g=F_m$. Comme $g$ est
de degr\'e $1$, il vient $\l_1(f)=\l_2(f)$. 
Lorsque $kod(X)=-\infty$, on peut \'egalement consid\'erer la fibration d'Albanese:
au lieu de consid\'erer des $2$-formes holomorphes, on consid\`ere des
$1$-formes holomorphes. L'application $f$ induit \'egalement 
un isomorphisme lin\'eaire sur $H^0(X,\Omega_X^1)$ qui pr\'eserve
le r\'eseau $H_1(X,\Z)$ (modulo torsion) et induit donc
une application lin\'eaire sur la vari\'et\'e d'Albanese,
$$
\text{Alb}(X):=H^0(X,\Omega_X^1) / H_1(X,\Z) \text{ mod torsion}.
$$
Lorsque $kod(X)=-\infty$, on peut supposer --\`a changement birationnel de coordonn\'ees
pr\`es-- que $X=\P^1 \times B$ est une surface r\'egl\'ee au dessus d'une courbe $B$.
Il y a beaucoup d'exemples int\'eressants, de tous types, lorsque $B=\P^1$. 
Lorsque $genre(B) \geq 1$, la fibration d'Albanese co\"{\i}ncide
avec la projection $\pi:X=\P^1 \times B \rightarrow B$ sur le second facteur (cf Lemme V.18
et Proposition V.15 dans [Bea]). On peut donc appliquer \`a nouveau le Lemme 2.15.
Il vient $\deg(g)=1$, donc $\l_1(f)=\l_2(f)$ lorsque $B$ est hyperbolique.
Lorsque $B$ est elliptique, on obtient $\l_2(f) \geq \l_1(f)$ et il est facile
de construire des exemples avec in\'egalit\'e stricte.
\end{preuve}

\begin{lem}
Soit $X$ une surface k\"ahl\'erienne compacte et $B$ une surface de 
Riemann compacte. On suppose qu'il existe des applications m\'eromorphes
$f,g,\pi$ telles que le diagramme suivant est commutatif
$$
\begin{array}{ccc}
X & \stackrel{f}{\rightarrow} & X \\
\pi \downarrow & \text{ } & \downarrow \pi \\
B & \stackrel{g}{\rightarrow } & B
\end{array}
$$
Alors $\l_1(f) \leq \l_2(f)$ avec \'egalit\'e si $\deg(g)=1$.
\end{lem}

\begin{preuve}
Remarquons qu'on peut r\'esoudre les singularit\'es de $\pi$ sans changer les donn\'ees
du probl\`eme: cela revient \`a remplacer $X$ par un mod\`ele bim\'ero-morphe $\tilde{X}$,
et $f$ par l'endomorphisme $\tilde{f}$ induit par $f$ sur $\tilde{X}$.
Comme $f$ et $\tilde{f}$ ont les m\^emes degr\'es dynamiques (corollaire 2.7), on s'est
ramen\'e au cas o\`u $\pi$ est holomorphe, ce que nous supposerons dans la suite.

Soit $d$ le degr\'e de l'application $g$. Observons que si $F$ est une fibre
g\'en\'erique de $\pi$, alors $F^2=0$ et $f^*F \simeq d F$. Par ailleurs le degr\'e topologique
de $f$ satisfait $\l_2(f)=d \cdot m$, o\`u $m$ est le degr\'e de
l'application $f_{|F}$.

Si $p \in I_f$, son image $f(p)$ est incluse dans la fibre $F_{g(a)}$ telle que $p \in F_a$.
Il s'ensuit que pour toute classe de cohomologie $\theta \in H^{1,1}(X,\R)$,
on a
\begin{equation}
\langle f^* \theta,f^*F\rangle =\l_2(f) \langle \theta,F\rangle .
\end{equation}
Il s'agit d'une version simple de la ``formule d'aller-retour'' (Proposition 4.8).

Soit en particulier $\theta \in H^{1,1}_{nef}(X,\R)$ une classe nef telle
que $f^* \theta=r_1(f) \theta$ (voir section 4.2). Si $\theta=F$, on obtient
$r_1(f)=\deg(g) \leq \l_2(f)$, donc $\l_1(f) \leq \l_2(f)$. De plus
si $\deg(g)=1$ on obtient
$r_1(f)=1$ donc $\l_1(f)=1$, d'o\`u $\l_2(f)=\l_1(f)=1$ (cf Th\'eor\`eme 2.4.a).

Si $\theta$ n'est pas proportionnelle \`a $F$,
il r\'esulte du Th\'eor\`eme de l'indice de Hodge
que $\theta \cdot F>0$. La formule d'intersection ci-dessus donne ainsi
$r_1(f) \cdot d=\l_2(f)$, d'o\`u le r\'esultat.
\end{preuve}

\vskip.2cm
Nous donnons \`a pr\'esent des exemples d'endomorphismes m\'eromorphes 
cohomologiquement hyperboliques.

\subsection{$kod(X)=-\infty$}

Il y a une multitude d'exemples d'endomorphismes rationnels de l'espace
projectif complexe $\P^k$ --donc de toute vari\'et\'e rationnelle.
Nous en donnons quelques uns ci-dessous.
Lorsque $X$ est une vari\'et\'e {\it unirationnelle}, i.e. lorsqu'il existe
une application m\'eromorphe dominante $\Phi:\P^k \rightarrow X$, alors
$X$ admet \'egalement des endomorphismes m\'eromorphes induits par ceux de $\P^k$
(voir Exemple 3.9). Il serait int\'eressant de donner une classification grossi\`ere
des vari\'et\'es k\"ahl\'eriennes compactes telles que $kod(X)=-\infty$ et
qui admettent des endomorphismes cohomologiquement hyperboliques.
Le Th\'eor\`eme 2.14 ci-dessus indique que l'existence de tels endomorphismes
impose des contraintes particuli\`eres sur la g\'eom\'etrie de la vari\'et\'e.

\begin{exa}
Soit $f:(z,w) \in \C^2 \mapsto (P(w),Q(z)+R(w)) \in \C^2$, o\`u $P,Q,R$ sont des polyn\^omes
de degr\'e $p,q,r$ avec $r \geq \max(p,q)$. Alors $f$ induit un endomorphisme
m\'eromorphe de $\P^2$ (ou de toute surface rationnelle) tel que
$$
\l_1(f)=r \text{ et } \l_2(f)=pq.
$$
On peut donc, en fonction des valeurs respectives de $p,q,r$, obtenir aussi bien $\l_1(f)>\l_2(f)$
que $\l_1(f)<\l_2(f)$.

Lorsque $p=q=1$ et $r \geq 2$, $f$ d\'efinit un automorphisme polynomial de $\C^2$
d'entropie positive: c'est une {\it application de H\'enon complexe}.
Lorsque $p=q=r$, on obtient un endomorphisme {\it holomorphe} de $\P^2$.

On peut bien s\^ur consid\'erer des compositions de telles applications.
La dynamique de ces endomorphismes est \'etudi\'ee dans [G 1]. On y montre
en particulier que l'ensemble d'ind\'etermination $I_f$ est $f^{-1}$-attirant
lorsque $\l_2(f)=pq <\l_1(f)=r$. Si $\Theta$ est une mesure de probabilit\'e lisse
localis\'ee pr\`es de $I_f$, on a donc
$\Theta_n:=\l_2(f)^{-n}(f^n)^* \Theta \rightarrow \d_{I_f}$: c'est un choix malheureux de $\Theta$
pour la strat\'egie propos\'ee en 2.2.2. On montre \'egalement dans [G 1] que
$\C^2$ est la r\'eunion disjointe du bassin d'attraction de $I_f$ (sous it\'eration
de $f^{-1}$) et de l'ensemble ferm\'e $K^-$ des points d'orbite n\'egative born\'ee. Si
ce dernier est d'int\'erieur non-vide (par exemple lorsque $f$ admet 
un point p\'eriodique r\'epulsif), on peut choisir $\Theta$ localis\'ee dans ce bassin
et obtenir des valeurs d'adh\'erence pour $\Theta_n$ qui ne chargent pas $I_f$.
\end{exa}

Nous verrons dans la section 3.3.1 que les endomorphismes holomorphes non inversibles
des surfaces rationnelles sont induits par les endomorphismes holomorphes
d'un mod\`ele minimal. Il y en a donc peu, hormis sur $\P^2$.
Il reste le cas --assez mal compris-- des automorphismes
d'entropie positive des surfaces rationnelles. En voici un exemple 
qui s'inspire de la construction des exemples de Latt\`es. 
Il est analys\'e par S.Cantat et C.Favre dans [CF].

\begin{exa}
Consid\'erons le tore complexe de dimension deux $Y=E \times E$, o\`u $E=\C/\Z[\zeta]$
est la courbe elliptique associ\'ee \`a une racine primitive de l'unit\'e $\zeta$ d'ordre
$3,4$ ou $6$. La matrice 
$$
A=\left[ \begin{array}{cc} 2 & 1 \\ 1 & 1 
\end{array} \right] \in GL(2,\Z)
$$
pr\'eserve le r\'eseau $\Lambda=\Z[\zeta] \times \Z[\zeta]$; elle induit donc un endomorphisme
holomorphe $g:Y \rightarrow Y$ dont on v\'erifie que c'est un automorphisme
d'entropie positive
$$
h_{\text{top}}(f)=\log \l_1(f)=2\log \frac{3+\sqrt{5}}{2}>0.
$$
On peut v\'erifier que c'est un diff\'eomorphisme d'Anosov (voir [Gh]).

Observons que $g$ commute avec l'homoth\'etie 
$\sigma:[x,y] \in Y \mapsto [\zeta x, \zeta y] \in Y$. Soit $X$ la surface
obtenue en d\'esingularisant le quotient $Y/ \! \!\langle \sigma\rangle $. C'est une surface rationnelle
sur laquelle $g$ induit un automorphisme $f$ de m\^eme entropie.
E.Ghys et A.Verjovsky ont donn\'e [GV] une liste compl\`ete des tores complexes de dimension
deux qui admettent de tels diff\'eomorphismes holomorphes d'Anosov. On en d\'eduit la liste
des automorphismes d'entropie positive des surfaces rationnelles qui proviennent
d'un tel passage au quotient. Ces automorphismes h\'eritent tous d'une mesure m\'elangeante
\'equivalente \`a la mesure de Lebesgue (obtenue
en poussant la mesure de Lebesgue du tore).
\end{exa}

Il existe d'autres exemples.
E.Bedford et K.Kim ont montr\'e [BK 1,3]
que pour certaines valeurs des param\`etres $\a_0,\a_1,\a_2,\b_0,\b_2$, l'application
rationnelle de $\C^2$
$$
f(x,y)=\left(y,\frac{\a_0+\a_1x+\a_2y}{\b_0+\b_1x} \right)
$$
d\'efinit un automorphisme d'entropie positive sur un \'eclat\'e de $\P^2$
qui poss\`ede un domaine de Siegel
(C.McMullen a r\'ecemment g\'en\'eralis\'e cette construction [M 3]): 
sa mesure d'entropie maximale
(que nous construirons au chapitre 4) n'est donc pas \'equivalente
\`a la mesure de Lebesgue.

\subsection{$kod(X)=0$}

Lorsque $X$ est de dimension de Kodaira nulle, le diviseur canonique $K_X$
est un $\Q$-diviseur effectif ($K_X \geq 0$). Cela contraint les ramifications possibles
de $f$. Notons $R_f$ le diviseur de ramification. Il r\'esulte de la formule de Hurwitz
$f^*K_X+R_f=K_X$, que pour tout $n \in \N$,
$$
(f^*)^n K_X+\sum_{j=0}^{n-1} (f^*)^j R_f=K_X.
$$
Si $D$ est un diviseur ample, il vient donc
$$
0 \leq \sum_{j=0}^{n-1} \langle (f^*)^j R_f,D^{k-1}\rangle  \leq K_X \cdot D^{k-1}.
$$
On en d\'eduit que $(f^*)^j R_f=0$ pour $j \geq j_0$.
Cela implique que $f$ est essentiellement non ramifi\'e.

Supposons par exemple que $\dim_{\C} X=2$. Quitte \`a effectuer un changement
bim\'eromorphe de coordonn\'ees, on peut supposer que l'on travaille sur le mod\`ele
minimal. Dans ce cas $12 K_X=0$, donc $12 R_f=0$. Nous utiliserons cette observation
\`a plusieurs reprises, notamment dans la section 4 pour en d\'eduire qu'un tel
endomorphisme $f$ est 1-stable (voir Proposition 4.5).

\begin{exa}
Soit $X=\C^k / \Lambda$ un tore complexe de dimension $k \geq 1$ i.e. le quotient
de $\C^k$ par un r\'eseau $\Lambda$ de rang $2k$. Soit $f:X \rightarrow X$
un endomorphisme m\'eromorphe de $X$. 
Alors $f$ est en fait automatiquement holomorphe: sinon on consid\`ere
une d\'esingularisation 
$$
\begin{array}{ccccc}
\text{ } & \text{ } & \tilde{\Gamma}_f & \text{ } & \text{ } \\
\text{ } & \stackrel{\pi_1}{\swarrow} & \text{ } &
\stackrel{\pi_2}{\searrow} & \text{ } \\
X & \text{ } & \stackrel{f}{\longrightarrow} & \text{ } & X
\end{array}
$$
de $f$. On peut trouver un diviseur exceptionnel $E \simeq \P^{k-1}$ de $\pi_1$
qui est envoy\'e par $\pi_2$ sur un sous-ensemble analytique de dimension positive,
image par $f$ d'un hypoth\'etique point d'ind\'etermination. 
Or $\P^{k-1}$ est simplement connexe et le rev\^etement universel de $X$ est $\C^k$,
on peut donc relever l'application ${\pi_2}_{|E}$ en une application holomorphe
vers $\C^k$. L'image de $E$ est alors un sous-ensemble analytique
compact (Th\'eor\`eme de Remmert) et connexe dans $\C^k$ qui est une vari\'et\'e de
Stein, c'est donc un point. Il s'ensuit que $\pi_2$ contracte $E$ sur un point,
donc $f$ est holomorphe.

On montre d'une fa\c{c}on analogue que $f$ est induit par un endomorphisme affine
$F$ de $\C^k$, $F(z)=A \cdot z+v$, o\`u
$A \in GL(k,\C)$ et $v \in \C^k$ sont tels que $F \Lambda \subset \Lambda$.
Observons que $D_xf=A$ en tout point $x \in X$. La matrice $A$ repr\'esente
l'action de $f^*$ sur $H^{1,0}(X,\R)$ dans la base canonique 
d\'etermin\'ee par les 1-formes $dz_j$.
Notons $a_1,\ldots,a_k$ ses valeurs propres ordonn\'ees de telle sorte que
$|a_1| \geq \cdots \geq |a_k|$. Un calcul imm\'ediat donne
$$
\l_j(f)=r_j(f)=|a_1|^2 \cdots |a_j|^2,
$$
pour tout $1 \leq j \leq k$. En particulier
$\l_1(f)$ est \'egal au carr\'e du rayon spectral de la matrice $A$, et 
$\l_k(f)$ est \'egal au carr\'e du determinant de $A$.

La condition ``$f$ est cohomologiquement hyperbolique'' est donc
\'equivalente ici au fait que les valeurs propres de $A$ sont toutes de module
diff\'erent de $1$, c'est \`a dire que $f$ est uniform\'ement hyperbolique.

Observons que la contrainte $F \Lambda \subset \Lambda$ n'est pas facile \`a
satisfaire. Pour un r\'eseau g\'en\'erique $\Lambda$, les seules matrices
$A$ qui conviennent sont les homoth\'eties de rapport entier.
Dans ce cas $\l_k$ est le plus grand des degr\'es dynamiques. On peut cependant,
pour des choix particuliers de r\'eseau (par exemple pour un r\'eseau produit), 
obtenir des endomorphismes tels que $\l_l(f)$ domine tous les autres
degr\'es dynamiques, quel que soit $l$ fix\'e dans $[1,k]$. 
E.Ghys et A.Verjovsky [GV] ont donn\'e une liste des tores de dimension deux
qui admettent des diff\'eomorphismes d'Anosov, i.e. pour lesquels $\l_1(f)>\l_2(f)=1$.
\end{exa}

Lorsque $\dim_{\C} X=2$, $kod(X)=0$, il r\'esulte de la classification 
d'Enriques-Kodaira (voir [Bea]) que le mod\`ele minimal de $X$ est
\begin{itemize}
\item soit un tore;
\item soit le quotient d'un tore par un groupe fini d'automorphismes sans point fixe 
   (surface hyperelliptique);
 \item soit une surface $K3$;
 \item soit le quotient d'une surface $K3$ par une involution holomorphe sans point fixe
  (surface d'Enriques).
\end{itemize}
Les endomorphismes m\'eromorphes d'une surface d'Enriques se rel\`event en des endomorphismes
de sa surface $K3$. Ceux des surfaces hyperelliptiques
se rel\`event sur le tore et 
sont dynamiquement moins int\'eressants car ils pr\'eservent la fibration d'Albanese qui est
non triviale dans ce cas --on obtient donc des produits crois\'es de grand degr\'e topologique.
Il nous reste donc, en dimension 2, \`a donner des exemples d'endomorphismes m\'eromorphes
des surfaces $K3$. Il existe une vaste litt\'erature traitant du cas des automorphismes
d'entropie positive
(voir notamment [Ca 1], [M 1]). Nous donnons \`a pr\'esent
quelques exemples non inversibles sur une surface de Kummer.

\begin{exa}
Soit $S$ une surface de Riemann de genre 2 et $A$ la Jacobienne de $S$: c'est le tore complexe
projectif de dimension 2 d\'efini par les $1$-formes holomorphes sur $S$ (modulo p\'eriodes).
Soit $X$ la surface obtenue en d\'esingularisant le quotient de $A$ par l'involution
$\sigma:z \mapsto -z$: c'est la surface de Kummer de la courbe $S$ 
--un type particulier de surface $K3$ (voir Exemple V.10 dans [Bea]).
Elle poss\`ede 16 courbes nodales (des courbes rationnelles d'autointersection -2) qui proviennent des
\'eclatements r\'ealis\'es pour d\'esingulariser $A/ \! \! \langle \sigma\rangle $. 
Tout endomorphisme $g$ de $A$ commute bien s\^ur avec $\sigma$ et induit
un endomorphisme --m\'eromorphe si $\l_2(g) \geq 2$-- de la vari\'et\'e $X$. Notons que
ceux-ci doivent permuter les 16 courbes nodales.

Lorsque $S$ est g\'en\'erique, J.Keum a construit [Ke] des automorphismes qui ne 
permutent pas les courbes nodales et ne proviennent donc pas d'auto-morphismes de $A$.
Soit $\p$ un tel automorphisme, on peut le supposer d'entropie positive --i.e. $\l_1(\p)>1$.
Soit $h$ l'endomorphisme m\'eromorphe de $X$ 
de degr\'e topologique $\l_2(h) \geq 2$ induit par une homoth\'etie de rapport entier sur $A$.
Alors
$$
f:=\p^q \circ h^p
$$
d\'efinit un endomorphisme m\'eromorphe de $X$ qui ne provient pas d'un endomorphisme de $A$
si $q \geq 1$. Observons que
$$
\l_2(f)=\l_2(h)^p \; \text{ et } \; \l_1(f) \leq \l_1(\p)^q \l_1(h)^p.
$$
On peut donc obtenir des exemples tels que $\l_2(f) >\l_1(f)$
en prenant $p>>q$, car $\l_1(h)=\l_2(h)^{1/2}$. A l'inverse, si on fixe $p$ et
on choisit $q$ tr\`es grand, on obtient $\l_1(f)>\l_2(f)$. En effet, soit
$\theta$ une $(1,1)$-classe nef non nulle telle que $\p^* \theta=\l_1(\p) \theta$
(cf section 4.2). Alors
$$
\l_1(\p)^q \lesssim \l_1(\p)^q \int_X (h^p)^* \theta \wedge \om
=\int_X f^* \theta \wedge \om \lesssim \l_1(f).
$$
On utilise ici de fa\c{c}on essentielle le fait que $kod(X)=0$: comme les endomorphismes
$f,h$ ne sont pas ramifi\'es, les actions lin\'eaires induites en cohomologie
sont compatibles avec la composition des applications m\'eromorphes, ce qui
justifie ici l'\'egalit\'e $f^* \theta=(h^p)^*( (\p^q)^* \theta)$ et $r_1=\l_1$
(voir Proposition 4.5).
Observons enfin que $(h^p)^* \theta$ est une classe nef non-nulle, comme on
peut le voir en utilisant la formule d'aller-retour (Proposition 4.8).
\end{exa}

Le lecteur int\'eress\'e trouvera d'autres exemples dans [Ca 3].

\chapter{Grand degr\'e topologique}

Dans ce chapitre nous consid\'erons le cas o\`u le degr\'e topologique $\l_k(f)$
domine strictement tous les autres degr\'es dynamiques.
Nous construisons une mesure $\mu_f$ m\'elangeante d'entropie maximale
dans la {\it section 3.1}, puis 
d\'emontrons dans la {\it section 3.2} la conjecture \'enonc\'ee dans l'introduction.
Nous donnons des exemples d'endomorphismes v\'erifiant nos hypoth\`eses
dans la {\it section 3.3} puis mentionnons quelques r\'esultats suppl\'ementaires
sur la dynamique de ces applications ({\it section 3.4}).
Nous pr\'esentons enfin  ({\it section 3.5}) quelques \'el\'ements du travail de T.C.Dinh et N.Sibony
sur les applications d'allure polynomiale: ce sont des perturbations (\'eventuellement
transcendantes) d'endomorphismes {\it polynomiaux} de grand degr\'e topologique.

\section{La mesure canonique}

Nous allons construire une mesure invariante canonique $\mu_f$ pour les endomorphismes
m\'eromorphes de grand degr\'e topologique. Lorsque
$X=\P^k$ et $f$ est {\it holomorphe}, cette mesure a \'et\'e construite
par J.-E.Fornaess et  N.Sibony dans [FS 2,3,4], et 
par J.-H.Hubbard et  P.Papadopol dans [HP 1].
Lorsque $f$ est m\'eromorphe, A.Russakovskii et B.Shiffman 
ont construit la mesure $\mu_f$ [RS] en observant une propri\'et\'e
d'\'equidistribution des pr\'eimages de points. Leur construction
-- d\'elicate -- ne leur a pas permis d'\'etablir des propri\'et\'es
dynamiques de la mesure $\mu_f$, 
mais elle a motiv\'e de nombreux
travaux (voir par exemple [FaG],[G 1,3,5], [DS 2,8], [HP 2]).
Nous avons obtenu dans [G 3,5] une construction simple de cette
mesure --dite de Russakovskii--Shiffman-- qui nous a permis d'\'etablir le
r\'esultat suivant:

\begin{thm}
Supposons que $\l:=\l_k(f)>\max_{j \neq k} \l_j(f)$. Alors il existe une mesure de
probabilit\'e $\mu_f$ telle que
$$
\frac{1}{\l^n}(f^n)^* \Theta \longrightarrow \mu_f,
$$
pour toute mesure de probabilit\'e lisse $\Theta$ sur $X$. De plus
\begin{enumerate}
\item toute fonction qpsh est int\'egrable par rapport \`a $\mu_f$; en particulier $\mu_f$
ne charge pas les hypersurfaces et $\log^+ ||Df^{\pm1}|| \in L^1(\mu_f)$;
\item la mesure $\mu_f$ est invariante et de jacobien constant $f^* \mu_f=\l \mu_f$, donc
d'entropie maximale
$$
h_{\mu_f}(f)=h_{top}(f)=\log \l>0.
$$
\item la mesure $\mu_f$ est m\'elangeante; plus pr\'ecis\'ement il existe $C>0$
telle que pour toutes fonctions test $\chi,\p$,
$$
\left| \int_X \p \chi \circ f^n d\mu_f -\int_X \p d\mu_f \int_X \chi d\mu_f \right|
\leq C ||\p||_{{\mathcal C}^2} ||\chi||_{L^{\infty}} \frac{\d_{k-1}(f^n,\om)}{\l^n}.
$$
\end{enumerate}
\end{thm}

\begin{preuve}
Soit $a \in X$ un point qui n'est ni une valeur critique de $f$,
ni un point d'ind\'etermination. Alors $f$ est localement inversible
pr\`es de $a$. Soit $\Theta$ une mesure lisse de probabilit\'e dont
le support est concentr\'e pr\`es du point $a$. Alors
$f^* \Theta$ est une mesure positive lisse de masse $\l$.
Comme $X$ est k\"ahl\'erienne, le $dd^c$-lemma assure
l'existence d'une forme lisse $S$ de bidegr\'e $(k-1,k-1)$
telle que
\begin{equation}
\frac{1}{\l}f^* \Theta=\Theta+dd^c S.
\end{equation}
Quitte \`a translater $S$ en lui ajoutant un multiple de $\om^{k-1}$,
on peut supposer que $0 \leq S \leq C \om^{k-1}$. Comme $f$ est $k$-stable,
on peut prendre l'image inverse de (4.1) par $f^n$ et utiliser
les identit\'es $(f^n)^*=(f^*)^n$ pour obtenir
$$
\frac{1}{\l^n}(f^n)^* \Theta=\Theta+dd^c S_n, \text{ o\`u } 
S_n:=\sum_{j=0}^{n-1} \frac{1}{\l^j}(f^j)^* S.
$$
Les courants $S_n$ forment une suite croissante de courants positifs
-car $S \geq 0$ donc $(f^j)^*S \geq 0$- dont la masse
est born\'ee par $\sum_{j \geq 0} \d_{k-1}(f^j,\om)/\l^j$ qui converge
car $\l>\l_{k-1}(f)$. On en d\'eduit la convergence
$$
\frac{1}{\l^n}(f^n)^* \Theta \longrightarrow \mu_f:=\Theta+dd^c S_{\infty}, \text{ o\`u } 
S_{\infty}:=\sum_{j \geq 0} \frac{1}{\l^j}(f^j)^* S.
$$

Si $\Theta'$ est une autre mesure de probabilit\'e lisse, alors
il existe une forme lisse $R$ de bidegr\'e $(k-1,k-1)$ telle que
$\Theta'=\Theta+dd^c R$. On en d\'eduit que
$\l^{-n}(f^n)^* \Theta' \rightarrow \mu_f$, car 
les courants $(f^n)^* R$ sont de masse major\'ee par $\d_{k-1}(f^n,\om)$.

Soit $\f$ une fonction quasiplurisousharmonique (qpsh) sur $X$.
C'est une fonction major\'ee dont la courbure est minor\'ee
par une forme lisse. Quitte \`a translater et dilater $\f$, on peut 
supposer $\f \leq 0$ et $dd^c \f \geq -\om$. On a alors
$$
0 \leq \int_X (-\f) d\mu_f =\int_X (-\f) \Theta+\int_X (-\f)dd^c S_{\infty}
\leq \int_X (-\f) \Theta+\int_X  \om \wedge S_{\infty},
$$
en int\'egrant par parties dans la deuxi\`eme int\'egrale, et en observant
que $-dd^c \f \wedge S_{\infty} \leq \om \wedge S_{\infty}$ car
le courant $S_{\infty}$ est positif. On peut justifier ces int\'egrations
par parties en approximant $\f$ par une suite d\'ecroissante de fonctions
lisses qpsh (voir [G 5], [GZ 1]). L'assertion (1) en r\'esulte.

En particulier $\mu_f$ ne charge pas l'ensemble d'ind\'etermination, ni
les valeurs critiques de $f$. Comme elle est limite de mesures
$\mu_n$ qui v\'erifient $f^*\mu_n=\l \mu_{n+1}$, on en d\'eduit
que $\mu_f$ est invariante et de jacobien constant. Il r\'esulte
de la Proposition 2.9 et de la majoration 
$h_{top}(f) \leq \max_j \log \l_j(f)$ que
$\mu_f$ est d'entropie maximale,
$$
h_{\mu_f}(f)=h_{top}(f)=\log \l >0.
$$

Il reste \`a montrer le m\'elange et \`a estimer la d\'ecroissance des
corr\'elations. C'est la m\^eme id\'ee que dans la preuve du Th\'eor\`eme 1.5,
mais les d\'etails techniques sont l\'eg\`erement plus \'elabor\'es.
Soit $\chi,\p$ des fonctions test. Posons
$c_{\chi}=\int \chi d\mu_f$ et $c_{\p}=\int \p d\mu_f$.
Quitte \`a translater et dilater $\chi$ on peut supposer
$\chi \geq 0$ et $c_{\chi}=1$. La mesure $\chi \mu_f$ est donc
une mesure de probabilit\'e et il s'agit de montrer
que $\chi \circ f^n \mu_f=\l^{-n}(f^n)^* (\chi \mu_f)$ converge
vers $\mu_f$, et d'estimer \`a quelle vitesse. L'id\'ee est la suivante:
comme $\mu_f$ et $\chi \mu_f$ sont cohomologues, on peut trouver
$R_{\chi}$, un courant de bidegr\'e $(k-1,k-1)$ qui d\'epend contin\^ument
de $\chi$ tel que $\chi \mu_f=\mu_f+dd^c R_{\chi}$. En prenant l'image inverse
par $f^n$ et en int\'egrant contre $\p$, on obtient ainsi
$$
I_n(\chi,\p)=\left| \langle dd^c \p, \l^{-n}(f^n)^* R_{\chi}\rangle \right| \leq
||\p||_{{\mathcal C}^2} ||R_{\chi}|| \d_{k-1}(f^n,\om) \l^{-n},
$$
o\`u $I_n(\chi,\p):=\left| \int \p \chi \circ f^n d\mu_f -c_{\chi} c_{\p}\right|$.

Il nous reste \`a expliciter un tel courant $R_{\chi}$, car celui-ci n'est
pas du tout unique, contrairement au cas de la dimension 1.
La th\'eorie de Hodge fournit des op\'erateurs 
$\partial^*,\overline{\partial}^*,G$ qui permettent de r\'esoudre 
canoniquement l'op\'erateur
$dd^c$ (voir [GH]). Dans la pratique, il s'agit d'int\'egrer contre
un noyau $K(x,y)$ tel que $dd^c K=[\Delta]-\theta$, o\`u
$[\Delta]$ d\'esigne le courant d'int\'egration sur la diagonale de $X^2$,
et $\theta$ est une forme diff\'erentielle lisse ferm\'ee cohomologue 
\`a $[\Delta]$. Nous renvoyons le lecteur \`a [BGS] pour une \'etude
syst\'ematique de ce type de noyau, et \`a la preuve de la Proposition 2.1 de [DS 4]
pour une application dans un contexte dynamique: il y est montr\'e
notamment qu'on peut d\'ecomposer le noyau $K=K_1-K_2$ en une diff\'erence de
formes positives, $K_i \geq 0$. Le courant $R_{\chi}$ est obtenu
en int\'egrant contre $K$,
$$
R_{\chi}(x)=\int_{y \in X} K(x,y) \wedge [\chi(y) \mu_f(y)-\mu_f(y)].
$$
Il se d\'ecompose en $R_{\chi}=R_1-R_2$ en suivant la d\'ecomposition $K=K_1-K_2$.
On obtient ainsi l'encadrement
$$
-\int K_i(x,y) \wedge \mu_f(y) \leq R_i(x) 
\leq ||\chi||_{L^{\infty}} \int K_i(x,y) \wedge \mu_f(y).
$$
L'estimation de la vitesse de m\'elange r\'esulte alors de
$|| R_{\chi}|| \leq C ||\chi||_{L^{\infty}}$.
\end{preuve}

\begin{rqes}
Lorsque $f$ est un endomorphisme holomorphe de $\P^k$, ce r\'esultat
est d\^u  \`a J.-E.Fornaess et N.Sibony [FS 3].
La m\'ethode que nous proposons est int\'eressante \'egalement dans ce cas:
elle \'evite le recours \`a des propri\'et\'es fines 
et techniques de th\'eorie du pluripotentiel.

L'estimation $lov(f) \leq \max \log \l_j(f)$ est 
d\'emontr\'ee dans [Fr] lorsque $f$ est holomorphe,
dans [G 5] lorsque $f$ est m\'eromorphe et $\dim_{\C} X \leq 3$
(ou lorsque $X$ est homog\`ene, par exemple pour $X=\P^k$);
le cas g\'en\'eral est d\'emontr\'e par 
T.-C.Dinh et N.Sibony dans [DS 3,4].

T.-C.Dinh et N.Sibony donnent dans [DS 1,8,9] une preuve
diff\'erente de l'estimation de la vitesse de m\'elange.
Ce contr\^ole pr\'ecis de la d\'ecroissance des coefficients de
corr\'elation implique, via le Th\'eor\`eme de Gordin-Liverani,
un Th\'eor\`eme central limite comme l'ont observ\'e
S.Cantat et S.Leborgne [CL] (lorsque $f$ est holomorphe),
et T.C.Dinh et N.Sibony [DS 9] (cas m\'eromorphe).

Notons enfin qu'on peut montrer par les m\^emes m\'ethodes que
la mesure $\mu_f$ m\'elange \`a tout ordre, et qu'elle est m\^eme
K-m\'elangeante (voir [CFS] pour les d\'efinitions et [DS 1]
pour une preuve). Il est probable que $(f,X)$ soit en fait
conjugu\'e au d\'ecalage de Bernoulli sur $\l$ symboles.
Cela a \'et\'e d\'emontr\'e par D.Heicklen, C.Hoffman [HH]
et J.-Y.Briend [Br] lorsque $f$ est un endomorphisme 
holomorphe de $\P^k$ (voir \'egalement [Buz]).
\end{rqes}

\section{Propri\'et\'es ergodiques de $\mu_f$}

Soit, comme pr\'ec\'edemment, $f:X \rightarrow X$ un endomorphisme dont
le degr\'e topologique $\l_k(f)$ domine 
strictement tous les autres degr\'es dynamiques.

\begin{thm} Sous les hypoth\`eses pr\'ec\'edentes,
\begin{enumerate}
\item la mesure $\mu_f$ est l'unique mesure d'entropie maximale;
\item elle est hyperbolique, ses exposants de Lyapunov v\'erifient
$$
\chi_1 \geq \cdots \geq \chi_k \geq \frac{1}{2} \log \l_k(f)/\l_{k-1}(f)>0;
$$
\item les points p\'eriodiques r\'epulsifs s'\'equidistribuent selon $\mu_f$.
\end{enumerate}
\end{thm}

Autrement dit la conjecture est v\'erifi\'ee lorsque $l=k$.
Ce r\'esultat est d\^u \`a J.-Y.Briend et J.Duval lorsque $f$
est un endomorphisme holomorphe de $\P^k$ [BrD 1,2].
Le cas m\'eromorphe est trait\'e dans [G 5], [DS 3,4].

L'id\'ee centrale de la d\'emonstration est de 
construire et contr\^oler beaucoup de branches inverses
de $f^n$, en suivant une m\'ethode initi\'ee
par M.Lyubich en dimension 1 [Ly] 
(voir section 1.4) et g\'en\'eralis\'ee par
J.-Y.Briend et J.Duval [BrD 2] en dimension sup\'erieure. 
Le lemme clef est le suivant:

\begin{lem}
Soit $V_l=\cup_{j=1}^l f^j({\mathcal C}_f)$, o\`u ${\mathcal C}_f$ d\'esigne 
l'ensemble critique de $f$. Fixons $\e>0$ et $1 <\d <\l/\l_{k-1}(f)$.
Alors il existe $l \in \N$ et $C>1$ tels que pour toute boule $\overline{B}$
qui \'evite $V_l$, on peut construire au moins
$(1-\e) \l^n$ branches inverses $f_i^{-n}$ de $f^n$ sur $B$, $n \geq l$, avec
$$
\text{diam}(f_i^{-n} B) \leq C \d^{-n/2}.
$$
\end{lem}

La m\'ethode de d\'emonstration est tr\`es proche du cas des endomorphismes
holomorphes trait\'e par J.-Y.Briend et J.Duval. Nous l'avons
esquiss\'ee dans la section 1.4 dans le cas de la dimension 1
en utilisant des arguments qui passent en dimension sup\'erieure.
Nous renvoyons le lecteur \`a l'article original [BrD 2] ainsi
qu'\`a l'expos\'e de X.Buff au s\'eminaire Bourbaki [Buf 2] pour
plus de d\'etails sur le cas holomorphe.

La pr\'esence de points d'ind\'etermination rend les d\'etails
techniques plus compliqu\'es dans le cas m\'eromorphe. 
Il est par exemple d\'elicat d'estimer
certaines int\'egrales qui, dans le cas holomorphe, se calculent
directement en cohomologie
(voir Lemme 1.5 dans [G5] pour une illustration).  
Nous renvoyons le lecteur \`a [G 5], [DS 1,3,4] pour les d\'etails
et nous nous contentons de donner ici la preuve du Lemme 3.4.

\begin{preuve}
Fixons $\e>0$ et $1<\d<\l/\l_{k-1}(f)$. Fixons \'egalement $l=l_{\e} >>1$
(la valeur sera pr\'ecis\'ee plus loin) et $B=B(p,r)$ une boule t.q.
$\overline{B} \cap V_l=\emptyset$.

\noindent {\bf Construction des branches inverses.} 
Nous allons construire des branches inverses $f_i^{-n}$ de $f^n$ sur $B$
par r\'ecurrence. 

 {\it Etape $l$ (initialisation).} Comme $B$ ne rencontre pas les
valeurs critiques de $f^l$, il y a $\l^l$ branches inverses $f_i^{-l}$ bien
d\'efinies de $f^l$ sur $B$. On note $B_i^{-l}:=f_i^{-l}B$ leurs images.

 {\it Etape $l+1$.} Si $\overline{B}_i^{-l}$ ne rencontre pas les valeurs critiques $V_1$
de $f$, on peut d\'efinir $\l$ nouvelles branches inverses de $f$ sur $B_i^{-l}$.
Lorsque $\overline{B}_i^{-l}$ rencontre $V_1$, il faut mesurer de quelle fa\c{c}on:
on diminue pour cela l\'eg\`erement le rayon de $B$ et on poursuit la
contruction des branches inverses si $f_i^{-l}(\r_{l+1} \overline{B}) \cap V_1=\emptyset$,
o\`u $\r_n:=1-\sum_{j=l}^{n} j^{-2}$. Lorsque
$f_i^{-l}(\r_{l+1} \overline{B}) \cap V_1 \neq \emptyset$, il vient 
$f^l(B_i^{-l} \cap V_1) \cap \r_{l+1}B \neq \emptyset$, c'est \`a dire que l'ensemble
analytique $Z_l:=f^l(B_i^{-l} \cap V_1)$ p\'en\`etre de fa\c{c}on consistante \`a
l'int\'erieur de la boule $B$. Fixons $z_l \in Z_l$ tel que $B(z_l,l^{-2} r) \subset B$.
Alors $B(z_l,l^{-2} r)  \cap Z_l$ est un sous-ensemble analytique sans bord de
dimension $s=\dim_{\C} V_1$ de 
$B(z_l,l^{-2} r)$. Comme $Z_l$ a un nombre de Lelong positif au point $z_l$, il vient
$$
C_1 (r l^{-2})^{2s} \leq \int_{B(z_l,l^{-2} r)} [Z_l] \wedge \om^s 
\leq \int [f^l B_i^{-l} \cap V_1] \wedge \om^s,
$$
pour une constante $C_1>0$ uniforme. Or
$$
\sum_i \int [f^l B_i^{-l} \cap V_1] \wedge \om^s \leq \int (f^l)_* [V_1] \wedge \om^s
\leq C_2 [\l_s(f)+\e']^l,
$$
o\`u $\e>0$ peut \^etre choisi arbitrairement petit et $C_2>0$ est une constante
qui ne d\'epend que de $\e'>0$. Observons que $\l_s(f) \leq \l_{k-1}(f)$ (cf Th\'eor\`eme 2.4.a).
Dans la suite on ajuste la valeur de $\e'$ de sorte que 
$\d=[\l_{k-1}(f)+\e']/\l >1$. il r\'esulte des in\'egalit\'es pr\'ec\'edentes que
$$
\sharp \{ i \, / \, f_i^{-l} (\r_{l+1} \overline{B}) \cap V_1 \neq \emptyset \}
\leq C_3 l^{4(k-1)} \d^{-l} \l^l.
$$
Si l est choisi suffisamment grand, le nombre des pr\'eimages $B_i^{-l}$ qui rencontrent
$V_1$ de fa\c{c}on consistante est donc n\'egligeable par rapport au nombre
total $\l^l$ de pr\'eimages. On a donc fabriqu\'e
$\l^{l+1} [1-C_3 l^{4(k-1)} \d^{-l}]$ pr\'eimages de $f^{l+1}$ sur $\r_{l+1} B$.

 {\it Etape n.} On proc\`ede ensuite par r\'ecurrence pour construire des pr\'eimages 
$f_i^{-n}$ de $f^n$ sur $\r_n B$, i.e. en diminuant \`a chaque \'etape 
l\'eg\`erement le rayon
de la boule $B(p,r)$  de d\'epart. A l'\'etape $n$, on a donc construit au moins
$$
d_n:=\l^n \left[ 1-C_3 \sum_{j=l}^{n-1} l^{4(k-1)} \d^{-l} \right] \geq \l^n (1-\e/2)
$$
branches inverses $f_i^{-n}$ qui sont toutes bien d\'efinies sur la boule
$B'=B(p,r')$ de rayon
$$
0<r'=r\left[1-\sum_{j \geq l} j^{-2} \right] \leq \r_n r,
$$
arbitrairement proche de $r$, si on choisit initialement $l$
suffisamment grand.
\vskip.1cm

\noindent {\bf Diam\`etre des images.} 
Nous souhaitons \`a pr\'esent estimer le diam\`etre des pr\'eimages de la boule
$B'=B(p,r')$. Pour cela nous allons trancher $B'$ par des disques holomorphes locaux
$\Delta_{\theta}$ passant par $p$ et estimer les diam\`etres 
$diam f_i^{-n}(\Delta_{\theta})$ pour beaucoup de $\theta$.

Soit $\om'$ une $(k-1,k-1)$-forme positive ferm\'ee \`a coefficients dans $L^1(X)$
qui est lisse dans $X \setminus \{p\}$ et telle que 
$$
\om'=\int_{\theta \in \P^{k-1}} [\Delta_{\theta}] d\nu(\theta)
\; \text{ dans } B',
$$
o\`u $\nu$ d\'esigne la mesure de Fubini-Study sur l'ensemble $\P^{k-1}$ des
droites locales passant par $p$. On peut construire une telle forme
$\om'$ en consid\'erant l'\'eclatement de $X$ au point $p$
ou, \`a la main, en posant
$$
\om'=\left[ A \om_1+dd^c (\chi \log dist(\cdot,p) \right]^{k-1}.
$$
Ici $\chi$ d\'esigne une fonction de troncature telle que 
$\chi \equiv 1$ dans $2B'$ et 
$\chi \log dist(\cdot,p) \in {\mathcal C}^{\infty}(X \setminus \{p\})$, et
$\om_1 \geq 0$ est une $(1,1)$-forme lisse ferm\'ee telle que
$\om_1 \equiv 0$ dans $B'$ et $\om_1>0$ hors de $\overline{B'}$. On choisit
$A>>1$ de sorte que la positivit\'e de $A\om_1$ compense la n\'egativit\'e
de $dd^c [ \chi \log dist(\cdot,p) ]$. Observons que
$$
0 \leq \sum_i \int_{B'} (f_i^{-n})_* \om' \wedge \om \leq
\int_X (f^n)^* \om' \wedge \om \leq C_4 \d^{-n} \l^n.
$$
Il s'ensuit qu'au moins $(1-\e) \l^n$ des branches inverses $f_i^{-n}$
sont telles que
$$
\int \text{Aire}(f_i^{-n} \Delta_{\theta}) d\nu(\theta)=
\int_{B'} (f_i^{-n})_* \om' \wedge \om \leq \frac{2C_4}{\e} \d^{-n}.
$$
On note $I_{\e}^n$ l'ensemble de ces branches inverses et, pour $i \in I_{\e}^n$,
$$
A_i^n:=\left\{ \theta \in \P^{k-1} \, / \, 
\text{Aire}(f_i^{-n} \Delta_{\theta}) \leq \frac{4C_4}{\e} \d^{-n} \right\}.
$$
J.Y.Briend et J.Duval montrent dans [BrD 2] que l'estimation d'aire implique,
quitte \`a l\'eg\`erement diminuer la taille de la boule $B'$,
un contr\^ole de diam\`etre,
$$
\text{diam}(f_i^{-n} \Delta_{\theta}) \leq C_5 \d^{-n/2}, \, \theta \in A_i^n,
$$
pour une constante $C_5>0$ ind\'ependante de $n$.
Or $\nu(A_i^n) \geq 1/2$, en particulier $A_i^n$ est de capacit\'e projective $\geq 1/2$.
Un r\'esultat de N.Sibony et P.M.Wong [SW] implique alors
$$
\text{diam}\left( f_i^{-n} \frac{1}{2}\Delta_{\theta}\right) \leq C_5 \d^{-n/2}, \, 
\text{ pour tout } \theta \in \P^{k-1}.
$$
L'estimation du diam\`etre de $f_i^{-n} B''$ en r\'esulte, $B''=B'/2$.
Nous renvoyons le lecteur \`a [GZ 1] (Th\'eor\`eme 6.3) pour la d\'efinition
et quelques propri\'et\'es 
de la capacit\'e projective.
\end{preuve}

\begin{rqe}
La preuve ci-dessus est une version l\'eg\`erement modifi\'ee de celle du Lemme 3.3 
de [G 5]: dans cette derni\`ere on supposait la vari\'et\'e ambiante projective,
ce qui permet, comme dans [BrD 2],
l'utilisation du Th\'eor\`eme de Bezout. On peut en fait
s'affranchir de l'hypoth\`ese de projectivit\'e comme 
nous l'avons indiqu\'e ci-dessus. Cela avait d\'ej\`a
\'et\'e observ\'e par T.C.Dinh et N.Sibony dans le cadre des applications d'allure
polynomiale [DS 1].
\end{rqe}

\section{Exemples}

\subsection{Endomorphismes holomorphes}

Il y a beaucoup d'endomorphismes {\it holomorphes} $f$ 
sur l'espace projectif complexe $\P^k$. Ils s'\'ecrivent en coordonn\'ees
homog\`enes $f=[P_0:\cdots:P_k]$, o\`u les $P_j$ sont des polyn\^omes homog\`enes
premiers entre eux de m\^eme degr\'e $d \in \N^*$, tels que
$\cap_{j=0}^k (P_j=0)=\{0\}$. Ils v\'erifient $\l_j(f)=d^j$, $0 \leq j \leq k$ et
sont donc de grand degr\'e topologique d\'es que $d \geq 2$.

Le but de cette section est de montrer que ce sont essentiellement 
les seuls exemples ``int\'eressants''
d'endomorphismes holomorphes de grand degr\'e topologique
sur les surfaces k\"ahl\'eriennes compactes.

\begin{thm}
Soit $X$ une surface k\"ahl\'erienne compacte, et
$f:X \rightarrow X$ un endomorphisme holomorphe tel que $\l_1(f) \neq \l_2(f)$.
Alors on est dans l'une des situations suivantes:
\begin{itemize}
\item soit $X=\P^2$;
\item soit $X$ est un tore;
\item soit $f$ est un automorphisme d'entropie positive;
\item soit $f$ pr\'eserve une fibration et $\l_2(f)>\l_1(f)$.
\end{itemize}
\end{thm}

La preuve de ce r\'esultat nous permettra d'\^etre beaucoup plus pr\'ecis.
Nous montrerons par exemple que dans le dernier cas, 
\begin{itemize}
\item soit $X$ est une surface hyperelliptique
et $f$ pr\'eserve la fibration elliptique d'Albanese (de plus que $f$ se 
rel\`eve en un endomorphisme holomorphe sur le produit des courbes elliptiques d\'efinissant
$X$);
\item soit $X$ est rationnelle et $f$ provient d'un endomorphisme sur un mod\`ele minimal; dans ce
cas on obtient $\l_2(f)=\l_1(f)^2$ si $X \neq \P^1 \times \P^1$;
\item soit $X$ est une surface r\'egl\'ee au dessus d'une courbe elliptique et $f$ pr\'eserve
la fibration rationnelle.
\end{itemize}

\begin{cor}
Les seuls endomorphismes holomorphes
tels que $\l_1(f) >\l_2(f) \geq 2$ sont ceux des tores.
\end{cor}

On peut par contre fabriquer de nombreux exemples d'endomorphismes m\'eromorphes tels
$\l_1(f) \neq \l_2(f)$ sur les surfaces rationnelles et les surfaces de dimension
de Kodaira nulle.

\begin{rqes}
Les automorphismes d'entropie positive 
des surfaces ont \'et\'e partiellement classifi\'es par S.Cantat [Ca 1]
(voir \'egalement [M 1,3], [BK 1] pour des exemples dynamiquement int\'eressants).

Les endomorphismes holomorphes
non inversibles des surfaces ont \'et\'e \'etudi\'es par de nombreux auteurs
(voir par exemple [Fu], [Na]). Sur les surfaces rationnelles minimales, on peut les d\'ecrire
\`a l'aide de coordonn\'ees (bi)homog\`enes (voir Proposition 3.3 dans [G 1]). 
Ils v\'erifient tous $\l_2(f)=\l_1(f)^2$, \`a l'exception d'un produit direct
de deux fractions rationnelles de degr\'es distincts sur $\P^1 \times\P^1$.
Si le mod\`ele minimal est diff\'erent de $\P^2$, l'endomorphisme (ou son carr\'e) pr\'eserve
une fibration rationnelle: c'est donc un produit crois\'e, dont la dynamique est
une version fibr\'ee de la dynamique unidimensionnelle. Il en est de m\^eme d'un endomorphisme
dans un \'eclat\'e $X$ de $\P^2$: pour que le relev\'e de $f:\P^2 \rightarrow \P^2$
soit un endomorphisme de $X$, il faut que $f$ pr\'eserve le pinceau de droites issues
du point auquel on \'eclate. En conclusion, les seuls endomorphismes 
holomorphes non inversibles
v\'eritablement int\'eressants sont ceux de $\P^2$ !
\end{rqes}

\begin{preuve.3.6}
Nous savons par le Th\'eor\`eme 2.14 qu'il suffit de s'int\'eresser aux cas
$kod(X)=0$ et $X$ rationnelle.
Le cas des automorphismes d'entropie positive ($1=\l_2(f)<\l_1(f)$) a \'et\'e \'etudi\'e
par S.Cantat, nous renvoyons le lecteur \`a son article [Ca 1].
Nous supposons donc dans la suite $\l_2(f) \geq 2$. Nous allons montrer que
\begin{itemize}
\item si $kod(X)=0$, alors $X$ est un tore;
\item si $X$ est rationnelle et $\pi:X \rightarrow \Sigma$ 
est un morphisme birationnel holomorphe vers
le mod\`ele minimal $\Sigma$ de $X$, alors il existe $N \in \N^*$ et $g:\Sigma \rightarrow \Sigma$
un endomorphisme holomorphe tel que $g \circ \pi=\pi \circ f^N$.
\end{itemize}
L'id\'ee consiste \`a analyser l'action de $f$ sur les courbes d'autointersection n\'egative
et d'utiliser la formule de Hurwitz $f^*K_X+R_f=K_X$, ainsi que la formule
d'aller-retour qui prend ici une forme particuli\`erement simple
$$
f_*f^*{\mathcal C}=\l_2(f) {\mathcal C}.
$$
Il r\'esulte par exemple de cette formule que les op\'erateurs $f_*,f^*$ sont des isomorphismes
sur l'espace de N\'eron-Severi r\'eel $NS(X,\R)=NS(X) \otimes \R$.

\noindent {\bf Etape 1}. On va montrer qu'il existe sur $X$ un nombre fini de courbes irr\'eductibles
d'autointersection n\'egative, et qu'elles sont totalement invariantes par un it\'er\'e de $f$.
Soit ${\mathcal C}$ une telle courbe et ${\mathcal C}'=f({\mathcal C})$.
Alors $f_*{\mathcal C}=a{\mathcal C}'$ pour un entier $a \in\N^*$. Supposons qu'une
autre courbe $\tilde{{\mathcal C}}$ est envoy\'ee sur ${\mathcal C}'$ par $f$. Alors
$f_*\tilde{{\mathcal C}}=\tilde{a}{\mathcal C}'$, donc
$$
a \tilde{{\mathcal C}} =\tilde{a} {\mathcal C}=:D
$$
puisque $f_*$ est un isomorphisme. Le fibr\'e en droites associ\'e au diviseur effectif $D$
\'etant d'autointersection n\'egative, il n'admet qu'une section holomorphe, donc
${\mathcal C}=\tilde{{\mathcal C}}$. On en d\'eduit que $f^* {\mathcal C}'=b{\mathcal C}$
pour un entier $b \in \N^*$ tel que $ab=\l_2(f)$. Notons
que $({\mathcal C}')^2=b{\mathcal C}^2/a<0$. 

Observons que $b \geq 2$ dans la formule pr\'ec\'edente si et seulement si
${\mathcal C}$ est incluse dans le lieu de ramification de $f$. Lorsque $b=1$, il vient
$a=\l_2(f)$ et donc
$$
\left| ({\mathcal C}')^2 \right|=\frac{1}{\l_2(f)} \left| {\mathcal C}^2 \right| 
<\left| {\mathcal C}^2 \right|.
$$
Il s'ensuit que $f^{m_0}({\mathcal C}) \subset R_f$ pour $m_0$ assez grand.

On en d\'eduit qu'il existe un nombre fini de courbes d'autointersection n\'egative. En effet
une telle courbe est pr\'ep\'eriodique, donc p\'eriodique puisque $f_*$ agit injectivement.
Notons
$$
m_{{\mathcal C}}=\inf \{m \geq 1 \, / \, f^m({\mathcal C})={\mathcal C} \}
\text{ et }
N:=\Pi_{{\mathcal C}^2<0, \, {\mathcal C} \subset R_f} m_{{\mathcal C}}.
$$
Alors $f^{N}{\mathcal C}={\mathcal C}$ pour toute courbe irr\'eductible ${\mathcal C} \subset R_f$
telle que ${\mathcal C}^2<0$. Si ${\mathcal C}$ est une courbe irr\'eductible d'autointersection
n\'egative qui n'est pas incluse dans $R_f$, on a $f^m {\mathcal C} \subset R_f$
pour un entier $m<N$, donc $f^{N}{f^m\mathcal C}=f^m{\mathcal C}$, donc
$f^{N} {\mathcal C}={\mathcal C}$ par injectivit\'e de $f_*$.
On en d\'eduit que l'ensemble des courbes irr\'eductibles d'autointersection n\'egative
est fini, contenu dans le support du diviseur $\cup_{j=0}^{N-1}f^j(R_f)$.
Remarquons que $f_*^N$ est l'identit\'e sur cet ensemble.
\vskip.1cm

\noindent {\bf Etape 2}. Quitte \`a changer $f$ en $f^N$, on peut donc supposer que chaque
courbe d'autointersection n\'egative est totalement invariante. Si une courbe irr\'eductible
lisse ${\mathcal C} \simeq \P^1$
est d'autointersection $-1$, on peut la contracter. L'endomorphisme $f$ descend alors en un endomorphisme
holomorphe puisque ${\mathcal C}$ est totalement invariante. On se ram\`ene ainsi
\`a travailler sur le mod\`ele minimal de $X$. Inversement, pour pouvoir
\'eclater un point et obtenir un endomorphisme holomorphe dans l'\'eclat\'e, il faut
que le point soit totalement invariant.

Lorsque $kod(X)=0$, il est facile de v\'erifier qu'il n'y a pas de point 
totalement invariant et que
les endomorphismes holomorphes vivent donc sur le mod\`ele minimal. En effet,
supposons le contraire et soit $X$
un tel mod\`ele minimal, $f:X \rightarrow X$
un endomorphisme holomorphe de degr\'e topologique $\l_2(f) \geq 2$,
$\pi:Y \rightarrow X$ l'\'eclatement de $X$ en un point $p$, et
$g:Y \rightarrow Y$ l'endomorphisme induit par $f$ sur $Y$. 
On note $E=\pi^{-1}(p)$ le diviseur exceptionnel.
Rappelons
que $12 K_X=0$ puisque $kod(X)=0$ (voir [Bea])
et que $f$ est non ramifi\'e (voir section 2.4.2). Les formules de Hurwitz 
$12E=12R_{\pi}=12K_Y$, $12g^*K_Y+12R_g=12K_Y$ combin\'ees \`a la relation de commutation
$\pi \circ g=f \circ \pi$ et \`a la relation d'invariance 
$g^* E=\sqrt{\l_2(f)} E$ donnent donc
$$
R_g+(\sqrt{\l_2(f)}-1) E=0,
$$
ce qui est contradictoire.

Comme les surfaces $K3$ sont simplement connexes, elles n'admettent aucun endomorphisme
holomorphe de degr\'e topologique $\geq 2$. Tout endomorphisme holomorphe d'une surface
d'Enriques se rel\`eve en un endomorphisme holomorphe de la surface $K3$ associ\'ee, et
est donc un automorphisme. Finalement, les  endomorphismes holomorphes
non inversibles des surfaces de dimension de Kodaira nulle ne peuvent exister que sur les
tores et les surfaces hyperelliptiques. Il est facile de construire de tels exemples
dans les deux cas. Dans le cas des surfaces hyperelliptiques, l'endomorphisme
pr\'eserve la fibration d'Albanese: c'est un produit crois\'e qui v\'erifie
donc $\l_2(f) \geq \l_1(f)$ (cf Lemme 2.15). Les tores admettent toujours
des endomorphismes de grand degr\'e topologique (via une homoth\'etie de rapport entier),
et parfois des endomorphismes de petit degr\'e topologique (voir Exemple 2.17).
\vskip.1cm

\noindent {\bf Etape 3}.
Il reste \`a traiter le cas des surfaces rationnelles. D'apr\`es ce qui pr\'ec\`ede, on peut,
lorsque $\l_2(f) \geq 2$ se ramener au cas d'un endomorphisme sur un mod\`ele minimal, i.e.
$\P^2$ ou une surface de Hirzebruch $\F_n$, $n \in \N \setminus \{1\}$. Le cas de
$\P^2$ est bien connu, on obtient $\l_2(f)=\l_1(f)^2>\l_1(f)$.
Celui de $\F_n$, $n \geq 2$ est d\'ecrit dans [G 1]: on obtient \'egalement
$\l_2(f)=\l_1(f)^2$; dans ce cas $f$ pr\'eserve la fibration rationnelle. Le cas
de $\F_0=\P^1 \times \P^1$ est d\'ecrit dans [FaG]: soit $f$ est un produit direct
d'endomorphismes $f_1,f_2$ de $\P^1$ de degr\'e $d_1,d_2$ auquel cas
$\l_1(f)=\max(d_1,d_2)$ et $\l_2(f)=d_1 d_2 \geq \l_1(f)$ (on ne peut alors
obtenir un endomorphisme holomorphe dans un \'eclat\'e que si $d_1=d_2$,
donc $\l_2(f)=\l_1(f)^2$); soit
$f$ est la compos\'ee d'un endomorphisme du type pr\'ec\'edent
et de l'involution $(x,y) \mapsto (y,x)$, dans ce cas
$\l_1(f)=\sqrt{d_1d_2}$ et $\l_2(f)=d_1d_2=\l_1(f)^2$.
\end{preuve.3.6}
\vskip.2cm

L'\'etude des endomorphismes holomorphes non inversibles fait l'objet
de plusieurs travaux r\'ecents (voir [PS], [ARV], [Bea 2], [Fu], [A], [AC 2]). 
On s'attend \`a ce qu'il y en ait peu qui soient cohomologiquement hyperboliques,
hormis ceux de l'espace projectif complexe.
Nous renvoyons le lecteur int\'eress\'e
\`a l'article de S.Cantat [Ca 2] (ainsi qu'aux r\'ef\'erences qu'il contient),
pour une justification de cette attente dans le cas des vari\'et\'es homog\`enes.
Nous \'etudierons la dynamique des endomorphismes holomorphes {\it inversibles}
(automorphismes) dans la section 5.2.

\subsection{Endomorphismes polynomiaux de $\C^k$}

Les endomorphismes polynomiaux de $\C^k$ de la forme
$$
f(z_1,\ldots,z_k)=(P_1(z_1),P_2(z_1,z_2),\ldots,P_k(z_1,\ldots,z_k)),
$$
o\`u les $P_j$ sont des polyn\^omes de degr\'e $d_j \geq 2$ en $z_j$ s'\'etendent
en des endomorphismes m\'eromorphes de $\P^k$ tels que
$$
\l_l(f)=\max_{i_1<\cdots<i_l} d_{i_1} \cdots d_{i_l};
\; \text{ en particulier } \l_k(f)>\max_{j \leq k-1} \l_j(f).
$$
Ces endomorphismes ont \'et\'e \'etudi\'es par de nombreux auteurs (voir
par exemple [ He], [Se], [J 1,2], [FaG]), notamment en dimension deux:
l'analyse de leur dynamique est en effet simplifi\'ee par le fait
que ces {\it produits crois\'es} pr\'eservent les feuilletages lin\'eaires
$\{ z_1=c_1,\ldots, z_j=c_j\}$. 

Nous montrons en particulier dans [FaG] que la mesure d'entropie maximale
$\mu_f$ peut \^etre exprim\'ee comme un produit ext\'erieur
$$
\mu_f=dd^c G_1 \wedge \cdots \wedge dd^c G_k,
$$
o\`u $G_1$ est une fonction plurisousharmonique continue dans $\C^k$ telle
que $G_1 \circ f=\l_1 G_1$ (la fonction de Green dynamique de $f$, introduite
par J.E.Fornaess et N.Sibony dans [FS 2,3,4]), $G_2$ est une fonction
de Green partielle d\'efinie uniquement sur le support de $dd^c G_1$ qui
v\'erifie $G_2 \circ f=\frac{\l_2}{\l_1} G_2$, etc.
L'int\'er\^et de cette construction r\'eside dans le fait que l'on sait 
contr\^oler le module de continuit\'e des fonctions $G_i$
(dans l'esprit de la Proposition 1.2, voir [DG]) et que cela donne 
par exemple des
informations sur la dimension de Hausdorff du support de la mesure $\mu_f$.
Lorsque $f:\C^k \rightarrow \C^k$ s'\'etend en un endomorphisme
holomorphe de $\P^k$, on obtient en fait $G_1=\cdots =G_k$ et la construction
est due \`a J.E.Fornaess, N.Sibony [FS 2] et J.H.Hubbard, P.Papadopol [HP 1].
Lorsque l'extension de $f$ a des points d'ind\'etermination \`a l'infini,
on obtient $(dd^c G_1)^k=0$ dans $\C^k$ et il devient n\'ecessaire de construire des 
fonctions de Green partielles. 

La construction de fonctions de Green partielles a \'et\'e \'etendue
dans [G 1] \`a certains endomorphismes polynomiaux de $\C^2$ qui ne sont pas des
produits crois\'es (voir Exemple 4.1),
et cette approche a \'et\'e g\'en\'eralis\'ee en dimension sup\'erieure
dans [DS 2].

Les applications concern\'ees sont cependant d'un type 
particulier: il faut un contr\^ole pr\'ecis de la croissance de $f$
sur le support du courant $dd^c G_1$, ce qui n'est pas toujours possible
(voir [DiDS]). On s'en rend d\'ej\`a compte en consid\'erant les endomorphismes
polynomiaux $f=(P_1,P_2)$ {\it quadratiques} de $\C^2$, i.e. ceux
tels que $\max(\deg P_1,\deg P_2)=2$.
Ils sont classifi\'es dans [G 3], \`a conjugaison pr\`es, en fonction de leurs
degr\'es dynamiques $(\l_1,\l_2)$: on obtient cinq familles
telles que $\l_2>\l_1$
et quatre possibilit\'es pour les couples de degr\'es, 
$$
(\sqrt{2},2), \; \left( \frac{1+\sqrt{5}}{2},2 \right), \;
(2,3), \; \text{ et } (2,4).
$$
Nous renvoyons le lecteur au Th\'eor\`eme 2.1 de [G 3] pour plus de d\'etails.

\subsection{La m\'ethode de Newton}

Si l'on souhaite trouver les racines du syst\`eme d'\'equations polynomiales
$$
w=P(z) \; \; \text{ et } \; \; z=Q(w),
$$
o\`u $P,Q$ sont des polyn\^omes de degr\'e $p,q \geq 2$, il est naturel d'utiliser
la m\'ethode de Newton qui consiste \`a it\'erer l'application rationnelle
$f$ d\'efinie pour $(z,w) \in \C^2$ par
$$
f \left( \begin{array}{c} z \\ w \end{array} \right)=
\left( \begin{array}{c} z \\ w \end{array} \right)
-\left[ \begin{array}{cc} P'(z) & -1 \\ -1 & Q'(w) \end{array} \right]^{-1}
\left( \begin{array}{c} P(z)-w \\ Q(w)-z \end{array} \right).
$$
C'est un probl\`eme mentionn\'e par J.E.Fornaess et N.Sibony dans [FS 4],
qui a \'et\'e \'etudi\'e par J.-H.Hubbard et P.Papadopol [HP 2].
Ces derniers montrent que $f$ induit un endomorphisme 1-stable de $\P^2$
tel que
$$
\l_1(f)=r_1(f)=p+q-1 \; \; \text{ et } \; \; 
\l_2(f)=pq>\l_1(f),
$$
et entament une \'etude syst\'ematique de la dynamique de $f$.

Le Th\'eor\`eme 3.1 fournit une mesure canonique $\mu_f$ d'entropie maximale
$=\log \l_2(f)$ qui ne charge pas les ensembles pluripolaires.
En particulier $\mu_f$ ne charge pas les points d'ind\'etermination: cela
r\'epond \`a une question d'A.Russsakovskii et B.Shiffman [RS].
On trouvera dans [HP 2] une preuve g\'eom\'etrique tr\`es \'el\'egante
de ce dernier point, due \`a A.Douady.

\subsection{Vari\'et\'es non rationnelles}

Rappelons qu'une vari\'et\'e projective $X$ de dimension $k$ est dite
unirationnelle s'il existe une application m\'eromorphe dominante
$\Phi:\P^k \rightarrow X$.
Lorsque $\dim_{\C} X \leq 2$, toutes les vari\'et\'es unirationnelles sont en fait
rationnelles (Th\'eor\`eme de Castelnuovo, voir e.g. corollaire V.5 dans [Bea]).
Ce n'est plus vrai en dimension $\geq 3$: V.A.Iskovskikh et Yu.A.Manin [IM] ont montr\'e
que les quartiques lisses
de $\P^4$ sont des vari\'et\'es unirationnelles dont le groupe 
d'endomorphismes birationnels est fini. Elles ne sont donc pas rationnelles. 
On peut v\'erifier qu'elles n'admettent pas d'endomorphisme holomorphe non inversible 
(voir[ARV]), bien qu'elles admettent de nombreux endomorphismes m\'eromorphes 
dynamiquement int\'eressants, comme le montre l'exemple suivant, d\^u \`a F.Campana, 
qui est mentionn\'e dans [DS 1].

\begin{exa}
Soit $X$ une vari\'et\'e projective: on peut la plonger dans un 
espace projectif $\P^N$ puis projeter sur un $k$-plan g\'en\'erique: autrement dit
on dispose toujours d'une application m\'eromorphe dominante
$p:X \rightarrow \P^k$. 
Supposons que $X$ est unirationnelle, on fixe
$\Phi:\P^k \rightarrow X$ une application m\'eromorphe dominante.
Soit alors $g:\P^k \rightarrow \P^k$
un endomorphisme rationnel. Il induit des endomorphismes m\'eromorphes sur $X$,
$$
f:=\Phi \circ g^j \circ p, \; \text{ o\`u } j \in \N.
$$
Il r\'esulte de la Proposition 2.6 qu'il existe $C_X \geq 1$ telle que pour tout $j \in \N$,
et pour tout $1 \leq l \leq k$,
$$
\l_l(f) \leq C_X^2 r_l(\Phi) r_l(p) r_l(g^j).
$$
Or pour $l=k$, on obtient bien s\^ur $\l_k(f)=\l_k(\Phi) \l_k(p) \l_k(g)^j$.
On en d\'eduit que si $g$ a un grand degr\'e topologique --i.e. si $\l_k(g)$
domine tous les autres degr\'es topologiques--, alors $\l_k(f)$ domine tous les
autres degr\'es dynamiques de $f$, pourvu que $j$ soit choisi assez grand.

\end{exa}

\section{Autres d\'eveloppements}

Nous concentrons nos efforts dans ce m\'emoire sur 
la construction d'une mesure d'entropie
maximale, et cherchons \`a \'etablir quelques unes
de ses propri\'et\'es ergodiques, en relation avec la conjecture
\'enonc\'ee dans l'introduction. 
Il y a bien entendu de nombreuses autres directions de recherche:
on peut s'int\'eresser \`a des propri\'et\'es plus fines
de cette m\^eme mesure (dimension, r\'egularit\'e),
ou \'etudier d'autres mesures qui rendent mieux compte de la g\'eom\'etrie
des ensembles de Julia (mesures conformes).
On peut \'egalement s'int\'eresser \`a la fa\c{c}on dont tous ces objets d\'ependent
d'un param\`etre. Nous passons en revue dans cette section
quelques unes des directions qui ont \'et\'e en partie explor\'ees.

\subsection{Dimension et exposants de Lyapunov}

Nous avons mentionn\'e \`a plusieurs reprises qu'il existe un lien entre la 
r\'egularit\'e des potentiels de la mesure $\mu_f$
et la dimension ponctuelle de celle-ci. Rappelons que la dimension
de Hausdorff de $\mu_f$ est
$$
\dim (\mu_f):=\inf \{ \dim_H(Y) \, / \, Y \text{ bor\'elien tel que } \mu_f(Y)=1 \},
$$
o\`u $\dim_H (Y)$ d\'esigne la dimension de Hausdorff du bor\'elien $Y$.
La dimension ponctuelle sup\'erieure de $\mu_f$ au point $p$ est
$$
\overline{\dim}(\mu_f,p):=\limsup_{r \rightarrow 0} \left[
\frac{\log \mu_f B(p,r)}{\log r} \right].
$$
On d\'efinit de fa\c{c}on analogue la dimension ponctuelle inf\'erieure
$\underline{\dim}(\mu_f,p)$. Nous renvoyons le lecteur au livre de Y.Pesin [Pe]
pour plus de renseignements sur ces notions. Il r\'esulte de la Proposition 1.2
et du corollaire 1.4 que, lorsque $k=1$,
$$
\underline{\dim}(\mu_f,p) \geq \frac{\log \l}{\chi_{top}(f)},
\text{ pour tout point } p \in X ;
$$
en particulier $\dim (\mu_f) \geq \log \l/\chi_{top}(f)$.
C'est le ``principe de distribution de masse'' (voir [Pe] p43).
On peut en fait raffiner ces estimations et obtenir la formule
de F.Ledrappier [L], A.Manning [Ma] et R.Ma\~ne [Mn];
$$
\dim (\mu_f)=\frac{\log \l}{\chi(\mu_f)}.
$$
Cette formule a \'et\'e partiellement g\'en\'eralis\'ee en
dimension sup\'erieure, dans le cadre des endomorphismes
de grand degr\'e topologique (voir [BiDeM], [DiD]). 
T.C.Dinh et C.Dupont [DiD] obtiennent notamment l'encadrement
\begin{equation}
\frac{\log \l_k(f)}{\chi_1(\mu_f)} \leq \dim (\mu_f) \leq
2k-\frac{\sum_{j=1}^k \chi_j(\mu_f)-\log \l_k(f)}{\chi_1(\mu_f)},
\end{equation}
o\`u $\chi_1(\mu_f) \geq \cdots \geq \chi_k(\mu_f)$ d\'esignent 
les exposants de Lyapunov de $\mu_f$.
\vskip.2cm

Nous avons vu (Th\'eor\`eme 3.3.2) que ceux-ci sont tous au moins \'egaux
\`a $\frac{1}{2} \log \l_k(f)/\l_{k-1}(f)$.
Il est naturel de s'int\'eresser au cas extr\'emal o\`u les 
exposants sont minimaux, i.e. lorsque
$$
\chi_1(\mu_f)=\cdots=\chi_k(\mu_f)=
\frac{1}{2} \log [ \l_k(f) / \l_{k-1}(f)].
$$
C'est un probl\`eme qui a \'et\'e \'etudi\'e par plusieurs 
auteurs (voir par exemple [L], [Z], [BeL], [BeD], [Ca 3]).
Les r\'esultats obtenus jusqu'ici nous incitent \`a poser
la question suivante:

\begin{ques}
Soit $f:X \rightarrow X$ un endomorphisme m\'eromorphe de grand degr\'e topologique,
$\l_k(f)>\max_{j \leq k-1} \l_j(f)$. Les propri\'et\'es suivantes sont-elles \'equivalentes ?
\begin{enumerate}
\item Les exposants de Lyapunov de $\mu_f$ sont tous \'egaux \`a 
$\frac{1}{2} \log \l_k(f)/\l_{k-1}(f)$.
\item La dimension de Hausdorff de $\mu_f$ est \'egale \`a $2k$ et
les degr\'es dynamiques v\'erifient $\l_j(f)=\l_1(f)^j$.
\item La mesure $\mu_f$ est absolument continue par rapport \`a la mesure de Lebesgue
et les degr\'es dynamiques v\'erifient $\l_j(f)=\l_1(f)^j$.
\item $f$ est un exemple de Latt\`es.
\end{enumerate}
\end{ques}

Un endomorphisme est un exemple de Latt\`es s'il existe un tore complexe
compact $A$ de dimension $k$, un groupe fini $G$ d'automorphismes de $A$ 
et une dilatation affine $F$ de $A$ tels que
\begin{enumerate}
\item $F$ passe au quotient en un endomorphisme de la vari\'et\'e $A/G$;
\item il existe une application birationnelle $\pi:X \rightarrow A/G$
t.q. $\pi \circ f=F \circ \pi$.
\end{enumerate}
Le lecteur trouvera dans [Mi 2] une classification compl\`ete des exemples
de Latt\`es en dimension 1. A.Zdunik [Z] a montr\'e dans ce contexte que
la r\'eponse \`a la question 3.10 est positive.
C.Dupont donne dans [Dup] une classification partielle des exemples de Latt\`es
en dimension deux. C.Dupont, F.Berteloot et J.-J.Loeb ont donn\'e une
r\'eponse positive \`a la question 3.10 lorsque $f$ est un endomorphisme
{\it holomorphe} de l'espace projectif complexe $X=\P^k$ (voir [BeL], [BeD], [Dup 2]).
Leur travail a \'et\'e g\'en\'eralis\'e par S.Cantat [Ca 3] au cas des
endomorphismes m\'eromorphes des {\it surfaces}.
La question reste donc ouverte pour les endomorphismes m\'eromorphes
en dimension $\geq 3$.
Nous donnons ci-dessous quelques indications qui permettent de justifier la partie
la plus simple de ces \'equivalences.

\begin{elements}
Supposons pour commencer que les exposants de Lyapunov sont minimaux.
Il r\'esulte alors de l'in\'egalit\'e de Ruelle-Margulis que,
\begin{equation}
k \log (\l_k/\l_{k-1})=2 \sum_{j=1}^k \chi_j \geq \log \l_k=h_{\mu_f}(f).
\end{equation}
Ainsi $\log \l_{k-1} \leq (1-1/k) \log \l_k$, ce que l'on peut r\'e\'ecrire
\begin{equation}
\log \l_{k-1}=\log \l_{(1-1/k)\cdot k+1/k\cdot 1}
\leq \left(1-\frac{1}{k} \right) \log \l_k+\frac{1}{k} \log \l_0,
\end{equation}
puisque $\l_0(f)=1$.
Or l'application $j \mapsto \log \l_j(f)$ est concave (Th\'eor\`eme 2.4.a). Il y a 
donc \'egalit\'e dans toutes les in\'egalit\'es pr\'ec\'edentes:
l'\'egalit\'e dans (3.4) implique que $\l_j(f)=\l_1(f)^j$ pour tout $0 \leq j \leq k$,
et l'\'egalit\'e dans (3.3) implique -- par les travaux de F.Ledrappier [L] -- que
$\mu_f$ est absolument continue par rapport \`a la mesure de Lebesgue.
Nous avons donc obtenu l'implication $(1) \Rightarrow (3)$ de la question 3.10.

L'implication $(3) \Rightarrow (2)$ est \'evidente. Nous montrons \`a pr\'esent
que $(2) \Rightarrow (1)$. Comme $\mu_f$ est de dimension maximale $2k$,
la majoration dans (3.2) montre que $ 2\sum_{j=1}^k \chi_j=\log \l_k$, donc
$\mu_f$ est absolument continue par rapport \`a la mesure de Lebesgue.
La relation alg\'ebrique $\l_j(f)=\l_1(f)^j$ montre que les exposants
de Lyapunov sont tous minor\'es par $\frac{1}{2} \log \l_1(f)$, donc
$2\sum_{j=1}^k \chi_j \geq k \log \l_1=\log \l_k$. Toutes ces in\'egalit\'es
sont donc des \'egalit\'es, en particulier les exposants de Lyapunov
sont tous \'egaux \`a $\frac{1}{2} \log \l_1=\frac{1}{2} \log \l_k/\l_{k-1}$.

Lorsque $f$ est un exemple de Latt\`es, les assertions (1),(2),(3)
sont trivialement v\'erifi\'ees. La partie difficile des travaux cit\'es
ci-dessus consiste \`a montrer la r\'eciproque \`a cette derni\`ere
implication. Nous renvoyons le lecteur aux articles originaux.
\end{elements}

\begin{rqe}
La preuve de S.Cantat [Ca 3] s'appuie
malheureusement sur la classification de Kodaira des surfaces et
ne passe donc pas en dimension sup\'erieure.
Notons que l'hypoth\`ese alg\'ebrique 
$\l_j(f)=\l_1(f)^j$ faite en (2), (3) est essentielle: S.Cantat donne
un exemple d'endomorphisme m\'eromorphe sur une surface $K3$ 
dont la mesure $\mu_f$ est lisse, mais qui ne provient pas
d'un endomorphisme affine d'un tore (voir Th\'eor\`eme C dans [Ca 3]). 
\end{rqe}

Rappelons que la dimension de Hausdorff de $\mu_f$ est en g\'en\'eral 
diff\'erente de la dimension de Hausdorff de son support: pour
un polyn\^ome $f(z)=z^{\l}+a_{\l-2}z^{\l-2}+\cdots+a_0$ de 
degr\'e $\l \geq 2$, on obtient
$$
\dim (\mu_f)=\frac{\log \l}{\chi(\mu_f)}
=\frac{\log \l}{\log \l+\sum_{f'(c)=0} G_f(c)} \leq 1,
$$
alors que l'ensemble de Julia $J_f=\text{Supp} \mu_f$ peut \^ etre
de dimension de Hausdorff \'egale \`a $2$ (voir [Shi]).
Mentionnons \`a ce sujet que X.Buff et A.Ch\'eritat [BuC] ont exhib\'e
de nombreux exemples de polyn\^omes quadratiques pour lesquels
$J_f$ est de mesure de Lebesgue strictement positive.
\vskip.2cm

Si l'on s'int\'eresse \`a des propri\'et\'es g\'eom\'etriques fines
du support de $\mu_f$, il faut donc \'etudier d'autres mesures invariantes que $\mu_f$:
c'est un des objets du formalisme thermodynamique et
nous renvoyons le lecteur \`a [Zi], [DPU] pour quelques r\'esultats dans cette direction
(en dimension 1).

\subsection{Ensemble exceptionnel}

Il r\'esulte du Lemme 3.4 que si un point $a$ n'appartient pas \`a 
l'ensemble postcritique $PC(f)=\cup_{j \geq 1} f^j({\mathcal C}_f)$, alors
\begin{equation}
\frac{1}{\l^n}(f^n)^* \e_a=\frac{1}{\l^n} \sum_{f^n(p)=a} \e_p \longrightarrow \mu_f,
\end{equation}
o\`u $\e_p$ d\'esigne la masse de Dirac au point $p$. On appelle
{\it ensemble exceptionnel} ${\mathcal E}_f$ l'ensemble des points $a$
pour lesquels la convergence (3.5) n'a pas lieu.

Lorsque $k=1$, il r\'esulte du Lemme 1.9 et du Th\'eor\`eme 1.10
que ${\mathcal E}_f$ est un ensemble fini constitu\'e d'au plus deux points.
Lorsque $k \geq 2$, l'ensemble ${\mathcal E}_f$ n'est plus n\'ecessairement 
alg\'ebrique car il contient l'orbite des points sur lesquels une courbe est
contract\'ee (voir Remarque 6.5 dans [G 3]). Plusieurs auteurs se sont
cependant int\'eress\'es \`a la caract\'erisation de l'ensemble
${\mathcal E}_f$ lorsque $f:\P^k \rightarrow \P^k$ est un endomorphisme
{\it holomorphe}:
\begin{itemize}
\item J.E.Fornaess et N.Sibony ont montr\'e [FS 3] que ${\mathcal E}_f$ est un ensemble
pluripolaire;
\item J.Y.Briend et J.Duval ont montr\'e [BrD 2] que ${\mathcal E}_f$ est un ensemble 
alg\'ebrique, et que c'est le plus grand sous-ensemble alg\'ebrique totalement invariant;
\item J.Y.Briend, S.Cantat et M.Shishikura montrent dans [BCS] que ${\mathcal E}_f$
est constitu\'e d'une union finie de sous-espaces lin\'eaires de $\P^k$,
g\'en\'eralisant plusieurs r\'esultats partiels dans cette direction 
(voir les r\'ef\'erences dans [BCS]).
\end{itemize}

Une \'etude fine de l'ensemble ${\mathcal E}_f$ a \'et\'e entreprise par
C.Favre et M.Jonsson dans [FaJ 1] lorsque $k=2$. On s'attend notamment 
\`a ce que le nombre de points totalement invariants
(i.e. les composantes ${\mathcal E}_f^0$ de dimension $0$ de ${\mathcal E}_f$)
soit au plus \'egal \`a $3$. E.Amerik et F.Campana ont obtenu dans
[AC 1] la majoration $\sharp \; {\mathcal E}_f^0 \leq 9$.
T.C.Dinh et N.Sibony montrent dans [DS 1], corollaire 3.2.8,
que $\sharp \; {\mathcal E}_f^0 \leq 3$ lorsque $f:\P^2 \rightarrow \P^2$
provient d'un endomorphisme polynomial de $\C^2$.

\subsection{Espace des param\`etres}

Il est int\'eressant d'\'etudier la fa\c{c}on dont les objets
construits jusqu'ici (mesure $\mu_f$, exposants $\chi_j(\mu_f)$,
ensemble de Julia, etc) d\'ependent de l'endomorphisme $f$.

Lorsque $f=f_t:\P^1 \rightarrow \P^1$ est une famille de fractions
rationnelles de  degr\'e $\l \geq 2$ qui d\'epend holomorphiquement
d'un param\`etre $t \in M$ -- $M$ une vari\'et\'e complexe de dimension $m$ --,
nous avons observ\'e dans la section 1.3 que l'exposant de Lyapunov 
$\chi(\mu_{f_t})$ est une fonction plurisousharmonique et 
H\"older-continue du param\`etre $t$.
L'\'etude des mesures de bifurcation $(dd^c \chi)^m$ 
fait l'objet de plusieurs travaux r\'ecents
(voir [DeM], [BaBe], [DuF]). Lorsque $f_t(z)=z^2+t$
est la famille des polyn\^omes quadratiques ($m=1$),
la mesure $dd^c \chi$ est pr\'ecis\'ement la mesure
d'\'equilibre de l'ensemble de Mandelbrot
qui a \'et\'e abondamment \'etudi\'ee
(voir [CG], [Mi 1]).

\section{G\'en\'eralisations}

La dynamique des endomorphismes  transcendants  de $\C^k$
est un sujet 
difficile qui n\'ecessite l'utilisation de techniques diff\'erentes
de celles pr\'esent\'ees dans ce m\'emoire.
Il y a cependant une situation interm\'ediaire
entre le monde des endomorphismes polynomiaux et celui des
endomorphismes transcendants de $\C^k$ de grand degr\'e topologique: 
c'est celui des applications d'allure polynomiale \'etudi\'ees par
T.C.Dinh et N.Sibony dans [DS 1].

\begin{defi}
Soit $V$ un ouvert connexe de $\C^k$ (ou d'une vari\'et\'e de Stein)
et $U \subset \subset V$ un ouvert relativement compact.
On appelle application d'allure polynomiale toute application
$$
f:U \rightarrow V \; \text{ holomorphe, propre},
$$
de degr\'e topologique $\l \geq 2$.
\end{defi}

Ces applications ont \'et\'e introduites en 
dimension $1$ par A.Douady et J.H.Hubbard [DH].
Dans ce cas elles sont conjugu\'ees, dans un voisinage de l'ensemble 
de Julia rempli
$$
K_f:=\cap_{n \geq 0} f^{-n} (V),
$$
\`a un polyn\^ome de degr\'e $\l$. 
Ce n'est plus vrai en dimension sup\'erieure comme l'ont observ\'e 
T.C.Dinh et N.Sibony (voir Exemple 3.4.11 dans [DS 1]).
Ces derniers ont men\'e \`a bien une \'etude syst\'ematique de la dynamique
de ces applications. Ils obtiennent notamment:

\begin{thm}
Soit $f:U \rightarrow V$ une application d'allure polynomiale.
Il existe une mesure de probabilit\'e invariante $\mu_f$ 
\`a support dans $\partial K_f$ telle que
$$
\frac{1}{\l^n}(f^n)^*\nu  \longrightarrow \mu_f,
$$
pour toute mesure de probabilit\'e lisse $\nu$ dans $V$.

La mesure $\mu_f$ est m\'elangeante, d'entropie maximale $=\log \l$.
\end{thm}

Notons que l'image inverse $f^* \nu$ d'une mesure de probabilit\'e
$\nu$ (pas n\'ecessai-rement lisse) est bien d\'efinie par dualit\'e
car $f$ est une application {\it propre} qui est un rev\^etement ramifi\'e:
si $\chi$ est une fonction continue \`a support compact dans $U$, alors
$$
f_* \chi(x):=\sum_{f(y)=x} \chi(y)
$$
est une fonction continue \`a support compact dans 
$f^{-1}U \subset U \subset \subset V$. 
On compte ici les pr\'eimages avec multiplicit\'e.

\begin{esquisse}
Soit $\nu$ une mesure de probabilit\'e lisse dans $V$.
La convergence de $\l^{-n}(f^n)^* \nu$ est \'equivalente \`a 
la convergence  de $\l^{-n}(f^n)_* \chi$ dans $L^1(\nu)$, 
pour toute fonction test
$\chi$ (fonction lisse \`a support compact).
Or une telle fonction peut s'\'ecrire comme diff\'erence
de deux fonctions psh born\'ees.
Il suffit donc d'\'etablir la convergence de $\l^{-n}(f^n)_* \f$,
pour toute fonction $\f$ psh born\'ee (born\'ee au voisinage
de $K_f$ suffit bien s\^ur).

L'observation de T.C.Dinh et N.Sibony qui est la clef de ce r\'esultat
(cf Lemme 3.2.2, [DS 1]) est que le principe du maximum garantit une
telle convergence. Plus pr\'ecis\'ement, si $\p$ est 
une fonction psh dans $V$ 
telle que $\l^{-1}f_* \p \geq \p$, alors
$$
\sup_U \p \geq \sup_V \frac{1}{\l} f_* \p \geq \sup_V \p,
$$
donc $\p$ est constante.

Appliquons ceci \`a la convergence de la suite $\f_n:=\l^{-n}(f^n)_* \f$.
Comme $\f$ est born\'ee, la suite $(\f_n)$ est uniform\'ement born\'ee
par $||\f||_{L^{\infty}}$.
On peut donc extraire une sous-suite convergente $(\f_{n_j})$
qui converge dans $L^1(\nu)$ vers une valeur d'adh\'erence $\p$.
Montrons que $\p$ est constante, \'egale \`a 
$c_{\f}=(\limsup \f_n)^*$, o\`u $(\cdot)^*$ d\'esigne la r\'egularisation
sup\'erieure, ce qui assurera $\f_n \rightarrow c_{\f}$ dans $L^1(\nu)$.
Posons $\p_0:=(\limsup \f_n)^*$. Alors $\p_0=c_{\f}$ est constante,
par le principe du maximum car $\l^{-1}f_* \p_0 \geq \p_0$.
Si $\p \neq c_{\f}$, alors il existe $\e>0$ tel que $\p \leq c_{\f}-2\e$ sur $U$,
donc $\f_{n_j} \leq c_{\f}-\e$ sur $f^{-1}U$, par le lemme de Hartogs.
Mais alors $\f_{p+n_j} \leq c_{\f}-\e$ pour tout $p \geq 1$, donc
$\p_0 \leq c_{\f}-\e$, contradiction.
Notons que l'ensemble $\{x \, / \, \limsup \f_n(x)<c_{\f} \}$ est pluripolaire,
il y a donc convergence presque partout (pour la la mesure de
Lebesgue) de la suite $\f_n(x)$ vers $c_{\f}$.
Nous verrons plus loin (Lemme 3.14) que la suite $(\f_n)$ converge
vers $c_{\f}$ dans $L^2(\mu_f)$.

Il r\'esulte de l'analyse pr\'ec\'edente que si $\chi$ est une fonction test
quelconque, alors $\chi_n:=\l^{-n}(f^n)_*\chi$ converge presque partout et
dans $L^1(\nu)$ vers une limite $c_{\chi}$ qui d\'epend lin\'eairement de
$\chi$. Il s'ensuit que 
$$
\mu_f: \chi \mapsto c_{\chi}=\lim \langle \l^{-n} (f^n)^* \nu, \chi \rangle
$$
est une mesure de probabilit\'e bien d\'efinie.
Elle ne d\'epend pas du choix de $\nu$ puisque $\chi_n$ converge 
presque partout vers $c_{\chi}$.
Si on part d'une mesure $\nu$ \`a support compact dans $V \setminus U$,
on obtient que $\text{Supp} \, \mu_f \subset \partial K_f$.

Observons que $f^* \mu_f=\l \mu_f$.
Pour montrer le m\'elange, il faut v\'erifier que si
$\chi,\p$ sont des fonctions test, alors
$$
I_n:=\langle \chi \circ f^n \p,\mu_f \rangle \longrightarrow c_{\chi} c_{\p}, \; \; 
\text{ o\`u } c_{\chi}=\int \chi d\mu_f, \; c_{\p}=\int \p d\mu_f.
$$
La relation fonctionnelle $f^* \mu_f=\l \mu_f$ permet de r\'e\'ecrire
$I_n$ sous la forme 
$$
I_n=\langle \chi \mu_f, \p_n \rangle , \; \text{ avec } \p_n:=\frac{1}{\l^n} (f^n)_* \p.
$$
Il r\'esulte du lemme de Hartogs (voir Lemme 3.14 ci-dessous)
que $\p_n$ converge vers $c_{\p}$ dans $L^2(\mu_f)$.
On en d\'eduit la convergence annonc\'ee.

Il reste \`a v\'erifier que $\mu_f$ est d'entropie maximale.
On introduit les degr\'es dynamiques
$$
\l_l(f)=\limsup_{n \rightarrow +\infty} [ \d_l(f^n,\om)]^{1/n},
\; \text{ o\`u } \d_l(f,\om):=\int_U f_* \om^{k-l} \wedge \om^l,
$$
et $\l_k(f)=\l$ est le degr\'e topologique de $f$.
La majoration de M.Gromov (voir Th\'eor\`eme 2.8) et
le principe variationnel donnent
$
h_{\mu_f}(f) \leq h_{top}(f) \leq \max_{ 1 \leq j \leq k} \log \l_j(f).
$
Or $\l=\l_k(f)$ domine tous les autres degr\'es dynamiques.
En effet comme $V$ est Stein on peut choisir une forme de K\"ahler $\om$ 
qui s'\'ecrit $\om=dd^c \Phi$, o\`u $\Phi$ est une fonction lisse
strictement psh. Comme $\l^{-n}(f^n)_* \Phi$ converge vers $c_{\Phi}$,
il vient 
$$
\d_{k-1}(f^n,\om)=\int_U (f^n)_* dd^c \Phi \wedge \om^{k-1}=o(\l^n),
$$
car $\l^{-n}(f^n)_* dd^c \Phi$ converge faiblement vers $0$.
Les autres degr\'es dynamiques se majorent de fa\c{c}on similaire
(voir corollaire 3.3.4, [DS 1]). Ainsi
$$
h_{\mu_f}(f) \leq h_{top}(f) \leq \log \l.
$$
La minoration $h_{\mu_f}(f) \geq \log \l$ est une cons\'equence
de ce que $\mu_f$ est de jacobien constant $=\l$ et ne charge pas les
valeurs critiques de $f$ (voir Proposition 2.9,
ainsi que le Th\'eor\`eme 2.3.1.(2) dans [DS 1]).
\end{esquisse}

\begin{lem}
Soit $\f$ une fonction psh born\'ee dans $U$.
Alors $\f_n:=\l^{-n} (f^n)_* \f$ converge dans $L^2(\mu_f)$
vers $c_{\f}=\int \f d\mu_f$.
\end{lem}

\begin{preuve}
Fixons $\e>0$ et posons $\Omega_n^{\e}:=\{ |\f_n-c_{\f}|>\e \}$.
Le lemme de Hartogs assure que $\f_n \leq c_{\f}+\e^2$
pour $n$ assez grand, d'o\`u $\Omega_n^{\e}:=\{ \f_n < c_{\f}-\e \}$.
On a donc pour $n \geq n_{\e}>>1$,
$$
\mu_f(\Omega_{n}^{\e}) \leq \frac{1}{\e} \int |\f_n-c_{\f}| d\mu_f
\leq \e+\frac{1}{\e}\int [c_{\f}+\e^2-\f_n] d\mu_f = 2\e,
$$
car $\int \f_n d\mu_f=c_{\f}$.
Comme la suite $(\f_n)$ est uniform\'ement born\'ee
($||\f_n||_{L^{\infty}} \leq ||\f||_{L^{\infty}})$, on en d\'eduit
$$
\int |\f_n-c_{\f}|^2 d\mu_f \leq \e^2+2\e ||\f||_{L^{\infty}},
$$
d'o\`u le r\'esultat.
\end{preuve}

\begin{rqe}
T.C.Dinh et N.Sibony montrent que la convergence de $\l^{-n}(f^n)_* \f$ vers $c_{\f}$
dans $L^2(\mu_f)$ implique que $\mu_f$ est $K$-m\'elangeante
(voir Proposition 2.2.2, [DS 1]).
\end{rqe}

Une difficult\'e importante, comme dans le cas des endomorphismes m\'ero-morphes,
est de montrer que $\mu_f$ n'est pas trop singuli\`ere,
par exemple qu'elle int\`egre les fonctions psh. Ce n'est pas toujours le cas
comme le montre l'exemple suivant (Exemple 3.10.3 dans [DS 1]).

\begin{exa}
Consid\'erons  $f:(z,w) \in \C^2 \mapsto (3z,Q(w)) \in \C^2$, o\`u 
$Q$ est un polyn\^ome de degr\'e $\l \geq 2$. Soit
$$
V_R:=\left\{ (z,w) \in \C^2 \, / \, |z|,|w| <R \right\}.
$$
Alors $f:U_R:=f^{-1}V_R \rightarrow V_R$ d\'efinit, pour $R>>1$
assez grand, une application d'allure polynomiale telle que 
$\mu_f$ est port\'ee par le sous-ensemble analytique $(z=0)$:
en particulier $\log |z| \notin L^1(\mu_f)$.
\end{exa}

Dans cet exemple, on a $\l_{k-1}(f)=\l$.
T.C.Dinh et N.Sibony montrent (corollaire 3.9.9, [DS 1])
que l'in\'egalit\'e $\l_{k-1}(f)>\l=\l_k(f)$
est une condition n\'ecessaire pour que $\mu_f$ int\`egre
les fonctions psh (``$\mu_f$ est PLB'') et que la condition 
``$\mu_f$ est PLB'' est stable par petites perturbations
(corollaire 3.9.7, [DS 1).
En particulier les petites perturbations $f_{\e}$
d'un endomorphisme polynomial $f:\C^k \rightarrow \C^k$
qui s'\'etend holomorphiquement \`a $\P^k$
d\'efinissent des applications d'allure polynomiale
$f_{\e}: f_{\e}^{-1} V_R \rightarrow V_R$ pour $R >>1$
telles que $\mu_{f_{\e}}$ est PLB.

Sous une telle condition T.C.Dinh et N.Sibony \'etablissent
les principales propri\'et\'es ergodiques de la mesure $\mu_f$:
\begin{itemize}
\item ses exposants de Lyapunov sont strictement positifs, $\geq \frac{1}{2} \log \l/\l_{k-1}$;
\item les points p\'eriodiques r\'epulsifs s'\'equidistribuent selon $\mu_f$;
\item il y a d\'ecroissance exponentielle des coefficients de corr\'elation.
\end{itemize}
Nous renvoyons le lecteur \`a [DS 1] pour la preuve de ces r\'esultats.
\vskip.5cm

Ces auteurs ont \'egalement \'etudi\'e [DS 8] des probl\`emes 
d'\'equidistribution pour les {\it correspondances m\'eromorphes}.
Une correspondance m\'eromorphe 
$$
\begin{array}{ccccc}
\text{ } & \text{ } & \G_f & \text{ } & \text{ } \\
\text{ } & \stackrel{\pi_1}{\swarrow} & \text{ } &
\stackrel{\pi_2}{\searrow} & \text{ } \\
X & \text{ } & \stackrel{f}{\longrightarrow} & \text{ } & X
\end{array}
$$
est la donn\'ee d'un graphe $\Gamma_f$, sous-ensemble analytique
de $X^2$, et de deux projections holomorphes surjectives
$\pi_1,\pi_2$. Lorsque $\pi_1$ est g\'en\'eriquement injective
(i.e. injective hors d'un sous-ensemble analytique propre),
on retrouve la notion d'endomorphisme m\'eromorphe.
On peut d\'efinir des degr\'es dynamiques associ\'es \`a 
une telle correspondance et montrer
l'existence d'une mesure canonique $\mu_f$ lorsque la correspondance
a un grand degr\'e topologique
(voir corollaire 5.3, [DS 8]).

T.C.Dinh et N.Sibony obtiennent aussi des r\'esultats d'\'equidistribution
g\'en\'eraux qui s'appliquent aussi bien au cas des 
{\it transformations m\'eromorphes} qu'\`a celui de l'it\'eration al\'eatoire
(voir Th\'eor\`eme 1.2, [DS 8]).

\chapter{Petit degr\'e topologique en dimension deux}

Dans ce chapitre nous supposons que $X$ est une {\it surface} --i.e. de dimension
complexe $2$-- et $f:X \rightarrow X$ est un endomorphisme m\'eromorphe
de petit degr\'e topologique, i.e. $\l_1(f)>\l_2(f)$.

\section{Bon mod\`ele}

Le premier probl\`eme auquel nous sommes confront\'es est celui
du calcul effectif du degr\'e dynamique $\l_1(f)$,
qui est d\'efini par une limite --contrairement \`a $\l_2(f)$.
Il se peut tr\`es bien par exemple que
$r_1(f)$ soit plus grand que le degr\'e topologique $\l_2(f)$, bien que
$\l_2(f)>\l_1(f)$. Plus s\'erieusement, nous n'avons pas seulement besoin de
conna\^itre la valeur de $\l_1(f)$ pour la comparer \`a celle de $\l_2(f)$,
mais, lorsque $\l_1(f)>\l_2(f)$, il est important de contr\^oler la croissance
de la norme $||(f^n)^*||$ sur $H^{1,1}(X,\R)$ (cf section 2.3): est-ce que la croissance
est de l'ordre de $\l_1(f)^n$, ou bien fait-elle appara\^itre
\'egalement des termes polynomiaux $n^s \l_1(f)^n$ ?
Nous nous proposons dans cette section d'apporter des \'el\'ements de r\'eponse
\`a ces questions.

\begin{defi}
On dit que $f:X \rightarrow X$ est $1$-stable si l'action lin\'eaire
induite par $f^*$ sur $H^{1,1}(X,\R)$ est compatible avec la dynamique, i.e.
si pour tout $n \in \N$, on a $(f^n)^*=(f^*)^n$.
\end{defi}

Consid\'erons l'endomorphisme polynomial de $\C^2$,
$$
f:(z,w) \in \C^2 \mapsto (zw+c_1,z+c_2) \in \C^2.
$$
C'est un endomorphisme {\it birationnel} de $\C^2$, i.e. son degr\'e topologique
$\l_2(f)$ est \'egal \`a 1. On peut calculer \`a la main $\l_1(f)=\frac{1+\sqrt{5}}{2}$,
mais nous allons retrouver cette valeur d'une autre fa\c{c}on.
On peut consid\'erer l'extension m\'eromorphe de l'endomorphisme $f$ \`a n'importe quelle
compactification lisse de $\C^2$. Notons $g:\P^2 \rightarrow \P^2$
son extension \`a $\P^2$. Elle induit une action lin\'eaire
$g^*$ sur l'espace $H^{1,1}(\P^2,\C) \simeq \C$ qui est engendr\'e par la classe d'une droite
projective $L$. On v\'erifie que $g^* L \sim 2 L$ et $(g^2)^*L \sim 3L$,
donc $g$ n'est pas 1-stable.

Consid\'erons \`a pr\'esent $h:\P^1 \times \P^1 \rightarrow \P^1 \times \P^1$ l'extension
m\'eromorphe de $f$ \`a $\P^1 \times \P^1$. Elle induit une action lin\'eaire
$h^*$ sur l'espace $H^{1,1}(\P^1 \times \P^1,\C) \simeq \C^2$ qui est engendr\'e
par la classe d'une droite verticale $L_z=(z=0)$ et 
celle d'une droite horizontale $L_w=(w=0)$.
On obtient
$$
h^*L_z \sim L_z+L_w \; \text{ et } \; h^*L_w \sim L_z,
$$
donc l'action de $h^*$ est donn\'ee, dans la base $(L_z,L_w)$ par la matrice
$A_h=\left[ \begin{array}{cc} 1 & 1 \\ 1 & 0 \end{array} \right]$.
L'application $h$ est 1-stable (voir le crit\`ere 4.4), ce qui assure
$$
\l_1(f)=\l_1(g)=\l_1(h)=r_1(h)=\text{rayon spectral de } A_h=\frac{1+\sqrt{5}}{2}.
$$
Cela r\'esulte de l'observation \'el\'ementaire suivante.

\begin{obs}
Si $f:X \rightarrow X$ est 1-stable, alors $r_1(f)=\l_1(f)$.
\end{obs}

Comme $\l_1(f)$ est invariant
par conjugaison bim\'eromorphe, il est naturel de se demander si l'on peut
toujours effectuer un changement bim\'eromorphe de coordonn\'ees -- et donc
remplacer $X$ par une vari\'et\'e bim\'eromorphiquement \'equivalente $\tilde{X}$ --
de telle sorte que l'endomorphisme $\tilde{f}$ induit par $f$ sur
$\tilde{X}$ soit 1-stable.
Ce n'est pas toujours possible, sans condition sur les degr\'es dynamiques,
comme l'a montr\'e C.Favre
dans [Fa 2] en analysant le cas des endomorphismes toriques de dimension 2:
il existe des endomorphismes toriques qui ne peuvent jamais \^etre rendus 1-stables.
Les endomorphismes concern\'es v\'erifient $\l_2(f)=\l_1(f)^2>\l_1(f)$,
ce qui motive la

\begin{ques}
Soit $f:X \rightarrow X$ un endomorphisme tel que $\l_1(f)>\l_2(f)$.
Peut-on effectuer un changement bim\'eromorphe de coordonn\'ees
$\pi:\tilde{X} \rightarrow X$ de sorte que
l'endomorphisme $\tilde{f}:=\pi^{-1} \circ f \circ \pi:\tilde{X} \rightarrow \tilde{X}$
soit 1-stable ?
\end{ques}

On dispose d'un crit\`ere alg\'ebrique simple pour d\'etecter
la 1-stabilit\'e, comme l'ont observ\'e J.E.Fornaess et N.Sibony lorsque
$X$ est un espace projectif complexe [FS 4].

\begin{pro}
Un endomorphisme $f:X \rightarrow X$ est 1-stable si et seulement si aucune
courbe ${\mathcal C}$ n'est contract\'ee par un it\'er\'e
$f^n$ sur un point d'ind\'etermination.
\end{pro}

\begin{preuve}
Supposons qu'une courbe ${\mathcal C}$ est contract\'ee par $f^n$ sur un point
$p \in I_f$. On peut supposer,
quitte \`a changer $f$ en $f^n$, que $n=1$.
Soit $\om$ une forme de K\"ahler. Alors $f^*\om$ est un courant positif ferm\'e
de bidegr\'e $(1,1)$ sur $X$ qui a un nombre de Lelong positif au point $p$.
Il s'ensuit que $f^* (f^* \om)$ est un courant qui charge la courbe ${\mathcal C}$.
Or $(f^2)^* \om$ est par d\'efinition une forme diff\'erentielle \`a coefficients
localement int\'egrables: elle co\"{\i}ncide presque partout
avec $f^*(f^*\om)$ mais ne charge pas
la courbe ${\mathcal C}$. On en d\'eduit que
les classes de cohomologie $f^*(f^*\{\om\})$ et $(f^2)^*\{ \om \}$ sont diff\'erentes.

R\'eciproquement supposons qu'aucune courbe n'est contract\'ee sur un point 
d'ind\'etermination.
Alors les ensembles $f^{-n}(I_f)$ sont tous de codimension $\geq 2$.
Soit $\theta$ une forme diff\'erentielle lisse ferm\'ee
de bidegr\'e $(1,1)$. Par d\'efinition,
$f^* (f^* \theta)$ et $(f^2)^* \theta$ co\"{\i}ncident hors
de $I_f \cup f^{-1}I_f$. Comme les ensembles de codimension $\geq 2$ sont n\'egligeables
pour les courants de bidegr\'e $(1,1)$, on en d\'eduit que ces formes co\"{\i}ncident
partout. Ainsi $(f^*)^2=(f^2)^*$ et on montre de m\^eme que
$(f^*)^n=(f^n)^*$ pour tout $n$.
\end{preuve}

\vskip.2cm

Dans l'exemple polynomial $f(z,w)=(zw+c_1,z+c_2)$ vu plus haut, les points
d'ind\'etermination sont situ\'es dans le diviseur \`a l'infini.
Les seules courbes qui peuvent \^etre contract\'ees sur un point d'ind\'etermination
sont les composantes irr\'eductibles de ce diviseur. On v\'erifie ainsi que
$h(w=\infty)=(z=\infty)$ et $h(z=\infty)=(\infty,\infty)$, or le point
$(\infty,\infty)$ est un point fixe qui n'est pas un point d'ind\'etermination pour $h$,
donc $h$ est 1-stable. 

Ce crit\`ere est \'egalement facile \`a v\'erifier sur
les surfaces de dimension de Kodaira nulle. On a vu \`a la section 2.4.2
que si $X$ est une surface minimale de dimension de Kodaira nulle,
alors $f$ est non-ramifi\'e. Or les seules courbes susceptibles d'\^etre
contract\'ees sur un point d'ind\'etermination sont les courbes de l'ensemble
critique. Comme celui-ci est vide, on en d\'eduit l'assertion suivante:

\begin{pro}
Si $kod(X)=0$, alors on peut supposer que $f$ est 1-stable en travaillant sur le
mod\`ele minimal de $X$.
\end{pro}

Le cas des surfaces rationnelles est plus d\'elicat.
Nous disposons toutefois du r\'esultat suivant de 
J.Diller et C.Favre [DF]:

\begin{thm}
Si $\l_2(f)=1$, alors on peut, quitte \`a conjuguer par une application 
bim\'eromorphe, supposer que $f$ est 1-stable.
\end{thm}

La preuve (voir Th\'eor\`eme 0.1 dans [DF]) utilise de fa\c{c}on cruciale le fait
que $f$ est birationnel: dans ce cas $f$ peut s'\'ecrire comme une composition
d'\'eclatements et de contractions. Le cas g\'en\'eral, $2 \leq \l_2(f) <\l_1(f)$,
semble beaucoup plus difficile. Notons qu'une avanc\'ee significative a \'et\'e
r\'ealis\'ee par C.Favre et M.Jonsson lorsque $f$ provient
d'un endomorphisme polynomial de $\C^2$ (voir [FaJ 3]).

\section{Analyse spectrale}

Dans toute la suite du chapitre 4, nous supposons que $f:X \rightarrow X$ est 
un endomorphisme 1-stable tel que $\l_1(f)>\l_2(f)$.
Nous allons en tirer des cons\'equences sur le spectre
des actions lin\'eaires induites par $f^*,f_*$ sur 
$H^{1,1}(X,\R)$. On note
dans la suite $\l:=\l_1(f)$. Comme $f$ est 1-stable, il r\'esulte de
la section pr\'ec\'edente que $\l=r_1(f)$ est \'egal au rayon spectral de
l'action lin\'eaire induite par $f^*$ sur $H^{1,1}(X,\R)$.

L'action 
$f_*:H^{1,1}(X,\R) \rightarrow H^{1,1}(X,\R)$ est duale de 
celle de $f^*$, car on est en dimension deux.
En particulier $(f_*)^n=(f^n)_*$ pour tout $n \in \N$.
Comme le c\^one $H^{1,1}_{nef}(X,\R)$ est dual du c\^one $H^{1,1}_{psef}(X,\R)$,
qui est pr\'eserv\'e par $f^*,f_*$, on
en d\'eduit que le c\^one nef est \'egalement pr\'eserv\'e par $f^*$ et $f_*$. 

Lorsque $X=\P^1 \times \P^1$, il s'ensuit que les actions lin\'eaires $f^*,f_*$
sont donn\'ees par des matrices \`a coefficients entiers naturels dans la base
$\{z=0\}$, $\{w=0\}$ de $H^{1,1}(X,\R)$. Il r\'esulte alors du Th\'eor\`eme
de Perron-Frob\'enius que $\l=r_1(f)$ est une valeur propre de $f^*,f_*$
et admet un  vecteur propre dans $H^{1,1}_{nef}(X,\R)$, comme cela a \'et\'e
observ\'e dans [FaG]. C'est en fait vrai en toute g\'en\'eralit\'e comme
l'ont montr\'e J.Diller et C.Favre dans [DF].

\begin{thm}
Soit $f:X \rightarrow X$ un endomorphisme $1$-stable tel que $\l:=\l_1(f)>\l_2(f)$.
Alors  il existe
une $(1,1)$-classe nef non nulle $\a^+$ (resp. $\a_-$) telle que
$f^* \a^+=\l \a_+$ (resp. $f_*\a_-=\l \a_-$); en particulier
le rayon spectral $\l:=r_1(f)=\l_1(f)$ est une valeur propre de $f^*,f_*$.
De plus
\begin{enumerate}
\item la valeur propre $\l$ est une racine simple du polyn\^ome caract\'eristique
de $f^*$. Les autres valeurs propres sont de module 
$\leq \sqrt{\l_2(f)}<\l$;
\item les vecteurs propres nef $\a^+$, $\a^-$ v\'erifient $\a^+ \cdot \a^->0$.
\end{enumerate}
\end{thm}

\begin{esquisse} 
Observons que $H^{1,1}_{nef}(X,\R)$ est un c\^one strict,
c'est \`a dire
$H^{1,1}_{nef}(X,\R) \cap -H^{1,1}_{nef}(X,\R)=\emptyset$.
Comme l'int\'erieur de $H^{1,1}_{nef}(X,\R)$ est non-vide
(c'est le c\^one de K\"ahler), il r\'esulte du Th\'eor\`eme
de Perron-Frob\'enius que le rayon spectral 
$\l=\l_1(f)=r_1(f)$ de $f^*$ sur $H^{1,1}(X,\R)$ est une valeur propre
qui poss\`ede un vecteur propre $\a^+$ dans $H^{1,1}_{nef}(X,\R)$.

Supposons dans un premier temps que $f$ est holomorphe.
Dans ce cas les op\'erateurs $f_*$ et $f^*$ v\'erifient
$f_*f^*(\cdot)=\l_2(f) (\cdot)$. On en d\'eduit que pour toutes (1,1)-classes
$\a,\b$, on a 
\begin{equation}
\langle f^*\a,f^*\b\rangle =\l_2(f) \langle \a,\b\rangle .
\end{equation}
En appliquant cette relation \`a $\a=\b=\a^+$, on obtient
$(\a^+)^2=0$, puisque $\l_2(f)<\l_1(f)<\l_1(f)^2$. Supposons
qu'il existe une classe $\b$ non nulle telle que
$f^* \b=\l \beta+c \a^+$. En appliquant (4.1) aux classes $\a=\a^+$ et
$\b$, on obtient de m\^eme $\langle \a^+,\b\rangle =0$. Une nouvelle application
de (4.1) avec $\a=\b$ donne enfin $\b^2=0$. La forme
d'intersection est donc non-n\'egative sur l'espace
vectoriel engendr\'e par $\a^+$ et $\b$. Il r\'esulte alors du Th\'eor\`eme
de l'indice de Hodge (voir [BPV]) que $\b$ est n\'ecessairement proportionnelle 
\`a $\a^+$: cela d\'emontre que $\l$ est une valeur propre simple du 
polyn\^ome caract\'eristique de $f^*$. Par dualit\'e il en est
de m\^eme pour $f_*$.

Soit \`a pr\'esent $t$ une valeur propre (complexe) de $f^*$ diff\'erente de $\l$
et  $\b$ un vecteur propre associ\'e. En appliquant \`a nouveau (4.1)
successivement \`a $\a=\b$, puis \`a $\a=\a^+$ et $\b$, on obtient
que soit $t^2=\l_2(f)$, soit $t=\l_2(f)/\l$, soit $\b^2=\langle \a^+,\b\rangle =0$.
Mais ce dernier cas n'est pas admissible car il implique $\b$ proportionnelle \`a 
$\a^+$ --par le Th\'eor\`eme de l'indice de Hodge--, contredisant
le fait que $t$ est distincte de $\l$. On en d\'eduit que
$|t| \leq \sqrt{\l_2(f)}$.

Il r\'esulte \`a nouveau du Th\'eor\`eme de l'indice de Hodge que les 
classes $\a^+$ et $\a^-$ s'intersectent positivement.

Lorsque $f$ est seulement m\'eromorphe, il faut remplacer l'identit\'e
(4.1) par une identit\'e plus compliqu\'ee, mais qui va jouer le m\^eme
r\^ole. C'est la formule d'aller-retour que nous rappelons ci-dessous.
L'analyse des spectres de $f^*,f_*$ se fait de fa\c{c}on analogue
en utilisant la positivit\'e de la forme quadratique
hermitienne $Q$. Nous renvoyons le lecteur \`a [DF] pour les d\'etails techniques.
\end{esquisse}

\begin{pro}[Formule d'aller-retour]
Il existe une forme quadratique hermitienne $Q \geq 0$ sur $H^{1,1}(X,\C)$
telle que
$$
\langle f^* \a,f^* \b\rangle =\l_2(f) \langle \a,\b\rangle +Q(\a,\b)
$$
pour toutes les classes $\a,\b \in H^{1,1}(X,\C)$. De plus 
$Q(\a,\a)=0$ si et seulement si $\langle \a,f(p)\rangle =0$ pour tout point 
d'ind\'etermination $p \in I_f$.
\end{pro}

\begin{esquisse}
La formule r\'esulte d'un calcul explicite de $f_*f^*\a$.
En passant par le graphe 
$X \stackrel{\pi_1}{\leftarrow} \Gamma_f \stackrel{\pi_2}{\rightarrow} X$
de l'application $f$, on d\'ecompose le calcul en deux temps:
la projection $\pi_2$ \'etant holomorphe, on a $(\pi_2)_*\pi_2^* \g=\l_2(f) \g$
pour toute classe $\g \in H^{1,1}(X,\C)$. 
Pour comprendre l'action de la projection
$\pi_1$ qui est une compos\'ee d'\'eclatements de points, il suffit
d'observer que $\pi^*\pi_* \g=\g+\langle \g,E \rangle \{E\}$
pour un \'eclatement $\pi$ dont on a not\'e $E$ le diviseur exceptionnel.
Nous renvoyons \`a [DF] (corollaire 3.4) pour plus de d\'etails.
\end{esquisse}
\vskip.2cm

Il r\'esulte de cette analyse que la trace des op\'erateurs $(f^n)^*$
cro\^it comme $\l_1(f)^n$ --\`a une constante multiplicative pr\`es. On obtient
donc une majoration pr\'ecise du nombre de points p\'eriodiques (voir
section 2.3). Il nous faut maintenant \'etablir une minoration de ce
nombre en cr\'eant beaucoup --environ $\l_1(f)^n$-- de points 
p\'eriodiques selles d'ordre $n$. La strat\'egie est \`a pr\'esent
de construire un $(1,1)$-courant ferm\'e positif $T^+$ (resp. $T^-$)
tel que $\{T^+\}=\a^+$ (resp. $\{T^-\}=\a^-$) et 
$f^* T^+=\l T^+$ (resp. $f_* T^-=\l T^-$). Comme les classes
$\a^+$ et $\a^-$ s'intersectent positivement, on peut esp\'erer
donner un sens au produit d'intersection $\mu_f=T^+ \wedge T^-$ qui produira
alors une mesure dynamiquement int\'eressante.

\section{Courants invariants}

Rappelons que dans toute la suite du chapitre 4, $f:X \rightarrow X$ est un
endomorphisme m\'eromorphe 1-stable tel que $\l:=\l_1(f)>\l_2(f)$.
Soit $\om$ une forme de K\"ahler sur $X$.
Soit $\a^+ \in H^{1,1}_{nef}(X,\R)$ (resp. $\a^-$) une classe non-nulle telle que
$f^* \a^+=\l \a^+$ (resp. $f_*\a^-=\l \a^-$). 
Il r\'esulte du Th\'eor\`eme 4.7 que les classes $\a^+, \a^-$ sont
uniques \`a une constante multiplicative pr\`es.
On normalise les classes $\a^+,\a^-,\{\om\}$ en imposant 
$$
\a^+ \cdot \a^-=\a^+ \cdot \{\om\}=\a^- \cdot \{\om\}=1.
$$

\subsection{Le courant $f^*$-invariant canonique}

\begin{thm}
Il existe un unique courant positif ferm\'e $T^+$ de bidegr\'e $(1,1)$ sur $X$
tel que $f^*T^+=\l T^+$, $\{T^+\}=\a^+$ et
$$
\frac{1}{\l^n}(f^n)^* \om \longrightarrow T^+.
$$
\end{thm}

\begin{esquisse}
Soit $\theta^+$ une forme lisse ferm\'ee de bidegr\'e $(1,1)$ sur $X$ dont la
classe de cohomologie est \'egale \`a $\a^+$. L'invariance
$f^* \a^+=\l \a^+$ se traduit, gr\^ace au $dd^c$-lemma de la th\'eorie de Hodge
(voir [GH] p 149), par
\begin{equation}
\frac{1}{\l}f^*\theta^+=\theta^++dd^c \g^+,
\end{equation}
o\`u $\g^+ \in L^1(X,\R)$ est une fonction lisse hors de $I_f$,
qui peut avoir des singularit\'es logarithmiques aux points d'ind\'etermination.

Comme $f$ est 1-stable, on peut prendre le pull-back de cette \'equation
par les it\'er\'es $f^n$ de $f$. La relation de compatibilit\'e
$(f^n)^*=(f^*)^n$ donne alors
$$
\frac{1}{\l^n}(f^n)^*\theta^+=\theta^+ +dd^c g_n^+,
\text{ o\`u } g_n^+:=\sum_{j=0}^{n-1} \frac{1}{\l^j} \g^+ \circ f^j.
$$
La difficult\'e consiste \`a montrer que la suite
$(g_n^+)$ converge dans $L^1(X)$. 

Nous affirmons que la fonction $\g^+$ est major\'ee sur $X$. Lorsque la forme
$\theta^+$ est semi-positive, la fonction $\g^+$ est quasiplurisousharmonique (qpsh),
i.e. localement diff\'erence d'une fonction psh (un potentiel local du courant
positif $\l^{-1}f^* \theta^+$) et d'une fonction lisse (potentiel local
de $\theta^+$); elle est donc bien major\'ee dans ce cas.
La semi-positivit\'e de $\theta^+$ est une hypoth\`ese forte sur la positivit\'e 
de la classe $\a^+$ que l'on peut toutefois v\'erifier lorsque
$X$ est une surface rationnelle minimale. Dans le cas g\'en\'eral, nous allons
montrer que la fonction $\g^+ \circ \pi_1$ est qpsh sur le graphe
d\'esingularis\'e $\tilde{\Gamma}_f$.
Si l'on tire en arri\`ere par $\pi_1$ l'\'equation (4.2), il vient
\begin{equation}
dd^c (\g^+ \circ \pi_1)=\frac{1}{\l} \pi_1^*((\pi_1)_* \pi_2^* \theta^+)-\pi_1^* \theta^+
=\frac{1}{\l} [\pi_2^*\theta^+ +R]-\pi_1^*\theta^+,
\end{equation}
o\`u $R=\pi_1^* (\pi_1)_* \eta-\eta$ et $\eta=\pi_2^* \theta^+$.
Comme $\pi_1$ est une composition d'\'eclatements, on v\'erifie 
que $R$ est un courant positif, support\'e par le diviseur exceptionnel de $\pi_1$:
c'est une formule de ``pull-push'' (appliquer la formule d'aller-retour 
\`a $\pi_1^{-1}$). Nous renvoyons le lecteur
\`a [DG] pour plus de d\'etails.
L'\'egalit\'e (3.3) exprime donc
$dd^c (\g^+ \circ \pi_1)$ comme la somme d'un courant positif et d'une forme lisse.
La fonction $\g^+ \circ \pi_1$ est ainsi qpsh, donc major\'ee sur $\tilde{\Gamma}_f$,
d'o\`u $\g^+$ est major\'ee sur $X$.

Quitte \`a retrancher une constante, on peut donc supposer que $\sup_X \g^+=0$. Cela
permet d'exprimer la suite $g_n^+$ comme une suite d\'ecroissante de fonctions
presque-qpsh. On montre alors par des arguments classiques que
la limite n'est pas identiquement $-\infty$, notons la $g^+$. Ainsi
$$
\frac{1}{\l^n}(f^n)^* \theta^+ \longrightarrow T^+:=\theta^++dd^c g^+.
$$
Le reste des assertions s'ensuit ais\'ement. Nous renvoyons le lecteur 
aux articles cit\'es ci-dessous pour plus de d\'etails.
\end{esquisse}
\vskip.2cm

\begin{rqe}
La construction du courant $T^+$ est due \`a N.Sibony pour un endomorphisme
rationnel quelconque de $X=\P^2$ [S] (certains cas particuliers avaient \'et\'e
trait\'es par J.E.Fornaess et N.Sibony dans [FS 2,3,4]). Nous avons g\'en\'eralis\'e
son approche en l'\'etendant au cas des surfaces de Hirzebruch dans [G 1]. 
J.Diller et C.Favre ont construit [DF] le courant $T^+$ lorsque $f$ est birationnel. 
Le cas g\'en\'eral \'enonc\'e ici est d\'emontr\'e dans [DG],
en collaboration avec J.Diller. 
\end{rqe}

Observons que le courant $T^+$ est canonique en ce sens que si $\theta$ est n'importe
quelle forme lisse ferm\'ee de bidegr\'e $(1,1)$, alors
$$
\frac{1}{\l^n}(f^n)^* \theta \rightarrow c T^+, \; \text{ avec } c=\{ \theta \} \cdot \a^-.
$$
Plus g\'en\'eralement on peut s'int\'eresser \`a la convergence de
$\l^{-n}(f^n)^* \theta$ lorsque $\theta$ est un (1,1)-courant positif ferm\'e
(voir la discussion qui pr\'ec\`ede le Th\'eor\`eme 4.12).
Nous traitons ici un cas simple d'\'equidistribution
(et renvoyons le lecteur \`a [RS], [S] pour des r\'esultats 
d'\'equidistribution plus g\'en\'eraux).

\begin{pro}
Supposons $X=\P^2$. Alors
$$
\frac{1}{\l^n}(f^n)^* [L] \longrightarrow T^+,
$$
pour presque toute droite projective $L \in (\P^2)^*$.
\end{pro}

\begin{preuve}
Notons $\nu^*$ la mesure de Fubini-Study sur l'ensemble $(\P^2)^*$
des droites projectives de $\P^2$. La formule de Crofton (voir [Ci])
permet d'exprimer la forme de K\"ahler $\om$ de Fubini-Study,
$$
\om=\int_{a \in (\P^2)^*} [L_a] d\nu^*(a)
$$
comme une moyenne de courants d'int\'egration sur les droites projectives
$L_a \subset \P^2$. Cela se traduit, au niveau des potentiels, par
$[L_a]=\om+dd^c \f_a$, o\`u
$$
\f_a[z]=\log \left[ \frac{|z \wedge a|}{||z|| \cdot ||a||} \right] \leq 0, \;
\text{ avec } \int_{a \in (\P^2)^*} \f_a[z]d\nu^*(a)=-1.
$$
Pour montrer la convergence de $\l^{-n}(f^n)^* [L_a]$ vers $T^+$,
il suffit donc de montrer que $\l^{-n} \f_a \circ f^n$ converge
vers $0$ dans $L^1(\P^2)$. Posons
$$
\p(a):=\sum_{n \geq 0} \p_n(a), \; \; 
\p_n(a)=\int_{[z] \in \P^2} -\frac{1}{\l^n} \f_a \circ f^n[z] d\nu [z],
$$
o\`u $\nu$ d\'esigne la forme volume de Fubini-Study sur $\P^2$.
Alors $\p$ est int\'egrable. En effet,
c'est une limite croissante de fonctions positives
($\f_a \leq 0$) telle que
\begin{eqnarray*}
\int_{a \in (\P^2)^*} \p(a)d\nu(a) &=& \sum_{n \geq 0}
\int_{[z] \in \P^2} \int_{a \in (\P^2)^*}
-\l^{-n} \f_a \circ f^n[z] d\nu^*(a) d\nu[z]\\
&=& \sum_{n \geq 0} \l^{-n} <+\infty.
\end{eqnarray*}
En particulier $\p$ est finie presque partout, donc
$\p_n(a) \rightarrow 0$ presque partout, 
d'o\`u $\l^{-n}(f^n)^* [L_a] \rightarrow T^+$, $\nu^*$-presque partout.
\end{preuve}

\subsection{Propri\'et\'es d'extr\'emalit\'e}

Il est naturel de se demander s'il existe d'autres courants invariants: peut-on
d\'ecrire le c\^one des courants positifs ferm\'es $S$ de bidegr\'e $(1,1)$
tels que $f^*S=\l S$ ? C'est l'analogue, dans le contexte des endomorphismes
de petit degr\'e topologique, du probl\`eme de la caract\'erisation des
points exceptionnels pour les endomorphismes de grand degr\'e topologique.

L'analyse du spectre de $f^*$ (Th\'eor\`eme 4.7) montre
que la classe $\{S\}$ d'un tel courant est proportionnelle \`a $\a^+$. On a donc
$S=c \theta^++dd^c \f_S$, pour une constante $c \geq 0$ et une fonction
qpsh $\f_S \in L^1(X)$ uniquement d\'etermin\'ee si on impose
$\sup_X \f_S=0$. Ainsi
$$
\frac{1}{\l^n}(f^n)^*S=c \frac{1}{\l^n}(f^n)^*\theta^+
+dd^c \left( \frac{1}{\l^n} \f_S \circ f^n \right).
$$
Il s'agit donc de savoir \`a quelle condition la suite $\l^{-n}\f_S \circ f^n$
converge vers 0 dans $L^1(X)$ pour pouvoir conclure (ou non) \`a la
convergence de la suite $\l^{-n}(f^n)^*S$ vers $c T^+$.
Plusieurs auteurs se sont int\'eress\'es \`a cette question 
(voir e.g. [BS 2], [FS 1,2,3,4], [HP 1], [Dil 1], [RS], [S], [FaJ 1]).
Nous avons donn\'e, en collaboration avec C.Favre, 
une solution compl\`ete \`a ce probl\`eme dans le cas
des endomorphismes birationnels (voir [FaG]). 
La  principale technique -- d\'evelopp\'ee par C.Favre dans sa th\`ese --
consiste \`a mettre au point des estim\'ees de volume 
via un contr\^ole asymptotique des multiplicit\'e locales (voir Th\'eor\`eme 1.7
et Lemma 1.9),
qui permettent de contr\^oler la convergence de la suite
$\l^{-n} \f_S \circ f^n$ dans l'esprit de la preuve du Th\'eor\`eme 1.10.

Nous n'explicitons pas ces r\'esultats ici et renvoyons le lecteur aux articles originaux.
Nous consid\'erons uniquement le cas, plus simple,
o\`u le courant $S$ est domin\'e par $T^+$.

Notons ${\mathcal T}$ le c\^one des courants positifs ferm\'es de
bidegr\'e $(1,1)$ sur $X$, et ${\mathcal T}_{f^*}$ le sous-c\^one
des courants $S \in {\mathcal T}$ tels que $f^*S=\l S$.

\begin{thm}
Si $0 \leq S \leq T^+$, alors $\l^{-n}(f^n)^* S \rightarrow c T^+$,
$c=\{ S\} \cdot \a^-$.

\noindent En particulier $T^+$ est extr\'emal dans le c\^one 
${\mathcal T}_{f^*}$ des courants $f^*$-invariants.

Si $\l_2(f)=1$, alors $T^+$ est un point extr\'emal du c\^one ${\mathcal T}$.
\end{thm}

Nous verrons un peu plus loin qu'il s'agit l\`a de propri\'et\'es de nature ergodiques
du courant $T^+$.

\begin{esquisse}
Lorsque le courant $S$ est domin\'e
par $T^+$, ses potentiels sont au moins aussi r\'eguliers que ceux de $T^+$. Cela se traduit
par l'in\'egalit\'e $\f_{T^+} \leq \f_S \leq 0$.
Or $\f_{T^+}$ diff\`ere de $g^+$ d'une constante additive,
donc par construction 
$$
\l^{-n} \f_{T^+} \circ f^n \simeq \l^{-n} g^+ \circ f^n \longrightarrow 0
\; \text{ dans } L^1(X).
$$
On en d\'eduit que $\l^{-n} \f_S \circ f^n \rightarrow 0$, et donc
$\l^{-n}(f^n)^*S \rightarrow c T^+$. 

En particulier si $S \in {\mathcal T}_{f^*}$ est un courant invariant tel que $0 \leq S \leq T^+$,
alors $S= cT^+$ donc
le courant $T^+$ est un point extr\'emal du c\^one ${\mathcal T}_{f^*}$. 

Lorsque $\l_2(f)=1$, on peut raffiner le r\'esultat de convergence ci-dessus pour le
rendre uniforme par rapport \`a $S$. Plus pr\'ecis\'ement soit $S \in {\mathcal T}$
tel que $0 \leq S \leq T^+$. Le courant $T^+$ v\'erifie $f^*T^+=\l T^+$
donc \'egalement $f_*T^+=(f^{-1})^*T^+=\l^{-1}T^+$ dans $X \setminus {\mathcal C}_{f^{-1}}$,
o\`u ${\mathcal C}_{f^{-1}}$ d\'esigne l'ensemble critique de $f^{-1}$.
Les courants positifs ferm\'es $S_j:=\l^j (f^j)_*S$ sont donc bien d\'efinis
et domin\'es par $T^+$ dans $X \setminus {\mathcal C}_{f^{-1}}$. 
On note encore $S_j$ leur extension
triviale \`a travers ${\mathcal C}_{f^{-1}}$: ce sont \`a nouveau des \'el\'ements de 
${\mathcal T}$ (Th\'eor\`eme d'extension de Skoda [Sk]) domin\'es par $T^+$.
Ils v\'erifient $\l^{-j}(f^j)^*S_j \equiv S$ dans $X$, donc
$S_j$ est cohomologue \`a $c \theta^+$, $0 \leq c \leq 1$, et on peut
appliquer le raisonnement ci-dessus pour conclure
\`a $S=\l^{-j}(f^j)^* S_j \rightarrow c T^+$, donc $S=cT^+$,
c'est \`a dire que $T^+$ est un point extr\'emal du
c\^one ${\mathcal T}$.
\end{esquisse}

\begin{rqe}
Ce r\'esultat est d\^u \`a J.E.Fornaess et N.Sibony [FS 3] lorsque $f$
est une application de H\'enon complexe (voir Th\'eor\`eme 2.4 dans [G 1] 
pour le cas g\'en\'eral).
\end{rqe}

\subsection{Laminarit\'e}

Lorsque $f$ est un diff\'eomorphisme holomorphe d'Anosov
(i.e. lorsque $X$ est un tore, cf [Gh]), $T^+$ est
un {\it cycle feuillet\'e} associ\'e au feuilletage holomorphe stable.
Dans notre contexte de non hyperbolicit\'e uniforme,
il faut remplacer la notion de cycle feuillet\'e par une notion
plus faible, introduite
par E.Bedford, M.Lyubich et J.Smillie dans [BLS 1].

\begin{defi}
Un courant $T$ est dit {\bf uniform\'ement laminaire} s'il peut
\^etre localement d\'ecrit par une moyenne de courants d'int\'egration
sur des graphes disjoints.
Un courant $T$ est dit {\bf laminaire} s'il peut \^etre approch\'e par
une suite croissante de courants localement uniform\'ement laminaires.
\end{defi}

Ces auteurs ont d\'emontr\'e que $T^+$ est laminaire
lorsque $f$ est une application de H\'enon complexe [BLS 1], [BS 3].
Leur approche a \'et\'e g\'en\'eralis\'ee au cas des automorphismes
des surfaces projectives par S.Cantat [Ca 1], puis \`a celui des
endomorphismes m\'eromorphes des surfaces projectives
par R.Dujardin [Du 1]. 
Enfin H.deThelin [DeT 1] a trait\'e le cas des surfaces k\"ahl\'eriennes 
non projectives:

\begin{thm}
Le courant $T^+$ est laminaire.
\end{thm}

\begin{esquisse}
On va supposer pour simplifier que $X=\P^2$ et $\l_2(f)=1$.
Soit $L \subset \P^2$ une droite projective g\'en\'erique.
Il r\'esulte de la Proposition 4.11 que $\l^{-n} (f^n)^* [L] \rightarrow T^+$.
Posons ${\mathcal C}_n:=f^{-n}L$. Ce sont des courbes singuli\`eres
de degr\'e $d_n=\l^n$ et de genre $1$ car
$f^{-n}_{|L}:L \simeq \P^1 \rightarrow {\mathcal C}_n$ est une r\'esolution (non minimale)
des singularit\'es de ${\mathcal C}_n$.

L'id\'ee est que si l'on projette ${\mathcal C}_n$ sur une 
droite g\'en\'erique $\P^1$ de $\P^2$, alors on peut esp\'erer \'ecrire
localement ${\mathcal C}_n$  comme un graphe au dessus de $\P^1$.
Il faut bien s\^ur se placer hors des singularit\'es et \'eviter
\'egalement les points o\`u ${\mathcal C}_n$ est tangente \`a
la projection. Nous allons d\'enombrer ceux-ci.

Commen\c{c}ons par les points singuliers de ${\mathcal C}_n$. Notons
$n_x({\mathcal C}_n)$ le nombre de composantes irr\'eductibles
locales de ${\mathcal C}_n$ au point $x$. Soit 
$\pi:\hat{{\mathcal C}_n} \rightarrow {\mathcal C}_n$ une r\'esolution minimale
des singularit\'es de ${\mathcal C}_n$. Alors
\begin{eqnarray*}
\sum_{x \in Sing({\mathcal C}_n)} n_x({\mathcal C}_n)
&=& \text{ nombre de points de } \pi^{-1}(Sing \, {\mathcal C}_n)  \\
&\leq& \sharp f^{n}(Sing \, {\mathcal C}_n) \cap L
\leq {\mathcal Crit}(f^{-n}) \cdot L,
\end{eqnarray*}
si l'on choisit $L$ de sorte qu'elle \'evite l'ensemble d\'enombrable
$\cup_{j  \geq 0} I_{f^{-j}}$.
Comme l'ensemble critique ${\mathcal Crit}(f^{-n})$ est de degr\'e 
$3\l^n-3$, on en d\'eduit
$$
\sum_{x \in Sing({\mathcal C}_n)} n_x({\mathcal C}_n)
\leq 3\l^n-3=O(d_n).
$$

Soit $\pi_p:\P^2 \setminus \{p\}: \rightarrow \P^1$ la projection holomorphe
sur une droite $\P^1$,
issue d'un point g\'en\'erique $p$. On d\'ecoupe $\P^1$ en morceaux $Q$
(``carr\'es'') de petit diam\`etre $\lesssim r$ (donc d'aire $\lesssim r^2$),
formant une subdivison ${\mathcal Q}$ de $\P^1$.
Soit $\Gamma$ une composante connexe de ${\mathcal C}_n \cap \pi^{-1}(Q)$.
On dit que $\Gamma$ est une {\it bonne composante} si c'est un graphe
au dessus de $Q$ (i.e. $\pi_p:\Gamma \rightarrow Q$ est un hom\'eomorphisme),
et que c'est une {\it mauvaise composante} sinon.
Observons que les mauvaises composantes sont celles qui contiennent
un point singulier de ${\mathcal C}_n$ ou un point $x$ en lequel
la droite $(px)$ est tangente \`a ${\mathcal C}_n$.
On peut contr\^oler le nombre de ces derniers en utilisant la formule
de Hurwitz, appliqu\'ee \`a l'application 
$\pi:=\pi_p \circ p: \hat{\mathcal C}_n \rightarrow \P^1$, o\`u 
$p: \hat{\mathcal C}_n \rightarrow {\mathcal C}_n$ est une r\'esolution
des singularit\'es de ${\mathcal C}_n$. On obtient
(voir Proposition 3.3 de [Du 1]),
$$
\sharp\{\text{mauvaises composantes}\}
\leq 4g_n+4d_n+
\sum_{x \in Sing({\mathcal C}_n)} n_x({\mathcal C}_n)=O(d_n),
$$
o\`u $g_n$ d\'esigne le genre de $\hat{\mathcal C}_n$ (il vaut 1 ici).
Consid\'erons alors
$$
T_{{\mathcal Q},n}:=\frac{1}{d_n} \sum_{Q \in {\mathcal Q}} 
\sum_{\Gamma b.c.} [\Gamma] \leq T_n:=\frac{1}{d_n}[{\mathcal C}_n].
$$
Le courant $T_{{\mathcal Q},n}$ est uniform\'ement laminaire (car les
bonnes composantes b.c. sont des graphes) et presque \'egal \`a $T_n$
puisque
\begin{equation}
0 \leq \langle T_n-T_{{\mathcal Q},n},\pi_p^* \om_{\P^1} \rangle
\leq C r^2,
\end{equation}
o\`u $\om_{\P^1}$ d\'esigne la forme de Fubini-Study sur la droite
de projection.

Soit $T_{{\mathcal Q}}$ une valeur d'adh\'erence de la suite
$T_{{\mathcal Q},n}$. Un argument de familles normales permet
de montrer que $T_{{\mathcal Q}}$ est uniform\'ement laminaire.
En raffinant la subdivision ${\mathcal Q}$ (donc en faisant tendre
$r$ vers $0$), on obtient une suite croissante de courants
uniform\'ement laminaires de masses uniform\'ement born\'ees.
Soit $T_{\infty}$ la limite croissante de ces courants.
Il r\'esulte de (4.4) que
$$
\langle T_{\infty}, \pi_p^* \om_{\P^1} \rangle=
\langle T^+, \pi_p^* \om_{\P^1} \rangle.
$$
Pour un choix g\'en\'erique de $p$, cela implique
$T^+=T_{\infty}$. Nous renvoyons le lecteur \`a [Du 1]
pour plus de d\'etails.
\end{esquisse}

\begin{rqe}
Le r\'esultat \'enonc\'e par R.Dujardin dans [Du 1] fait intervenir une
hypoth\`ese technique, la condition (H), qui est superflue (voir [DDG]).

Lorsque $X$ n'est pas projective (ainsi $kod(X)=0$ d'apr\`es le Th\'eor\`eme 2.14),
on ne peut pas utiliser de projection
sur une courbe comme nous l'avons esquiss\'e. Il faut faire une analyse
locale. Celle-ci a \'et\'e men\'ee \`a bien par H.deThelin dans [DeT 1],
o\`u la th\'eorie d'Ahlfors remplace la formule de Hurwitz.
Cependant l'estimation (4.4) de la masse de $T^+-T_{{\mathcal Q}}$ est moins
pr\'ecise que dans le cas projectif. Cela pose
probl\`eme lorsque l'on souhaite intersecter g\'eom\'etriquement
deux tels courants (voir paragraphe 4.4.3)..
\end{rqe}

\subsection{Courants $f_*$-invariants}

Nous souhaitons \`a pr\'esent construire un courant $f_*$-invariant $T^-$ 
qui va jouer un r\^ole analogue \`a celui de $T^+$.
Lorsque $\l_2(f)=1$, il suffit d'appliquer la construction pr\'ec\'edente
\`a $f^{-1}$. Lorsque $\l_2(f) \geq 2$, 
certaines difficult\'es techniques suppl\'ementaires apparaissent,
r\'esolues dans [DDG]:

\begin{thm}
Il existe un unique courant positif ferm\'e $T^-$ de bidegr\'e $(1,1)$ sur $X$
tel que $f_*T^-=\l T^-$, $\{T^-\}=\a^-$ et
$$
\frac{1}{\l^n}(f^n)_* \om \longrightarrow T^-.
$$
\end{thm}

\begin{esquisse}
La construction est tout \`a fait analogue \`a celle de $T^+$.
Le point de d\'epart est encore le $dd^c$-lemma, qui
donne ici
$$
\frac{1}{\l^n} (f^n)_* \theta^-=\theta^-+dd^c g_n^-,
\text{ avec } g_n^-:=\sum_{j=0}^{n-1} \frac{1}{\l^j} (f^j)_* \g_-,
$$
pour une fonction $\g^- \in L^1(X,\R)$.
Comme dans la preuve du Th\'eor\`eme 4.9, la difficult\'e est de montrer
que la fonction $\g^-$
est essentiellement major\'ee, ce que
nous obtenons dans [DDG] en montrant que
la classe $\a^-$ peut-\^etre repr\'esent\'ee par un courant positif
ferm\'e qui admet un potentiel born\'e.
On peut alors supposer que $\g^-$ est n\'egative,
et on obtient une suite d\'ecroissante de potentiels $(g_n^-)$
qui converge dans $L^1(X)$.

Certains cas particuliers sont plus faciles \`a traiter:
lorsque $X$ est une surface rationnelle
minimale, la classe $\a^-$ est semi-positive et la construction
de $T^-$ est obtenue dans [G 1].
Lorsque $kod(X)=0$, on peut contourner la difficult\'e en utilisant
une forme volume invariante.
On travaille dans ce cas sur un mod\`ele minimal. Comme $12K_X=0$,
il existe une 2-forme holomorphe $\eta$ sur $X$
--essentiellement unique-- qui ne s'annule nulle part.
La forme $f^* \eta$ est une 2-forme holomorphe dans $X \setminus I_f$
qui se prolonge de fa\c{c}on holomorphe \`a travers $I_f$,
car $codim_{\C} I_f \geq 2$. On a donc $f^* \eta= \zeta \eta$,
avec $\zeta \in \C$. Soit alors $\nu:=i \eta \wedge \overline{\eta}$:
c'est une forme volume sur $X$ telle que $f^* \nu=|\zeta|^2 \nu$.
On a donc $|\zeta|^2=\l_2(f)$ et $f^* \nu=\l_2(f) \nu$.
On en d\'eduit la convergence de $g_n^-$ dans $L^1(\nu)$, car
$$
\sum_{j \geq 0} \frac{1}{\l^j} \int_X |(f^j)_* \g^-| d\nu
=\sum_{j \geq 0} \left( \frac{\l_2(f)}{\l} \right)^j \int |\g^-| d\nu <+\infty,
$$
puisque $\l=\l_1(f)>\l_2(f)$.
\end{esquisse}

\vskip.2cm
Il serait int\'eressant de caract\'eriser l'ensemble des courants $f_*$-invariants,
mais nous ne disposons pas d'estim\'ees de volume pour les images r\'eciproques
lorsque $\l_2(f) \geq 2$. 
Lorsque $\l_2(f)=1$, $f_*=(f^{-1})^*$ et on peut appliquer les r\'esultats
de la section 4.3.2.
Le courant $T^-$ jouit cependant d'une propri\'et\'e
ergodique forte, qui s'apparente au m\'elange [DDG]:

\begin{thm}
Le courant $T^-$ est extr\'emal dans le c\^one 
${\mathcal T}_{f_*}$ des courants $f_*$-invariants.
Si $kod(X)=0$ et $(\a^-)^2=0$,
alors $T^-$ est un point extr\'emal du c\^one ${\mathcal T}$.
\end{thm}

\begin{preuve}
La preuve de l'extr\'emalit\'e de $T^-$ dans le c\^one 
${\mathcal T}_{f_*}$ des $(1,1)$-courants positifs ferm\'es
$S$ tels que $f_*S=\l S$ est analogue \`a celle du Th\'eor\`eme 4.12.

Supposons \`a pr\'esent que $X$ est de dimension de Kodaira nulle
et que $(\a^-)^2=0$.
Soit $S$ un courant positif ferm\'e de bidegr\'e $(1,1)$ tel que
$0 \leq S \leq T^-$. Il s'agit de montrer que $S$ est proportionnel \`a $T^-$.
L'hypoth\`ese $(\a^-)^2=0$ assure que $\a^-$ est une classe extr\'emale
dans $H^{1,1}_{nef}(X,\R)$ (donc $S$ est cohomologue \`a $c \theta^-$,
$0 \leq c \leq 1$), et que c'est un vecteur propre pour l'op\'erateur $f^*$,
$$
f^* \a^-=\frac{\l_2}{\l} \a^- \; \text{ et } \;
\frac{1}{\l_2}f_*f^* \a^-=\a^-.
$$
Cela r\'esulte de la formule d'aller-retour et nous renvoyons le lecteur
\`a [DDG] pour plus de d\'etails.
Posons
$$
T_j^-:=\frac{\l^j}{\l_2^j} (f^j)^* T^- \; \; 
\text{ et }
\; \; S_j:=\frac{\l^j}{\l_2^j} (f^j)^* S.
$$
Les courants $T_j^-$ et $S_j$ sont encore cohomologues \`a $\theta^-$, $c\theta^-$
et v\'erifient
$$
\frac{1}{\l^j} (f^j)_* T_j^- \equiv T^- \;
\text{ et } \; 
\frac{1}{\l^j} (f^j)_* S_j \equiv S
$$
Nous allons montrer que la suite de courants positifs $(\l^{-j}(f^j)_*S_j)$ converge
vers $cT^-$, ce qui assurera $S=cT^-$.

Posons $R_j:=T_j^--S_j \geq 0$ et fixons $u_j,v_j,w_j \in L^1(X)$ telles que
$$
T_j^-=\theta^-+dd^c u_j, \;
S_j=c\theta^-+dd^c v_j, \;
\text{ et }
R_j=(1-c)\theta^-+dd^c w_j.
$$
On normalise les fonctions $u_j,v_j,w_j$ par
$\int_X u_j d\nu=\int_X v_j d\nu=\int_X w_j d\nu=0$.
Comme la normalisation est lin\'eaire, on a $u_j=v_j+w_j$. Cette normalisation implique \'egalement
que les suites $(u_j),(v_j),(w_j)$ sont relativement compactes dans $L^1(X)$
(voir [GZ 1]), en particulier ces suites sont toutes uniform\'ement major\'ees sur $X$.

Il r\'esulte de l'\'egalit\'e $\l^{-j}(f^j)_*T_j^-=T^-$ 
et de l'invariance de $\nu$ que
$$
\tilde{u_j}:=
\frac{1}{\l^j}(f^j)_* u_j=g^--g_j^-
 \longrightarrow 0 \text{ dans } L^1(X).
$$
L'in\'egalit\'e $u_j-C \leq v_j \leq C$ donne ainsi
$$
\tilde{u}_j-C \frac{\l_2^j}{\l^j} \leq \frac{1}{\l^j}(f^j)_* v_j \leq C \frac{\l_2^j}{\l^j},
$$
d'o\`u $\tilde{v_j}:=\l^{-j}(f^j)_* v_j \rightarrow 0$.
Il s'ensuit que
$$
S=\l^{-j}(f^j)_* S_j=
\lim_{j \rightarrow +\infty} \left[ c \l^{-j}(f^j)_*\theta^-+ dd^c ( \tilde{v_j}) \right]
=c T^-.
$$
\end{preuve}
\vskip.2cm

Remarquons que la condition $(\a^-)^2=0$ est 
v\'erifi\'ee lorsque
$f$ est holomorphe (voir preuve du Th\'eor\`eme 4.7).

\section{La mesure canonique}

Dans toute la suite du chapitre 4, nous supposons que $f:X \rightarrow X$
est un endomorphisme m\'eromorphe 1-stable 
tel que $\l:=\l_1(f)>\l_2(f)$.
Nous avons contruit deux courants positifs canoniques invariants 
$T^+,T^-$ qui repr\'esentent les classes de cohomologie $\a^+,\a^-$.
Rappelons que nous avons normalis\'e ces classes et impos\'e
$$
\a^+ \cdot \a^-=\{ T^+ \} \cdot \{ T^- \}=1.
$$

\subsection{D\'efinition pluripotentialiste}

Nous souhaitons \`a pr\'esent d\'efinir une mesure invariante $\mu_f$ en consid\'erant
le produit d'intersection $\mu_f:=T^+ \wedge T^-$.
On ne peut pas toujours donner un sens \`a un tel produit:
les courants $T^{\pm}$ sont des formes diff\'erentielles \`a coefficients distributions,
et il n'est pas toujours possible de d\'efinir le produit de deux distributions.
On y arrive  cependant lorsque les potentiels $g^{\pm}$ de $T^{\pm}$ ne sont
pas trop singuliers (voir [BT], [De]). Lorsque $g^+$ est int\'egrable par rapport \`a la
mesure trace de $T^-$, on peut ainsi d\'efinir le courant $g^+ T^-$ et poser
$$
\mu_f:=\theta^+ \wedge T^-+dd^c (g^+ T^-).
$$
Notons que cette d\'efinition est sym\'etrique au sens o\`u $g^+ \in L^1(T^-)$
si et seulement si $g^- \in L^1(T^+)$, comme on le v\'erifie en
int\'egrant par parties,
$$
\int_X g^+ T^- \wedge \om=
\int_X g^- T^+ \wedge \om+\int_X [g^+ \theta^--g^-\theta^+] \wedge \om.
$$
Nous ne connaissons \`a l'heure actuelle aucun exemple pour lequel cette condition n'est pas v\'erifi\'ee.
Cela justifie la

\begin{ques}
A-t-on toujours  $g^+ \in L^1(T^- \wedge \om)$ ? 
\end{ques}

En supposant cette condition satisfaite, nous montrons, en collaboration
avec J.Diller dans [DG], que la mesure
$\mu_f$ est une mesure de probabilit\'e
qui ne charge pas les points d'ind\'etermination: c'est donc une mesure invariante, 
$f_* \mu_f=\mu_f$, comme on peut le voir en approximant $\mu_f$
par les mesures \`a densit\'e $\mu_n:=\l^{-2n}(f^n)^* \theta^+ \wedge (f^n)_* \theta^-$.

Il reste \`a d\'eterminer quand cette condition est effectivement satisfaite. C'est 
une motivation de l'article [DG] qui fait le point 
sur les propri\'et\'es de r\'egularit\'e des fonctions de Green dynamiques $g^{\pm}$.
Lorsque $f$ est holomorphe, ces fonctions sont h\"olderiennes (m\^eme preuve
que pour la Proposition 1.2) et la condition est donc satisfaite. 
Nous avons vu
cependant au paragraphe 3.3.1 que de tels endomorphismes
sont rares (voir corollaire 3.7).

Lorsque $f$ est seulement m\'eromorphe, les points d'ind\'etermination ont tendance \`a cr\'eer
des singularit\'es logarithmiques. Or il peut y en avoir beaucoup (Exemple 4.35). Il faut
donc s'int\'eresser \`a des notions de r\'egularit\'e plus faibles. 
Si les fonctions $g^{\pm}$ sont de gradient $L^2$, alors $g^+ \in L^1(T^- \wedge \om)$. 
En effet, on peut supposer $g^+ \leq 0$ et l'in\'egalit\'e de
Cauchy-Schwarz implique
\begin{eqnarray}
0 &\leq&  \int_X (-g^+) T^- \wedge \om =O(1)+\int_X d g^+ \wedge d^c g^- \wedge \om \\
&\lesssim & \left( \int dg^+ \wedge d^c g^+ \wedge \om \right)^{1/2}
\left(  \int d g^- \wedge d^c g^- \wedge \om \right)^{1/2}<+\infty.
\end{eqnarray}
Nous analysons cette condition dans [DG].
Elle n'est pas toujours satisfaite lorsque $X$ est rationnelle (voir Exemple 4.36).
Lorsque  $kod(X)=0$, on a le r\'esultat suivant [DDG].

\begin{pro}
Si $kod(X)=0$ alors $\nabla g^{-} \in L^2(X)$.
\end{pro}

\begin{preuve}
Rappelons que nous travaillons ici sur le mod\`ele minimal.
Dans ce cas $f$ est non-ramifi\'e, donc 1-stable (Proposition 4.5).
De plus l'op\'erateur $f_*$ pr\'eserve l'ensemble des fonctions
continues. Il s'ensuit que la fonction $\g^-$ de
la preuve du Th\'eor\`eme 4.17 est en r\'ealit\'e continue,
et que la suite de fonctions continues 
$(g_n^-)$ converge uniform\'ement vers $g^-$.

On v\'erifie ais\'ement que le gradient d'une fonction quasiplurisousharmonique
born\'ee (en particulier continue) est dans $L^2(X)$.
\end{preuve}

\vskip.2cm
Dans un article r\'ecent [BDi 2], E.Bedford et J.Diller ont \'etudi\'e une condition
l\'eg\`erement plus forte qui leur permet de mener \`a bien une partie de l'\'etude
dynamique des endomorphismes birationnels.

\begin{defi}
L'endomorphisme $f$ v\'erifie la condition (BD) si la fonction
$g^-$ est finie en tout point d'ind\'etermination de $f$.
\end{defi}

Il s'av\`ere que cette condition est sym\'etrique en $f/f^{-1}$. Elle implique que
les fonctions
$g^{\pm}$ sont de gradient $L^2$, et m\^eme un peu plus (voir [BDi 2], [DG]).
Les diff\'erentes conditions rencontr\'ees jusqu'\`a pr\'esent peuvent 
s'interpr\'eter en termes de propri\'et\'es de r\'ecurrence des
points d'ind\'etermination:
\vskip.2cm

\begin{itemize}
\item $f$ est 1-stable ssi $\log \text{dist}(I_{f^{-1}},f^{-n}I_f) >-\infty$
pour tout $n \geq 1$;
\vskip.2cm

\item $\nabla g^{+} \in L^2$ ssi $\sum_{n \geq 1} \l^{-2n} \log \text{dist}(I_{f^{-1}},f^{-n}I_f) >-\infty$;
\vskip.2cm

\item $f$ satisfait la condition (BD) ssi
$\sum_{n \geq 1} \l^{-n} \log \text{dist}(I_{f^{-1}},f^{-n}I_f) >-\infty$.
\end{itemize}
\vskip.2cm

\noindent
Lorsque $f$ est 1-stable, l'ensemble d'ind\'etermination it\'er\'e
$I_{f^n}$ est \'egal \`a $I_f \cup f^{-1} I_f \cup \cdots f^{-n+1} I_f$.
Les conditions de r\'egularit\'e ci-dessus mesurent donc \`a quelle vitesse
les ensembles $I_{f^n}$ et $I_{f^{-n}}$ se rapprochent l'un de l'autre.
Nous renvoyons le lecteur \`a [BDi 2], [DG], [DDG] pour la justification 
de ces \'equivalences et la preuve du r\'esultat suivant:

\begin{pro} Si $f$ v\'erifie la condition (BD) alors $\log ||Df|| \in L^1(\mu_f)$.

Si $kod(X)=0$ alors $f$ v\'erifie la condition (BD).

Si $f:\P^2 \rightarrow \P^2$ est birationnel et v\'erifie (BD), alors
$f_A:=A \circ f$ est birationnel et v\'erifie (BD) pour toute matrice
$A \in PGL(3,\C)$ hors d'un ensemble pluripolaire de $PGL(3,\C)$.
\end{pro}

Nous verrons plus loin (Exemple 4.36) des exemples d'endomorphismes
birationnels 1-stables de $\P^2$ qui ne satisfont pas la condition (BD).

\subsection{Ergodicit\'e}

Dans toute la suite nous supposons que $g^+ \in L^1(T^- \wedge \om)$.

\begin{thm}
La mesure $\mu_f$ est ergodique.

Si $T^-$ est un courant extr\'emal, alors $\mu_f$ m\'elange.
\end{thm}

\begin{rqe}
Ce r\'esultat a \'et\'e \'etabli pour les applications de H\'enon complexes par 
E.Bedford et J.Smillie [BS 2].
La preuve indiqu\'ee ci-dessous reprend les simplifications apport\'ees par
J.E.Fornaess et N.Sibony dans [FS 3] et par l'auteur et N.Sibony dans [GS].

Notons que $\mu_f$ est m\'elangeante lorsque $\l_2(f)=1$.
Nous indiquerons au paragraphe 5.1.1 comment montrer
que $\mu_f$ est toujours faiblement m\'elangeante.
\end{rqe}

\begin{esquisse}
Supposons $T^-$ extr\'emal et montrons que $\mu_f$ m\'elange.
Soit $\chi,\p$ des fonctions test. On ne perd rien en supposant que $0 \leq \chi \leq 1$.
Il s'agit de montrer que 
$$
\int_X \chi \p \circ f^n d\mu_f \longrightarrow c_{\chi}c_{\p},
$$
o\`u $c_{\chi}=\int \chi d\mu_f$, $c_{\p}=\int \p d\mu_f$. On observe que 
$$
\int \chi \p \circ f^n d\mu_f=\langle  \l^{-n}(f^n)^*(\p T^+), \chi T^-\rangle 
=\langle \p T^+,\l^{-n}(f^n)_*(\chi T^-)\rangle .
$$
Il s'agit donc de  montrer que la suite de mesures 
$(\mu_n)$ converge faiblement vers la mesure $c_{\chi} \mu_f$, o\`u 
$$
\mu_n= R_n \wedge T^+, \; \;    \text{ avec }
R_n:=\frac{1}{\l^n}(f^n)_*(\chi T^-).
$$
Les courants $R_n$ sont des courants positifs non-ferm\'es qui sont domin\'es par $T^-$.
On montre que leur bord tend en masse vers 0 (voir Lemme 4.25 ci-dessous). Il s'ensuit
que toute valeur d'adh\'erence $R=\lim R_{n_i}$ de la suite $(R_n)$ est un courant
positif ferm\'e domin\'e par $T^-$.
Comme $T^-$ est extr\'emal, il vient $R=c_R T^-$. Or on calcule
$$
c_R=\langle R,\theta^+\rangle =\lim_{n_i \rightarrow +\infty} \langle R_{n_i},\theta^+\rangle =c_{\chi},
$$
d'o\`u $c_R$ est ind\'ependante de $R$. C'est donc que la suite
$(R_n)$ converge en fait vers le courant $c_{\chi}T^-$.

Si $T^+$ \'etait une forme lisse, on en d\'eduirait imm\'ediatement que
$\mu_n=R_n \wedge T^+$ converge vers $c_{\chi} T^- \wedge T^+$.
Ce n'est pas le cas, mais $T^+$ est tr\`es bien approch\'e par les
formes  $\theta_j^+:=\l^{-j}(f^j)^* \theta^+$: $j$ \'etant 
fix\'e, on obtient $R_n \wedge \theta_j^+ \rightarrow c_{\chi}T^- \wedge \theta_j^+$
car $\theta_j^+$ est lisse (hors d'un nombre fini de points),
puis on montre que $R_n \wedge (T^+-\theta_j^+)$ tend vers $0$
uniform\'ement par rapport \`a $n$, lorsque $j$ tend vers $+\infty$
(voir la preuve du Th\'eor\`eme 4.1 dans [GS] pour les d\'etails).

Lorsque $T^-$ est seulement extr\'emal dans le c\^one 
${\mathcal T}_{f_*}$ des courants $f_*$-invariants, on peut
reprendre les arguments ci-dessus en rempla\c{c}ant partout
$\mu_n$ par $\mu_n':=n^{-1}\sum_{j=0}^{n-1} \mu_j$ et
$R_n$ par $R_n':=n^{-1} \sum_{j=0}^{n-1} R_j$.
Les valeurs d'adh\'erences de la suite $(R_n')$ sont des courants
positifs ferm\'es {\it invariants}, car $f_*R_j=\l R_{j+1}$.
On en d\'eduit la convergence de $R_n'$ vers $c_{\chi} T^-$, puis celle
de $\mu_n'$ vers $c_{\chi} \mu_f$, ce qui montre que $\mu_f$ est ergodique.
\end{esquisse}
\vskip.2cm

Pour faire fonctionner la m\'ethode pr\'ec\'edente, il est  n\'ecessaire de
contr\^oler la masse des courants $dR_n$ et $dd^c R_n$. C'est l'objet du lemme suivant
d\'emontr\'e par E.Bedford et J.Smillie lorsque $f$ est une application de H\'enon
complexe (voir Lemme 1.2 dans [BS 2]) et qui se g\'en\'eralise ais\'ement
\`a notre cadre.

\begin{lem} 
On a $||dR_n||=O(\l^{-n/2})$ et $||dd^c R_n||=O(\l^{-n})$.
\end{lem}

\begin{rqe}
L'estimation de la vitesse de m\'elange est plus d\'elicate
que dans le cas des endomorphismes de grand degr\'e topologique.
Elle a \'et\'e obtenue par T.C.Dinh [Di 2] pour les applications de H\'enon complexes.
Elle n'est pas connue dans le cadre plus g\'en\'eral des endomorphismes birationnels qui v\'erifient
la condition (BD).
\end{rqe}

\subsection{D\'efinition g\'eom\'etrique}

Pour pousser plus avant l'\'etude des propri\'et\'es ergodiques de la
mesure $\mu_f$, il est n\'ecessaire de donner une interpr\'etation
g\'eom\'etrico-dynamique de celle-ci.
Nous supposons pour simplifier dans la suite que $\l_2(f)=1$.
Dans ce cas les courants $T^+$ et $T^-$ sont laminaires
(Th\'eor\`eme 4.15).

Si $S^+,S^-$ sont deux courants uniform\'ement laminaires,
localement d\'ecrits par une moyenne sur des disques 
$$
S^{\pm}=\int [\Delta_a^{\pm}] d\nu^{\pm}(a),
$$
on peut d\'efinit leur {\it intersection g\'eom\'etrique} par
$$
S^+ \dot{\wedge} S^-:=
\int \int [\Delta^+_a \cap \Delta^-_b] d\nu^+(a) d\nu^-(b).
$$
On v\'erifie que cette d\'efinition co\"{\i}ncide avec la d\'efinition
pluripotentialiste (voir Th\'eor\`eme 3.1 dans [Du 2]).

Lorsque les courants $S^+,S^-$ sont seulement laminaires,
on d\'efinit leur intersection g\'eom\'etrique
$S^+ \dot{\wedge} S^-$ comme limite des intersections
g\'eom\'etriques de leurs approximants uniform\'ements 
laminaires. Malheureusement cette d\'efinition ne co\"{\i}ncide pas
toujours avec la d\'efinition pluripotentialiste (voir
le second paragraphe de [Du 2] pour des exemples).
Lorsqu'on a un contr\^ole quantitatif de l'approximation
de $S^{\pm}$ par des courants uniform\'ement laminaires
(comme c'est le cas dans le Th\'eor\`eme 4.15, cf
(4.4)), on peut esp\'erer montrer que les deux
notions de produit ext\'erieur co\"{\i}ncident.

R.Dujardin justifie cette attente dans [Du 2,4,5] 
en d\'emontrant:

\begin{thm}
Supposons $X$ projective.
Si $f$ v\'erifie la condition (BD) alors les courants 
$T^+$ et $T^-$ s'intersectent g\'eom\'etriquement.
\end{thm}

\begin{esquisse}
On approxime $T^{\pm}$ par des courants $T_r^{\pm}$ 
qui sont uniform\'ement laminaires associ\'es \`a une
subdivison ${\mathcal Q}$ en cubes $Q$ de taille $r$, de
telle sorte que la masse des diff\'erences satisfait
$$
|| T^{\pm}-T_r^{\pm}|| \leq C r^2,
\text{ o\`u } T_r^{\pm}:=\sum_{Q \in {\mathcal Q}} T_Q^{\pm}.
$$
Il s'agit d'estimer la masse de la diff\'erence entre
$T^+ \wedge T^-$ (d\'efinition pluripotentialiste) et
$\sum_Q T_Q^+ \wedge T_Q^-$ (intersection g\'eom\'etrique).

Quitte \`a l\'eg\`erement bouger la subdivision ${\mathcal Q}$,
on peut supposer que la mesure $T^+ \wedge T^-$ ne charge pas le
bord des cubes $Q$. Soit $\chi$ une fonction test qui
s'annule pr\`es du bord des cubes et est
identique \`a 1 dans la plus grande partie de $Q$ (les d\'eriv\'ees
d'ordre 1 de $\chi$ explosent donc comme $C/r$, lorsque $r \rightarrow 0$).
Il nous faut montrer que la
quantit\'e
$$
\int_X \chi (T^+ \wedge T^- -T_r^+ \wedge T_r^-)=
\sum_Q \int_Q \chi (T^+ \wedge T^- -T_Q^+ \wedge T_Q^-)
$$
tend vers $0$ lorsque $r \rightarrow 0$.
Comme $\chi$ est \`a support compact dans chaque cube $Q$, on peut
traiter les int\'egrales s\'epar\'ement.
Par sym\'etrie il suffit de contr\^oler
la masse de $T^+ \wedge (T^- -T_Q^-)$ dans $Q$. Or $T^+=dd^c G^+$
dans $Q$, pour un potentiel local $G^+$ qui est de gradient
$L^2$ par rapport \`a $T^-$ (c'est la propri\'et\'e (BD)). 
L'in\'egalit\'e de Cauchy-Schwarz implique alors
\begin{eqnarray*}
\lefteqn{
\int \chi dd^c G^+ \wedge (T^- -T_Q^-)
= -\int d\chi \wedge d^c G^+ \wedge (T^- -T_Q^-) }\\
&\leq& \left( \int dG^+ \wedge d^c G^+ \wedge (T^- -T_Q^-)\right)^{1/2}
 \left( \int d\chi \wedge d^c \chi \wedge (T^- -T_Q^-)\right)^{1/2} \\
&\leq& C \left( \int dG^+ \wedge d^c G^+ \wedge (T^- -T_Q^-)\right)^{1/2}.
\end{eqnarray*}
Lorsque $r \rightarrow 0$, $T_Q^-$ cro\^it vers $T^-$ et
cette derni\`ere int\'egrale converge vers $0$
(il suffit d'approximer $dG^+$ par une suite 
de formes lisses dans $L^2(T^-)$) .
\end{esquisse}

\section{Propri\'et\'es ergodiques de $\mu_f$}

\noindent Nous testons ici la conjecture dans l'une des deux situations
suivantes:
\begin{enumerate}
\item {\it Hypoth\`eses:} $kod(X)=-\infty$, $\l_2(f)=1$ et condition (BD).
Ainsi,
 \begin{itemize}
 \item on peut supposer que $f$ est 1-stable (Th\'eor\`eme 4.6);
 \item on sait construire les courants invariants $T^+,T^-$ (Th\'eor\`eme 4.9);
 \item ils sont extr\'emaux et laminaires (Th\'eor\`emes 4.12 et 4.15);
 \item la mesure $\mu_f:=T^+ \wedge T^-$ est bien d\'efinie (voir (4.6)),
 m\'elangeante (Th\'eor\`eme 4.23) et g\'eom\'etrique (Th\'eor\`eme 4.27);
 \item les exposants de Lyapunov de $\mu_f$ sont bien d\'efinis (Proposition 4.22).
 \end{itemize}
La conjecture a \'et\'e d\'emontr\'ee dans ce contexte par E.Bedford, M.Lyubich et J.Smillie
lorsque $f$ est une application de H\'enon complexe [BLS 1,2] et partiellement
d\'emontr\'ee par E.Bedford, J.Diller [BDi 2], et R.Dujardin  [Du 5] sous ces hypoth\`eses.

\item {\it Hypoth\`eses:} $kod(X)=0$ et $X$ est projective. 
Ainsi
  \begin{itemize}
  \item $f$ est 1-stable si on travaille sur le mod\`ele minimal (Proposition 4.5);
  \item on sait construire les courants $T^+$, $T^-$ (Th\'eor\`emes 4.9 et 4.17);
  \item $T^+$ est laminaire,
   $T^-$ est extr\'emal si $(\a^-)^2=0$ (Th. 4.15 et 4.18);
  \item la condition (BD) est toujours satisfaite (Proposition 4.22), donc $\mu_f=T^+ \wedge T^-$
   est bien d\'efinie, ergodique (et m\^eme m\'elangeante si $(\a^-)^2=0$, Th\'eor\`eme 4.23),
   et g\'eom\'etrique (Th\'eor\`eme 4.27);
  \item les exposants de Lyapunov de $\mu_f$ sont bien d\'efinis (Proposition 4.22).
  \end{itemize}
La conjecture a \'et\'e d\'emontr\'ee dans ce contexte par S.Cantat [Ca 1] lorsque
$\l_2(f)=1$ -- c'est \`a dire lorsque $f$ et un automorphisme --, et partiellement 
d\'emontr\'ee dans [DDG] lorsque $\l_2(f) \geq 2$.
\end{enumerate}

\subsection{Lamination stable}

Tout courant laminaire $T$ (voir d\'efinition 4.14) est localement repr\'esentable
\begin{equation}
T=\int_{\a \in \mathcal A} [\Delta_{\a}] d\nu(\a)
\end{equation}
par une moyenne sur une famille de disques holomorphes $\Delta_{\a}$.
En g\'en\'eral on ne peut pas lui associer une structure laminaire convenable
(voir [Du 1,2,4] pour des exemples), mais c'est possible lorque le courant
$T$ est {\it fortement approximable} (i.e. redevable de l'estim\'ee quantitative (4.4)):
c'est le r\'esultat principal de [Du 4] qui assure que si
${\mathcal L}$ est une r\'eunion de disques $\Delta_{\a}$ (on dit que 
${\mathcal L}$ est une boite de flot), alors ${\mathcal L}$ a une structure
de lamination plong\'ee et $T_{|{\mathcal L}}$ est uniform\'ement laminaire.
Il s'ensuit que $T_{|{\mathcal L}}$ induit une mesure transverse invariante
par holonomie qui est ergodique lorsque $T$ est extr\'emal
(voir Th\'eor\`eme 2.4 dans [Du 5]).

Lorsque $T=T^+$ est le courant canonique $f^*$-invariant, il r\'esulte du
Th\'eor\`eme 4.15 que $T^+$ est fortement approximable. Nous montrons \`a 
pr\'esent que les disques $\Delta_{\a}^+$ intervenant dans la repr\'esentation
laminaire de $T^+$ (on parle de disques subordonn\'es \`a $T^+$) sont
des morceaux de vari\'et\'es stables.

Soit $\tau$ une union finie de disques holomorphes transverses \`a une boite de flot
${\mathcal L}$. La th\'eorie du slicing assure que la mesure 
$T^+_{|{\mathcal L}} \wedge [\tau]$ est bien d\'efinie pour un choix g\'en\'erique
de $\tau$. On note $T^+ \wedge \tau$ la mesure induite sur la transversale
$\tau$ (elle est invariante par holonomie).
Le lemme suivant, d\^u \`a R.Dujardin [Du 5], est l'analogue dans ce contexte
du lemme de Briend-Duval-Lyubich (Lemme 3.4).

\begin{lem}
Soit ${\mathcal L}=\{ \Delta_t^+, t \in \tau \}$ une boite de flot
subordonn\'ee \`a $T^+$ et $\tau$ transverse \`a ${\mathcal L}$. 
Pour tout $\e>0$, il existe $C_{\e} >0$ et une
transversale $\tau_{\e} \subset \tau$ tels que
$||T^+ \wedge \tau_{\e} || \geq (1-\e) ||T^+ \wedge \tau||$ et
$$
\text{Aire}(f^n \Delta_t^+) \leq \frac{C_{\e} n^2}{\l^n}, \; \; \; 
\forall t \in \tau_{\e}, \, \forall n \in \N.
$$
\end{lem}

\begin{esquisse}
Quitte \`a l\'eg\`erement bouger $\tau$, on peut supposer que
${\mathcal L}$ ne rencontre pas ${\mathcal C}_f \cup I_f$.
Dans ce cas $f({\mathcal L})$ est une boite de flot pour $T^+$
et $f \tau$ est transverse \`a $f({\mathcal L})$, car $f$ est 
localement biholomorphe et $f_*T^+=\l^{-1}T^+$ dans
$X \setminus ( {\mathcal C}_f \cup I_f)$. On en d\'eduit
$$
f_*(T^+ \wedge \tau)=f_* T^+ \wedge f \tau=\l^{-1} T^+ \wedge f\tau.
$$

Quitte \`a enlever un ensemble de $\tau$-mesure nulle, on peut supposer
que $\tau$ ne rencontre pas $\cup_{n \geq 0} ( I_{f^n} \cup {\mathcal C}_{f^n})$.
On a donc comme pr\'ec\'edemment
$(f^n)_* (T^+ \wedge \tau)=\l^{-n} T^+ \wedge f^n \tau$.
Deux ``disques'' $f^n \Delta_{f^{-n}t_i}$ ne peuvent se rencontrer qu'en un nombre fini
de points lorsque $t_1 \neq t_2$ varient dans $f^n \tau$.
Il s'ensuit que pour tout $n \in \N$,
$$
\int_t \text{Aire} f^n (\Delta_{f^{-n}t}) d (T^+ \wedge f^n \tau)(t) \leq ||T^+||=1.
$$
Comme $T^+ \wedge f^n \tau=\l^{n} (f^n)_* (T^+ \wedge \tau)$
et $(f^n)^* (\text{Aire} f^n(\Delta_{f^{-n}t}))=\text{Aire}(f^n \Delta_t)$,
on en d\'eduit
$$
\int_t \text{Aire}(f^n \Delta_t) d(T^+ \wedge \tau)(t) \leq \l^{-n}.
$$
La plupart des disques $\Delta_t$ ont donc une aire qui d\'ecroit
en $\l^{-n}$ sous it\'eration de $f$.
Nous renvoyons au Lemme 3.2 de [Du 5], pour 
plus de d\'etails.
\end{esquisse}

\subsection{Exposants de Lyapunov}

Il r\'esulte du lemme pr\'ec\'edent que
le plus grand exposant de Lyapunov $\chi^+(\mu_f)$
de $\mu_f$ v\'erifie l'estimation attendue (voir Th\'eor\`eme 4.6 de [Du 5]),
$$
\chi^+(\mu_f) \geq \frac{1}{2} \log \l >0.
$$
Lorsque $f$ est birationnel, on en d\'eduit bien s\^ur la majoration
$\chi^-(\mu_f) \leq -\frac{1}{2} \log \l<0$, en intervertissant les r\^oles
de $f,f^{-1}$. Lorsque $2 \leq \l_2(f) <\l_1(f)=\l$, il faut donner un argument 
suppl\'ementaire pour contr\^oler $\chi^-(\mu_f)$.
En voici un lorsque $X$ est de dimension de Kodaira nulle [DDG].

\begin{pro}
Si $kod(X)=0$ alors
$$
\chi^+(\mu_f)+\chi^-(\mu_f)=\frac{1}{2} \log \l_2(f).
$$
En particulier, lorsque $X$ est projective, 
$$
\chi^-(\mu_f) \leq -\frac{1}{2} \log [ \l_1(f) / \l_2(f)] <0<
\frac{1}{2} \log \l_1(f) \leq \chi^+(\mu_f).
$$
\end{pro}

\begin{esquisse}
Rappelons que lorsque $X$ est de dimension de Kodaira nulle, il existe
une forme volume $\nu$ de jacobien constant, $f^* \nu=\l_2(f) \nu$.
Elle est d\'efinie \`a l'aide d'une 2-forme holomorphe $\eta$ telle que
$f^* \eta=\zeta \eta$, o\`u $\zeta \in \C$ est tel que $|\zeta|^2=\l_2(f)$
(voir la preuve du Th\'eor\`eme 4.17).
Soit $x \in X$ un point $\mu_f$-g\'en\'erique. Il r\'esulte du Th\'eor\`eme
de r\'ecurrence de Poincar\'e que $f^n(x)$ est arbitrairement proche de $x$
pour une infinit\'e de $n \in \N$. L'invariance de la forme $\eta$
assure que pour ces indices,
$$
\frac{1}{n}\log |Jac(f^n)(x)| \simeq \log |\zeta|.
$$

Or le Th\'eor\`eme de r\'eduction d'Oseledec-Pesin (voir [KH]) assure
que pour un point $\mu_f$-g\'en\'erique, on a 
$$
\chi^+(\mu_f)+\chi^-(\mu_f) \simeq \frac{1}{n}\log |Jac(f^n)(x)|.
$$
L'\'egalit\'e $\chi^+(\mu_f)+\chi^-(\mu_f)=\log |\zeta|=\frac{1}{2}\log |\l_2(f)|$
s'ensuit. 

Lorsque $X$ est projective, il r\'esulte de l'estimation 
$\chi^+ \geq \frac{1}{2} \log \l$ que le plus petit exposant de
Lyapunov $\chi^-(\mu_f)$ est major\'e par
$-\frac{1}{2} \log \l/\l_2(f)$.
\end{esquisse}

\begin{rqe}
La d\'emonstration montre en fait que la formule
$\chi^+(\mu)+\chi^-(\mu)=\frac{1}{2} \log \l_2(f)$ est valable pour 
toute mesure invariante ergodique $\mu$. Lorsqu'on l'applique \`a
une combinaison lin\'eaire de masses de Dirac 
\'equidistri-bu\'ees selon un cycle p\'eriodique,
on en d\'eduit que le produit des valeurs propres est \'egal
\`a la racine carr\'ee du degr\'e topologique: il ne peut donc
y avoir ni cycle attractif, ni domaine de Siegel lorsque
le degr\'e topologique est $\geq 2$ (comparer avec [M 1]).
Il s'ensuit \'egalement que l'ensemble de Fatou est vide dans une telle situation.
\end{rqe}

\subsection{Entropie et point selles}

Le Lemme 4.28 permet d'interpr\'eter les disques $\Delta_t^{\pm}$ subordonn\'es
\`a $T^{\pm}$ comme des disques stables/instables.
Il r\'esulte du Th\'eor\`eme 4.27 que $\mu_f$ a une structure de produit
local par rapport aux vari\'et\'es stables/instables.
Des arguments classiques permettent alors de montrer que $\mu_f$ est d'entropie
maximale $\log \l$ (voir Th\'eor\`eme 4.6 dans [Du 5]) et que les points
p\'eriodiques selles s'\'equidistribuent selon $\mu_f$
(voir [BLS 2] et Th\'eor\`eme 5.4 dans [Du 5]).

\section{Exemples}

\subsection{Endomorphismes polynomiaux de $\C^2$}

Les automorphismes polynomiaux de $\C^2$ ont \'et\'e classifi\'es
par S.Friedland et J.Milnor dans [FrM].

\begin{thm}
Soit $f:\C^2 \rightarrow \C^2$ un automorphisme polynomial de $\C^2$
qui est {\it cohomologiquement hyperbolique}. Alors $f$ est conjugu\'e
\`a une application de H\'enon complexe, i.e. une composition finie
d'automorphismes de la forme
$$
h:(z,w) \in \C^2 \mapsto (P(z)-aw,z) \in \C^2, \; a \neq 0,
$$
o\`u $P$ est un polyn\^ome de degr\'e $d \geq 2$.
\end{thm}

La preuve s'appuie sur le Th\'eor\`eme de Jung qui d\'ecrit la structure
du groupe des automorphismes polynomiaux de $\C^2$. Une preuve
simple et \'el\'egante en a \'et\'e donn\'ee par S.Lamy dans [La 1].
La dynamique des applications de H\'enon complexes a \'et\'e abondamment
\'etudi\'ee, notamment par E.Bedford, J.Smillie [BS 1-3], M.Lyubich [BLS 1,2],
J.H.Hubbard  et ses coauteurs [Hu], [HOV 1,2], [HPV], N.Sibony et J.E.Fornaess [FS 1-4].
\vskip.1cm

Il est naturel de vouloir consid\'erer d'autres classes d'exemples.
On peut s'int\'eresser aux endomorphismes polynomiaux {\it birationnels}
(l'inverse n'est pas n\'ecessairement polynomial) comme dans [FaG],
mais la structure de ce semi-groupe est encore mal comprise (voir [Da], [SY]).
On peut au contraire vouloir conserver le caract\`ere propre des automorphismes,
mais autoriser un degr\'e topologique plus grand. S.Lamy a 
partiellement classifi\'e dans [La 2] les endomorphismes polynomiaux propres de
$\C^2$ de degr\'e topologique 2: ils sont conjugu\'es \`a $\tau \circ h$,
o\`u $h$ est un automorphisme polynomial de $\C^2$ et $\tau(z,w)=(z^2,w)$.
Il serait int\'eressant de pousser plus avant cette classification et
d'\'etudier la dynamique de ces applications.
\vskip.1cm

Nous avons \'etudi\'e dans [G 3] la dynamique des endomorphismes polynomiaux 
{\it quadratiques} de $\C^2$ --i.e. ceux qui induisent un endomorphisme m\'eromorphe
de $\P^2$ tel que $r_1(f)=2$. Nous y montrons en particulier le

\begin{thm}
Soit $f:(z,w) \in \C^2 \mapsto (P(z,w),Q(z,w)) \in \C^2$ un endomorphisme polynomial
quadratique de $\C^2$ (i.e. $P,Q$ sont des polyn\^omes de degr\'e $\leq 2$)
tel que $\l_1(f)>\l_2(f)$.
Alors $f$ s'\'ecrit, apr\`es conjugaison par un automorphisme polynomial de $\C^2$,
sous l'une des deux formes suivantes:
\begin{enumerate}
\item $f(z,w)=(w+c,zw+c')$, o\`u $c,c' \in \C$. Dans ce cas 
$f$ est 1-stable dans $\P^1 \times \P^1$ et on obtient $\l_1(f)=\frac{1+\sqrt{5}}{2}$,
$\l_2(f)=1$;
\item $f(z,w)=(w+c,w[w-az]+bz+c')$,
o\`u $a,b,c,c' \in \C$ avec $(a,b) \neq (0,0)$. Dans ce cas 
$f$ est 1-stable dans $\P^2$ et on obtient $\l_1(f)=2$, $\l_2(f)=1$.
\end{enumerate}
\end{thm}

Observons que ces endomorphismes sont birationnels 
et qu'ils v\'erifient la condition (BD) sur $\P^1 \times \P^1$
(classe 1) ou sur $\P^2$ (classe 2). 
Plus g\'en\'eralement, si $f:\C^2 \rightarrow \C^2$
est un endomorphisme polynomial tel que $\l_1(f)>\l_2(f)$, il est naturel
de se demander si $f$ v\'erifie la condition (BD) dans une compactification ad\'equate.
Nous proposons une question plus pr\'ecise.

\begin{ques}
Soit $f:\C^2 \rightarrow \C^2$ un endomorphisme polynomial tel que $\l_1(f)>\l_2(f)$.
Peut-on trouver une compactification $X=\C^2 \cup D_{\infty}$ de $\C^2$
telle que $f$ s'\'etende en un endomorphisme de $X$ avec
$$
\! \!\! \! \! \!\! \! \! \!\! \! 
(\dag) \hskip1cm
f^m(X \setminus I_{f^m})=q_{\infty} \text{ est un point fixe superattractif de } f ?
$$
Ici $m$ d\'esigne le nombre de composantes irr\'eductibles du diviseur $D_{\infty}$.
\end{ques}

Il r\'esulte imm\'ediatement de $(\dag)$ que $f^m$ est 1-stable dans $X$ et
v\'erifie la condition (BD), puisque $f^{-m}(I_{f^m}) \subset I_{f^m}$.
Nous avons v\'erifi\'e $(\dag)$ dans [FaG] pour les endomorphismes birationnels
de $\C^2$ qui sont engendr\'es par les automorphismes polynomiaux de $\C^2$
et par l'application $(x,y) \mapsto (x,xy)$.
Il y a certes des endomorphismes birationnels polynomiaux qui ne sont pas de ce type
(cf [Da]), mais tous les exemples que nous connaissons (voir par exemple ceux 
explicit\'es dans [SY]) v\'erifient \'egalement $(\dag)$ dans une compactification 
ad\'equate. Voici une observation qui motive \'egalement
la question 4.33.

\begin{pro}
Soit $f:\C^2 \rightarrow \C^2$ un endomorphisme polynomial tel que $\l_1(f)>\l_2(f)$.
Supposons que $f$ admette une extension m\'eromorphe 1-stable \`a 
$\P^2=\C^2 \cup L_{\infty}$. Alors la droite \`a l'infini $L_{\infty}$ est contract\'ee
par $f$ sur un point fixe superattractif $q_{\infty}$.
\end{pro}

\begin{preuve}
Notons $f=(P,Q)$ les composantes polynomiales de $f$ dans $\C^2$, avec 
$\l=\max (\deg P,\deg Q)$: c'est le rayon spectral $r_1(f)$ de
l'extension m\'eromorphe de $f$ \`a $\P^2$.
L'ensemble d'ind\'etermination $I_f$ de $f$ est localis\'e  
\`a l'infini. La seule droite qui puisse \^etre
contract\'ee par $f$ sur un point d'ind\'etermination est la droite $L_{\infty}$.
Puisque $f$ est 1-stable dans $\P^2$, soit celle-ci est contract\'ee sur
un point $q_{\infty}$ qui n'est pas d'ind\'etermination (dans ce cas c'est un
point fixe superattractif), soit $L_{\infty}$ n'est pas contract\'ee: nous allons
voir que cela entra\^ine $\l_2(f) \geq \l=\l_1(f)$.

Supposons donc que $L_{\infty}$ n'est pas contract\'ee. Alors
les parties homog\`enes de plus haut degr\'e $P_{\l},Q_{\l}$
de $P,Q$ se d\'ecomposent en $P_{\l}=R \cdot A$ et $Q_{\l}=R \cdot B$,
o\`u $A,B$ sont des polyn\^omes homog\`enes premiers entre eux de
degr\'e $d \geq 1$, et $(R=0) \cap L_{\infty}=I_f$. Nous allons
montrer que $\l_2(f) \geq \l d$.

Soit $L,L'$ deux droites projectives g\'en\'eriques. Alors $L \cdot L'$ est un
point g\'en\'erique de $\C^2$ et $\l_2(f)$ est le nombre de pr\'eimages de ce point,
c'est \`a dire le cardinal de $f^{-1}(L) \cap f^{-1}(L')$ dans $\C^2$.
On peut estimer ce cardinal \`a l'aide du produit d'intersection
$f^*[L] \wedge f^* [L']$: celui-ci est de masse totale $\l^2$ dans $\P^2$,
qui se r\'epartit en $\l_2(f)$ masses de Dirac en des points de $\C^2$, et
$\l^2-\l_2(f)$ masses de Dirac aux points d'ind\'etermination de $f$ (localis\'es
\`a l'infini). Comme l'application
$f_{\infty}:=f_{|L_{\infty}}$ est de degr\'e $d \geq 1$, la courbe
$f^{-1}(L)$ intersecte $L_{\infty}$ en $d$ points qui ne sont pas
des points d'ind\'etermination. Si on choisit pour $L'$ une petite perturbation de
$L_{\infty}$, on d\'eduit de la relation d'invariance 
$f^*[L_{\infty}]=\l [L_{\infty}]$ que les courbes
$f^{-1}(L)$ et $f^{-1}(L')$ s'intersectent en au moins $\l d$ points dans $\C^2$.
\end{preuve}

\vskip.2cm
La proposition ci-dessus indique \'egalement une marche \`a suivre
pour trouver une compactification 1-stable d'un endomorphisme polynomial
$f$ de $\C^2$: on commence par consid\'erer son extension m\'eromorphe
\`a $\P^2=\C^2 \cup L_{\infty}$. Si $f$ ne contracte pas la droite
$L_{\infty}$, alors $f$ est 1-stable dans $\P^2$ et $\l_2(f) \geq \l_1(f)$.
Si $f$ contracte $L_{\infty}$ sur un point $q_{\infty}$ qui n'est pas un point
d'ind\'etermination, alors $f$ est 1-stable dans $\P^2$ et v\'erifie la condition
(BD). Sinon $f(L_{\infty} \setminus I_f)=m \in I_f$: on \'eclate
alors au point $m$ et on analyse le comportement de l'endomorphisme
induit dans l'\'eclat\'e...

C.Favre et M.Jonsson ont appliqu\'e leur analyse de l'arbre des valuations [FaJ 2]
dans le cadre d\'ecrit ici et ont r\'ealis\'e des progr\`es substantiels
en direction d'une r\'eponse positive \`a la question 4.33
(voir [FaJ 3]).

\subsection{Endomorphismes rationnels pathologiques}

Nous pr\'esentons \`a pr\'esent quelques exemples d'endomorphismes rationnels 
qui illustrent la difficult\'e qu'il peut y avoir \`a contr\^oler la dynamique
pr\`es des points d'ind\'etermination.

\begin{exa}
Consid\'erons, comme dans l'Exemple 2.17, le tore $Y=E \times E$, o\`u
$E=\C/\Z[\zeta]$ est la courbe elliptique associ\'ee \`a une racine primitive de l'unit\'e
$\zeta$ d'ordre $3,4$ ou $6$. La matrice
$$
A=\left[
\begin{array}{cc} 1 & d+2 \\ 1 & d \end{array} \right],
\; d \geq 1,
$$
pr\'eserve le r\'eseau $\Lambda=\Z[\zeta] \times \Z[\zeta]$. Elle induit donc un
endomorphisme holomorphe $g:Y \rightarrow Y$ tel que
$$
\l_2(g)=4 \;  \text{ et } \;
\l_1(g)=\left(\frac{d+1+\sqrt{(d+1)^2+8}}{2} \right)^2 > \l_2(g).
$$
Soit $\sigma:Y \rightarrow Y$ l'homoth\'etie de rapport $\zeta$, et
soit $f:X \rightarrow X$ l'endomorphisme induit par $g$ sur la surface rationnelle $X$ obtenue
en d\'esingularisant le quotient $Y/\langle \sigma\rangle $, i.e. en \'eclatant aux points fixes de $\sigma$.
Soit $a$ un tel point fixe. Comme $g$ est de degr\'e topologique $\l_2(g) \geq 2$,
$g^{-1}(a)$ contient d'autres pr\'eimages que les points fixes de $\sigma$. Chaque
point de $g^{-1}(a) \setminus Fix(\sigma)$ correspond, dans $X$, \`a un point d'ind\'etermination
de $f$. Comme les $g$-pr\'eimages de tout point s'\'equidistribuent selon la mesure de Lebesgue 
$\nu_Y$ du tore $Y$,
l'ensemble $(g^{-n}(a))_{n \in \N}$ est dense dans $Y$.
On en d\'eduit que l'ensemble
$$
I_f^{\infty}:=\bigcup_{n \in \N} f^{-n} I_{f^n}
\text{ est dense dans } X.
$$

Observons que $f$ est 1-stable: comme $g$ ne contracte aucune courbe,
il en est de m\^eme de $f$, qui v\'erifie donc le crit\`ere 
de la Proposition 4.4.
\end{exa}

C.Favre a donn\'e dans [Fa 1]  des exemples d'endomorphismes birationnels 1-stables
de $\P^2$ qui ne v\'erifient pas la condition (BD)
(voir \'egalement Exemple 5 dans [B]): 
ce sont des exemples qui admettent une droite
invariante $L$ sur laquelle $f$ est conjugu\'ee \`a une rotation
tr\`es irrationnelle. Nous g\'en\'eralisons ici
cette construction, en suivant [DG], pour obtenir des exemples similaires
(trois droites invariantes) pour lesquels nous sommes tout de m\^eme
en mesure de v\'erifier la condition $g^+ \in L^1(T^- \wedge \om)$.

\begin{exa}
Soit $(a,b,c) \in \C^3$ et $f=f_{abc}:\P^2 \rightarrow \P^2$ d\'efini par
$$
f[x:y:z]=[bcx(-cx+acy+z):
acy(x-ay+abz):abz(bcx+y-bz)].
$$
On v\'erifie que $f$ est un endomorphisme birationnel de $\P^2$
tel que $f^{-1}=f_{a^{-1}b^{-1}c^{-1}}$. Observons que les droites
$(x=0),(y=0),(z=0)$ sont invariantes par $f$,$f^{-1}$
et contiennent tous les points d'ind\'etermination des $f^n$, $n \in \Z$.
On calcule par exemple
$$
I_f=\left\{ [a:1:0];[0:b:1];[1:0:c] \right\}.
$$

Fixons $b=-2e^{2i\pi\theta}$, $c=1/2$ et $a=i$. On v\'erifie dans
ce cas que $f$ est 1-stable sur $\P^2$ et 
$g^-$ est finie aux deux points d'ind\'etermination 
$[0:b:1],[1:0:c]$. L'action de $f$ sur la droite $(z=0)$ est celle
d'une rotation d'angle $\theta$. Or les points $[i:1:0] \in I_f$ et
$[-i:1:0] \in I_{f^{-1}}$ sont sur un m\^eme cercle invariant.
Lorsque $\theta$ est choisi judicieusement, l'orbite n\'egative
$\{ f^{-n}[i:1:0], \, n \in \N \}$ du point $[i:1:0]$ passe
parfois tr\`es pr\`es de $[-i:1:0]$, de sorte que
$g^-([i:1:0])=-\infty$. La condition (BD) n'est donc pas satisfaite.
Nous montrons cependant dans [DG] que $g^- \in L^1(T^+)$.
\end{exa}

\subsection{Exemples provenant de la Physique}

De nombreuses applications birationnelles interviennent en
Physique en relation avec la th\'eorie des syst\`emes int\'egrables
(voir [BTR], [GNR] pour un survol r\'ecent).
La famille 
$$
f_a(x,y)=\left( y \frac{x+a}{x-1},x+a-1 \right), \; a \in \R,
$$
est constitu\'ee d'applications birationnelles du plan $\R^2$
qui pr\'eservent l'aire (elle pr\'eservent une 2-forme m\'eromorphe)
et sont r\'eversibles ($f_a$ est conjug\'ee \`a $f_a^{-1}$
par une involution). Elles ont \'et\'e \'etudi\'ees num\'eriquement
par de nombreux auteurs (voir [Ab], [BoM]).
C.Favre et J.Diller [DF], puis E.Bedford et J.Diller [BDi 1] ont
montr\'e comment les m\'ethodes d'analyse et g\'eom\'etrie complexe
permettent d'\'etudier la dynamique de telles applications:
elles v\'erifient notamment la condition (BD),
lorsqu'on consid\`ere leur complexification et compactification
\`a $\P^1 \times \P^1$. On dispose donc
d'une mesure d'entropie maximale $\mu_f$.
E.Bedford et J.Diller montrent [BDi 1], lorsque $a<0$, $a \neq -1$,
que la mesure $\mu_f$ est \`a support dans $\R^2$,
et que $(f,\mu_f)$ est conjugu\'e \`a $(\sigma,\nu)$,
o\`u $\sigma$ est le sous-d\'ecalage du nombre d'or
et $\nu$ d\'esigne son unique mesure d'entropie maximale.

Le lecteur trouvera dans [BDi 3] une \'etude analogue pour une famille
\`a deux param\`etres d'applications birationnelles.

Notons \'egalement que de nombreux endomorphismes rationnels interviennent
(comme op\'erateurs de renormalisation) dans l'analyse spectrale d'op\'erateurs
diff\'erentiels (Laplace, Schr\"odinger,...) sur des structures mod\`eles
self-similaires. On pourra consulter [Sa] \`a ce sujet.

\subsection{Applications monomiales}

Soit $A=\left[ \begin{array}{cc} a & b \\ c & d \end{array} \right]
\in {\mathcal M}(2,\Z)$ une matrice {\it inversible} \`a coefficients
entiers relatifs. On lui associe l'endomorphisme rationnel de $\C^2$
$$
f_A:(z,w) \in \C^2 \mapsto (z^aw^b,z^cw^b) \in \C^2.
$$
C'est une application {\it monomiale} dont les it\'er\'es s'expriment simplement
en fonction des puissances de la matrice $A$,
$(f_A)^n=f_{A^n}$. On en d\'eduit que
$$
\l_2(f_A)=|\det A| \; \; \text{ et } \;
\l_1(f_A)=\text{rayon spectral de } A.
$$
Ces applications ont \'et\'e consid\'er\'ees dans [FaG] lorsque 
$A \in {\mathcal M}(2,\N)$: ce sont dans ce cas des endomorphismes
polynomiaux de $\C^2$ qui induisent un endomorphisme 1-stable 
de $\P^1 \times \P^1$.

Le cas g\'en\'eral a \'et\'e syst\'ematiquement \'etudi\'e 
par C.Favre dans [Fa 3].
On peut trouver une surface $X$ rationnelle sur laquelle $f_A$ induit
un endomorphisme 1-stable lorsque le spectre de $A$ est r\'eel
(en particulier lorsque $\l_2(f)<\l_1(f)^2$).
On ne peut pas rendre $f_A$ 1-stable lorsque les valeurs
propres de $A$ sont de la forme $\rho e^{2i\pi\theta}$
avec $\theta \notin \Q$ (dans ce cas $\l_2(f)=\l_1(f)^2>1$).

La mesure de Lebesgue $\mu_f$ sur le tore r\'eel
$S^1 \times S^1=(|z|=|w|=1)$ est l'unique mesure d'entropie maximale lorsque 
$\l_1(f) \neq \l_2(f)$. Ses exposants de Lyapunov sont
$$
\chi_1(\mu_f)=\log \l_1(f) \; \; 
\text{ et } \chi_2(\mu_f)=\log \l_2(f)/\l_1(f).
$$
On v\'erifie ais\'ement les propri\'et\'es d'\'equidistribution des 
cycles p\'eriodiques qui sont majoritairement situ\'es 
sur le tore invariant $S^1 \times S^1$.
\vskip.2cm

On peut s'int\'eresser \`a ce m\^eme type d'applications,
$f=f_A$, $A \in {\mathcal M}(k,\Z)$, en dimension sup\'erieure ($k \geq 3$).
On v\'erifie qu'on a encore $\l_k(f_A)=|\det A|$ et
$\l_1(f_A)=$rayon spectral de $A$.
Il est probable que les autres degr\'es dynamiques s'expriment
\'egalement en fonction des valeurs propres $|a_1| \geq \ldots \geq |a_k|$
de la matrice $A$, mais le calcul s'av\`ere un peu plus d\'elicat:

\vskip.2cm
\noindent {\bf Probl\`eme.}
{\it Est-ce que les degr\'es dynamiques $\l_j(f)$ coincident avec les produits $|a_1| \cdots |a_j|$ 
pour tout $1 \leq j \leq k$ ?}
\vskip.2cm

\noindent On pourra consulter [HaPr] pour quelques informations sur cette question.

\chapter{Dimension sup\'erieure}

Il est peut-\^etre utile de pr\'eciser la strat\'egie que nous proposons 
pour d\'emontrer la conjecture \'enonc\'ee dans l'introduction.

\begin{defi} Soit $l \in [1,k]$.
On dit que $f:X \rightarrow X$ est l-stable si l'action lin\'eaire induite
par $f^*$ sur $H^{l,l}(X,\R)$ est compatible avec la dynamique, i.e. si
$(f^n)^*=(f^*)^n$ pour tout $n \in \N$.
\end{defi}

Supposons que le degr\'e dynamique $\l:=\l_l(f)$ domine 
strictement tous les autres degr\'es
dynamiques. On peut alors proc\'eder de la fa\c{c}on suivante:
\begin{itemize}
\item on cherche un mod\`ele $\tilde{X}$ birationnellement \'equivalent \`a $X$,
sur lequel $f$ induit un endomorphisme l-stable. Sur ce mod\`ele 
$\l=r_l(f)=r_{k-l}(f_*)$;
\item on essaie de montrer que $\l$ est une valeur propre ``distingu\'ee'' pour les
actions lin\'eaires induites $f^*:H^{l,l}(X,\R) \rightarrow H^{l,l}(X,\R)$
et $f_*:H^{k-l,k-l}(X,\R) \rightarrow H^{k-l,k-l}(X,\R)$
(voir d\'efinition 5.3). On note 
$\a^+ \in H^{l,l}_{psef}(X,\R)$ et $\a^- \in H^{k-l,k-l}_{psef}(X,\R)$
des vecteurs propres associ\'es \`a $\l$, lorsque ces c\^ones sont invariants;
\item on cherche \`a construire des courants invariants 
canoniques $T_l^+,T_{k-l}^-$ avec
$\{T_l^+\}=\a^+$, $\{T_{k-l}^-\}=\a^-$, puis \`a \'etablir des propri\'et\'es
d'extr\'emalit\'e et de laminarit\'e de ces courants;
\item on cherche enfin \`a d\'efinir $\mu_f=T_l^+ \wedge T_{k-l}^-$
et \`a \'etablir ses principales propri\'et\'es ergodiques.
\end{itemize}
\vskip.2cm

Lorsque $l=k$, toute application est k-stable, sur n'importe quel mod\`ele. 
L'analyse spectrale est tr\`es simple puisque
$H^{k,k}(X,\R) \simeq \R$. 
De plus il n'est pas n\'ecessaire d'intersecter de courants,
car $T_k^+$ est d\'ej\`a de bidegr\'e maximal.
Ceci explique en partie pourquoi le cas de grand degr\'e topologique,
trait\'e au chapitre 3,
est un peu plus simple.

Lorsque $l \leq k-1$, il y a plusieurs difficult\'es \`a surmonter,
comme nous l'avons vu au chapitre 4, lorsque la vari\'et\'e $X$
est de dimension deux. Nous nous int\'eressons \`a pr\'esent
aux endomorphismes de petit degr\'e topologique en dimension
sup\'erieure. C'est une source de difficult\'es nouvelles et 
il semble raisonnable, dans un premier temps, de limiter le champ d'investigations
\`a des classes  significatives d'exemples.
Nous en consid\'erons dans ce chapitre deux grandes familles:
\begin{itemize}
\item les {\it automorphismes} cohomologiquement hyperboliques: ils sont
automatiquement l-stables et la construction de courants invariants est simplifi\'ee
par l'absence de points d'ind\'etermination. C'est une direction d\'efrich\'ee
r\'ecemment par T.-C.Dinh et N.Sibony dans [DS 5,6]; 
\item les endomorphismes polynomiaux de $\C^k$ qui induisent un endomorphisme l-stables sur
$X=\P^k$: la g\'eom\'etrie est plus simple et les points d'ind\'etermination
sont confin\'es dans l'hyperplan \`a l'infini. C'est la direction de 
recherche adopt\'ee
dans [BP],[S],[GS],[CoG], [G 6], [DS 2,7].
\end{itemize}

\section{Pourquoi la dimension sup\'erieure ?}

Vues les difficult\'es que nous avons rencontr\'ees au chapitre 4
dans l'\'etude dynamique des endomorphismes des surfaces, on peut
\^etre sceptique \`a l'id\'ee d'attaquer cette question en dimension
sup\'erieure. Il y a pourtant de nombreuses raisons qui justifient
cette \'etude, comme nous l'expliquons \`a pr\'esent.

\subsection{Endomorphismes produits}

C'est une astuce classique dans le monde des syst\`emes dynamiques
d'\'etablir des propri\'et\'es ergodiques fines d'un endomorphisme
$f:X \rightarrow X$ en \'etudiant des propri\'et\'es -- a priori --
plus grossi\`eres de l'endomorphisme produit
$$
f \otimes f:(x,y) \in X^2 \mapsto (f(x),f(y)) \in X^2.
$$
Si $\mu_f$ est une mesure de probabilit\'e invariante pour $f$,
alors $\nu_f(x,y):=\mu_f(x) \otimes \mu_f(y)$ est une mesure
de probabilit\'e invariante pour $f \otimes f$, et on a par exemple
le r\'esultat suivant:

\begin{pro}
la mesure $\nu_f$ est ergodique si et seulement si $\mu_f$ est faiblement m\'elangeante.
\end{pro}

Rappelons que $\mu_f$ est dite faiblement m\'elangeante si pour toute paire
de bor\'eliens $A,B \subset X$, il existe un ensemble $J=J(A,B) \subset \N$
de densit\'e nulle dans $\N$ tel que
$
\lim_{J(A,B) \ni \hspace{-.5em} / \, n \rightarrow +\infty} \mu( f^{-n} A \cap B)=\mu(A) \mu(B).
$
Lorsque $J(A,B)=\emptyset$, on retrouve la propri\'et\'e de m\'elange fort.
Nous renvoyons le lecteur \`a [W] (Th\'eor\`eme 1.24 p 46), 
pour une preuve de ce r\'esultat \'el\'ementaire.

Observons \`a pr\'esent que la conjecture \'enonc\'ee dans l'introduction est
stable par passage au produit $f \otimes f$.
Cela r\'esulte notamment du Th\'eor\`eme 2.4.b: si
$\l_l(f)>\max_{j \neq l} \l_j(f)$, alors
$
\l_{2l}(f \otimes f)=[\l_l(f)]^2>\max_{j \neq 2l} \l_j(f \otimes f).
$
Il r\'esulte de plus des in\'egalit\'es de concavit\'e 2.4.a que
$$
\l_{2l-1}(f \otimes f)=\l_l(f) \l_{l-1}(f),
\; \text{ donc } \;
\frac{\l_{2l}(f \otimes f)}{\l_{2l-1}(f \otimes f)}=
\frac{\l_{l}(f)}{\l_{l-1}(f)},
$$
et de m\^eme $\l_{2l}(f \otimes f)/\l_{2l+1}(f \otimes f)=\l_{l}(f)/\l_{l+1}(f)$.
Nous laissons le soin au lecteur de v\'erifier comment les autres invariants
(exposants de Lyapunov, entropie, courants invariants, etc) se comportent par passage
au produit.

Lorsque l'on sait construire une mesure canonique invariante $\mu_f$ pour $f$, 
Il suffit ainsi, d'apr\`es la Proposition 5.2, de montrer l'ergodicit\'e de
$(f \otimes f,\mu_f \otimes \mu_f)$ pour obtenir le m\'elange faible.
\vskip.1cm

Dans un esprit similaire,  T.C.Dinh a r\'ecemment obtenu
[Di 2] une estimation de la d\'ecroissance des corr\'elations pour les applications de
H\'enon complexes $f:\C^2 \rightarrow \C^2$ en \'etudiant les propri\'et\'es
ergodiques de l'endomorphisme produit $(f,f^{-1}):\C^4 \rightarrow \C^4$.

\subsection{R\'eduction au cas polynomial}

Soit $f:\P^k \rightarrow \P^k$ un endomorphisme rationnel. On peut lui associer
canoniquement un endomorphisme $F: \P^{k+1} \rightarrow \P^{k+1}$
qui a les m\^emes caract\'eristiques dynamiques que $f$, et dont la
restriction \`a $\C^{k+1}$ est polynomiale. Plus pr\'ecis\'ement, 
rappelons que $f$ s'\'ecrit en coordonn\'ees homog\`enes
$
f[z_0:\cdots:z_k]=[P_0(z):\cdots:P_k(z)],
$
o\`u les $P_j$ sont des polyn\^omes homog\`enes premiers entre eux de degr\'e $r:=r_1(f)$.
Ils sont uniquement d\'etermin\'es \`a une constante multiplicative pr\`es.
L'endomorphisme polynomial $F:z \in \C^{k+1} \mapsto (P_0(z),\ldots,P_k(z)) \in \C^{k+1}$
rend le diagramme suivant commutatif,
$$
\begin{array}{ccc}
\P^k & \stackrel{f}{\longrightarrow} & \P^k \\
\pi \uparrow & \text{ } & \uparrow \pi \\
(\C^{k+1})^* & \stackrel{F}{\longrightarrow} & (\C^{k+1})^*
\end{array}
$$
Il s'\'etend en un endomorphisme rationnel de $\P^{k+1}$,
$
F[z:t]=[P_0(z):\cdots:P_k(z):t^{r}],
$
o\`u l'on a not\'e $(t=0)$ l'hyperplan \`a l'infini, $\P^{k+1}=\C^{k+1} \cup (t=0)$.
L'ensemble d'ind\'etermination $I_F$ est localis\'e dans l'hyperplan \`a l'infini,
et v\'erifie $I_F=I_f$ si l'on consid\`ere que $f$ agit sur
$(t=0) \simeq \P^k$: l'endomorphisme $f$ peut ainsi \^etre consid\'er\'e
comme la restriction de $F$ \`a $(t=0)$.
Observons que $F$ d\'efinit ainsi un endomorphisme 1-stable de $\P^{k+1}$
tel que $\l_1(F)=r_1(F)=r$. On suppose bien s\^ur $r \geq 2$. Plus g\'en\'eralement,
$$
\l_{j+1}(F)=r \l_j(f), \; \text{ pour }
0 \leq j \leq k. 
$$
Cela se v\'erifie ais\'ement
gr\^ ace au Th\'eor\`eme 2.4, si l'on consid\`ere l'endomorphisme produit
induit par $F$ sur $\P^k \times \P^1$,
$
([z],t) \in \P^k \times \P^1 \mapsto (f[z],t^r) \in \P^k \times \P^1.
$
Il r\'esulte en effet de 2.4.b que $\l_{j+1}(F)=\max (r \l_j(f); \l_{j+1}(f))$. Mais
les in\'egalit\'es de concavit\'e 2.4.a montrent par ailleurs
que $\l_{j+1}(f) \leq \l_1(f) \l_j(f)$ $\leq r \l_j(f)$, d'o\`u 
$\l_{j+1}(F)=r \l_j(f)$.
Il s'ensuit que $\l_l(f)$ domine strictement tous les autres degr\'es dynamiques
de $f$ si et seulement si $\l_{l+1}(F)$ domine strictement tous les degr\'es
$\l_j(F)$, $j \neq l+1$. 

Une mesure F-invariante
$\mu_F$ induit une mesure f-invariante, $\mu_f:=\pi_* \mu_F$.
Tout cycle $f$-p\'eriodique d'ordre $n$ ayant $l$ valeurs propres de module plus grand que $1$
et $k-l$ valeurs propres de module plus petit que $1$ [un cycle selle de type
$(l,k-l)$] induit $r$ cycles selles d'ordre $n$ pour $F$ et de type $(l+1,k-l)$:
il suffit de relever un tel cycle dans $\C^{k+1}$ et de faire agir
dessus la multiplication par une racine primitive de l'\'equation $\zeta^r=1$.

Les exposants de Lyapunov de $(f,\mu_f)$ et $(F,\mu_F)$
sont \'egalement reli\'es les uns aux autres de fa\c{c}on tr\`es simple: l'un des exposants
$\chi_j(\mu_F)$ est \'egal \`a $\log r$, les autres prennent exactement les m\^ emes
valeurs que les $\chi_i(\mu_f)$, en particulier
$
\sum_{j=0}^{k+1} \chi_j(\mu_F)=\log r+\sum_{j=1}^k \chi_j(\mu_f).
$
\vskip.1cm

\noindent {\bf Conclusion.} Pour comprendre la dynamique des applications rationnelles 
$f:\P^k \rightarrow \P^k$, notamment pour \'etablir la conjecture, il suffit
-- au prix d'une augmentation artificielle de la dimension --
de consid\'erer ceux qui sont induits par un endomorphisme polynomial de $\C^k$.

\subsection{Dynamique \`a param\`etre}

Nous avons mentionn\'e (voir paragraphe 3.4.3) qu'il est int\'eressant d'\'etudier
la dynamique de famille d'endomorphismes $(f_t)_{t \in M}$ qui d\'ependent
holomorphiquement d'un param\`etre $t \in M$.
Lorsque la vari\'et\'e $M$ est un ouvert de Zariski d'une vari\'et\'e projective,
on peut \'etudier en famille la dynamique des applications $f_t$ en consid\'erant
l'endomorphisme $F(x,t)=(f_t(x),t)$.

Consid\'erons par exemple le cas d'endomorphismes polynomiaux $f_t:\C^k \rightarrow \C^k$
qui d\'ependent polynomialement d'un param\`etre $t \in \C^m$.
Alors
$$
F:(z,t) \in \C^{k+m} \mapsto (f_t(z),t) \in \C^{k+m}
$$
est un endomorphisme polynomial de $\C^{k+m}$ qui pr\'eserve les feuilletages
$(t_i=constante)$: c'est un produit crois\'e.
On peut consid\'erer son extension m\'eromorphe -- not\'ee encore $F$ -- \`a $\P^k$
ou \`a toute compactification lisse de $\C^k$.
L'endomorphisme $F$ n'est pas cohomologiquement hyperbolique
mais on peut cependant \'etudier les courants dynamiques qui interviennent en bidegr\'e
$(p,p)$, $p \leq k$, et en d\'eduire des informations sur la dynamique de $f_t$.
Par exemple la fonction de Green dynamique de $F$,
$$
G_F(z,t):=\lim_{n \rightarrow +\infty} \frac{1}{\l_1(F)^n} \log^+ ||F^n(z,t)||
$$
co\"{\i}ncide avec la fonction de Green dynamique de $f_t$,
$$
G_{f_t}(z):=\lim_{n \rightarrow +\infty} \frac{1}{\l_1(f_t)^n} \log^+ ||f_t^n(z)||.
$$
Des r\'esultats de r\'egularit\'e (resp. plurisousharmonicit\'e)
de $(z,t) \mapsto G_F(z,t)$ (dans l'esprit de la 
Proposition 1.2, voir [DG]) donnent donc des informations sur la continuit\'e
(resp. plurisousharmonicit\'e)
de $t \mapsto G_{f_t}$.

\section{Automorphismes}

Nous testons ici la strat\'egie propos\'ee ci-dessus dans le cas
o\`u $f:X \rightarrow X$ est un automorphisme
cohomologiquement hyperbolique,
i.e. $I_f= \emptyset$, $\l_k(f)=1$ et
on note $\l:=\l_l(f) > \max_{j \neq l} \l_j(f)$.
Comme $f$ est holomorphe, $f$ est $l$-stable sur $X$.
Nous montrons que $\l$ est une valeur propre distingu\'ee
de l'op\'erateur $f^*:H^{l,l}(X,\C) \rightarrow H^{l,l}(X,\C)$
lorsque $\dim_{\C} X \leq 3$.

Nous en d\'eduisons la construction de courants invariants canoniques
$T_l^+$ et $T_{k-l}^-$, et d'une mesure de probabilit\'e invariante
$\mu_f=T_l^+ \wedge T_{k-l}^-$ dont nous montrons qu'elle est
m\'elangeante et d'entropie maximale lorsque $X$ est projective.

Ce probl\`eme a \'et\'e \'etudi\'e par  T.C.Dinh et
N.Sibony dans [DS 6], en suivant une m\'ethode originale,
qui fournit des informations sur la dimension de Hausdorff
du support de $\mu_f$. 
Nous esquissons leur approche 
dans le paragraphe 5.2.3

\subsection{Analyse spectrale}

Soit $f:X \rightarrow X$ un endomorphisme holomorphe.
Nous allons tester la condition suivante:

\begin{defi}
L'endomorphisme $f$ satisfait la condition {\bf Spec(f,l)} si 
$\l_l(f)$ est une valeur propre simple de l'op\'erateur
$f^*:H^{l,l}(X,\C) \rightarrow H^{l,l}(X,\C)$ qui domine strictement
toutes les autres valeurs propres de cet op\'erateur.
\end{defi}

Notons que la condition {\bf Spec(f,k)} est trivialement v\'erifi\'ee.
Nous avons vu que la condition {\bf Spec(f,1)} est v\'erifi\'ee
en dimension deux lorsque $\l_2(f)<\l_1(f)^2$ (Th\'eor\`eme 4.7).
C'est encore vrai en dimension sup\'erieure:

\begin{pro}
Soit $f:X \rightarrow X$ un endomorphisme holomorphe tel que
$\l_2(f)<\l_1(f)^2$. Alors $f$ v\'erifie la condition 
{\bf Spec(f,1)}
et toute valeur propre $\zeta \neq \l_1(f)$ de l'op\'erateur 
$f^*:H^{1,1}(X,\C) \rightarrow H^{1,1}(X,\C)$ est telle que
$|\zeta| \leq \sqrt{\l_2(f)}<\l_1(f)$.
\end{pro}

\begin{preuve}
Comme $f$ est holomorphe, $\l_1(f)=r_1(f)$ est le rayon spectral de l'action
lin\'eaire $f^*:H^{1,1}(X,\C)\rightarrow H^{1,1}(X,\C)$. Le c\^one $H_{nef}^{1,1}(X,\R)$ \'etant
pr\'eserv\'e, il existe une classe nef $\a^+ \neq 0$ telle que
$f^* \a^+=\l_1(f) \a^+$.

Observons que $\a^+ \wedge \a^+=0$, sinon $f^* (\a^+ \wedge \a^+)=\l_1(f)^2 \a^+ \wedge \a^+$
et donc $\l_1(f)^2 \leq \l_2(f)$, contredisant notre hypoth\`ese.
Supposons que $f^* \beta=\l_1(f)\b+\e \a^+$ pour une
classe $\b \in H^{1,1}(X,\C)$ et $\e \in \C$.
Alors $f^*(\b \wedge \a^+)=\l_1(f)^2 \b \wedge \a^+$, donc $\b \wedge \a^+=0$.
De m\^eme $\b \wedge \b =0$.
Il r\'esulte alors des relations bilin\'eaires de Hodge-Riemann (voir [GH]) que
$\b$ est proportionnelle \`a $\a^+$.
Cela montre que $\l_1(f)$ est une valeur propre simple.

Soit $\zeta \in \C$ une valeur propre de $f^*:H^{1,1}(X,\C) \rightarrow H^{1,1}(X,\C)$
et $\b \neq 0$ un vecteur propre associ\'e.
Si $\b \wedge \a^+ \neq 0$ alors $f^*(\beta \wedge \a^+)=\zeta \l_1(f) \beta \wedge \a^+$
implique $|\zeta| \leq \l_2(f)/\l_1(f) \leq \sqrt{\l_2(f)}$.
Supposons \`a pr\'esent $\b \wedge \a^+=0$.
Si $\b \wedge \b=0$, il r\'esulte de l'analyse pr\'ec\'edente que
$\b$ est proportionnel \`a $\a^+$ et donc $\zeta=\l_1(f)$.
Lorsque $\zeta \neq \l_1(f)$ on a donc
$\b \wedge \b \neq 0$; dans ce cas $f^* \beta \wedge  \beta = \zeta^2 \beta \wedge \beta$
implique $|\zeta| \leq \sqrt{\l_2(f)}$.
\end{preuve}
\vskip.2cm

En appliquant l'analyse pr\'ec\'edente \`a $f^{-1}$, lorsque $f$ est inversible,
on obtient que  $f$ v\'erifie la condition
{\bf Spec(f,k-1)} si  $\l_{k-2}(f)< \l_{k-1}(f)^2$.
Cela r\`egle le cas des automorphismes
des vari\'et\'es de dimension trois:

\begin{pro}
Tout automorphisme $f$ cohomologiquement hyperbolique d'une vari\'et\'e
de dimension 3 v\'erifie la condition {\bf  Spec(f,l)}, \,
$l$ tel que $\l_l(f) >\max_{j \neq l} \l_j(f)$.
\end{pro}

Nous pensons qu'un r\'esultat similaire a lieu lorsque $\dim_{\C} X \geq 4$.
Notons qu'il se peut que $f$ ne v\'erifie pas la condition {\bf Spec(f,j)}, 
$j \leq l-1$, mais c'est bien la condition {\bf Spec(f,l)} qui importe.

\begin{exa}
Soit $X=\C^k/\Lambda$ un tore complexe compact de dimension $k \geq 1$.
Soit $f=f_A:X \rightarrow X$ un endomorphisme induit
par $z \in \C^k \mapsto A \cdot z \in \C^k$, o\`u $A \in GL(k,\C)$
pr\'eserve le r\'eseau $\Lambda$. 

On note $\{a_1,\ldots,a_k\}=Spec(A)$ les valeurs propres de la matrice
$A$ rang\'ees par ordre d\'ecroissant, $|a_1| \geq \ldots \geq |a_k|$.
L'endomorphisme $f$ est cohomologiquement hyperbolique si et seulement si
aucune valeur propre $a_i$ n'est de module 1 (voir Exemple 2.18).
Soit $l \in [1,k]$ tel que $|a_l|>1>|a_{l+1}|$ (avec la convention
$a_{k+1}=0$). Alors
$$
\l_j(f)=\prod_{i=1}^j |a_i|^2 ,
\text{ donc } \l_l(f)>\max_{j \neq l} \l_j(f).
$$
L'endomorphisme $f$ v\'erifie la condition {\bf Spec(f,l)}: soit
$(\eta_i)$ une base de $1$-formes holomorphes sur $X$
trigonalisant l'action de $f^*$ sur $H^{1,0}(X,\C)=\C^k$,
$$
f^* \eta_i=a_i \eta_i+\e_i \eta_{i-1} , \; \e_i \in \{ 0,1\}, \; 1 \leq i \leq k.
$$
Alors $\l_l(f)$ est une valeur propre simple de $f^*:H^{l,l}(X,\C) \rightarrow H^{l,l}(X,\C)$
associ\'ee au vecteur propre 
$\eta_1 \wedge \overline{\eta_1} \wedge \ldots \wedge \eta_l \wedge \overline{\eta_l}$,
et toute autre valeur propre 
est de module strictement plus petit.

Notons que l'endomorphisme $f$ ne v\'erifie pas n\'ecessairement la condition
Spec(f,j), $j \leq l-1$: on peut avoir par exemple $|a_1|=|a_2|>1$.
\end{exa}

Supposons \`a pr\'esent que $f$ v\'erifie la condition {\bf Spec(f,l)},
$l$ tel que $\l_l(f)>\max_{j \neq l} \l_j(f)$.
Comme $f^*$ pr\'eserve le c\^one $H^{l,l}_{psef}(X,\C)$,
on en d\'eduit l'existence d'une classe $\a^+ \in H^{l,l}_{psef}(X,\C)$
-- unique \`a constante multiplicative pr\`es -- telle que $f^* \a^+=\l_l(f) \a^+$.
Par dualit\'e de Serre, on a un r\'esultat analogue pour l'action duale
$f_*:H^{k-l,k-l}(X,\C) \rightarrow H^{k-l,k-l}(X,\C)$. On normalise alors
le choix des classes invariantes $\a^+ \in H^{l,l}_{psef}(X,\C)$,
$\a^- \in H^{k-l,k-l}_{nef}(X,\C)$ et d'une classe de 
K\"ahler $\{\om\} \in H^{1,1}(X,\R)$ en imposant
$$
\a^+ \cdot \a^-=\a^+ \cdot \{\om^{k-l}\} = \a^- \cdot \{\om^l\}=1.
$$
Le point important ici est que la canonicit\'e de $\a^+,\a^-$
assure que le produit $\a^+ \cdot \a^-$ est strictement positif:
il suffit de v\'erifier que $\l^{-n}(f^n)^* \{\om^l\}$
converge vers un multiple strictement positif de $\a^+$.

Une autre cons\'equence int\'eressante de la condition {\bf Spec(f,l)}
est li\'ee \`a l'extr\'emalit\'e des courants invariants que
nous construisons un peu plus loin.

\begin{pro}
Si $f$ v\'erifie la condition {\bf Spec(f,l)}, alors
$\a^+$ est une classe extr\'emale dans le c\^one
$H^{l,l}_{psef}(X,\R)$.
\end{pro}

\begin{preuve}
Soit $\eta \in H^{l,l}_{psef}(X,\R)$ une classe pseudoeffective telle que
$0 \leq \eta \leq \a^+$. Nous devons montrer que $\eta$ est proportionnelle
\`a $\a^+$. Une telle classe se d\'ecompose en
$\eta=c \a^+ +\sum_{i=1}^s \tau_i$, o\`u $0 \leq c \leq 1$ et
les classes $\tau_i$ appartiennent
aux sous-espaces caract\'eristiques de $f^*:H^{l,l}(X,\C) \rightarrow H^{l,l}(X,\C)$
associ\'es \`a la valeur propre $\zeta_i$, $|\zeta_i| <\l_l(f)$.
Comme $f_*f^* \tau=\l_k(f) \tau$, il vient
$$
||(f^n)^* \tau_i|| \simeq n^{m_i} \zeta_i^n ||\tau_i|| \text{ et }
||(f^n)_* \tau_i|| \simeq n^{m_i} \left( \frac{\l_k(f)}{\zeta_i} \right)^n ||\tau_i||.
$$
Or la suite de classes pseudoeffectives  $\l^n \l_k^{-n} (f^n)_* \eta$ est
domin\'ee par $\a^+$, c'est donc que $\tau_i=0$ pour tout $i$.
\end{preuve}

\subsection{La mesure canonique}

Nous construisons ici une mesure canonique invariante lorsque
$f:X \rightarrow X$ est un automorphisme cohomologiquement hyperbolique
qui v\'erifie la condition {\bf Spec(f,l)}, $l$ tel que
$\l_l(f)>\max_{j \neq l} \l_j(f)$ (par exemple lorsque
$\dim_{\C} X=3$).
On fixe dans la suite $\a^+ \in H^{l,l}_{nef}(X,\R)$, $\a^- \in H^{k-l,k-l}_{nef}(X,\R)$
et $\om$ une forme de K\"ahler telles que
$$
(f^{\pm})^* \a^{\pm}=\l \a^{\pm} \,
\text{ et  } \, 
\a^+ \cdot \{\om^{k-l} \}=\a^- \cdot \{ \om^{l} \}=\a^+ \cdot \a^-=1.
$$

\noindent {\bf Construction des courants invariants.}
Soit $\theta^+$ une forme lisse ferm\'ee de bidegr\'e $(l,l)$ qui repr\'esente
$\a^+$. Soit $R^+$ une forme lisse de bidegr\'e $(l-1,l-1)$ telle que
$\l^{-1}f^* \theta^+=\theta^+ +dd^c R^+$.
Quitte \`a changer $R^+$ en $R^+ +C \om^{l-1}$, on peut supposer
que $0 \leq R^+ \leq C \om^{l-1}$.
La positivit\'e des formes et des courants peut \^etre ici interpr\'et\'ee
au sens faible ou fort, cela n'a pas d'importance (voir [De] pour les d\'efinitions
et propri\'et\'es de base des courants positifs). On it\`ere alors cette
\'equation fonctionnelle en prenant son image inverse par $f^n$, ce qui donne
$$
\frac{1}{\l^n}(f^n)^* \theta^+=\theta^+ +dd^c R_n^+,
\; \; R_n^+:=\sum_{j=0}^{n-1} \frac{1}{\l^j} (f^j)^* R^+.
$$
La suite $(R_n)^+$ est une suite croissante de formes diff\'erentielles positives
qui converge au sens des courants car sa masse est uniform\'ement major\'ee:
$$
||R_n^+||:=\int_X R_n^+ \wedge \om^{k-l+1} \leq 
C \sum_{j \geq 0} \frac{1}{\l^j} \d_{l-1}(f^j,\om) <+\infty,
$$
car $\l=\l_l(f)>\l_{l-1}(f)$. 
On en d\'eduit que la suite $\l^{-n}(f^n)^* \theta^+$ converge au sens des courants
vers un courant ferm\'e de bidegr\'e $(l,l)$,
$$
T_l^+:=\theta^+ +dd^c R_{\infty}^+,
\text{ avec } R_{\infty}^+=\sum_{j \geq 0} \l^{-j} (f^j)^* R^+ \geq 0.
$$
Notons que le courant $T_l^+$ est {\it invariant}, $f^* T_l^+=\l T_l^+$, et {\it positif}
car limite de la suite de formes positives $\l^{-n}(f^n)^* \om^l$, comme on le v\'erifie
en d\'ecomposant $\{ \om^l\}$ dans  une base qui trigonalise l'action de 
$f^*$ sur $H^{l,l}(X,\C)$: on utilise ici de fa\c{c}on d\'ecisive la condition
{\bf Spec(f,l)}. Plus g\'en\'eralement, si $\Theta$ est une 
$(l,l)$-forme lisse ferm\'ee, alors
$$
\frac{1}{\l^n}(f^n)^* \Theta \longrightarrow c T_l^+
\; \text{ avec } \; c=\{ \Theta \} \cdot \a^-.
$$
On construit de m\^eme un courant positif ferm\'e $T_{k-l}^-$
de bidegr\'e $(k-l,k-l)$ tel que $f_* T^-_{k-l}=\l T^-_{k-l}$,
$$
T^-_{k-l}=\theta^-+dd^c R_{\infty}^-, \; \; 
R_{\infty}^-:= \sum_{j \geq 0} \frac{1}{\l^j} (f^{-j})^* R^- \geq 0,
$$
o\`u $\theta^-$ est une $(k-l,k-l)$-forme lisse ferm\'ee repr\'esentant
$\a^-$ et $R^- \geq 0$ est une $(k-l-1,k-l-1)$-forme lisse telle que
$\l^{-1}(f^{-1})^* \theta^-=\theta^-+dd^c R^-$:
on utilise ici l'hypoth\`ese $\l=\l_l(f)>\l_{l+1}(f)$.
Rappelons que $\l=\l_l(f)$ domine strictement tous les autres degr\'es dynamiques
si et seulement si $\l_l(f)>\max[ \l_{l-1}(f), \l_{l+1}(f)]$.
Cela r\'esulte des in\'egalit\'es de concavit\'e 2.4.a.
\vskip.2cm

L'existence de potentiels positifs $R_{\infty}^{\pm} \geq 0$ assure, 
comme dans le Th\'eor\`eme 3.1,
de bonnes propri\'et\'es d'int\'egrabilit\'e des courants $T_l^+,T_{k-l}^-$:

\begin{lem}
Soit $T$ un courant positif ferm\'e de bidegr\'e $(p,p)$ sur $X$ 
qui admet la d\'ecomposition $T=\theta+dd^c R$, o\`u $\theta$ est une 
$(p,p)$-forme lisse et $R$ est un $(p-1,p-1)$-courant {\bf positif}.
Alors toute fonction quasiplurisousharmonique est int\'egrable par rapport
\`a la mesure trace $T \wedge \om^{k-p}$.
\end{lem}

\begin{preuve}
C'est l'observation d\'eja faite au Th\'eor\`eme 3.1.1: soit $\f$ une fonction qpsh
sur $X$; quitte \`a translater et dilater $\f$, on peut supposer que $\f \leq 0$
et $dd^c \f \geq -\om$. Une int\'egration par parties donne alors
\begin{eqnarray*}
0 \leq \int_X (-\f) T \wedge \om^{k-p}
&= &\int_X (-\f) \theta \wedge \om^{k-p} +\int_X R \wedge (-dd^c \f) \wedge \om^{k-p} \\
&\leq & \int_X (-\f) \theta \wedge \om^{k-p} +\int_X R \wedge \om^{k-p+1}<+\infty,
\end{eqnarray*}
car les fonctions qpsh sont int\'egrables par rapport aux mesures lisses.
\end{preuve}
\vskip.2cm

Voici une cons\'equence importante de la condition {\bf Spec(f,l)}:

\begin{thm}
Le courant $T_l^+$ (resp. $T_{k-l}^-$) est un point extr\'emal du c\^one convexe
des courants positifs ferm\'es de bidegr\'e $(l,l)$ (resp. $(k-l,k-l)$).
\end{thm}

La preuve est dans le m\^eme esprit que celle du Th\'eor\`eme 4.18.
Il faut commencer par s'assurer que tout courant positif
ferm\'e domin\'e par $T_l^+$ est cohomologue \`a un multiple de $\a^+$
(c'est le contenu de la Proposition 5.7), puis mettre en place un
r\'esultat de convergence uniforme des images inverses (resp. directes)
normalis\'ees. Ce dernier point est plus d\'elicat lorsque $l \geq 2$,
car les potentiels des courants de bidegr\'e $(l,l)$ ne sont pas 
uniques. Nous renvoyons le lecteur au Th\'eor\`eme 4.1 de [DS 6]
pour une preuve (voir \'egalement [G 6], [DS 7] pour un contexte
proche).

\vskip.3cm
\noindent {\bf Construction de $\mu_f$.}
Nous d\'efinissons \`a pr\'esent la mesure invariante canonique
$\mu_f=T_l^+ \wedge T_{k-l}^-$: il est possible de donner un sens \`a ce
produit d'intersection car les courants $T_l^+,T_{k-l}^-$ sont tr\`es bien
approxim\'es par les formes lisses $\l^{-n}(f^{\pm n})^* \theta^{\pm}$.

\begin{pro}
Les suites de mesures 
$$
\l^{-n}(f^n)^* \theta^+ \wedge T_{k-l}^-, \; T_l^+ \wedge \l^{-n}(f^{-n})^* \theta^- \;
\text{ et  } \;  \l^{-2n}(f^n)^* \om^l \wedge (f^{-n})^* \om^{k-l}
$$
convergent toutes vers une m\^eme mesure de probabilit\'e $\mu_f$.
\end{pro}

\begin{preuve}
Posons $\theta_n^{\pm}=\l^{-n}(f^{\pm n})^* \theta^{\pm}$,
$\mu_n=\theta_n^+ \wedge T_{k-l}^-$ et
$\nu_n=\theta_n^+ \wedge \theta_n^-$.
Observons que
$\mu_n=\theta^+ \wedge T_{k-l}^-+dd^c (R_n^+ \wedge T_{k-l}^- )$. Or
$R_n \wedge T_{k-l}^-$ est une suite croissante de courants positifs
dont la masse est uniform\'ement born\'ee,
$$
||R_n^+ \wedge T_{k-l}^-|| \leq C \sum_{j=0}^{n-1} \frac{1}{\l^j}
\int_X (f^j)^* \om^{l-1} \wedge T_{k-l}^- \wedge \om 
\leq C' \sum_{j \geq 0} \frac{\d_{l-1}(f^j,\om)}{\l^j} <+\infty.
$$
On en d\'eduit que $\mu_n$ converge vers une mesure $\mu_f$ telle que
$$
\mu_f=\theta^+ \wedge T_{k-l}^-+dd^c (S),
\text{ avec }
S:=\sum_{j \geq 0} 
\frac{1}{\l^j} (f^j)^* R_+ \wedge T_{k-l}^- \geq 0.
$$
Observons que $\mu_n-\nu_n= dd^c (\theta_n^+ \wedge \sum_{j \geq n} \l^{-j} (f^j)_* R^-)$,
or
$$
|| \theta_n^+ \wedge \sum_{j \geq n} \l^{-j} (f^j)_* R^-)||
\leq C' \sum_{j \geq n} \frac{\d_{l+1}(f^j,\om)}{\l^j} \longrightarrow 0,
$$
donc les suites $\nu_n$ et $\mu_n$ ont la m\^eme limite.
Les autres suites se traitent de fa\c{c}on similaire.
En particulier la mesure limite $\mu_f$ est une mesure de probabilit\'e, 
car limite des mesures positives lisses 
$\l^{-2n}(f^n)^*\om^l \wedge (f^{-n})^*\om^{k-l}$ dont la masse tend vers 1.
\end{preuve}

\begin{thm}
La mesure $\mu_f$ est une mesure invariante 
m\'elangeante qui int\`egre
les fonctions qpsh. Elle est d'entropie maximale
si $X$ est projective,
$$
h_{top}(f)=h_{\mu_f}(f)=\log \l_l(f).
$$
\end{thm}

\begin{preuve}
Observons que $f_* \mu_n=\mu_{n+1}$, donc la mesure $\mu_f$ est invariante.
Elle se d\'ecompose en 
$$
\mu_f=\theta^+ \wedge \theta^-+ dd^c( [C \om^l+\theta^+] \wedge R_{\infty}^- +S)
-C \om^l \wedge T_{k-l}^-+C \om^l \wedge \theta^-,
$$ 
o\`u $S$ est un courant positif et $C \geq 0$ est choisie de sorte que 
la forme $C\om^l+\theta^+$ soit positive.
Il r\'esulte alors d'une double application du
Lemme 5.8 que les mesures $\theta^+ \wedge T_{k-l}^-$ 
et $[C \om^l+\theta^+] \wedge \theta^-+dd^c( [C \om^l\theta^+] \wedge R_{\infty}^- +S)$
int\`egrent les fonctions
qpsh. Il en est donc de m\^eme de $\mu_f$.

Le m\'elange, comme dans le Th\'eor\`eme 4.23, 
est une cons\'equence de l'extr\'emalit\'e des courants $T_l^+,T_{k-l}^-$
(Th\'eor\`eme 5.9).

L'entropie de $\mu_f$ est major\'ee par $h_{top}(f)=\log \l_l(f)$
d'apr\`es le principe variationnel et le Th\'eor\`eme 2.8.
La minoration s'appuie sur les travaux de Y.Yomdin (voir [Sm], [BS 2]
pour une utilisation de ces travaux dans le contexte des applications de 
H\'enon complexes, et le Th\'eor\`eme 3.2 de [G 6] pour la dimension
sup\'erieure).
\end{preuve}
\vskip.2cm

\noindent {\bf Bilan.} Soit $f:X \rightarrow X$ un automorphisme cohomologiquement
hyperbolique sur une vari\'et\'e projective de dimension trois.
Quitte \`a changer $f$ en $f^{-1}$, nous pouvons supposer que $\l:=\l_1(f)$ est
le degr\'e dynamique dominant.
Alors $f$ v\'erifie la condition {\bf Spec(f,1)} (Proposition 5.4), donc
la demi-droite canonique $\R^+ \a^+$ est extr\'emale (Proposition 5.7).
On sait donc construire des courants invariants canoniques $T_1^+$,
$T_2^-$ extr\'emaux (Th\'eor\`eme 5.9), et une mesure
de probabilit\'e invariante canonique qui est m\'elangeante et
d'entropie maximale (Th\'eor\`eme 5.11). Il reste \`a estimer
ses exposants de Lyapunov (indiquons le tout r\'ecent travail
[DeT 5] qui donne cette estimation lorsque l'on sait calculer
l'entropie m\'etrique)
et pr\'eciser la nature des points selles de $f$.
Cela passe vraisemblablement par une meilleure compr\'ehension de la nature
g\'eom\'etrique des courants $T_1^+,T_2^-$ dans l'esprit de ce qui a \'et\'e expos\'e
en dimension deux (voir paragraphe 4.3.3). Les travaux r\'ecents de T.C.Dinh
[Di 1] et H.deThelin [DeT 4] vont dans ce sens.
Notons enfin que le lemme 4.28 de R.Dujardin est valable en toute dimension.

\subsection{L'approche de Dinh-Sibony}

Soit $f:X \rightarrow X$ un automorphisme cohomologiquement hyperbolique
d'une vari\'et\'e k\"ahl\'erienne compacte. Nous supposons pour simplifier 
l'exposition que $X$ est de dimension trois.
Quitte \`a changer $f$ en $f^{-1}$, on peut donc supposer que
$$
\l:=\l_1(f^{-1})=\l_2(f)>\l_1(f)>1.
$$
Il r\'esulte de la Proposition 5.4 que $\l$ est une valeur propre
distingu\'ee de l'op\'erateur $f_*:H^{1,1}(X,\C) \rightarrow H^{1,1}(X,\C)$.
Nous avons construit pr\'ec\'edemment des courants invariants canoniques
$T_1^-$ et $T_2^+$ de bidegr\'es respectifs $(1,1),(2,2)$, et montr\'e
que la mesure $\mu_f:=T_1^- \wedge T_2^+$ est dynamiquement int\'eressante.

Nous indiquons \`a pr\'esent quelques \'el\'ements de l'approche de
T.C.Dinh et N.Sibony [DS 6] qui permet de montrer que $\mu_f$ ne charge
pas les ensembles de petite dimension de Hausdorff.
Pr\'ecisons que cette approche fonctionne en toute dimension.
\vskip.2cm

L'id\'ee de la m\'ethode est de construire par r\'ecurrence 
sur $p$, des courants positifs ferm\'es 
invariants $\hat{T}=T \wedge S_T$ de bidegr\'e $(p+1,p+1)$,
en partant d'un courant positif ferm\'e invariant $T$ de bidegr\'e $(p,p)$.
A chaque \'etape, le courant $S_T$ est 
invariant, de bidegr\'e $(1,1)$, ferm\'e, \`a
potentiels h\"old\'eriens, et tel que $T \wedge S_T$ est positif.
Le point de d\'epart ($p=1$) est le courant $T_1^-$, celui d'arriv\'ee
($p=k-1$) est la mesure $\mu_f$.

Il s'agit donc de faire de l'analyse sur un courant 
positif ferm\'e invariant $T$. Cela n\'ecessite de contr\^oler
l'action des op\'erateurs $f^*,f_*$ sur les espaces de cohomologie
de tels courants: on introduit
$$
N^{1,1}(T,\R):=\left\{ \a \in H^{1,1}(X,\R) \, / \, 
\a \wedge \{T\}=0 \text{ dans } H^{p+1,p+1}(X,\R) \right\}.
$$
On introduit \'egalement $N_{\nu}^{1,1}(T,\R)$ le sous-espace constitu\'e
des classes $\a=\{S\}$ qui peuvent \^etre repr\'esent\'ees par un courant
(non n\'ecessairement positif) ayant un potentiel $\nu$-h\"old\'erien.
On consid\`ere alors les espaces
$$
H^{1,1}(T,\R):=H^{1,1}(X,\R) / N^{1,1}(T,\R)
\text{ et }
H^{1,1}_{\nu}(T,\R):=H^{1,1}(X,\R) / N^{1,1}_{\nu}(T,\R).
$$

Lorsque $T$ est invariant, $f^*T=\l_T T$, $\l_T>0$, l'op\'erateur
$f^*$ pr\'eserve $N^{1,1}(T,\R)$ et $N^{1,1}_{\nu}(T,\R)$ et 
induit donc un op\'erateur sur les espaces
$H^{1,1}(T,\R)$, $H^{1,1}_{\nu}(T,\R)$.

\begin{defi}
Soit $T$ un courant positif ferm\'e de bidegr\'e $(p,p)$
qui est $f$-invariant, $f^*T=\l_T T$, $\l_T>0$.
Soit $l$ un entier tel que $0 \leq l+p \leq k$.

Le $l^{\text{i\`eme}}$ degr\'e dynamique de $T$ est
$$
\l_l(f,T):=\liminf_{n \rightarrow +\infty} \left[
\int_X T \wedge (f^n)^* \om^l \wedge \om^{k-p-l} \right]^{1/n}.
$$
\end{defi}

Observons que $\l_l(f,T)=\l_l(f)$ lorsque $T$
est de bidegr\'e $(0,0)$. En g\'en\'eral on a seulement une 
in\'egalit\'e: on montre comme au Th\'eor\`eme 2.4 que
\begin{enumerate}
\item $\l_l(f,T) \leq \l_1(f)$;
\item $\l_l(f,T) \leq \l_1(f,T)^l$;
\item $\l_{k-p}(f,T)=\l_T^{-1}$.
\end{enumerate}

Lorsque $T=T_1^-$, $p=1$ et $k=3$, il vient $\l_T=\l_1(f^{-1})^{-1}$
et on obtient
$$
\l_1(f,T)^2 \geq \l_2(f,T)=\l_T^{-1}=\l_1(f^{-1})>1,
$$
donc $\l_1(f,T_1^-)>1$: c'est la condition qui va permettre 
(cf Th\'eor\`eme 5.13) de construire un courant 
$S_{T_1^-}$ courant de bidegr\'e $(1,1)$ 
\`a potentiels h\"old\'eriens, tel que 
$$
T':=T_1^- \wedge S_{T_1^-} \geq 0 \text{ et }
T_1^- \wedge f^*S_{T_1^-}=\l_1(f,T_1^-) T_1^- \wedge S_{T_1^-}.
$$
Autrement dit, le courant $S_{T_1^-}$ est positif et $f$-invariant
uniquement par ``restriction'' au courant $T_1^-$.
Observons que $T'$ est un courant positif ferm\'e 
invariant de bidegr\'e $(2,2)$, avec
$$
\l_{T'}=\frac{\l_1(f,T_1^-)}{\l_1(f^{-1})} \leq \frac{\l_1(f)}{\l_1(f^{-1})} <1,
$$
car $f$ est cohomologiquement hyperbolique. On en d\'eduit
$$
\l_1(f,T')=\l_{T'}^{-1}>1.
$$

Une nouvelle application du Th\'eor\`eme 5.13 permet alors d'obtenir
la mesure de probabilit\'e invariante
$$
\mu_f:=T' \wedge S_{T'}=T_1^- \wedge S_{T_1^-} \wedge S_{T'} \geq 0.
$$
Celle-ci ne charge pas les ensembles de petite dimension  de Hausdorff,
car chacun des courants $T_1^-,S_{T_1^-},S_{T'}$ est \`a potentiels
h\"old\'eriens (voir corollaire 1.4).
La canonicit\'e des constructions permet de montrer qu'il s'agit bien
de la mesure construite en 5.2.2.
\vskip.2cm

Il reste \`a \'etablir l'existence des courants $S_T$, ce que nous 
faisons \`a pr\'esent. On fixe $\om$ une forme de K\"ahler sur $X$.

\begin{thm}
Soit $T$ un courant positif ferm\'e de bidegr\'e $(p,p)$ sur $X$,
$1 \leq p \leq k-1$, qui est invariant, $f^*T=\l_T T$ avec $\l_T>0$.

Si $\l_1(f,T)>1$ alors il existe un entier $l_T \in \N$ et un courant ferm\'e
$S_T$ de bidegr\'e $(1,1)$ \`a potentiels h\"old\'eriens tels que
$$
T \wedge \left( \frac{1}{n} \sum_{j=1}^{n} \frac{1}{j^{l_T} \l_1(f,T)^j}
(f^j)^* \om \right) \longrightarrow \hat{T}:=T \wedge S_T,
$$
o\`u $\hat{T}$ est un courant positif ferm\'e non nul de bidegr\'e $(p+1,p+1)$.
Le courant $S_T$ est invariant sur $T$, au sens o\`u 
$$
T \wedge f^* S_T=T \wedge \l_1(f,T) S_T.
$$
\end{thm}

Notons que le courant $S_T$ n'est pas n\'ecessairement positif, mais
que $\hat{T}=T \wedge S_T$ l'est. Le produit d'intersection est bien
d\'efini car $S$ est \`a potentiels $\nu$-H\"older, donc born\'es [BT].
L'exposant $\nu>0$ est, comme dans la Proposition 1.2, contr\^ol\'e
par le rapport $\log \l_1(f,T)/\chi_{top}(f)$. Observons enfin
que le courant $\hat{T}$ est lui aussi invariant,
$$
f^* \hat{T}=\l_T \l_1(f,T) \hat{T}.
$$
Lorsque $p=k-1$, $\hat{T}$ est donc une mesure invariante car 
$\l_1(f,T)=\l_{k-p}(f,T)=\l_T^{-1}$ (\'egalit\'e 3. de la page pr\'ec\'edente).

\begin{esquisse}
Soit $l$ le plus grand entier tel que les classes 
$\{\om\},f^* \{\om\},$
$\ldots,(f^l)^* \{\om\}$ soient lin\'eairement ind\'ependantes
dans $H^{1,1}_{\nu}(T,\R)$.
On note $E$ le sous-espace engendr\'e par ces classes: il est stable sous l'action
de l'op\'erateur $f^*$, et le rayon spectral de $f^*:E \rightarrow E$ est \'egal
\`a $\l_1(f,T)$. L'entier $l_T$ de l'\'enonc\'e correspond \`a la non-diagonalisibilit\'e
de cette action (la norme de $(f^n)^*_{|E}$ cro\^it comme $n^{l_T}\l_1(f,T)^n$).
Nous supposons pour simplifier l'exposition que $l=0$, donc $l_T=0$
et $f^* \{\om\}=\l_1(f,T) \om$.

Il existe donc $R$ un $(1,1)$-courant positif ferm\'e \`a potentiels h\"old\'eriens et
$u:X \rightarrow \R$ une fonction h\"old\'erienne, tels que
$$
\frac{1}{\l_1(f,T)}f^* \om=\om+dd^c u+R,
\; \text{ avec } \; T \wedge R=0.
$$
En appliquant l'op\'erateur $(f^{n-1})^*$ \`a cette equation fonctionnelle
on obtient
$$
\frac{1}{\l_1(f,T)^n}(f^n)^* \om=\om+dd^c u_n+R_n,
\; \; T \wedge R_n=0, \; 
\text{ et } \;
u_n:=\sum_{j=0}^{n-1} \frac{1}{\l_1(f,T)^j} u \circ f^j.
$$
On montre, comme dans la Proposition 1.2 que la suite $(u_n)$ converge 
uniform\'ement vers une fonction h\"old\'erienne $u_{\infty}$, et on
pose
$$
S_T:=\om+dd^c u_{\infty}.
$$
Observons que $S_T$ n'est a priori ni positif, ni invariant (on n'a aucun contr\^ole
sur les termes d'erreurs $R_n$), mais
$\hat{T}:=T \wedge S_T$ est positif, comme limite de courants positifs,
et invariant (le terme d'erreur s'annulant sur $T$).

Lorsque $l_T \geq 1$, il est n\'ecessaire de consid\'erer des moyennes de
C\'esaro. La construction de $S_T$ proc\`ede de la m\^eme id\'ee, mais
les d\'etails techniques sont plus d\'elicats et nous renvoyons le lecteur
\`a la Proposition 2.4 et au Th\'eor\`eme 3.1 de [DS 6] pour plus de d\'etails.
\end{esquisse}

\section{Endomorphismes polynomiaux de $\C^k$}

Nous consid\'erons
ici le cas des endomorphismes rationnels qui sont $l$-stables sur
l'espace projectif complexe $\P^k$. Quitte \`a travailler dans $\P^{k+1}$,
il suffit de consid\'erer le cas d'endomorphismes qui sont
polynomiaux dans une carte affine (voir paragraphe 5.1.2).
Nous commen\c{c}ons par consid\'erer le cas des automorphismes.
Les techniques mises en jeu sont proches de celles esquiss\'ees
pr\'ec\'edemment, nous nous contentons donc d'\'enoncer 
les principaux r\'esultats et donnons quelques exemples.

\subsection{Automorphismes r\'eguliers}

N.Sibony a introduit dans [S] une classe d'exemples
dynamiquement int\'eressants qui contient les automorphismes
\'etudi\'es par E.Bedford et V.Pam-buccian dans [BP].

\begin{defi}
Soit $f:\C^k \rightarrow \C^k$ un automorphisme polynomial non affine.
On dit que $f$ est r\'egulier lorsque 
l'extension m\'eromorphe de $f$ \`a $\P^k$ v\'erifie
$I_f \cap I_{f^{-1}}=\emptyset$.
\end{defi}

Il s'agit d'une g\'en\'eralisation pluridimensionnelle des applications de 
H\'enon complexes: celles-ci sont pr\'ecis\'ement les automorphismes r\'eguliers
lorsque $k=2$. On v\'erifie (voir [S]) que
\begin{itemize}
\item $f$ induit un endomorphisme $l$-stable de $\P^k$, avec $l=\dim_{\C} I_{f^{-1}}+1$;
\item $\l_j(f)=\l_1(f)^j$ si $1 \leq j \leq l$ et $\l_j(f)=\l_{k-1}(f)^{k-j}$
si $l \leq j \leq k$. En particulier le degr\'e $\l_l(f)=\l_1(f)^l=\l_{k-1}(f)^{k-l}$
domine strictement tous les autres degr\'es dynamiques;
\item il existe des courants invariants canoniques
$T_l^+$, $T_{k-l}^-$ qui admettent des potentiels suffisamment r\'eguliers pour
pouvoir d\'efinir leur intersection potentialiste $\mu_f=T_l^+ \wedge T_{k-l}^-$.
\end{itemize}

\begin{exa}
Consid\'erons l'endomorphisme de $\C^k$ d\'efini par
$$
f(z_1,\ldots,z_k)=(P_{k-1}(z_1,\ldots,z_{k-1})+a_kz_k,\ldots,P_1(z_1)+a_2z_2,a_1z_1),
$$
o\`u $a_j \in \C^*$ et les $P_j$ sont des polyn\^omes de degr\'e $d \geq 2$ 
tels que $\deg_{z_j} P_j=d$.
On v\'erifie que $f$ est un automorphisme de $\C^k$.
Son extension m\'eromorphe \`a $\P^k$ -- encore not\'ee $f$ -- est telle que
$I_f=(z_1=\cdots=z_{k-1}=t=0)$ est r\'eduit \`a un point qui n'appartient pas \`a 
$I_{f^{-1}}=\{z_k=t=0\}$. Ici $(t=0)$ d\'esigne l'hyperplan \`a l'infini.
Ainsi $f$ est r\'egulier avec dans ce cas $l=k-1$.

En intervertissant les r\^oles de $f,f^{-1}$, on obtient des exemples
tels que $l=k-1$. Pour fabriquer des exemples tels que $2 \leq l \leq k-2$, lorsque
$k \geq 4$, on peut utiliser l'observation suivante: soit
$f_i:\C^{k_i} \rightarrow \C^{k_i}$, $i=1,2$, des automorphismes 
r\'eguliers de $\C^{k_i}$ avec 
$\l_1(f_1)=\l_1(f_2)$. Posons
$l_i=1+\dim_{\C} I_{f_i^{-1}}$. Alors le produit
direct $f=f_1 \times f_2$ d\'efinit un automorphisme polynomial r\'egulier
de $\C^k$, $k=k_1+k_2$, tel que $\dim_{\C} I_{f^{-1}}=l_1+l_2-1$.
\end{exa}

Nous avons montr\'e dans [G 6] le r\'esultat suivant:

\begin{thm}
Soit $f:\C^k \rightarrow \C^k$ un automorphisme polynomial r\'egulier.
Alors la mesure $\mu_f:=T_l^+ \wedge T_{k-l}^-$ est m\'elangeante 
et d'entropie maximale,
$$
h_{\mu_f}(f)=h_{top}(f)=\log \l_l(f)>0.
$$
\end{thm}

Le m\'elange d\'ecoule comme pour le Th\'eor\`eme 4.23 de l'extr\'emalit\'e
des courants $T_l^+$, $T_{k-l}^-$ que nous \'etablissons dans [G 6].
Notons que cette propri\'et\'e d'extr\'emalit\'e a \'et\'e
d\'emontr\'ee ind\'ependamment par T.-C.Dinh et N.Sibony [DS 7].
\vskip.2cm

Il est naturel de se demander si tous les automorphismes polynomiaux
cohomologiquement hyperboliques sont -conjugu\'es \`a- des automorphismes r\'eguliers.
C'est le cas en dimension 2 (voir [FrM]). Ce n'est plus vrai en dimension sup\'erieure.
On peut s'en rendre compte en analysant les automorphismes quadratiques de $\C^3$
qui ont \'et\'e partiellement classifi\'es par J.-E.Fornaess et H.Wu [FW] 
(voir \'egalement [Mae]).
Nous avons pr\'ecis\'e cette classification dans [CoG],
en collaboration avec D.Coman (voir Proposition 5.23 ci-apr\`es).

\subsection{Automorphismes faiblement r\'eguliers}

Soit $f:\C^k \rightarrow \C^k$ un automorphisme polynomial. On note 
$$
X_f:=\overline{f^k((t=0) \setminus I_{f^k})},
$$
o\`u $(t=0)$ d\'esigne l'hyperplan \`a l'infini, $\P^k=\C^k \cup (t=0)$.
Nous avons introduit, en collaboration avec N.Sibony [GS],
la notion suivante:

\begin{defi}
On dit que $f$ est faiblement r\'egulier si $I_f \cap X_f=\emptyset$.
\end{defi}

Tout automorphisme r\'egulier est faiblement r\'egulier
car dans ce cas $X_f=I_{f^-1}$
(voir Proposition 2.5.3 dans [S]). Un produit direct
d'automorphismes r\'eguliers $f_1,f_2$ (resp. faiblement r\'eguliers) est un
automorphisme faiblement r\'egulier (mais il n'est r\'egulier que
lorsque $\l_1(f_1)=\l_1(f_2)$).

\noindent On v\'erifie dans ce cas (voir section 2 dans [GS]) que
\begin{itemize}
\item $f$ induit un endomorphisme $l$-stable sur $\P^k$ avec $l=\dim_{\C} X_f+1$;
\item l'ensemble $X_f$ est un attracteur;
\item il existe un courant invariant $T_l^+$ canonique tel que $f^* T_l^+=\l_l(f) T_l^+$;
\item $\dim_{\C}I_f=k-l-2$, donc $\l_l(f)=\l_1(f)^l$ mais $\l_{l+1}(f)<\l_1(f)^{l+1}$.
\end{itemize}

\begin{thm}
Soit $f$ un automorphisme polynomial faiblement r\'egulier de $\C^k$
tel que $\l_l(f) > \l_{l+1}(f)$, $l=\dim_{\C} X_f+1$.
Alors il existe un courant positif ferm\'e
$T_{k-l}^-$ de bidegr\'e $(k-l,k-l)$,
tel que $f_* T_{k-l}^-=\l T_{k-l}^-$ et
$$
\frac{1}{\l^n}(f^n)_* \om^{k-l} \longrightarrow T_{k-l}^-,
$$
o\`u $\om$ d\'esigne la forme de Fubini-Study sur $\P^k$.
\end{thm}

\begin{rqe}
Ce r\'esultat est d\'emontr\'e dans [GS] (Th\'eor\`eme 3.1) sous des hypoth\`eses
non-optimales. Il est d\'emontr\'e  dans
[CoG] (section 2.1) lorsque $l=1$, et dans [G 6] (section 2.2)
dans le cas g\'en\'eral. 

Notons que l'hypoth\`ese $\l_l(f)>\l_{l+1}(f)$
entra\^ine, gr\^ace aux in\'egalit\'es de concavit\'e 2.4.a, que $\l_l(f)$
domine strictement tous les autres degr\'es dynamiques.
Observons \'egalement que cette hypoth\`ese est toujours satisfaite pour les
automorphismes quadratiques de $\C^3$ (voir Proposition 5.23).
\end{rqe}

Comme $X_f$ est un attracteur, le courant $T_{k-l}^-$ admet de bons potentiels.
Plus pr\'ecis\'ement, il r\'esulte de sa construction qu'on peut l'\'ecrire
$$
T_{k-l}^-=\Theta+dd^c (\T_{\infty}^-),
$$
o\`u $\Theta$ est une forme lisse ferm\'ee de bidegr\'e $(k-l,k-l)$ qui
est cohomologue \`a $\om^{k-l}$ , et 
$\T_{\infty}^-$ est un courant de bidegr\'e $(k-l-1,k-l-1)$
qui peut \^etre choisi {\it positif} hors d'un voisinage arbitrairement petit
de l'attracteur $X_f$. 
C'est l'analogue de la construction faite au paragraphe 5.2.2: la condition
{\bf Spec(f,l)} est ici trivialement v\'erifi\'ee puisque
$H^{l,l}(X,\C) =\C$, mais le prix \`a payer de cette simplification cohomologique
est la pr\'esence de points d'ind\'etermination.

Nous montrons dans [G 6] que
le courant $T_{k-l}^-$ est extr\'emal
et que la mesure $\mu_f:=T_l^+ \wedge T_{k-l}^-$ est bien d\'efinie,
en utilisant une variante du Lemme 5.8.
Lorsque $I_f$ est $f^{-1}$-attirant, la mesure $\mu_f$ est \`a support compact
dans $\C^k$, ce qui simplifie son \'etude.
Nous obtenons ainsi [G 6]:

\begin{thm}
Soit $f:\C^k \rightarrow \C^k$ un automorphisme faiblement r\'egulier.
Supposons $\l_l(f)>\l_{l+1}(f)$ et
que l'ensemble $I_f$ est un attracteur pour $f^{-1}$.

Alors la mesure $\mu_f=T_l^+ \wedge T_{k-l}^-$ est une mesure de probabilit\'e invariante
qui ne charge pas les hypersurfaces. C'est une mesure m\'elangeante
d'entropie maximale,
$$
h_{\mu_f}(f)=h_{top}(f)=\log \l_l(f).
$$
\end{thm}

Nous renvoyons le lecteur aux Th\'eor\`emes 3.1 et 3.2 dans [G 6] pour une d\'emonstration
de ce r\'esultat. Nous pensons qu'il n'est pas n\'ecessaire de
supposer que l'ensemble
$I_f$ est un attracteur.

\begin{ques}
Est-ce que les r\'esultats du Th\'eor\`eme 5.20 subsistent 
lorsque $I_f$ n'est pas un attracteur ?
\end{ques}

Nous montrons dans le Th\'eor\`eme 4.1 de [CoG] que l'ensemble $I_f$ est 
{\it souvent} (mais pas toujours) $f^{-1}$-attirant pour les
automorphismes polynomiaux quadratiques de $\C^3$ dont 
les degr\'es dynamiques sont deux \`a deux distincts, i.e.
qui sont cohomologiquement hyperboliques.

\subsection{Le cas des endomorphismes}

La plupart des r\'esultats pr\'ec\'edents s'\'etendent au cas des endomorphismes
polynomiaux non inversibles de $\C^k$ qui sont faiblement r\'eguliers,
moyennant une hypoth\`ese technique sur l'ensemble critique (hypoth\`ese
(H4) dans [G 6]).

Nous avons d\'emontr\'e (Th\'eor\`eme 2.6 dans [G 6]) l'extr\'emalit\'e de
$T_{k-l}^-$ uniquement dans le c\^one des courants positifs invariants;
celle-ci entra\^ine l'ergodicit\'e de la mesure $\mu_f$
(c'est l'analogue du Th\'eor\`eme 4.23), et on peut \'etablir le m\'elange faible
en consid\'erant l'automorphisme produit $f \otimes f$.

Il est probable que $T_{k-l}^-$ soit fortement extr\'emal. Notre incapacit\'e 
\`a le prouver r\'esulte uniquement de la mauvaise compr\'ehension que
nous avons des courants de bidegr\'e $(k-l,k-l)$, lorsque
$k-l>1$.
Lorsque $k-l=1$ et $k=2$, le courant $T_{k-l}^-$ est effectivement
fortement extr\'emal comme nous l'avons montr\'e dans [G 1], Th\'eor\`eme 5.2: dans 
le cas de bidegr\'e $(1,1)$, les potentiels des courants sont des fonctions
(uniques \`a constante additive pr\`es), ce qui facilite \'enorm\'ement
leur analyse dynamique. 

Nous renvoyons le lecteur \`a [G 6] pour la formulation pr\'ecise des
r\'esultats.

\section{Exemples}

\subsection{Automorphismes cohomologiquement hyperboliques}

Le tore complexe $X=(\C/\Z[i])^k$, $k \geq 2$, admet de nombreux automorphismes
d'entropie positive, induits par une matrice $A \in GL(k,\Z)$.
Les automorphismes cohomologiquement hyperboliques sont ceux pour lesquels
la matrice $A$ n'a aucune valeur propre de module 1 (automorphismes
d'Anosov). 
On v\'erifie ais\'ement que $\l_l(f)>\max_{j \neq l} \l_j(f)$ est une valeur propre
distingu\'ee de l'op\'erateur $f^*:H^{l,l}(X,\C) \rightarrow H^{l,l}(X,\C)$
(condition {\bf Spec(f,l)}):
la strat\'egie propos\'ee pour construire $\mu_f$ fonctionne donc
bien dans ce cas.
Nous renvoyons le lecteur \`a l'Appendix de [GV] pour la classification
des automorphismes d'entropie positive sur les tores complexes compacts de 
dimension deux.
\vskip.1cm

B.Mazur donne dans [Maz] des exemples d'automorphisme d'entropie positive
sur certaines surfaces $K3$. Sa construction fonctionne en toute dimension
comme l'ont observ\'e T.C.Dinh et N.Sibony [DS 6].

\begin{exa}
Soit $P(z^0,\ldots,z^k)$ un polyn\^ome multihomog\`ene de multidegr\'e
$(2,\ldots,2)$ en $z^0=(x_0,y_0),\ldots,z^k=(x_k,y_k)$, i.e.
tel que 
$$
P(\l_0z^0,\ldots,\l_kz^k)=\l_0^2\cdots \l_k^2 P(z^0,\ldots,z^k),
\text{ pour tout } \l_0,\ldots,\l_k \in \C^*.
$$
Un tel polyn\^ome d\'efinit une hypersurface $X$ de dimension $k$ et de
multidegr\'e $(2,\ldots,2)$ dans l'espace $(\P^1)^{k+1}=\P^1 \times \cdots \times \P^1$.
Cette hypersurface est une vari\'et\'e de Calabi-Yau lisse pour un choix
g\'en\'erique de $P$. Consid\'erons
$$
\pi_i: X \subset (\P^1)^{k+1} \rightarrow (\P^1)^k
$$
la projection parall\`element \`a la  $i^e$ coordonn\'ee. C'est un rev\^etement holomorphe
de degr\'e $2$ sur $(\P^1)^k$, qui permet de d\'efinir $\sigma_i:X \rightarrow X$
l'involution holomorphe qui \'echange les deux pr\'eimages de $\pi_i$.
Soit enfin
$$
f:=\sigma_0 \circ \sigma_1 \circ \cdots \circ \sigma_k:X \rightarrow X.
$$
Alors $f$ est un automorphisme cohomologiquement hyperbolique.
On trouvera une preuve de ce fait dans [Ca 1], lorsque $k=2$.
\end{exa}

Le lecteur trouvera d'autres exemples dans [BK 1,2], [Ca 1], [Ke], [M 3].

\subsection{Automorphismes polynomiaux de $\C^3$}

Les automorphismes polynomiaux quadratiques de $\C^3$ ont \'et\'e 
classifi\'es par J.E.Fornaess et H.Wu dans [FW]
(voir \'egalement [Mae]). Nous avons, en collaboration avec D.Coman [CoG], 
pr\'ecis\'e cette classification en nous int\'eressant aux automorphismes
cohomologiquement hyperboliques.

\begin{pro}
Soit $f:\C^3 \rightarrow \C^3$ un automorphisme polynomial quadratique  
cohomologiquement hyperbolique.
Alors apr\`es conjugaison, et quitte \`a changer  $f$
en $f^{-1}$ ou $f^2$, on obtient que
\begin{itemize}
\item $f$ est $1$-stable dans $\P^3$ et $\l_1(f) >\max_{j \neq 1} \l_j(f)$;
\item $f$ est soit r\'egulier, soit faiblement r\'egulier.
\end{itemize}
\end{pro}

Notons que pour un automorphisme polynomial de $\C^3$, on a 
$\l_3(f)=1$ et $\l_2(f)=\l_1(f^{-1})$. Les automorphismes 
cohomologiquement hyperboliques sont donc ceux tels que $\l_1(f) \neq \l_1(f^{-1})$.

Nous renvoyons le lecteur au Th\'eor\`eme 4.1 de [CoG] pour une preuve ainsi qu'un \'enonc\'e 
plus pr\'ecis. Ce r\'esultat montre qu'il est n\'ecessaire de consid\'erer
la notion plus g\'en\'erale d'automorphisme faiblement r\'egulier et que les
automorphismes quadratiques de $\C^3$ v\'erifient les hypoth\`eses du
Th\'eor\`eme 5.18. La plupart v\'erifient \'egalement l'hypoth\`ese
suppl\'ementaire du Th\'eor\`eme 5.20, mais pas tous comme le montre l'exemple suivant.

\begin{exa}
Consid\'erons
$$
f(x,y,z) \in \C^3 \mapsto (xy+az,x^2+by,x) \in \C^3,
$$
o\`u $a,b \in \C^*$. C'est un automorphisme quadratique de $\C^3$ (membre de la famille 5
du Th\'eor\`eme 4.1 de [CoG]).

On v\'erifie ais\'ement que l'inverse $f^{-1}$ de $f$ satisfait les hypoth\`eses
du Th\'eor\`eme 5.18: $f^{-1}$ est faiblement r\'egulier avec $l=1$,
$\l_1(f^{-1})=3>\l_2(f^{-1})=\l_1(f)=2$. Cependant $I_{f^{-1}}$ n'est pas toujours
$f$-attirant. En effet la droite $(x=z=0)$ est invariante,
$$
f(0,y,0)=(0,by,0), \; \; \; f^{-1}(0,y,0)=(0,b^{-1}y,0).
$$
Elle rencontre l'hyperplan \`a l'infini en un point qui est 
\`a la fois d'ind\'etermina-tion pour $f$ et pour $f^{-1}$.
Il s'ensuit que $I_{f^{-1}}$ n'est pas $f$-attirant si $|b| \leq 1$.
Nous montrons dans [CoG] que $I_{f^{-1}}$ est $f$-attirant lorsque $|b|>1$.
\end{exa}

Le lecteur int\'eress\'e trouvera dans [BP], [S], [G 1,3,6], [GS], [CoG], [DS 2] d'autres exemples
d'automorphismes et d'endomorphismes polynomiaux de $\C^k$ qui sont cohomologiquement
hyperboliques et satisfont les hypoth\`eses suppl\'ementaires des r\'esultats
\'enonc\'es plus haut.

\vskip 3cm

Vincent Guedj

Laboratoire Emile Picard

UMR 5580, Universit\'e Paul Sabatier

118 route de Narbonne

31062 TOULOUSE Cedex 04 (FRANCE)

guedj@picard.ups-tlse.fr

\end{document}